\documentclass[10pt]{report}
\usepackage[latin1]{inputenc}
\usepackage{graphicx}
\usepackage{amsmath}
\usepackage{amssymb}
\usepackage{amsmath,amsthm,amsfonts,amssymb,amscd}
\usepackage{enumerate}
\usepackage[latin1]{inputenc}
\usepackage{graphicx}

\usepackage{makeidx}
\makeindex

\def\esimo{${}^{\text{\b o}}$}

\newcommand{\co}{\mathbb C}

\def\hot{\operatorname{{h.o.t.}}}
\def\SL{\operatorname{{SL}}}
\def\Res{\operatorname{{Res}}}

\def\Fix{\operatorname{{Fix}}}

\def\dim{\operatorname{{dim}}}
\def\Aut{\operatorname{{Aut}}}
\def\GL{\operatorname{{GL}}}
\def\Aff{\operatorname{{Aff}}}

\def\Hol{\operatorname{{Hol}}}

\def\Ker{\operatorname{{Ker}}}

\def\Dom{\operatorname{{Dom}}}
\def\Id{\operatorname{{Id}}}
\def\Tni{\operatorname{{Int}}}

\def\Diff{\operatorname{{Diff}}}
\def\sing{\operatorname{{sing}}}
\def\Sing{\operatorname{{Sing}}}
\def\codim{\operatorname{{codim}}}

\def\Det{\operatorname{{Det}}}
\def\deg{\operatorname{{deg}}}

\def\Spec{\operatorname{{Spec}}}

\def\sep{\operatorname{{sep}}}

\def\arg{\operatorname{{arg}}}
\def\virt{\operatorname{{virt}}}

\def\Mon{\operatorname{{Mon}}}

\def\L{{\mathcal{L}}}
\def\G{{\mathcal{G}}}
\def\fa{{\mathcal{F}}}
\def\O{{\mathcal{O}}}
\def\eR{{\mathcal{R}}}
\def\eS{{\mathcal{S}}}

\def\U{{\mathcal{U}}}

\def\E{{\mathcal{E}}}

\def\U{{\mathcal{U}}}

\def\L{{\mathcal{L}}}
\def\G{{\mathcal{G}}}
\def\fa{{\mathcal{F}}}
\def\O{{\mathcal{O}}}
\def\eR{{\mathcal{R}}}

\def\U{{\mathcal{U}}}

\def\E{{\mathcal{E}}}

\def\U{{\mathcal{U}}}

\def\po{{\partial}}

\def\te{{\theta}}

\def\om{{\omega}}
\def\Om{{\Omega}}
\def\vr{{\varphi}}
\def\ga{{\gamma}}
\def\Ga{{\Gamma}}
\def\la{{\lambda}}
\def\La{{\Lambda}}
\def\ov{\overline}
\def\al{{\alpha}}
\def\ve{{\varepsilon}}
\def\lg{{\langle}}
\def\rg{{\rangle}}

\def\be{{\beta}}

\def\bh{{\mathbb{H}}}
\def\bp{{\mathbb{P}}}

\def\re{{\mathbb{R}}}
\def\bz{{\mathbb{Z}}}
\def\bq{{\mathbb{Q}}}
\def\bd{{\mathbb{D}}}
\def\bc{{\mathbb{C}}}
\def\bn{{\mathbb{N}}}

\def\bh{{\mathbb{H}}}
\def\bp{{\mathbb{P}}}

\def\re{{\mathbb{R}}}
\def\bz{{\mathbb{Z}}}
\def\bq{{\mathbb{Q}}}
\def\bd{{\mathbb{D}}}
\def\bc{{\mathbb{C}}}
\def\bn{{\mathbb{N}}}

\def\bl{{\mathbb{L}}}

\def\Aut{\operatorname{{Aut}}}
\def\GL{\operatorname{{GL}}}
\def\Aff{\operatorname{{Aff}}}
\def\Hol{\operatorname{{Hol}}}

\def\Ker{\operatorname{{Ker}}}

\def\Dom{\operatorname{{Dom}}}
\def\Id{\operatorname{{Id}}}
\def\Tni{\operatorname{{Int}}}

\def\Diff{\operatorname{{Diff}}}
\def\sing{\operatorname{{sing}}}
\def\Sing{\operatorname{{Sing}}}
\def\codim{\operatorname{{codim}}}

\def\deg{\operatorname{{deg}}}

\def\sep{\operatorname{{sep}}}

\def\arg{\operatorname{{arg}}}
\def\virt{\operatorname{{virt}}}

\def\SL{\operatorname{{SL}}}
\def\Res{\operatorname{{Res}}}

\def\dim{\operatorname{{dim}}}
\def\Aut{\operatorname{{Aut}}}
\def\GL{\operatorname{{GL}}}
\def\Aff{\operatorname{{Aff}}}
\def\Hol{\operatorname{{Hol}}}

\def\Ker{\operatorname{{Ker}}}

\def\Dom{\operatorname{{Dom}}}
\def\Id{\operatorname{{Id}}}
\def\Tni{\operatorname{{Int}}}

\def\Diff{\operatorname{{Diff}}}
\def\sing{\operatorname{{sing}}}
\def\Sing{\operatorname{{Sing}}}
\def\codim{\operatorname{{codim}}}

\def\deg{\operatorname{{deg}}}

\def\sep{\operatorname{{sep}}}

\def\arg{\operatorname{{arg}}}
\def\Inv{\operatorname{{Inv}}}

\newtheorem{Theorem}{Theorem}[section]
\newtheorem{Corollary}[Theorem]{Corollary}
\newtheorem{Proposition}[Theorem]{Proposition}
\newtheorem{Lemma}[Theorem]{Lemma}

\newtheorem{Claim}[Theorem]{Claim}

\theoremstyle{Definition}
\newtheorem{Definition}[Theorem]{Definition}
\newtheorem{Example}[Theorem]{Example}
\newtheorem{Question}[Theorem]{Question}
\newtheorem{Problem}[Theorem]{Problem}

\newtheorem{Fact}[Theorem]{Fact}
\newtheorem{Exercise}[Theorem]{Exercise}
\theoremstyle{Remark}
\newtheorem{Remark}[Theorem]{Remark}
\numberwithin{section}{chapter} \numberwithin{equation}{chapter}
\numberwithin{figure}{chapter}

\def\SL{\operatorname{{SL}}}
\def\Res{\operatorname{{Res}}}

\def\dim{\operatorname{{dim}}}
\def\Aut{\operatorname{{Aut}}}
\def\GL{\operatorname{{GL}}}
\def\Aff{\operatorname{{Aff}}}
\def\Hol{\operatorname{{Hol}}}

\def\Ker{\operatorname{{Ker}}}

\def\Dom{\operatorname{{Dom}}}
\def\Id{\operatorname{{Id}}}
\def\Tni{\operatorname{{Int}}}

\def\Diff{\operatorname{{Diff}}}
\def\sing{\operatorname{{sing}}}
\def\Sing{\operatorname{{Sing}}}
\def\codim{\operatorname{{codim}}}

\def\deg{\operatorname{{deg}}}

\def\sep{\operatorname{{sep}}}

\def\arg{\operatorname{{arg}}}

\title{A swift  introduction to holomorphic foliations with singularities}

\author{Bruno Sc\'ardua}

\date{}

\makeindex

\begin{document}
\maketitle

\tableofcontents

\chapter*{Preface}
The theory of foliations is one of those subjects in mathematics
that gathers several distinct domains such as topology, dynamical
systems and  geometry, among others. Its origins go back to the
works of C. Ehresmann and Shih (\cite{E-R}, \cite{E-Sh}) and G. Reeb
(\cite{Reeb1,Reeb2}). It provides an interesting and valuable
approach to the qualitative study of dynamics and ordinary
differential equations on manifolds.

Although  its origins are in the classical framework of real
functions and manifolds,  the notion of foliation is also very
useful in the holomorphic world. Indeed, it has ancient origins in
the study of complex differential equations. From these first
problems, the introduction of singularities as an object of study is
a natural step. We mention the works of P. Painlev\'e
(\cite{Painleve1,Painleve2}) and Malmquist (\cite{Malmquist}).
With P. Painlev\'e the study of rational complex differential equations of the form
$\frac{dy}{dx} = \frac{P(x,y)}{Q(x,y)}$
has its first more specific methods and results.
After Painlev\'e many authors have contributed for the initial push
up of the theory, among  them are E. Picard, G. Darboux, H.
Poincar\'e, H. Dulac, Briot and Bouquet.

Complex differential equations appear naturally in  mathematics and
in natural sciences (\cite{Arnold1,Arnold,hille}). For instance, we
mention the theory of electrical circuits, valves and
electromagnetic  waves (\cite{Clayton Paul}). Another motivation is
the search for and study of new (classes of) transcendent functions,
as the Liouvillian functions (\cite{singer}).

With the advent of the geometric theory of foliations and the modern
results of Cartan, Oka, Nishino, Suzuki  and others, on the theory
of analytic functions of several complex variables and some from
algebraic and analytic geometry, this field of research became quite
active again. To these days it is one of the active branches of
modern research in mathematics.

These are the notes of a series of lectures delivered by the author
at the Graduate School of Mathematical Sciences of the University of
Tokyo, during the month of October 2015. They were meant to be as
self-contained as possible, taking into account time and space. The
basic idea was to introduce the concepts and some  of the basic
results in the theory of holomorphic foliations with singularities.
Another goal is to guide the reader to some of the recent questions
and problems in the field, providing in this way a motivating
introduction to those who are interested in studying a new subject.

I wish to express my gratitude to Professor Taro Asuke for his
personal effort in making  this project possible and for his warm
hospitality. I want to thank Professor T. Tsuboi for his kind
hospitality and support. I want to thank all those at the Graduate
School of Mathematical Sciences of the University of Tokyo for their
support and for making my stay in Tokyo such a pleasant and fruitful
period.

The first draft of this text was conceived  during a visit to the
Instituto de Matemática y Ciencias Afines - IMCA (Lima). It is my
honor to present these notes in the series Monografias del IMCA. I
am very much indebted with this prestigious Institute for all the
support during the last couple of decades.  My special thanks to
Professors Félix Escalante and Roger Metzger for their warm
hospitality. \vglue.2in

\rightline{Rio de Janeiro, November 2018}
 \rightline{Bruno
Sc\'ardua}

\chapter{The classical notions of foliations}
\label{chapter:classicalnotions} This chapter is intended to
introduce the classical notions of foliation in the real framework.
The reader which is already familiar with these notions may skip to
the next chapter. We refer to \cite{camacho-linsneto},
\cite{Godbillon}, \cite{Hector-Hirsch} or \cite{morales-scardua} for
a more complete exposition of  the theory of real foliations.

\section{Definition of foliation}
\label{section:conceptoffoliation}

There are some ways of motivating  the concept of foliation.
Probably, the very first  is given by a  submersion $f\colon M \to
N$ from a manifold $M$ into a manifold $N$. If $f$ is sufficiently
differentiable (usually of class $C^r$, $r\ge2$) then by the local
form of submersions, the level sets $f^{-1}(y)$, $y \in N$ are
embedded submanifolds of $M$. These fibers  are  {\it locally\/}
organized as the fibers of a projection $(x,y) \mapsto y$. This
local picture is not necessarily global, and the fibers may be
disconnected.

\par A second important example is given by a closed non-singular
$1$-form $\om$ on a manifold $M$.  Again, under sufficient
differentiability conditions, by the  integration lemma of Poincaré
we can write {\it locally\/} $\om = df$ for a submersion map $f$
taking values on $K$. Here  $K$ is the field of real numbers if
$\omega$ is differentiable and $M$ is a real differentiable
manifold. In case $M$ is a complex manifold and $\omega$ is  a
holomorphic $1$-form on  $M$ we have  $K=\bc$ the field of complex
numbers. In this later case $f$ is holomorphic. Any local function
$f$ as above, defined in an open subset in $M$, is called a {\it
first integral}\index{first integral} for $\omega$.
\newline Notice that two local {\it first integrals\/} $f$ and
$\widetilde f$ for $\om$ in a same connected subset of $M$ are
related by $\widetilde f = f +$ constant. Therefore, they share
level sets, these local sets can therefore be globalized as immersed
(locally closed) submanifolds of $M$, again  locally organized as
fibers of a projection.

The third and last basic example we shall  mention is the one
provided by a  differentiable, namely $C^r$ with $r \geq 1$,  vector
field $X$ on a manifold $M$. Given a non-singular point $p \in M$
which is not a singular point of $X$, the flow-box theorem gives a
conjugation between $(X,U)$, where $p \in U\subset M$ is an open
neighborhood, and a constant vector field on $\mathbb R^m$, where $m
= \dim M$. The orbits of $X$ in $U$ then follow the same geometrical
condition of the above examples. The above examples motivate the
classical definition of foliation below  as follows:

\begin{Definition}[foliation, \cite{camacho-linsneto},\cite{Godbillon},\cite{Hector-Hirsch},\cite{morales-scardua}]
\label{Definition:foliation0} {\rm Given a differentiable  manifold
$M$ of dimension $m$ and class $C^r, r \geq 0$; by a {\it
codimension $0 \leq n \leq m$ foliation of class $C^r$ of
$M$}\index{foliation}, we mean an atlas $\mathcal F= \{(U_j,
\vr_j)\}_{j \in J}$ of $M$, where each coordinate chart $\vr_j
\colon U_j \subset M \to \vr_j(U_j)\subset \mathbb R^{m-n} \times
\mathbb R^{n}$ is of class $C^r$ and we have the following
compatibility condition:

\par{\em For each non-empty intersection $U_i \cap U_j \ne \phi$, the corresponding change of coordinates

\[
\vr_j\circ\vr_i^{-1}\big|_{\vr_i(U_i\cap U_j)}{ \colon \vr_i(U_i
\cap U_j) \longrightarrow \vr_j(U_i \cap U_j)}
\]
preserves the natural horizontal fibration $y=const.$  of $\mathbb
R^{m-n}\times\mathbb R^n \ni (x,y)$. }}
\end{Definition}
This is equivalent to say that, in  coordinates $(x,y) \in \mathbb
R^{m-n}\times\mathbb R^n$, we have
$$
\vr_j\circ\vr_i^{-1}(x,y) = \big(h_{ij}(x,y), g_{ij}(y)\big) \in
\mathbb R^{m-n}\times\mathbb R^n.
$$
The charts $\vr_j\colon U_j \to \vr_j(U_j)$ are called {\it
foliation charts, trivializing charts\/} or {\it distinguished
charts\/} of $\fa$. The local {\it plaques\/} of $\fa$ are the
fibers of a foliation chart in $\fa$. Given a diffeomorphism
$\vr\colon U\hookrightarrow \vr(U)\subset \mathbb R^m = \mathbb
R^{m-n} \times \mathbb R^{n}$ of class $C^r$, we say that $\vr$ is
{\it compatible} with the foliation $\fa$ if for any $j\in J$ such
that $U_j\cap U\ne \emptyset$, we also have
$$
\vr_j\circ\vr^{-1}(x,y) = \big(h(x,y), g(y)\big) \in \mathbb
R^{m-n}\times\re^n.
$$
In short, this is equivalent to say that $\fa \cup\{(U,\vr)\}$ is
still a foliation. Using this and Zorn's lemma, we may consider the
foliation atlas $\fa$ as {\it maximal}, in the sense that it
contains all the compatible charts of class $C^r$ of $M$.

In $M$ we consider the
equivalence relation induced by  the connected finite union of local
plaques. This means that two point $x, y \in M$ are equivalent
$x\thicksim y$ iff $x$ and $y$ lie in the same plaque of $\fa$ or
there is a finite number of plaques $P_1,...,P_r, \, r \geq 2$; of
$\fa$ such that $x \in P_1,\, y \in P_r$ and $P_i\cap P_{i+1}\ne
\emptyset$ for all $i=1,...,r-1$. Given a point $x \in M$ we call
the corresponding equivalence class $[x]\subset M$ is {\it the leaf
of $\fa$ through $x$}\index{leaf}. Usually we denote this leaf by
$\fa_x$ or by $L_x$. The leaf $L_x\subset M$ is an immersed $C^r$
submanifold, but not necessarily embedded. These leaves  then
decompose $M$ into disjoint immersed $C^r$ submanifolds. Each leaf
has dimension $m-n$ and meets a foliation chart domain along plaques
of the foliation. For instance, in the case of a submersion $f
\colon M \to N$, the leaves of the corresponding foliation are the
connected components of the level sets $f^{-1}(y), y \in N$. The
quotient space $M/\thicksim$ is the {\it leaf space}\index{leaf
space} of $\fa$, also denoted by $M/\fa$.

\section{Other  definitions of foliation}

According to the literature there are essentially three ways to
define foliations in the real differentiable manifolds (cf.
\cite{camacho-linsneto},\cite{Godbillon}, \cite{Hector-Hirsch}). In
addition to the one we just have given in
Definition~\ref{Definition:foliation0} above, we have the following.
Let $M$ be a $m$-dimensional manifold, $m\in {\mathbb N}$. Let $D^k$
be the open unit ball of $\mathbb R^k$ where $k\in {\mathbb N}$. Let
$0\leq n\leq m$ be fixed.

\begin{Definition}
\label{def1}{\rm  A  {\it foliation}\index{foliation} of $M$, of
codimension $n$ and class $C^r$, is  a partition $\fa$ of $M$
consisting of  pairwise disjoint immersed $C^r$ submanifolds
$L\subset M$ of dimension $m-n$, distributed as follows: for each
point $x\in M$ there is  a neighborhood $U$ of $x$, and a $C^r$
diffeomorphism $\vr\colon U \to D^{m-n}\times D^{n}$, such that for
each $y\in D^{n}$ there is $L\in \fa$ satisfying
$$
\vr^{-1}(D^{m-n}\times y)\subset L.
$$
The elements of the partition $\fa$ are the {\em leaves of
$\fa$}\index{leaf}. The element $L_x$ of $\fa$ containing $x\in M$
is  the {\it leaf of $\fa$ containing} $x$\index{leaf}.}
\end{Definition}

We observe that,  not every decomposition of $M$ into immersed
submanifolds with the same dimension is a foliation (see
\cite{morales-scardua}).

The third definition of foliation uses the notion of distinguished
maps. Let ${\mathcal F}=\{(U_j,\vr_j)\,, j \in J \}$ be a foliation
of a manifold $M$ in the sense of Definition
~\ref{Definition:foliation0}. Then $\forall i,j$ the transition map
$\vr_j\circ(\vr_i)^{-1}$ has the form
$$
\vr_j\circ(\vr_i)^{-1}(x,y)= (f_{i,j}(x,y),g_{i,j}(y)).
$$

The map $g_{i,j}$ is a local diffeomorphism in its domain of
definition. This follows from the fact that the derivative of the
transition map is given by  $D(\vr_j\circ(\vr_i)^{-1})(x,y)\cdot
(v,w)=(\partial_xf_{i,j}(x,y)\cdot v, Dg_{i,j}(y)\cdot w), \, (v,w)
\in \mathbb R^{m-n} \times \mathbb R^n$. We define for all $i$ the
map $ g_i=\Pi_2\circ \vr_i $, where $\Pi_2$ is the projection onto
the second coordinate: $\Pi_2 \colon D^{m-n}\times D^{n} \to D^{n},
\, (x,y)\mapsto y$. We claim that  $g_j=g_{i,j}\circ g_i$. Indeed,
we have $g_{i,j} \circ g_i= g_{i,j} \circ \Pi_2 ^i \circ \vr_i=
\Pi_2 ^j \circ (\vr_j  \circ \vr^{-1}_i)\circ \vr_i= \Pi_j ^j \circ
\vr_j = g_j$. Therefore, a $C^r$ foliation ${\mathcal F}$ of
codimension $n$ of a manifold $M^m$ is equipped with an open cover
$\{U_i\}_{i \in I}$ of $M$ and $C^r$ submersions $g_i\colon U_i\to
D^{n}$ such that for all $i,j$ there is a local diffeomorphism
$g_{i,j}:V_i \subset D^{n}\to V_j \subset D^{n}$ satisfying the
cocycle relations
$$
g_j=g_{i,j}\circ g_i, \,\,\,\,g_{i,i}=\Id.
$$
The $g_i$'s are the {\em distinguished maps}\index{distinguished
map} of ${\mathcal F}$.

Conversely, suppose that $M^m$ admits  an open cover
$M=\bigcup\limits_{i \in I} U_i$ such that for each $i\in I$ there
is a  $C^r$ submersion $g_i\colon U_i\to D^{n}$ such that for all
$i,j$ there is a  diffeomorphism $g_{i,j}:V_i \subset D^{n}\to V_j
\subset D^{n}$ satisfying the cocycle relations above. By the local
form of the submersions we can assume that for each $i\in I$ there
is a $C^r$ diffeomorphism $\vr_i\colon U_i\to D^{m-n}\times D^{n}$
such that
$$
g_i=\Pi_2\circ \vr_i.
$$
since
$$
\Pi_2\circ( \vr_j\circ (\vr_i)^{-1}) =g_j\circ (\vr_i)^{-1}=
g_{i,j}\circ g_i\circ (\vr_i)^{-1}= g_{i,j}\circ \Pi_2,
$$
we have that the atlas
$$
{\mathcal F}= \{(U_i,\vr_i)\}_{i \in I}
$$
defines a foliation of class $C^r$ and codimension $n$ of $M$. The
above suggests the following equivalent definition of foliation.

\begin{Definition}
\label{of3} {\rm A  {\it foliation}\index{foliation}   of $M^m$ of
class $C^r$ and of codimension $n$, is given by the following:

\begin{enumerate}
\item An open  cover
$\{U_i:i\in I\}$ of $M$.
\item A family of $C^r$ submersions
$g_i\colon U_i\to D^{n}, \forall i \in I$; with the  following
compatibility property: $\forall i,j\in I$ with $U_i \cap U_j \ne
\emptyset$,   there is a local diffeomorphism $g_{i,j}:V_i \subset
D^{n}\to V_j \subset D^{n}$ satisfying the cocycle relations
$$
g_j=g_{i,j}\circ g_i, \,\,\,\,g_{i,i}=\Id.
$$
\end{enumerate}
The submersions $g_i$'s are  the {\em distinguished maps} of the
foliation ${\mathcal F}$.}
\end{Definition}

This last definition leads to several interesting definitions. For
instance, a foliation $\fa$ of $M$ is said to be {\em transversely
holomorphic} or {\em transversely affine} depending on whether, for
some convenient choice,  its distinguished maps $g_{i,j}$ are
holomorphic  or affine maps. We shall resume this subject later on.
In order to distinguish foliations, we shall use the following
definition.

\begin{Definition}
\label{d2} {\rm Two foliations ${\mathcal F}$ and ${\mathcal F}'$ of
manifolds $M$ and $M'$ respectively are {\it
$C^r$-equivalent\/}\index{foliation! equivalence} if there is a
$C^r$-diffeomorphism $h\colon M \to M'$ ($h$ is a homeomorphism if
$r=0$), sending leaves of ${\mathcal F}$ into leaves of ${\mathcal
F}'$. In other words, if ${\mathcal F}_x$ denotes the leaf of $\fa$
that contains $x \in M$ and ${\mathcal F}'_y$ denotes the leaf of
$\mathcal F'$ that contains $y \in M'$ then we have:
$$
h({\mathcal F}_x) = {\mathcal F}'_{h(x)},\,\,\forall\, x \in M.
$$}
\end{Definition}

The above notion can be stated for the case of holomorphic objects,
in the obvious way. This relation defines an equivalence in the
space of foliations.

\section{Frobenius theorem} \label{subsection:Frobeniustheorem}

Let $X$, $Y$ two vector fields on a manifold $M$ and $p \in M$ be
fixed. Denote by $X_t$ the local flow of $X$ and similarly by  $Y_t$
the local flow of $Y$ assuming that $X,Y \in C^r$, \,\, $r \ge 2$.
Given $p \in M$ we define $X_t^*(Y)(p) = DX_{-t}(X_t(p))\cdot
Y(X_t(p))\in T_p(M)$. Note that $X_t^*(X)(p) = X(p)$, \,\,
$\forall\, t$. All this holds for $|t|$ small enough.

\begin{Definition} {\rm  The {\it Lie bracket\/}
of $X,Y$ is the vector field $[X,Y]$ on $M$ defined at each point
$p\in M$ by
$$
L_X(Y)(p) = [X,Y](p) = \frac{d}{dt}\big|_{t=0} (X_t^*(Y)(p))\quad
X,Y \in C^r, \,\, r \ge 2.
$$

In local coordinates $(x_1,...,x_m)\in M$, the Lie bracket  $[X,Y]$
has the following form: writing
$$
X = \sum\limits_{i=1}^m a_i\,\dfrac{\po}{\po x_i}, \,\,\,\, Y =
\sum\limits_{i=1}^m b_i\,\dfrac{\po}{\po x_i}
$$
one has
$$
[X,Y] = \sum_{i,j=1}^m \left(a_i \,\frac{\po b_j} {\po x_i} -
b_i\,\frac{\po a_j}{\po x_i}\right) \frac{\po}{\po x_j}\,.
$$

When $X$ and $Y$ are defined in an open set of $\re^m$, the formula
above yields
$$
[X,Y] = DY(p)\cdot X(p) - DX(p)\cdot Y(p).
$$
A vector field $X$ on $M$ is {\em tangent} to a plane field $P$ on
$M$ (denoted by $X\in P$) if $X(p)\in P(p)$ for all $p\in M$.}
\end{Definition}
\vglue .1in

\begin{Definition} {\rm  A plane field $P$ on $M$
is {\it involutive\/} if $X,Y\in P$ $\Rightarrow$ $[X,Y]\in P$.}
\end{Definition}

\vglue .1in

\begin{Lemma}
If $\fa$ is a foliation, then its associated plane field $T\fa$ is
involutive.
\end{Lemma}

\begin{proof} Let $X$, $Y$ be two vector fields  tangent to $T\fa$. By
using suitable local coordinates
$(x,y)=(x_1,...,x_{m-n},y_1,...,y_{n})$ such that $\fa$ is given in
these coordinates by $y=(y_1,...,y_n)=const.$, one can assume that
$X,Y$ are  of the form
\begin{align*}
X(x,y) &= (f(x,y),0), \,\,\, Y(x,y) = (g(x,y),0)\\
[X,Y] &=
\begin{pmatrix}
\po_xg &\po_yg\\
0 &0
\end{pmatrix}
\begin{pmatrix}
f\\
0
\end{pmatrix} -
\begin{pmatrix}
\po_xf &\po_yf\\
0 &0
\end{pmatrix}
\begin{pmatrix}
g\\
0
\end{pmatrix} \\ & = (f\cdot \po_xf - f\cdot \po_xf,0).
\end{align*}
Hence $[X,Y]\in T\fa$ and the proof follows.
\end{proof}

A plane field $P$ on $M$ of dimension $k$ is {\it completely
integrable} if there is a foliation $\fa$ of $M$ of codimension $k$
such that $T\fa=P$, {\it i.e.}, if for each $p \in M$ we have
$T_p(\fa)=P(p)$. In particular, given $p \in M$, the space $P(p)$ is
the tangent space of the leaf $L_p\in \fa$ of $\fa$ that contains
$p$. From the above lemma, a completely integrable plane field is
involutive.

The converse of the above fact is a well-known result in theory of
foliations:

\begin{Theorem}[Frobenius theorem \index{Theorem! of
Frobenius}\cite{camacho-linsneto},\cite{Candel-Conlon},
\cite{Godbillon}] \label{Theorem:Frobenius}

Involutive plane
fields are completely integrable.
\end{Theorem}

\begin{Example} [integrable systems of differential forms]
\label{Example:integrableformssystem} {\rm Let
$\omega_1,...,\omega_r$ be differential $1$-forms of class $C^r$ on
a manifold $M$ and assume that they are linearly independent at each
point $p\in M^n$. We call the set $\mathcal
S:=\{\omega_1,...,\omega_r\}$ a {\it system} of $1$-forms on $M$. We
may consider the corresponding distribution $P(\mathcal S)$ of
$(n-r)$-dimensional planes defined as follows: given $p \in M$ we
set  $P(\mathcal S)(p)\subset T_pM$ as
$$
P(\mathcal S)(p)=\{v\in T_p M, \,\omega_j(p)\cdot v=0,
\,j=1,...,r\}.
$$
The system $\{\omega_1,...,\omega_r\}$ is called {\it integrable} if
we have  $d\omega_j \wedge \omega_1\wedge ...\wedge \omega_r=0$ for
all $j=1,...,r$. In particular a distribution given by a 1-form
$\omega$ is integrable iff $\omega \wedge d\omega=0$. This occurs
for instance if we have a closed 1-form $\omega$ with $\omega(p)\ne
0, \forall p\in M$.

The system $\mathcal S$ is integrable if and only the distribution
$P(\mathcal S)$  is involutive. Therefore, according to Frobenius
theorem above, $\mathcal S$ is  integrable if and only if
$P(\mathcal S)$ is completely integrable.

In the case of an integrable non-singular one-form $\omega$ we have
a codimension one foliation $\fa$ of $M$ which is defined by the
{\em Pffafian equation} $\omega=0$. }

\end{Example}

\section{Holonomy}
\label{section:holonomy}

The concept of holonomy of a foliation is motivated by the concept
of {\it return map\/} or {\it Poincaré map\/} of a periodic orbit of
a vector field.

\begin{figure}[h]
\begin{center}
\includegraphics[scale=1.2]{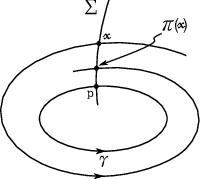}
\caption{Representation of a first return map $\pi \colon (\Sigma,p)
\to (\Sigma,p), \, x \mapsto \pi(x)$, of a periodic orbit $\gamma$
of a real vector field.}
\end{center}
\end{figure}

\begin{Remark}
{\rm Given a manifold $\Sigma$ of class $C^r$ and a point $p \in
\Sigma$ we denote by $\Diff^r(\Sigma,p)$ the group of germs of
diffeomorphisms of  class $C^r$ fixing $p \in \Sigma$. Each germ is
induced by a $C^r$ diffeomorphism $f\colon (U,p)\to (V,p)$, where
$U,V\subset \Sigma$ are neighborhoods of $p$, and $f(p)=p$. }
\end{Remark}

\noindent The ``return map"  $\pi\colon(\Sigma,p) \to (\Sigma,p)$,
$x \mapsto \pi(x)$ above illustrated is, always under suitable
differentiability conditions, a germ of  diffeomorphism. \newline
Identifying the transverse section $\Sigma$ with a disc in $\mathbb
R^{m-1}$ centered at the origin $0 \in \mathbb R^{m-1}$ ($m=$
dimension of the ambient manifold) we may consider $f$ as a (germ of
a) diffeomorphism $(\mathbb R^{m-1}{,\,0}) \to (\mathbb
R^{m-1}{,\,0})$ (fixing the origin). \vglue.1in

The same kind of idea above gives the concept of holonomy group of a
leaf of a foliation. Let us introduce it in a more formal way:

Let $\fa$ be a codimension $n$ foliation of class $C^r$ of a
manifold $M^m$. Given a leaf $L$ of $\fa$ we fix a point $p \in L$
and, by a local trivialization of $\fa$ around $p$, choose a
transverse section $p \in \Sigma_p \subset M$, diffeomorphic to a
disc in $\mathbb R^n$ and transverse to all the leaves of $\fa$.
Let now $\ga\colon [0,1] \to L$ be any closed path of
class $C^r$.

Then by suitable choice of the local foliation charts, we can cover
the image $\gamma([0,1])$ by domains of foliation charts $(\vr_j,
U_j), j=0,...,\ell$ with the following properties:

(a) There is a partition $0=t_0< t_1<\cdots < t_j < t_{j+1} < \ldots
< t_\ell =1$ of $[0,1]$ such that $\gamma([t_j,t_{j+1}])\subset
U_j$. In particular, we have $U_j\cap U_{j+1} \ne \emptyset$.

(b) For each $j=0,...,\ell-1$ the  union $U_j \cup U_{j+1}$
is contained in some foliation chart domain. Moreover, the same
holds for the union  $U_\ell \cup U_0$.

(c) There is a neighborhood $ p \in A \subset \Sigma$ such that for
each $y \in A$ there is a path of  plaques from the plaque $P^0 _y$
of $(\vr_0,U_0)$ through $y$, to the plaque $P(1) _y$ of
$(U_1,\vr_1)$ that meets $P^0_y$, then from the plaque $P(1)_y$ to
the plaque $P^2_y$ of $(U_2,\vr_2)$ that meets $P(1)_y$, and so on.
Thus we reach a plaque $P^\ell _ y$ of $(U_\ell, \vr_\ell)$  that
meets the plaque $P^{\ell -1} _y$. The intersection $P^\ell _y \cap
A$ is a single point $f_\gamma (y) \in \Sigma$ called the
$\gamma$-{\it holonomy image} of $y$.

\noindent This defines a map germ $f_\ga\colon (\Sigma,p) \to
(\Sigma,p)$ which has the following properties:

\noindent (i)\,\, $f_\ga$ has the same differentiability class as
that of $\fa$.

\noindent (ii)\, $f_\ga$ does not depend on the choice of the  cover
$\{(\vr_j,U_j)\}$ of $\ga([0,1])$\,\, (as a germ).

\noindent(iii) The map germ $f_\ga$ only depends on the homotopy
class $[\ga] \in \pi_1(L;p)$.

\noindent (iv)\, The map $\Hol_L: \pi_1(L;p) \to \Diff^r(\Sigma;p)$,
from the fundamental group of $L$ based at $p$, into the group of
germs of $C^r$ diffeomorphism of $\Sigma$ fixing $p$, and given by
given $[\ga]\mapsto f_\ga$, is a well-defined group homomorphism.
The image is denoted by $\Hol(\fa,L,\Sigma;p)\subset
\Diff^r(\Sigma;p)$.

\noindent (v)\, The groups $\Hol(\fa,L,\Sigma;p)$ depend on $\Sigma$
and $p$ by natural conjugation by diffeomorphisms.

Thus we may speak of the {\it holonomy group of the leaf\/} $L$
denoted by $\Hol(\fa,L)$ or simply by $\Hol(L)$, and identify it
with (a conjugacy class of) a subgroup of $\Diff^r(\mathbb R^n,0)$,
where $n$ is the codimension of the foliation $\fa$ of $M$.

We close this chapter illustrating the concepts above introduced by
means of a remarkable example.

\begin{Example}[foliations generated by closed $1$-forms\index{foliation! generated by closed $1$-form}]
{\rm In order to illustrate the above concepts of foliation and
holonomy we consider the class of  codimension one foliations
generated by closed differential $1$-forms. Indeed, we prove the
following:

\begin{Proposition} A codimension one smooth foliation $\fa$ tangent to a closed
non-singular $C^\infty$ $1$-form $\omega$ on a manifold $M$ has
trivial holonomy.
\end{Proposition}

\begin{proof}
Indeed, let us fix   a Riemannian metric $<,>$ on $M$. Let $X$ be the {\em gradient} of
$\omega$, {\it i.e.}, the smooth vector field  on $M$
defined by
$$
\omega_p(v_p)=\lg X(p),v_p\rg,
$$
for all $p\in M$ and $v_p\in T_pM$. Clearly $X$ is well-defined and
non-singular since $\omega$ is non-singular. In addition $X$ is
transverse to $\fa$. Let $L$ be a leaf of $\fa$ and $\gamma$ a
closed curve in $L$. We can assume that $\gamma:S^1\to L$ is an
immersion. Set $I=[-1,1]$ and define the map $ \phi:S^1\times I\to
S=\phi(S^1\times I) $ by
$$
\phi(\theta,t)=X_t(\gamma(\theta)),
$$
where $X_t(\gamma(\theta))$  stands for the integral curve of $X$
with initial point at $\gamma(\theta)$.  It is clear that $\phi$ is
an immersion of class $C^r, r \geq 2$. Then
$\omega^*=\phi^*(\omega)$ is a well-defined $1$-form on $S^1\times
I$. Because $d\omega^*=d\phi(\omega^*)= \phi^*(d\omega)=\phi^*(0)=0$
we have that $\omega^*$ is closed. Hence $\omega^*$ defines a
foliation $\fa^*$ on $S^1\times I$. Note that $\fa^*$ is conjugated
to $\fa\cap S$. It follows that the curves $\gamma^*=S^1\times 0$
and $\gamma$ have the same holonomy. Let us calculate the holonomy
of $\gamma^*$. Fix $(\theta^*,0)\in \gamma^*$ and
$\Sigma^*=\{\theta^*\}\times I$. Clearly $\Sigma^*$ is a transversal
of $\fa^*$. Let $f^*:\Dom(f^*)\subset \Sigma^*\to \Sigma^*$ be the
holonomy of $\gamma^*$, $p\in \Dom(f^*)$ and $q=f^*(p)$. Let
$\alpha$ be an arc in $\Sigma^*$ joining $p$ and $q$.

Let $l$ be a path in a leaf of $\fa^*$ joining $p,q$. Let $R$ be the
closed region bounded by the curves $\gamma^*$, $l$ and $\alpha^*$.
Because
$$
0=\int_Rd\omega^*= \int_{\po R}\omega^*=
\int_{l}\omega^*+\int_{\alpha}\omega^*= 0+\int_{\alpha}\omega^*
$$
one has
$$
\int_{\alpha}\omega^*=0.
$$
This equality implies that $\alpha$ is trivial and so $p=q=f^*(p)$.
We conclude that $\gamma^*$ has trivial holonomy. Hence $\gamma$ has
trivial holonomy and the proof follows.
\end{proof}

}
\end{Example}


\chapter{Some results from several complex variables}

In the course of the text we shall refer to some results from the
theory of several complex variables.  For the sake of clarity we
shall now state them separately.
\section[Some extension theorems]
{Some extension theorems from several complex variables}
\label{section:coupleextension}

This section is dedicated to some useful extension theorems from
several complex variables. We shall start with a statement due to
Riemman. We recall that a subset \, $X \subset M$ of a topological
space is {\it nowhere dense\/} if the closure $\overline{X}$ has
empty interior.

\begin{Theorem}[Riemann extension theorem,
\cite{Gunning 1}\index{Theorem! Riemann extension}]
\label{Theorem:riemannextension} Let $M$ be a connected complex
manifold and $X \subset M$  an analytic nowhere dense subset of $M$.
Then a holomorphic function $f$ is $M\setminus X$ which is bounded
near $X$ has a unique holomorphic extension to $M$.
\end{Theorem}

We shall now state local and simple versions of two powerful results
from the theory of functions and analytic sets in several complex
variables:

\begin{Theorem}[Hartogs' extension theorem, \cite{Gunning 1}\index{Theorem! Hartogs' extension}]
\label{Theorem:Hartogs} Let $U \subset \bc^n$ be a connected open
subset and $W\subset U$ be a codimension $\ge2$ analytic subset.
Then any holomorphic (respectively, meromorphic) function $f$
defined on the open subset $U\setminus W$ admits an unique
holomorphic (respectively, meromorphic) extension to $U$.
\end{Theorem}

\begin{Remark}{\rm  The same holds for (as an immediate consequence) vector fields and differential forms.
We recall that an {\it analytic subset\/} of a complex space is one given locally by set of common zeros of (local) holomorphic functions.
}\end{Remark} Remmert-Stein theorem gives conditions for the closure
of an analytic set to be analytic.

\begin{Theorem} [Theorem of Remmert-Stein, \cite{Gunning 2,gunning-rossi}\index{Theorem! Remmert-Stein}\index{Theorem! of Remmert-Stein}]
\label{Theorem:Remmert-Stein}
{\it Let $M$ be a complex manifold, $W \subset M$ an irreducible analytic subset of $M$ and $V$ an
irreducible analytic subset of $M\setminus W$; such that $\dim(V) > \dim(W)$. Then the closure $\overline{V} \subset M$
is an analytic subset of dimension $\dim(V)$}.
\end{Theorem}

\noindent Also we need:

\begin{Theorem}[Theorem of Chow, \cite{gunning-rossi}\index{Theorem! of Chow}]
\label{Theorem:Chow}{\it A {\rm(}closed{\rm)} analytic subvariety on a complex projective space is algebraic}.
\end{Theorem}

\noindent And the following immediate consequence:

\begin{Corollary} An irreducible closed analytic subset of pure codimension one of the complex
projective space $\bc P(n)$ is an algebraic hypersurface.
\end{Corollary}

\section{Levi's global extension theorem}
\label{section:leviextension}

Recall that given an open subset $U \subset \bc^n$, $n \ge 2$ and a
map $f\colon U \to \re$ of class $C^2$ we say that $f$ is {\it
plurisubmarhomic\/} ({\it plush\/} for short) or {\it strictly
plurisubmarhomic\/} ({\it s-plush\/} for short) if $\forall\, p_i
\in U$ and $\forall\, v \in \bc^n$ the restriction $u\colon z
\mapsto f(p_0+zv)$ is subharmonic or strictly subharmonic, in the
sense that the Laplacian satisfies
$\Delta u \ge 0$ or $\Delta u > 0$ respectively.

\noindent If we denote by $H_f(p_0)$ the complex $n \times n$ matrix
$$
H_f(p_0) = \left(\frac{\po^2f}{\po z_i \po
\overline{z_j}}\,(p_0)\right)_{i,j=1,\dots,n}
$$
then $H_f(p_0)$ is hermitian $\left(\overline{H_f(p_0)} =
H_f(p_0)^t\right)$. It well-known that:

\begin{Lemma} Let $f\colon U \subset \mathbb C^n \to \mathbb R$ be a $C^2$ map in the open subset $U$. Then:

\begin{enumerate}
\item $f$ is plush in $U \, \, \Leftrightarrow H_{f_{\forall\,p\in U}}$  is definite non negative.

\item $f$  is s-plush in $U\, \,  \Leftrightarrow
H_{f_{\forall\,p\in U}}$  is definite non positive.
\end{enumerate}
\end{Lemma}

As a consequence:
\begin{Lemma}
Given  $f \in C^2(U)$, $U \subset \bc^n, \, n\ge2$ as above, then we
have:
\begin{enumerate}
\item $f$ is plush in $U \, \, \Leftrightarrow$ for each holomorphic curve $\ga\colon V \subset \bc \to U$
the map $f\circ\ga\colon V \to \re$ is subharmonic.

\item $f$  is s-plush in $U\, \,  \Leftrightarrow
f\circ\ga\colon V \to \re$ is strictly subharmonic for every
holomorphic immersed curve $\ga\colon V \subset \bc \to U$.
\end{enumerate}
\end{Lemma}

Let us see how to extend this to complex manifolds. Given $f \in
C^2(U)$ as above we define the {\it Levi form\/}\index{Levi form} of
$f$ as the quadratic form
$$
L_f(p) := \sum_{j=1}^n \sum_{i=1}^n \frac{\po^2f}{\po z_i \po
\overline{z_j}}\,(p)\, dz_i d\overline{z_j}, \, p \in U.
$$
Thus we have a quadratic form at $T_p\,\bc^n$ defined by
$$
L_f(p)\cdot \omega =  \frac{\po^2f}{\po z_i \po
\overline{z_j}}\,(p)\cdot w_i\overline{w_j} = w\cdot H_f(b)\cdot
\overline{w^t}, \,\, \forall\, \omega =  (w_1,\dots,w_n) \in \bc^n.
$$
Using this form we can state (cf. \cite{Gunning 2},\cite{Range}):

\begin{enumerate}
\item $f$ is plush in $U \, \, \Leftrightarrow L_f(p) \ge 0, \forall\, p \in U.$

\item $f$  is s-plush in $U\, \,  \Leftrightarrow
L_f(p) > 0, \forall\, p \in U. $
\end{enumerate}

Given a holomorphic map $\phi\colon V \to U$, $V \subset \bc^m$ open
subset, we can use the Taylor expansion of order two for $f$ in
order to prove that $L_{f\circ\phi} = \phi^*(L_f)$ \,\, for the
Levi-forms of $f$ and of $f \circ \phi$ (\cite{Gunning
2},\cite{Range}).

\noindent Given now $f\colon M \to \re$ of class $C^2$, \, $M^n$ a
complex manifold and a chart $\phi\colon U \to \bc^n$ of $M$
we can consider the {\it Levi form\/}\index{Levi form} of $f$ at $p \in
M \subset U$ as the quadratic form on the (complex) tangent space
$T_p(M)$ defined by
$$
L_f(p)\cdot v := L_{f\circ\phi^{-1}} (\phi(p))\cdot(D\phi(p)\cdot
v), \,\, \forall\, v \in T_p(M).
$$
This is well defined according to the above remarks. Finally we
reach the following definition:

\begin{Definition}[\cite{[Ho],Range,Si}] {\rm Given $f\colon M \overset{C^2}{\longrightarrow} \re$,\, $M$
complex manifold, we say that $f$ is {\it
plurisubharmonic}\index{function! plurisubharmonic} ({\it plush\/}
for short),  respectively {\it strictly
plurisubharmonic}\index{function! strictly plurisubharmonic}
($s$-{\it plush} for short) if $L_f(p) \ge 0$, $\forall\, p \in M$
(respectively, $L_f(p) > 0$, $\forall\, p \in M$). Given $1 \le k
\le n$ we shall say that $f$ is $k$-({\it strictly
plurisubharmonic}) if $\forall\, p \in M,\, \exists$ a subspace $E
\subset T_p(M)$ of (complex) dimension $k$ such that
$L_f(p)\big\vert_E > 0$, i.e., $\forall\, v \in E\setminus\{0\}$
\,\,\ $L_f(p)\cdot v > 0$.

\vglue.1in Let $M$ be a differentiable manifold.  An {\it
exhaustion\/} of $M$ is a continuous function $g\colon M \to \re$
such that:
\begin{itemize}
\item[{\rm(a)}] $g$ is bounded from below $g \ge c$ in $M$
\item[{\rm(b)}] $g$ is proper: $\forall$\, sequence $\{p_n\} \subset M$ with no accumulation point in $M$,
the sequence $g(p_n)$ satisfies $g(p_n) \to +\infty$ as $n \to
+\infty$.
\end{itemize}

\vglue.1in A {\it Stein manifold\/}\index{Stein manifold} is a
complex manifold admitting a $C^\infty$ $s$-plush exhaustion
(Hörmander's theorem\cite{[Ho]}). }
\end{Definition}

\noindent On the affine space $\bc^n$ the function \,\, $f(z) =
||z||^2$ is a s-plush exhaustion. The important result below is due
to Levi:

\begin{Theorem} [Levi's global  extension theorem, \cite{Si}\index{Theorem! Levi global extension}]
\label{Theorem:Levi}

{\it Let $M$ be a complex manifold admitting a $k$-s-plush
exhaustion where $k\ge2$. If $K \subset M$ is compact subset such
that $M\setminus K$ is connected then any meromorphic $q$-form $\om$
on $M\setminus K$ admits an unique extension as a meromorphic
$q$-form on $M$.}
\end{Theorem}
We also need the following result:

\begin{Proposition} Let $X \subset \bc P(n)$
be an algebraic subset defined by $k$ homogeneous polynomials in
$\bc^{n+1}$. Then the open manifold $M = \bc P(n)\setminus X$ admits
a $\ell$-s-plush exhaustion where $\ell = n-k+1$. In particular if
$X$ is an algebraic hypersurface (codimension one) the $M = \bc
P(n)\setminus X$ is a Stein manifold.
\end{Proposition}

\begin{proof} If $X$ is defined by the homogeneous polynomials $f_1,\dots,f_k$ in $\bc^{n+1}$ then
we define $f\colon M \to \re$ by setting in homogeneous coordinates
$(z_1,\dots,z_{n+1}) \in \bc^{n+1}$
$$
f(z_1,\dots,z_{n+1}) := \ln \bigg(\frac{(\sum\limits_{j=1}^{n+1}
|z_j|^2)^q}{\sum\limits_{j=1}^k |f_j^{q_j}\,(z)|^2}\bigg)
$$
where if $d_j = \deg(f_j)$ then $q_1,\dots,q_k \in \bn$ are such
that $d_1q_1 =\cdots= d_kq_k = q \in \bn$.

\noindent Computing the Levi form of $f$ we can conclude by using
the following:

\begin{Lemma}\quad $\forall\, x \in \bc^n$ the quadratic form
$$
Q:=\sum_{j=1}^n |dx_j|^2 + \sum_{i < j} |x_idx_j-x_jdx_i|^2
$$
is definite positive.
\end{Lemma}
\end{proof}

\chapter{Holomorphic foliations: nonsingular case}

\section{Basic concepts}
The purpose of this section is to  introduce, in a formal way,
the concept of holomorphic foliation. In fact, a holomorphic foliation
is, in particular,  a foliation in the  classical sense.
Nevertheless, in this section we shall see some important examples
for the development of the theory, that illustrate the concept.
\vglue.1in

\begin{Definition}
\label{Definition: 1.I}  {\rm Let $ M $ be a complex manifold of
(complex) dimension $ n $. A \index{foliation! holomorphic} {\it
holomorphic foliation of $M$, of dimension} $ k $, or {\it
codimension} $n-k, \ 1 \leq k \leq n-1 $, is a decomposition $ \fa $
of $ M $ in pairwise disjoint immersed complex submanifolds (called
\index{leaves} {\it leaves}\index{holomorphic foliation! leaf} of
the foliation $ \fa $) of  dimension (complex) $ k $,  and having
the following properties:
\vglue.05in
\noindent {\rm(i)} $\forall p\in M$ there exists a unique
submanifold $L_p$ of the decomposition that passes by $p$ (called
the\index{leaf through $p$} {\it leaf through $p$}). \vglue.05in
\noindent {\rm(ii)} $\forall p\in M$, there exists a holomorphic chart of $M$
(called \index{distinguished chart}{\it distinguished chart of $\fa$}),
$(\vr,U), \, p\in U, \, \vr\colon U \to \vr(U)\subset \mathbb{C}^n$,
such that $\vr(U)= P\times Q$, where $P$ and $Q$ are open polydiscs
em $\mathbb{C}^k$ and $\mathbb{C} ^{n-k}$ respectively. \vglue.05in
\noindent {\rm(iii)}  If $L$ is a leaf
of $\fa$ such that $L\cap U\ne \emptyset$,
then $L\cap U= \bigcup\limits_{q\in D_{L,U}}
\vr ^{-1}(P\times \{q\})$, where $D_{L,U}$ is a countable subset of
$Q$. \vglue.05in The subsets of $U$ of the form $\vr ^{-1}(P\times
\{q\})$ are called of \index{plaques}{\it plaques} of the
distinguished chart $(\vr,U)$. \vglue.05in A foliation of dimension
one is also called  {\it foliation by curves}\index{foliation! by
curves}.
In this case, the
leaves are  Riemann surfaces. \vglue.05in Observe that {\rm(iii)}
also implies that the leaves are immersed submanifolds immersed in
$M$. Indeed, the intersection of a leaf with a distinguished chart
is a union of plaques pairwise disjoint. Later on we shall see
examples of foliations exhibiting leaves which are (immersed but)
not embedded submanifolds. }
\end{Definition}



\vglue.1in
\begin{Remark}\label{Remark:1.I}{\rm A dimension $k$  foliation $\fa$ of $M$,
induces on $M$  a  \index{distribution}{\it distribution} of planes
of dimension  $k$, denoted by $T \fa$, which is defined by $ T_p \fa
= T_{p} (L_{p})$, the tangent plane at $p$ of the leaf $L_p$ passing
through  $p$. From  {\rm(iii)}, this distribution is holomorphic. It
defines a holomorphic vector sub-bundle of the tangent bundle  $TM$,
which will also be denoted by  $T \fa$. \vglue.1in}
\end{Remark}
The most simple example of holomorphic foliation of dimension $k$ is
the following: \vglue.05in
\begin{Example}
\label{Exemplo:1.I} {\rm Given the affine space  $\mathbb{C}^n$ we may consider any decomposition  $\mathbb{C}^n = \mathbb{C}^k\times
\mathbb{C}^{n-k}$. Such a decomposition defines a foliation  $\fa$
of dimension $k$ in $\mathbb{C}^n$, whose leaves are the affine
subspaces $\mathbb{C}^k\times \{q\}, \, q\in \mathbb{C} ^{n-k}$.}
\end{Example}
\vglue.1in
Next we will see two ways to define foliations, equivalent to the
above, and which will be further used throughout the text.
\vglue.05in

\begin{Proposition}
\label{Proposition: 1.I} A holomorphic foliation $ \fa $ of
dimension
$ k $ on a complex manifold  $ M
$ can also be defined in the following equivalent ways:
\vglue.05in
\noindent {\rm(I)}  {\it Description by charts
distinguished:}\index{holomorphic foliation! description by charts}
\vglue.05in
$\fa$ is given by an atlas of $M$ {\rm(}also denoted
$\fa${\rm)}, $\{(\vr_\alpha, U_\alpha) / \alpha \in A\}$ where:
\vglue.05in
\noindent {\rm(I.1)} $\vr_\alpha(U_\alpha)= P_\alpha \times Q_\alpha$, where
$P_\alpha, Q_\alpha$ are polydiscs of dimensions $k$ and $n-k$
respectively. \vglue.05in
\noindent {\rm(I.2)} If $U_\alpha \cap U_\beta \ne \emptyset$ then
the change of charts $\vr_\beta \circ \vr_\alpha ^{-1}$
is locally of the form
$$\vr_\beta \circ \vr_\alpha ^{-1}(x_\alpha,y_\alpha) =
(h_{\alpha \beta} (x_a,y_\alpha),g_{\alpha \beta} (y_\alpha))$$
\vglue.05in In this case the plaques of $\fa$ in $U_{\alpha}$ are
the sets of the form $\vr_{\alpha}^{-1}(P_\alpha\times \{q\})$.
\vglue.05in
\noindent {\rm(II)}  {\it   Description by local submersions\index{holomorphic foliation! description by local submersions}:} \vglue.05in
$\fa$ is given by an open  cover
$M=\bigcup\limits_{\alpha \in A} U_\alpha$  by
collections $\{y_\alpha\}_{\alpha \in A}$ and $\{g_{\alpha \beta}
\}_{U_{\alpha} \cap U_{\beta} \ne \emptyset}$, that satisfy:
\vglue.05in
\noindent {\rm(II.1)} $\forall \alpha \in A$, $y_{\alpha}\colon
U_\alpha \to \mathbb{C}^{n-k}$ is a submersion.
\vglue.05in
\noindent {\rm(II.2)} If $U_\alpha \cap U_\beta \ne \emptyset$
then $y_\alpha = g_{\alpha \beta}(y_\beta)$ where $g_{\alpha \beta }
\colon y_\beta(U_\alpha\cap U_\beta)\subset \mathbb{C}^k \to
y_\alpha (U_\alpha\cap U_\beta)\subset \mathbb{C}^k$ is a
local holomorphic. \vglue.05in In this case the
plaques of $\fa$ in $U_{\alpha}$ are the sets of the form  $y_\alpha
^{-1}(q), \, q\in V_\alpha$. \vglue.05in
\end{Proposition}

\begin{proof} First we  show that (I) is equivalent to the
definition of foliation.  Let $ \fa $ be a dimension $k$
foliation  on $ M $. We shall  build a
holomorphic atlas
${\mathcal A} $ of  $ M $ satisfying  conditions (I.1) and (I.2)
above. From the original definition of foliation it follows that
there is a holomorphic atlas $ {\mathcal A} = \{(\vr_ \alpha, U_
\alpha); \alpha \in A \} $ of $ M $ such that all systems of
coordinates $ (\vr_ \alpha, U_ \alpha) $ of $ {\mathcal A} $ satisfy
{\rm (ii)} and {\rm (iii)},  and $ \vr_ \alpha (U_ \alpha) = P_
\alpha \times Q_ \alpha $, where $ P_ \alpha $ and $ Q_ \alpha $ are
polydiscs of dimension $ k $ and $ n - k$ respectively. Let us
consider the change of charts $\vr_\beta \circ \vr_\alpha ^{-1}
\colon \vr_\alpha (U_\alpha \cap U_\beta) \to \vr_\beta (U_\alpha
\cap U_\beta)$. We can write
$$\vr_\beta \circ \vr_\alpha ^{-1} (x_\alpha, y_\alpha) =
(h_{\alpha \beta} (x_a,y_\alpha),g_{\alpha \beta} (x_\alpha,
y_\alpha)) = (x_\beta,  y_\beta)$$ where $(x_a,y_\alpha) \in
P_\alpha \times Q_\alpha$. We claim that $g_{\alpha \beta}$ does not
depend   on $x_\alpha$. \vglue.05in
Indeed, a point $y_\alpha \in Q_\alpha$ defines the  plaque
$\vr_{\alpha}^{-1}(P_\alpha\times \{y_\alpha \})$ in $U_\alpha$,
which  is contained in a  leaf $L$ of $\fa$. On the other hand
$L \cap U_\beta$ consists of a  countable union of plaques of
$U_\beta$, of the form $\cup_i \vr_{\beta}^{-1}(P_\beta\times
\{y_{\beta}^{i} \})$. From this we obtain
$$\vr_\beta \circ \vr_\alpha ^{-1}
((P_\alpha \times \{ y_\alpha \})\cap U_\beta )\subset \vr_\beta (L
\cap U_\beta) = \cup_i ( P_\beta \times \{ y_{\beta}^{i} \}). $$
This last  implies the following
$$g_{\alpha \beta} (P_\alpha \times \{ y_\alpha \}) \subset \cup_i \{ y_{\beta}^{i} \}$$
\vglue.05in It is not difficult to see that the above implies that
$\partial g_{\alpha \beta} / \partial x_\alpha = 0$. This proves the
claim. \vglue.05in
Assume now that there exists a holomorphic atlas
${\mathcal F}$ of $M$ satisfying (I.1) and (I.2). Since $M$
is a manifold, we can assume that ${\mathcal F}$ is countable. Next
we shall define the  ``leaves of $\fa$", taking into account
(I.1)and (I.2). \vglue.05in

In $M$ we consider an equivalence relation
that identifies two points $p,q\in M$ if, and only if, there exists
a finite chain of plaques (as in (I)) $P_1,...,P_r$ of $\fa$ such
that $P_i\cap P_{i+1}\ne \emptyset, \forall i$, and $p\in P_1, q\in
P_r$. The leaves are the equivalence classes of $M$ by this
relation. Thus, two points $p,q\in M$ are in the same leaf if, and
only if, there exists a chain of plaques as above that contains
these points. Since the plaques are connected, it follows that the
leaves are connected. \vglue.05in In order to see that the leaves
are immersed submanifolds in $M$ it is necessary to endow each leaf
$L$
of $\fa$ with a structure
of a holomorphic variety, in such a way that  the inclusion map
$i\colon L \to M$ is an immersion. This structure, which
is called
\index{intrinsec structure}{\it intrinsec structure}, is defined as
follows: \vglue.05in Fix a leaf $L$ and consider the cover  of $L$
consisting of all the plaques contained in $L$. Given a plaque
$P_{\alpha}^{q} = \vr_{\alpha}^{-1}(P_\alpha\times \{q\}) \subset L$
we define the ``coordinate system"
$$\vr_{\alpha}^{q}= \pi_{1} \circ \vr \mid_{P_{\alpha}^{q}}
\colon P_{\alpha}^{q} \to \mathbb{C}^k$$ where $\pi_1 \colon
\mathbb{C}^k \times \mathbb{C}^{n-k} \to \mathbb{C}^k $ is the first
projection. Then we obtain an ``atlas",
$${\mathcal F}_L = \{(\vr_{\alpha}^{q},P_{\alpha}^{q})\ ;
\ P_{\alpha}^{q} \ \text{is a plaque contained in\, } L \}.$$
\vglue.05in In order to assure that ${\mathcal F}_L$ is a
holomorphic atlas of $L$ it is necessary to prove that the changes
of charts are biholomorphisms  between open subsets of
$\mathbb{C}^k$. This fact, left to reader as an exercise, follows
from  (I.2). We mention also without an explicit proof, that $L$,
with the above defined structure, is a Hausdorff space. \vglue.05in
From the definition of leaf given above, it follows that $L\cap
U_\alpha$ is the disjoint union of plaques of $U_\alpha$ of the form
$$(*)\ \ L\cap U_{\alpha}= \bigcup\limits_{q\in D_{L,\alpha}}
\vr ^{-1}(P_{\alpha}\times \{q\}),$$
where $ D_{L,\alpha}\subset \mathbb{C}^{n-k}$. Note that
each plaque of $L$ in $U_{\alpha}$
corresponds to a single  point in
$ D_{L,\alpha}$.
From this, from the
definition of leaf and from the fact that $\fa$ is countable, we
obtain that  $D_{L,\alpha}$ is countable. Therefore $L$ contains
only countably many  plaques. Therefore $L$ has countable basis of
open sets, therefore it is a manifold. Finally
observe that (*) shows that  $i \colon L \to M$ is an injective
immersion. \vglue.05in We leave the proof  that (II) is equivalent
to the  definition of foliation as an exercise to reader.\end{proof}

\section{Examples}
Now we explore some examples of holomorphic foliations. These
examples are illustrate  the basic constructions and also motivate
the concept of holomorphic foliation with singularities, to be
introduced later on. Special attention should be paid to the
dimension one and to the codimension one cases.

\begin{Example} \label{Example: 2.I} {\rm Let $ f \colon M \to N$
be a holomorphic submersion, where $ M $ and $ N $ are holomorphic
manifolds of dimensions $ n + k $ and $ k $ respectively. In this
case, by the local form of  holomorphic submersions, the level sets
$ \{F = c \}, c \in N $ are holomorphic submanifolds of codimension
$ k $ of $ M $. Definition (II) of Proposition~\ref{Proposition:
1.I} then ensures that there is a foliation on $ M $ whose leaves
are the connected components  level sets $ f$. We leave the proof of
this fact as exercise for the reader. } \end{Example}

\begin{Example} [Pull-back or inverse image of a foliation]
\label{Example: 3.I} {\rm  Let $ M $ and $ N $ complex manifolds,
$ f \colon M \to N $ a holomorphic map and $ \fa $ a foliation on $
N $ of codimension $ k $.}
\end{Example}

\begin{Definition} \label{Definition: 2.I} {\rm
We say that $f$ is \index{transverse map to a foliation} {\it
transverse} to $\fa $ if for every point $ q \in N $, the subspace $
df_q (T_q M) $ and $ T_ {p} \fa $ generate the tangent space $ T_
{p} N $, and $ p = f (q) $. \vglue.05in If this is the case then
there is a foliation of $M$, denoted by $ f ^ * (\fa) $, of the same
codimension $ k $, whose leaves are the connected components of the
inverse images by $ f $, $ f ^ {- 1} (L) $, of the leaves $ L $ of $
\fa $ in $ N $. The foliation $ f ^
* (\fa) $ is called \index{foliation! pull-back} {\it pull-back}
or {\it inverse image}\index{foliation! inverse image} of $ \fa $ by
$ F $.} \end{Definition}
The
foliation $f^*(\fa)$ is obtained using (II) in
Proposition~\ref{Proposition: 1.I}. Indeed, consider an open  cover
$\{U_\alpha \}_{\alpha \in A}$ of $N$ and collections
$\{y_\alpha \}_{\alpha \in A}$ and $\{g_{\alpha \beta} \}_{U_\alpha
\cap U_\beta \ne \emptyset}$ satisfying (II.1) and (II.2) in
Proposition~\ref{Proposition: 1.I}. Given $\alpha \in A$ let
$V_\alpha
= f^{-1}(U_\alpha)$ and $z_\alpha = y_\alpha \circ f \colon V_\alpha
\to \mathbb{C}^k$.
In this way we obtain  an open cover
$\{V_\alpha \}_{\alpha \in A}$ of
$M$ and a collection of
holomorphic maps $\{z_\alpha \}_{\alpha \in A}$. Note that $V_\alpha
\cap V_\beta = f^{-1}(U_\alpha \cap U_\beta)$,  so that $V_\alpha
\cap V_\beta \ne \emptyset$ if, and only if, $U_\alpha \cap U_\beta
\ne \emptyset$. Moreover,  if $V_\alpha \cap V_\beta \ne \emptyset$
then $z_\alpha = g_{\alpha \beta} \circ z_\beta$. Hence, in order to
check that that $f^*(\fa)$ is a foliation it is enough to prove that
$z_\alpha \colon V_\alpha \to \mathbb{C}^k$ is a submersion for
every $\alpha \in A$. This is  a consequence of the fact  that $f$
is transverse to $\fa$, as it can be verified  from the definition.
\vglue.05in

\begin{Example} [foliation  generated by a
holomorphic
vector field] \label{Example: 4.I} \vglue.05in {\rm Let $ M $ be a
complex manifold of dimension $ n $ and $ X $ a holomorphic vector
field not identically zero in $ M $. Let $ S = \{p \in M; X (p) = 0
\}$, the singular set of $ X $.  Then $ X $ generates a holomorphic
foliation $ \fa $ of  dimension $ 1 $ in the open $ N = M \setminus
S$. The leaves of $ \fa $ are the trajectories of $ X $ in $ N $.
The structure of foliation arises from the Flow box
theorem\index{Theorem! Flow box theorem} for holomorphic vector
fields, which can be stated as follows: \vglue.05in `` For every $ p
\in M$ such that $ X (p) \ne  0$,
there
is a holomorphic coordinate system $ (\phi = (z_1, ..., z_n), U) $,
where $ P \in U, \phi: U \to \phi (U) = A \times B \subset
\mathbb{C} \times \mathbb{C} ^ {n-1} $ and where $ X = \partial /
\partial z_1$. "\vglue.05in

Since the trajectories of $ X $ are the solutions of the
differential equation $ Dz / dt = X (z) $ and $ X \mid_U = \partial
/ \partial z_1 $, we conclude that the trajectories of $ X $ in $U$
are  the form $ \phi ^ {- 1} (A \times \lbrace w \rbrace)$ with $
W\in B $. We get from this
and from Definition (I) in
Proposition~\ref{Proposition: 1.I}, a dimension $1$
foliation, whose leaves are the trajectories
of $ X $. \vglue.05in In fact, every foliation of dimension one is
locally defined by vector fields, as we shall see in
the following
result, whose proof we leave to the reader as an exercise. }
\end{Example}

\vglue.05in

\begin{Proposition}
\label{Proposition: 2.I} Let $ M $ a complex manifold of   dimension
$ n \ge  2$ and $ \fa $ a dimension one foliation of $ M $. There
are collections    ${\mathcal X} = \{X_ \alpha \} _ {\alpha \in A}
$, $ {\mathcal U} = \{U_ \alpha \} _ {\alpha \in A} $ and $ \mathcal
G = \{g _ {\alpha \beta} \} _ {U_ \alpha \cap U_ \beta \ne
\emptyset} $ such that:

\noindent {\rm(i)} ${\mathcal U}$ is a cover of $M$ by open sets.

\noindent  {\rm(ii)} $X_\alpha$ is a holomorphic vector field
on $U_\alpha$, which is nonzero at each point.

\noindent {\rm(iii)} $g_{\alpha \beta}\in
\O^*(U_\alpha \cap U_\beta)$,
that is, is a holomorphic function that does not vanish
on $U_\alpha \cap U_\beta$.

\noindent {\rm(iv)} In $U_\alpha \cap U_\beta \ne \emptyset$
we have
$X_\alpha=g_{\alpha \beta}.X_\beta$.

\noindent {\rm(v)} If $p\in U_\alpha$,
then $T_p \fa = \mathbb{C}.X_\alpha(p)$,
the subspace of $T_p M$ generated by $X_\alpha(p)$.
\vglue.05in Conversely, given  collections
$\mathcal X$,
${\mathcal U}$ and $\mathcal G$ satisfying {\rm(i)}, {\rm(ii)},
{\rm(iii)} and {\rm(iv)},  there exists a foliation $\fa$ that
satisfies {\rm(v)}.
\end{Proposition}


The following example is a complex version of results found in
\S~\ref{subsection:Frobeniustheorem} (cf.
Example~\ref{Example:integrableformssystem}).

\begin{Example} [foliations generated by differential
$1$-forms\index{foliation! generated by  differential $1$-forms}]
\label{Example: 5.I} \vglue.05in {\rm Let $ M $ a complex manifold
of dimension $n$ and $ \om $ a holomorphic  $1$-form, not
identically zero in $ M $. Let $ S = \{p \in M; \Om_p \ne 0 \} $,
the singular set of $ \om $. In this case, $ \om $ induces a
distribution of hyperplanes $\Omega$ in the open set $ N = M
\setminus S $ defined by

$$
\Omega_{p}= \text{ker}(\om_{p})= \{ v \in T_{p}M \ ;\ \om_{p}(v)=0
\}
$$

\begin{Definition} \label{Definition: 3.I} {\rm We say that $ \om $ (or $ \Omega $) is
\index{form! integrable} {\it integrable}, if there is a holomorphic
foliation $ \fa $ in $ N $ such that $ T \fa = \Omega $. In other
words, the tangent space at $ p $ to the leaf of $ \fa $ passing
through $ P $, coincides with $ \Omega_p$. }
\end{Definition}

A well known fact is the following (see
\S~\ref{subsection:Frobeniustheorem} or \cite{Godbillon}):
\begin{Lemma}
The $1$-form $\om$  is integrable if, and only if $\om \wedge
d\om=0.$
\end{Lemma}

The above result is known as \index{Theorem! of
Frobenius} {\it Frobenius theorem}. It is often said that the
foliation $ \fa $ is defined by the differential equation $ \om = $
0 and that leaves of $ \fa $ are integral submanifolds of this
equation. \vglue.05in

It should be noted that if $ \eta $ is a $1$-form with $\eta = f \,
\om $, where $ f $ is a holomorphic function in N that does not
vanish, then the hyperplanes distribution induced by $ \eta $
coincides with  $\Omega $. In particular, $ \eta $ is also
integrable and the foliations defined by $ \eta = $ 0 and $\om = 0$
coincide.

}
\end{Example}

Codimension one foliations are locally defined by integrable
differential $1$-forms, as  it is stated below.
\vglue.1in

\begin{Proposition}
\label{Proposition:3.I} Let $M$ be a complex manifold of dimension
$n\ge 2$ and $\fa$ a codimension one foliation of $M$. There exist
collections $\mathcal W=\{\omega_\alpha \}_{\alpha\in A}$,
${\mathcal U}=\{U_\alpha\}_{\alpha\in A}$ and $\mathcal
G=\{g_{\alpha \beta} \}_{U_\alpha \cap U_\beta \ne \emptyset}$ such
that:

\noindent {\rm(i)} ${\mathcal U}$ is a cover
of $M$ by open sets.

\noindent  {\rm(ii)}
$\omega_\alpha$ is an integrable holomorphic  $1$-form
in $U_\alpha$, which does not vanish at any point.

\noindent {\rm(iii)} $g_{\alpha \beta}\in \O^*(U_\alpha \cap U_\beta)$.

\noindent {\rm(iv)} In $U_\alpha \cap U_\beta \ne \emptyset$ we have
$\omega_\alpha=g_{\alpha \beta}.\omega_\beta$.

\noindent {\rm(v)} If $p\in U_\alpha$, then $T_p \fa = \ker(\omega_\alpha(p))$.

\vglue.05in

Conversely, if there exist collections $\mathcal W$, ${\mathcal U}$
and $\mathcal G$ satisfying {\rm(i)}, {\rm(ii)}, {\rm(iii)} and
{\rm(iv)}, then there exists a foliation $\fa$ that
satisfies~{\rm(v)}.
\end{Proposition}

\vglue.05in The proof is similar to the proof of
Proposition~\ref{Proposition: 2.I} and is also left to the reader as
an exercise. \vglue.05in

\section{The Identity principle for holomorphic foliations}
A well-known property of holomorphic functions and differential
forms is the Identity principle  (cf.\cite{Gunning 1, Gunning 2,
gunning-rossi}). The next proposition is the statement of this same
principle for holomorphic foliations.
\begin{Proposition}  [identity principle
for holomorphic foliations] \label{Proposition:12.I} Let $M^n$ be a
connected complex manifold and let $\fa$, $\fa_1$, be two
holomorphic foliations of same  codimension $1\leq k \leq n-1$ in
$M$. Assume that  $\fa$ and $\fa_1$ coincide  in some nonempty open
subset  $U\subset M$. Then $\fa = \fa_1$ in $M$.
\end{Proposition}
\begin{proof}
As it can be easily checked, since $M^n$ is connected, it is enough
to prove the case where $M^n$ is an  open polydisc $M=\Delta^n$ in
$\mathbb C^n$. Moreover, for the same reason, we may assume that
$\fa$ and $\fa_1$ are given holomorphic submersions $f, f_1\colon
\Delta^n \to \mathbb C^{n-k}$. In this case, there is some
nontrivial open subset $U\subset \Delta^n$ where $f$ and $f_1$
satisfy: $df \wedge df_1 \equiv 0$ (use the local form of
submersions to one of these maps and then conclude it from the fact
that the other map is constant along the fibers of the first one).
Thanks then to the classical Identity Principle for holomorphic
maps, we conclude that $df \wedge df_1\equiv 0$ in $\Delta^n$. This
proves that $\fa$ and $\fa_1$ coincide in $\Delta^n$.
\end{proof}

\chapter{Holomorphic foliations with singularities}

In this chapter we introduce the concept of holomorphic foliation
with singularities, focusing  on the two cases: dimension one and
codimension one. Intuitively, as a first step, one may think of a
holomorphic foliation with singularities on a complex manifold $M$
as a
a pair $\fa := (\fa^*,S)$ where $S \subsetneqq M$ is an analytic subset and $\fa^*$ is a (classical)
holomorphic foliation on $M^* := M\setminus S$. We call $S$ the {\it singular set\/} of $\fa$, and
write $S = \sing(\fa)$. The {\it leaves\/} of $\fa$ are by definition the leaves of $\fa^*$ in $M^*$. All other concepts for
foliations (holonomy, plaques, etc) extend to this case by associating them to the non-singular foliation $\fa^*$.
Nevertheless, in the above mentioned cases we can say more and
make a more precise definition. This is done in the next sections.

\section{Linear vector fields on the plane}
\label{section:vectorfields} In this section we motivate the concept
of one-dimensional holomorphic foliation with the study of linear
plane vector fields. First of all, let us resume a  quite general
remark.

Let $M$ be a complex (always assumed to be connected) manifold of
complex dimension $m$. Let $X$ be a  holomorphic vector field on
$M$. We first assume that $X$ is {\it non-singular}, i.e., $X(p) \ne
0$, $\forall\,p \in M$. To $X$ we associate a holomorphic autonomous
ordinary  differential equation
$$
(*)\quad\begin{cases}
x &=\,\, X\big(x(t)\big), \,\,\, t \in \bc\\
x(0) &=\,\, x_0
\end{cases}
$$
The solutions of (*) define the {\it local flow\/} of $X$ and give a
(complex) dimension one holomorphic  foliation $\fa(X)$ of $M$,
obtained by gluing these local solutions

The leaves of this foliation $\fa(X)$ are the orbits/trajectories of
$X$ and they are immersed Riemann surfaces on $M$.

\begin{Example}{\rm  We consider a  linear diagonal vector field $X(x,y) = \la x\,\dfrac{\po}{\po x} +
\mu y\, \dfrac{\po}{\po y}, \, \mu,\la \in \bc\setminus\{0\}$ in the
complex plane $\mathbb C^2$. Then $X$ is holomorphic and
non-singular in the punctured plane $M = \bc^2\setminus\{(0,0)\}$.
From the above  considerations, to $X$ we can associate a dimension
one holomorphic foliation $\fa(X)$ of $M$. The leaves of $\fa(X)$ in
$M$ are parameterized by the (global) flow of $X\colon \bc \ni t
\mapsto \big(xe^{\la t}, ye^{\mu t}\big)$. The dynamical/geometrical
behavior of these leaves very much depends on the quotient $\la/\mu$
(whether it is real, rational positive,...).

\begin{Remark}{\rm  Even though the above example seems nice, actually we are not (usually in theory of foliations)
interested in the parametrization of the orbits; but on their
dynamical/geometrical behavior on the ambient manifold.
}\end{Remark} The point of view we want to propose is illustrated by
the following description of $\fa(X)$:

\noindent $\bullet$\,\, $\la/\mu = n/m$,\, $n/m \in \bq_+
\Rightarrow \fa(X)$ has a meromorphic first integral $f(x,y) =
y^n/x^m$.  The leaves are contained in the curves $ax^m + by^n =
0$,\, $a,b \in \bc$ and all of them accumulate only at the origin $0
\in \bc^2$.

\noindent $\bullet$\,\, $\la/\mu = -n/m$,\, $-n/m \in \bq_-
\Rightarrow \fa(X)$ has a holomorphic first integral $f = x^my^n$.

The leaves are closed in $\bc^2$ except for (those) the two axes;
the leaves are contained in the curves $s^my^n = c \ni \bc$.

We have a {\em real saddle-type behavior}.

\noindent $\bullet$ $\la/\mu \in \mathbb R_- \setminus \mathbb Q$:
the first integral $f(x,y)= y ^\lambda x^{\-mu}$ shows  that the
leaves are not closed, except for the axes. The closure of a typical
leaf is a $3$-dimensional real manifold given by $|f(x,y)|= |c| >0$,
i.e., by $|y|^\lambda |x|^{-\mu} = |c|$ (we may assume that
$\lambda, -\mu \in \mathbb R_+$ still with $\lambda / \mu \in
\mathbb R \setminus \mathbb Q$).

\noindent $\bullet$ $\lambda / \mu \in \mathbb C \setminus \mathbb
R$ (hyperbolic case). In this case again the leaves are not closed,
except for the axes, and all leaves accumulate on both axes. Indeed,
the leaves are closed outside of the the two axes (meaning that a
leaf accumulates only at the two axes).

\noindent $\bullet$ $\lambda / \mu \in \mathbb R_+$: exercise (think
of the radial vector field)!

}
\end{Example}

\section{One-dimensional foliations  with isolated singularities}
\label{section:foliationssingularities}

From what we have seen above, a holomorphic vector field $X$ with
singular set  $\sing(X) \subset M$ on a complex manifold $M$,
defines a dimension one holomorphic foliation $\fa(X)^*$ on the open
set $M^* = M\setminus\sing(X)$ (remark: $\sing(X) \subset M$ is an
analytic subset, assumed to be proper, so $M\setminus \sing(X)$ is
open and connected). The leaves of $\fa(X)^*$ are the orbits of $X$
on $M^*$, i.e., the non-singular orbits of $X$. A natural question
is whether any one-dimensional holomorphic foliation is induced by a
holomorphic vector field. The answer is {\em no}, indeed there are
manifolds which may be equipped with one-dimensional foliations, but
that do not admit non-trivial holomorphic vector fields. We shall
see this in a while. Let us now introduce one of the main concepts
in this text:

In what follows, as usual, $M$ is a connected complex manifold.

\begin{Definition}
\label{Definition:singularfoliation}{\rm A {\it holomorphic
foliation of dimension one of $M$ with
singularities\/}\index{holomorphic foliation! of dimension one with
singularities} is a pair $\fa := (\fa^*,S)$ where $S \subset M$ is a
proper analytic subset and $\fa^*$ is a (classical) one-dimensional
holomorphic foliation on $M^* := M\setminus S$. We call $S$  the
{\it singular set\/} of $\fa$, and write $S = \sing(\fa)$. The {\it
leaves\/} of $\fa$ are by definition the leaves of $\fa^*$ on $M^*$.
}
\end{Definition}

\begin{Remark}{\rm  Regarding the above concept:
\begin{enumerate}[{\rm(i)}]

\item We may assume that $S$ is minimal with the property ``there is a non-singular foliation  on
$M\setminus S$ that agrees with $\fa^*$ in some open set" (cf.
Propositions~\ref{Proposition:extensionfoliation} and
~\ref{Proposition:saturationfoliation}).

\item  Thanks to (i) above, in these notes we shall assume that
$\fa$ {\it has isolated singularities\/}\index{foliation! of
dimension one with isolated singularities},  {\it i.e.}, if $\dim
\fa = 1$ then  $\sing(\fa)$ is assumed to be a discrete subset of
$M$.
\end{enumerate}
}
\end{Remark}

The importance of the above concept also relies on the following
description of isolated singularities.

\begin{Lemma} Let $\fa$ be a one-dimensional holomorphic foliation with {\rm(}isolated{\rm)}
singularities on a complex manifold $M$. Then there is an open cover
$M = \bigcup\limits_{j\in J} U_j$  by connected subsets $U_j \subset
M$ such that:

\begin{enumerate}[{\rm(i)}]
\item $\sing(\fa) \cap U_j = \phi$ or $\sing(\fa) \cap U_j = \{p_j\}$;

\item On each $U_j$\,, $\fa$ is given (its plaques) by a holomorphic vector field $X_j$\,, with $\sing(X_j) = \sing(\fa) \cap U_j$\,.

\item  If $U_i \cap U_j \ne \phi$ then on $U_i \cap U_j$ we have $X_i = g_{ij}\cdot X_j$
for some non-vanishing holomorphic function $g_{ij}\colon U_i \cap
U_k \to \bc^* = \bc\setminus\{0\}$.  Conversely any such data
$\big(U_j, X_j, g_{ij}\big)$ defines a holomorphic, one-dimensional
foliation with isolated singularities on a manifold $M$.
\end{enumerate}

\end{Lemma}

\begin{proof} The proof is a consequence of the local case, i.e., we may assume that $M$ is a connected neighborhood of the origin $0 \in \bc^m$ and that $\sing(\fa)=\{0\}$. For sake of simplicity we assume $m=2$. Given now any point $p \in M\setminus\{0\}$ we consider the tangent space $T_p(\fa) := T_p(L_p)$ of the leaf $L_p$ of $\fa$ through $p$.

\noindent Regarded as a line through the origin $T_p(L_p)$
identifies with an element of the projective space
$\bc^2\setminus\{0\}\big/\bc^* = \bc P(1)$ is the Riemann sphere
$\overline{\bc} = \bc \cup \{\infty\}$, so that we have defined a
map $f\colon M\setminus\{0\} \to \overline{\bc}$. This map is
holomorphic (as a map into $\overline{\bc}$) and extends (by
Hartogs' extension theorem) to a meromorphic map $f\colon M \to
\overline{\bc}$. Indeed, the  local theory of analytic functions
says that $f(x,y)$ writes as a quotient $f(x,y) =
\dfrac{A(x,y)}{B(x,y)}$\,, where $A$, $B$ are holomorphic functions
on $M$ (for $M$ a small bidisc for instance).

Now we introduce the holomorphic vector field $X$ on $M$ by $X(x,y)
= B(x,y)\,\dfrac{\po}{\po x} + A(x,y)\,\dfrac{\po}{\po y}\,\cdot$
Then a local integral curve $\big(x(t),y(t)\big)$, $t \in D \subset
\bc$ of $X$ satisfies

$$\begin{cases}\dot x(t) &=\, B\big(x(t),y(t)\big)\\ \dot y(t) &=\, A\big(x(t),y(t)\big)\end{cases}$$

so that $\dfrac{y'(t)}{x'(t)} =
\dfrac{A\big(x(t),y(t)\big)}{B\big(x(t),y(t)\big)}$\,, that is
$\dfrac{dy}{dx} = \dfrac AB$  along the integral curves of $X$. This
shows that such integral curves are tangent to the leaves of $\fa$
and by same dimension we conclude that $\fa = \fa(X)$ in
$M\setminus\{(0,0)\}$.\end{proof}

\noindent Summarizing:
\begin{Proposition}
\label{Proposition:existsvectorfield}Every isolated singularity of a holomorphic foliation
of dimension one is defined by a holomorphic vector field.
\end{Proposition}

\noindent Thus, in a certain sense, {\sl the study of singular
points of holomorphic  foliations of dimension one, is the study of
singularities of holomorphic vector fields.}

The next example may require some further knowledge in theory of
singularities of projective varieties.
\begin{Example}[implicit complex ordinary differential equations]
\label{Example:implicit} {\rm  An {\it algebraic implicit ordinary
differential equation}
in $n\ge 2$ complex variables is given by expressions:

$$
\, \, \, \, \, \,  (**) \, \, f_j(x_1,...,x_n,x_j ^\prime)=0
$$

\noindent where $f_j(x_1,...,x_n,y)\in \co [x_1,...,x_n,y]$ are
polynomials $j=2,...,n$  and the $(x_1,...,x_n)\in \co^n$ are affine
coordinates.
Clearly, any polynomial vector field $X$ on $\co ^n$ defines
such an equation. In general $(**)$ defines a
one-dimensional singular foliation of some algebraic
variety of dimension $n$. For this we begin
by defining $F_j(x_1,...,x_n,y_2,...,y_n):=
f_j(x_1,...,x_n,y_j)\in \co[x_1,...,x_n,y_2,...,y_n]$ polynomials
in $n+(n-1)= 2n -1$ variables. Put also
$S_j := \{(x,y)\in \co ^n \times \co ^{n-1};
F_j(x,y)=o\} \simeq \{(x_1,...,x_n,y_j)\in
\co_x ^n \times \co _{y_j}; \, \, f_j(x_1,...,x_n,y_j) =0\}\times
\co ^{n-2} =: \Lambda_j \times \co ^{n-2} _{(y_2,...,\hat{y_j},...,y_n)}$.

We consider the projectivizations $\ov {S_j} \subset
\co P(2n-1)$ and the complete intersection subvariety
${S}:= \ov {S_2} \cap ... \cap  \ov{S_n}\subset
\co P(2n-1)$. Given by the differential forms
$\omega_j := y_j dx_1 - dx_j$ \, \,  $(j=2,...,n)$ on
$\co ^{n}\times \co^{n-1}$. Then
$\{ \omega_j=0, \, j=2,...,n\}$ defines an integrable system
on ${S}$. We say that the implicit differential equation $(*)$ is {\it
normal} if $S$ admits a normalization\index{normalization of a
curve} (desingularization\index{desingularization of a curve}) by
blow-ups $\sigma \colon \hat S \to S$. In particular we obtain in
general a singular foliation $\fa(**)$ of dimension one on the
algebraic $n$-dimensional subvariety $S \subset \co P(2n-1)$. Denote
by $f_1 \colon S\cap \co ^{n} \to \co ^1$ the projection in the
first coordinate $f_1(x_1,...,x_n,y_2,...,y_n)=x_1$, and extend it
to a holomorphic proper mapping $f_1 \colon \ov {S} \to \co P(1)$.
Assume now  that $S$ admits a normalization $\sigma \colon \hat
{{S}} \to {S}$. It is then possible to show that the foliation
$\fa(**)$ lifts to a foliation by curves $\hat{\fa}(**)$ on
$\hat{S}$ and $\hat f _1 = f_1 \circ \sigma$ defines a holomorphic
proper mapping from $\hat{\ov{S}}$ over $\co P(1)$. Finally,  using
Stein factorization theorem (\cite{G-R})\index{Theorem! Stein
Factorization theorem}  we can find a splitting $\hat f_1 \colon
\hat S \overset{\hat f}\to B \overset{\alpha} \to \co P(1)$ where
$\alpha \colon B \to \co P(1)$ is a finite ramified covering and
$\hat f \colon \hat S \to B$ is an extended holomorphic fibration
over the compact Riemann surface $B$ such that the following diagram
therefore commutes

$$
\begin{matrix}
\hat {S} \overset{\sigma}\longrightarrow \ov{S} \\
\, \, \hat f \downarrow  \qquad    \downarrow f_1 \\
\qquad  B \overset{\alpha}\longrightarrow \co P(1)
\end{matrix}
$$

\noindent for a map $\hat{f_1}\colon \hat S \to \co P(1)$.

}

\end{Example}

\section{Differential forms and vector fields}
\label{section:formsandvectorfields}

When in dimension two, the notions of codimension one and dimension
one, for a foliation, coincide.
In this case, the vector field viewpoint though quite valid is not the only. Indeed the differential
forms approach turns out to be quite useful when we look for integral of the foliation.
This is what we show through some simple examples in what follows.

In complex dimension two, to any holomorphic vector field $X(x,y) =
P(x,y)\,\dfrac{\po}{\po x} + Q(x,y)\,\dfrac{\po}{\po y}$ we can
associate the corresponding {\it dual} $1$-form
$$
\omega(x,y) = P(x,y)dy - Q(x,y)dx.
$$
The autonomous differential equation $\dot x = X(x)$ then is
equivalent to the {\it Pfaffian equation} $\om=0. $ \, The solutions
are the same but sometimes the Pfaffian viewpoint is more suitable.
Let us show this in a couple of examples:

\begin{Example}{\rm  Again we consider a linear diagonal vector field
$X_\la = x\,\dfrac{\po}{\po x} + \la y\,\dfrac{\po}{\po y}$\,, \,\, $\la \in \bc^*$,
the dual $1$-form is
$$
\omega =  x\,dy - \la y\,dx
$$
and Pfaffian equation is $x\,dy - \la\,dx=0$. Dividing by $xy$  we
obtain $\dfrac{dy}{y} - \la\,\dfrac{dx}{x} = 0 \Rightarrow d(\ell
n\,y - \la \ell n\,x) = 0 \Rightarrow d \ell n\,yx^{-\la} = 0$. This
shows that $f = yx^{-\la}$ is a (multivalued) first integral for
$\fa$.}
\end{Example}

\begin{Example} [Poincaré-Dulac normal form]
\label{Example:poincaredulac} {\rm For $ n \in \bn\setminus\{1\}, \,
a \in \bc^*$ we consider the vector field $X(x,y) =
(nx+ay^n)\,\dfrac{\po}{\po x} + y\,\dfrac{\po}{\po y}\,\cdot$ This
is called {\it Poincaré-Dulac normal form\/}\index{normal form!
Poincaré-Dulac} as we shall see later on. It is possible to see that
the only invariant curve through the origin is the axis $\{y=0\}$.

\noindent The dual $1$-form is
$$
\omega =  nx + ay^n)dy - y\,dx.
$$
Dividing by $y^{n+1}$ we obtain
\begin{align*}
\frac{\om}{y^{n+1}} &= \frac{(nx+ay^n)dy-y\,dx}{y^{n+1}} = \frac{mx\,dy}{y^{n+1}} - \frac{dx}{y^n} + \frac{a\,dy}{y}\\
&= -d\left(\frac{x}{y^n}\right) + d(\ell n\,y^a) =
d\left(-\frac{x}{y^n} + \ell n\,y^a\right).
\end{align*}
This gives us the first integral
$$
f = \exp\left(-\frac{x}{y^n} + \ell n\,y^a\right) = y^a \cdot
e^{-x/y^n}.
$$
As we can easily see we may change coordinates so that $a=1$ and
$f(x,y) = y\cdot e^{-x/y^n}$.  Since $f$ has a (line of) essential
singularity at $(y=0)$ we can conclude that {\it all leaves of,
except for the one contained in $(y=0)$, accumulate at $(y=0)$ and
are closed off $(y=0)$}.

}
\end{Example}

\section{Codimension one foliations with singularities}

As we have already seen
(Proposition~\ref{Proposition:existsvectorfield}),  any
one-dimensional holomorphic  foliation $\fa^*$, defined in a
punctured polydisc $\Delta\setminus\{0\}$, is given by a holomorphic
vector field $X$ defined on the polydisc $\Delta$, and with
$\sing(X)\subset \{0\}$. Moreover, by definition, $\fa^*$ extends to
$\Delta$ as a holomorphic foliation of dimension one, such that
$\sing(\fa)=\sing(X)$. Similarly to this one can prove for the case
of codimension one:

\begin{Proposition}
\label{Proposition:oneform} Given a polydisc $\Delta^n \subset
\mathbb C^n$, centered at the origin $0 \in \mathbb C^n$, denote by
$H\subset \mathbb C^n$ any codimension two plane through the origin
and put $\Delta^*=\Delta\setminus H$. Let $\fa^*$ be a
codimension-one holomorphic foliation defined in $\Delta^*$. Then
there exists a holomorphic $1$-form $\omega$ defined in the whole
polydisc $\Delta$ with the following properties:
\begin{enumerate}
\item $\omega$ is integrable, {\it i.e.}, $\omega \wedge d \omega =0$.
\item $\sing(\omega) \subset H$.
\item The restriction $\omega\big|_{\Delta^*}$  defines the foliation $\fa^*$ in the sense of Proposition~\ref{Proposition:3.I}.
\end{enumerate}
\end{Proposition}
\begin{proof}
Indeed, by definition there is an open cover $\mathcal U=\{U_j, \, j
\in J\}$ of $\Delta^*x$, such that on each open subset $U_j$ the
foliation is given by a holomorphic integrable $1$-form $\omega_j$.
Moreover, on each nonempty intersection $U_i \cap U_j \ne\emptyset$
there is a nonvanishing holomorphic function $g_{ij}\colon U_i \cap
U_j \to \mathbb C\setminus \{0\}$, such that $\omega_i = g_{ij}
\omega_j$, in this intersection. Then, the data $\{U_j, g_{ij}\}$
give a Cousin multiplicative problem in $\Delta^*$ (\cite{Gunning
2},\cite{Range}). We may assume that $n \ge 3$, otherwise we are in
dimension two, where we may use vector fields and the dual
$1$-forms. Then, from Cartan's theorem (\cite{gunning-rossi,Range})
on the solution of the second Cousin problem (also called Cousin
multiplicative problem), there are holomorphic functions $g_j \colon
U_j \to \mathbb C \setminus \{0\}$, such that on each intersection
$g_{ij}= g_i/g_j$. Then $\frac{1}{g_i} \omega_i = \frac{1}{g_j}
\omega_j$ on each intersection $U_i \cap U_j \ne \emptyset$. This
defines a $1$-form  $\omega$ in $\Delta^*$ by $\omega\big|_{U_j}=
\frac{1}{g_j} \omega_j$. This $1$-form is holomorphic, integrable
(each $\omega_j$ is integrable) and defines $\fa$ on each $U_j$.
Thus $\omega$ defines $\fa$ on $\Delta^*$. By Hartogs' extension
theorem, $\omega$ admits a holomorphic extension to the origin.

\end{proof}

The above remark also motivates the following definition:

\begin{Definition}[codimension one holomorphic foliation with singularities]
{\rm Let $M$ be a complex manifold.  A {\it singular holomorphic
foliation}\index{holomorphic foliation! of codimension one with
singularities} {\em of codimension one} $\fa$ of $M$ is given by an
open cover $M=\bigcup_{j\in J}U_j$ and holomorphic integrable
$1$-forms $\omega_j \in \bigwedge^1(U_j)$ such that if $U_j \cap U_j
\ne \emptyset$, then $\omega_i = g_{ij}\omega_j$ in $U_i\cap U_j$,
for some $g_{ij}\in \O^{*}(U_i\cap U_j)$.  We put $\sing(\fa)\cap
U_j =\{p\in U_j; \omega_j(p)=0\}$ to obtain $\sing(\fa) \subset M$,
a well-defined analytic subset of $M$, called singular set of $\fa$.
$M \setminus \sing(\fa)$ is foliated by a holomorphic codimension
one (regular) foliation $\fa_0$.   By definition the {\it
leaves}\index{leaf} of $\fa$ are the leaves of $\fa_0$.}
\end{Definition}

\begin{Remark}
{\rm We may always assume that $\sing(\fa)\subset M$ has codimension
$\geq 2$.  If $(f_j = 0)$ is a local equation of a codimension one
component of $\sing(\fa) \cap U_j$, then we get $\omega_j = f_j^n
\bar{\omega}_j$ where $\bar{\omega}_j$ is a holomorphic $1$-form and
$\sing(\bar{\omega}_j) $ does not contain $(f_j =0)$.}
\end{Remark}

\begin{Example}
[Darboux foliations] {\rm  Let $M$ be a complex manifold and let
$f_j\colon M\to{\co }$ be holomorphic  functions and $\la_j\in{\co
}^*$ complex numbers, $j=1,\dots,r$. The holomorphic $1$-form $\om
=(\prod\limits_{j=1}^r f_j). \sum\limits_{i=1}^{r} \la_i
\frac{df_i}{f_i}$ \, is integrable. Indeed, $\frac{1}{h}\omega$ is
closed and meromorphic, where $h=\prod\limits_{j=1}^r f_j$.
Therefore, $\omega$ is integrable in $M\setminus \{h=0\}$ and, by
the Identity Principle, $\omega$ is integrable on $M$. Thus $\omega$
defines a codimension one holomorphic foliation $\fa(\omega)$ in
$M$. Later on, we will see that we may assume that the singular set
of $\fa(\omega)$ has codimension $\geq 2$
(Proposition~\ref{Proposition:saturationfoliation}). We will also
see that the functions $f_j$ may be assumed to be meromorphic (cf.
Proposition~\ref{Proposition:extensionfoliation}). We call
$\fa(\omega)$ a  {\it Darboux  foliation}\index{foliation! Darboux}
of $M$.  The foliation $\fa(\omega)$ has $f = \prod\limits_{j=1}^{r}
f_j^{\la_j}$ as a {\it logarithmic} first integral. For this reason,
$\fa(\omega)$ is also called a {\it logarithmic
foliation}\index{foliation! logarithmic}. }
\end{Example}

\section{Analytic leaves}

In general the leaves of a foliation are not embedded submanifolds,
as we have already mentioned. Nevertheless when the
foliation has a properly embedded leaf, this leaf is an
analytic subset of the
ambient manifold. We shall see next a criteria  for a foliation,
defined in a  complex manifold $M$ by a  $1$-form, to have an
analytic leaf.  We consider the following situation: \vglue.05in
Let $\fa$ be a holomorphic foliation defined in
a connected manifold $M$, by a holomorphic integrable $1$-form
$\om$. Let also  $f\in \O(M)$ be a non-constant
holomorphic function, that vanishes at some point
of $M$, so that  the analytic subset
$(f=0)$ of $M$ is nonempty and has codimension one.
We shall say that the analytic set $(f=0)$ is {\it
invariant}\index{invariant analytic set} by $\fa$ if
is its
connected components connected are leaves of $\fa$. \vglue.05in
\begin{Proposition}
\label{Proposition:4.I} In the above situation, the analytic variety
$(f=0)$ is invariant by $\fa$ if, and only if,  there exists a
holomorphic 2-form $\theta$ on $M$ such that
$$(*)\ \ \om\wedge df= f\theta$$

\vglue.05in
\end{Proposition} \begin{proof} Suppose
that $(f=0)$
is  invariant by $\fa$. In this case, since
each connected component
of $(f=0)$ is a leaf of $\fa$, these are smooth and properly
embedded submanifolds of  $M$. Thus, given a point
$p$ such that $f(p)=0$,
we can choose a trivializing chart of $\fa$, $(\phi=(x,y), U)$, such
that
$p \in U$, $\phi (p)=0$,
$x \colon U \to \mathbb{C}^{n-1}$, $y \colon U \to \mathbb{C}$ and
the plaques of $\fa$ in $U$ are of the form $y^{-1}(q),q\in y(U)$.
Note that, since $(f=0)$ is embedded, we can assume that $(f=0)\cap
U = y^{-1}(0)$. We obtain  then that $f(x,0)\equiv 0$. From this it
follows that  $f(x,y) = y^{k}.u(x,y)$, where $k \geq 1$ and $u$ is
holomorphic and does not vanish in $U$.

\vglue.05in

On the other hand, since the plaques of $\fa$ in
$U$ are of the form $y=const.$, we can write $\om \mid_U = g.dy$,
where $g$ is holomorphic and does not vanish in $U$. This implies

$$
\om\wedge \frac{df}{f} = g.dy \wedge (k.\frac{dy}{y} +
\frac{du}{u})=g.dy \wedge \frac{du}{u}
$$
\vglue.05in

This proves that a 2-form $\theta=\om \wedge \frac{df}{f}$ is
holomorphic, completing this part.

\vglue.05in

Assume now that $\om\wedge df= f\theta$.
Let $L$ be an irreducible component of $(f=0)$ and take $p \in L$.
Computing  (*) in $p$, we obtain,
$$\om_{p} \wedge df_p=0\ \Rightarrow\ (**)\ df_p=\lambda(p).\om_p,\
\text{where}\ \lambda(p) \in \mathbb{C}$$ since $\om_p \ne 0$. We
have two cases to consider: (a)\, $df \not\equiv 0$ in $L$, (b)\,
$df \equiv 0$ in $L$. \vglue.05in

Let us consider  case (a). In this case the set $A= \{p\in L; df_p
\ne 0\}$ is open and dense in $L$ (see the Identity principle in
\cite{Gunning 1}). On the other hand, $(**)$ implies that, if
$p\in A$
then $\lambda(p)\ne 0$ and $T_p L=\ker(df_p)= \ker(\om_p)$. It
follows that $A$ is contained in a  leaf of $\fa$ and therefore its
$L$, is a leaf of $\fa$. \vglue.05in

Let us consider
case (b). Here we use the fact that  the
set of smooth points in $L$ is open and dense in $L$ (see
\cite{Gunning 2, gunning-rossi}). Given a smooth point $p$ of $L$, there exists a coordinate
system $(\phi=(x,y), U)$ such that
$p \in U$, $\phi (p)=0$, $x \colon U \to
\mathbb{C}^{n-1}$, $y \colon U \to
\mathbb{C}$ and $U\cap L= (y=0)$.
Since $f \mid_L \equiv 0$, we obtain that $f(x,y)= y^k . u(x,y)$,
where $k\geq 2$ and $u$ is holomorphic and does not vanish
in $U$.
On the other hand, we can write $\om \mid_U = b\,dy+
\sum_{i=1}^{n-1} a_i\,dx_i $, and therefore of (*), we obtain that
the 2-form below is holomorphic
$$\om \wedge \frac{df}{f}= \om\wedge\frac{du}{u} + k\, \sum_{i=1}^{n-1}a_i \, dx_i
\wedge \frac{dy}{y}.$$

\vglue.05in

Since $\frac{du}{u}$ is
holomorphic, we obtain  that $\sum_{i=1}^{n-1}a_i \,dx_i \wedge
\frac{dy}{y}$ is holomorphic. We conclude that $y$ divides $a_i$,
for every $i=1,...,n-1$, i.e., that we can write $\om\mid_U=y.\eta +
b.dy$ where $b$ and $\eta$ are holomorphic. This implies that
$(y=0)=L\cap U$ is invariant by $\fa$. Since the set of smooth
points  of $L$ is open and dense in $L$, we can conclude that $L$ is
invariant by $\fa$, {\it i.e.},  it is a leaf of $\fa$.
\end{proof}

\vglue.1in


\section{Two extension lemmas for holomorphic foliations}

This section is dedicated to the extension of codimension one
holomorphic foliations with singularities,  through codimension
$\geq 2$ analytic subsets. Moreover, we conclude that the singular
set of such foliation, may always be assumed to have codimension
$\geq 2$.
\begin{Proposition}
\label{Proposition:saturationfoliation}  Let $\fa$ be a codimension
one holomorphic foliation with  singularities on a connected complex
manifold $M$. Then there exists  a foliation $\fa_1$ of $M$ with the
following properties:
\begin{itemize}
\item[{\rm(a)}]  The irreducible components  of $\sing(\fa_1)$ have codimension $\ge 2$,
and  $\sing(\fa_1)\subset \sing(\fa)$.

\item[{\rm(b)}]   $\fa_1$ coincides with $\fa$ on $M\setminus \sing(\fa)$

\item[{\rm(c)}]  $\fa_1$ is {\it maximal} in the following sense: if  $\fa_2$ is another  foliation
of $M$ satisfying {\rm(a)} and {\rm(b)}, then $\fa_2=\fa_1$.

\end{itemize}
\end{Proposition}

\begin{proof} The proof is basically a consequence of Hartogs' extension theorem already mentioned in this text.
Let  $\{\om_\alpha\}_{\alpha\in A}$, $\{U_\alpha\}_{\alpha\in  A}$ and
$\{g_{\alpha \beta}\}_{U_\alpha \cap U_\beta \ne \phi}$,  be collections defining  $\fa$.
Assume that $\sing(\fa)$  has some codimension one irreducible components. Denote by  $W$ the union of these
components. We develop an elimination procedure for these
components. Given a point  $p\in W$  we choose local coordinates
$(x=(x_1,...,x_n),U_p)$ such that $p\in U_p$, $x\colon U_p \to
\Bbb{C}^{n}$, and  $x(U_p)$ is an open polydisc in $\Bbb{C}^n$, and
also
$W\cap U_p$ has only finitely many irreducible components, say  $W^{p}_{1},...,W^{p}_{r}$,
with corresponding irreducible equations $f_1,...,f_r$ respectively (see \cite{Gunning 2}).
We can assume that $U_p \subset U_\alpha$, for some $\alpha=\alpha(p) \in A$.

If $g$ is a holomorphic function in $U_p$ vanishing  in $W\cap U_p$,
then $g=f_{1}^{n_1}...f_{r}^{n_r}.h$, where $n_1,...,n_r \in \Bbb{ N
}$ and $h\in \O(U_p)$ (cf. \cite{Gunning 2,gunning-rossi}).  We can
also write $\om_\alpha \mid_{U_p} = \sum_{j=1}^{n} a_j\, dx_j$.
Since $\om_\alpha\mid_{W\cap U_p}\equiv 0$, we conclude that the
coefficients  $a_j$ of $\om_\alpha\mid_{U_p}$ vanish in $W\cap U_p$
and therefore $\om_\alpha =f_{1}^{n_1}...f_{r}^{n_r} . \om^\prime
_{p}$, where $n_1,..., n_r \in \Bbb{ N }$ and $\om^\prime _{p}$
is an integrable $1$-form on
$U_p$, having  singular set of codimension $\ge 2$.

On the other hand, if $p\notin W$ then we take  $U_p\subset
U_\alpha$, $\om_{p}^\prime= \om_{\alpha}\mid_{U_p}$, for some
$\alpha =\alpha(p)\in A$, in such a way that
$U_p \cap W=\phi$ and that $U_p$ is the domain of a   local chart $x=(x_1,...,x_n)$.
\vglue.05in In this way, we can define an open  cover $\{U_p\}_{p\in
M}$ and a collection  $\{\om^{\prime}_{p}\}_{p\in M}$, where
$U_p\subset U_{\alpha(p)}$,  $\om^{\prime}_{p}$ is a holomorphic
integrable $1$-form on $U_p$ such that \linebreak $\codim(\sing
X\om^{\prime}_{p}))\ge 2$ and
$\om^{\prime}_{p}$ generates $\fa$ in $U_p\setminus \sing(\om_{\alpha(p)})$
({\it i.e.}, if $q\in U_p\setminus \sing(\om_{\alpha(p)})$, then $T_q \fa=\ker(\om_\alpha(q)$).
We shall see next that there exists a collection $\{g_{p,q}\}_{U_p\cap U_q \ne \phi}$, where $g_{p,q}\in
\O^{*}(U_p\cap U_q)$,  such that
$\om^{\prime}_{p}=g_{p,q}.\om^{\prime}_{q}$ in $U_p\cap U_q \ne
\phi$. \vglue.05in Let $p,q\in M$ such that $U_p\cap U_q \ne \phi$ e
$\alpha=\alpha(p)$ and $\beta=\alpha(q)$. Let also
$x=(x_1,...,x_n)\colon U_p \to \Bbb{C}^n$ be a local system of
coordinates. We can write $\om^{\prime}_{p}=\sum_{j=1}^{n}a_j\,
dx_j$ and $\om^{\prime}_{q}\mid_{U_p\cap U_q}=\sum_{j=1}^{n}b_j\,
dx_j$. Observe that  $\om_\alpha= g_{\alpha \beta}.\om_\beta$
implies that
$$\frac{a_1}{b_1}=...=\frac{a_n}{b_n}= g_{p,q}\ \ \ \text{em}\ U_p\cap U_q.$$ This means that
$\om^{\prime}_{p}=g_{p,q}.\om^{\prime}_{q}$, where
$g_{p,q}$ is at first sight  meromorphic. It is then enough to prove
that $g_{p,q}$ extends to a function in $ \O^{*}(U_p\cap U_q)$.
\vglue.05in Indeed, first we  observe that the singular sets  $S_p$
and $S_q$, of $\om^{\prime}_{p}$ and of $\om^{\prime}_{q}$, are of
codimension $\ge 2$. Let us put  $Z=(S_p\cup S_q)\cap(U_p\cap U_q)$.
\vglue.05in Given $z_o\in (U_p\cap U_q)\setminus Z$, there exists an
index $j\in \{1,...,n\}$ such that $b_j(z)\neq 0$, for every $z$ in
a small neighborhood of $z_o$. Hence  $g_{p,q}=\frac{a_j}{b_j}\in
\O^(U_p\cap U_q \setminus Z)$. Since $Z$ is of codimension $\ge 2$,
it follows from Hartogs' extension theorem
(Theorem~\ref{Theorem:Hartogs}) that  $g_{p,q}$ extends to a
holomorphic function in $U_p\cap U_q$ (see \cite{Gunning 2,
gunning-rossi}). For the same reason, $\frac{1}{g_{p,q}}$ also
extends. Therefore the extension we just obtained does not vanish.
The proof of  (c) is a consequence of the following proposition.
\end{proof}

Similarly to above we obtain:

\begin{Proposition}
\label{Proposition:extensionfoliation} Let $M$ be a connected
complex manifold of dimension $\ge 2$, and let $V$ be an analytic
subset of $M$ of codimension $\ge 2$. Given $\fa$ a codimension one
holomorphic foliation of $U=M\setminus V$, there  exists a unique
foliation $\fa^{\prime}$ of $M$, whose restriction to $U$ coincides
with $\fa$.
\end{Proposition}

\chapter{Holomorphic foliations given by closed $1$-forms}

This chapter is dedicated to the study of an important class of
codimension one foliations with singularities. The class of
foliations given by closed $1$-forms. We study their holonomy and
some extension property. Starting with the holomorphic case we are
to consider the meromorphic case, making use of the extension
results Propositions~\ref{Proposition:saturationfoliation} and
~\ref{Proposition:extensionfoliation}. This study also is related to
the classification and normal forms of the isolated singularities of
foliations, as suggested for instance by
\S~\ref{section:poincaredulac} and
Example~\ref{Example:poincaredulac}.

\section{Foliations given  by closed  holomorphic $1$-forms}

{\rm  \noindent Let $M$ be a complex manifold of dimension $\ge 2$
and $\om$ be a closed holomorphic $1$-form  on $M$ (that is
$d\om=0$) assumed to be not identically zero. Then, $\om$ is
integrable ( $\omega\wedge d \omega=0$) and therefore it defines a
foliation $\fa$ of $M$. Classical  integration lemma of  Poincaré
assures that, given any simply-connected open subset  $U\subset M$,
there exists a holomorphic function $f\colon U \to \mathbb{C}$, such
that $\om \mid_{U} =df$. Observe that if  $g\colon V \to \mathbb{C}$
is a function such that $dg=\om$, where $U\cap V$ is connected and
non-empty, then $g-f$ is constant in  $U\cap V$. Therefore, as we
have already observed before, the corresponding  foliation $\fa$ can
be locally defined by holomorphic functions as follows: there exist
collections ${\mathcal U}=\{U_\alpha\}_{\alpha\in A}$,
$F=\{f_\alpha\}_{\alpha \in A}$ and $C=\{c_{\alpha \beta}
\}_{U_\alpha \cap U_\beta \ne \emptyset}$, such that:
\noindent{\rm{\rm(i)} ${\mathcal U}$ is a cover of $M$ by
simply-connected open sets.}

\noindent{\rm{\rm(ii)} If $\alpha \in A$,
then $f_\alpha$ is a
holomorphic function  in $U_\alpha$ such that
$df_\alpha=\om\mid_{U_\alpha}$.}

\noindent{\rm{\rm(iii)} If
$U_\alpha \cap U_\beta \ne \emptyset$, then
$U_\alpha \cap U_\beta$ is connected,
$c_{\alpha \beta} \in \mathbb{C}$ and
$f_\alpha = f_\beta + c_{\alpha \beta}$ in $U_\alpha \cap U_\beta$.}
\vglue.05in Observe that if $\om$ is free of
singularities, then the
functions $f_\alpha$ are submersions  and $\fa$ is a non-singular
foliation. In this case, if we denote by
$g_{\alpha \beta}$ the translation
$g_{\alpha \beta}(z)=z+c_{\alpha \beta}$, then
$f_\alpha=g_{\alpha \beta}\circ f_\beta$, and then
$\fa$ can be defined by
local submersions as in (II) of
Proposition~\ref{Proposition: 1.I}, in which case the  $g_{\alpha
\beta}$ are translations. We say then that $\fa$ has  an {\it
additive transverse structure}\index{transverse structure!
additive}. In case $\sing(\om)\ne \emptyset$, then  $\fa$ has an
additive transverse structure
in $M\setminus \sing(\fa)$. \vglue.05in

Conversely, if  $\fa$ is a foliation with an additive transverse
structure  in $M\setminus \sing (\fa)$ and such that $\codim(\sing
(\fa))\ge 2$, then $\fa$ can be defined by a closed holomorphic
$1$-form  $\omega$ on $M$.  Indeed, let be given ${\mathcal
U}=\{U_\alpha\}_{\alpha\in A}$, $F=\{f_\alpha\}_{\alpha \in A}$ and
also $C=\{c_{\alpha \beta} \}_{U_\alpha \cap U_\beta \ne \phi}$
collections satisfying {\rm(ii)} and {\rm(iii)}, where ${\mathcal
U}$ is an open cover of $M\setminus \sing(\fa)$. From {\rm(iii)} we
obtain that, if $U_\alpha \cap U_\beta \ne \emptyset$, then
$df_\alpha=df_\beta$ in $U_\alpha \cap U_\beta$. This implies that
there exists the $1$-form holomorphic $\om$ on $M\setminus
\sing(\fa)$ such that $\om\mid_{U_\alpha}=df_\alpha$. It is not
difficult to see that the form $\om$ is closed and defines (induces)
the foliation $\fa$ in $M\setminus \sing(\fa)$. On the other hand,
since $\codim(\sing(\fa))\ge 2$, the classical extension theorem of
Hartogs then implies that $\om$ extends to a holomorphic form on
$M$, which is also closed and defines a foliation $\fa$.

We can then state the following result:

\begin{Proposition}
\label{Proposition:14.I} Let $M$ be a  holomorphic manifold and
$\fa$ a foliation of $M$, with singular set of codimension $\ge 2$.
Then $\fa$ can be defined by a closed $1$-form if, and only if,
$\fa$ has an additive transverse  structure  in $M\setminus
\sing(\fa)$.
\end{Proposition}

\subsection{Holonomy of foliations defined by closed holomorphic $1$-forms}

Let $M$ be a complex manifold of dimension $n\ge 2$ and $\om$ a
closed holomorphic $1$-form not identically zero in  $M$. Let $\fa$
be the singular foliation of codimension one defined by $\om$ in
$M$. Our purpose is to prove:
\begin{Proposition}
\label{Proposition:trivialholonomy} The foliation $\fa$ is without
holonomy. More precisely, if $L \subset M\setminus \sing(\om)$ is a
leaf of $\fa$,  then $L$ has trivial holonomy.
\end{Proposition}
\begin{proof}
We make use of the existence of an additive transverse structure for the foliation.
Fix a  closed regular curve $\gamma \colon I \to L$ with
$\gamma(0)=\gamma(1)=p_o$. Given $q\in \gamma(I)$ there exists a
local chart $(x,y)\colon U \to \mathbb{C}^{n-1}\times \mathbb{C}$
such that $U\cap L=(y=0)$ and $\om \mid_U =dy$ (Poincaré integration
lemma). We can then obtain a collection $\mathcal
C=\{((x_j,y_j),U_j) \}_{j=1}^{k}$ of such charts and a partition
$\{0=t_0 <t_1<...<t_k=1 \}$ of $I$, such that:
\begin{itemize}
\item[{\rm(i)}] $\cup_{j=1}^{k}U_j =\gamma (I)$.
\item[{\rm(ii)}] $\gamma([t_{j-1},t_{j}])\subset U_j,\ \ \forall j=1,...,k.$
\item[{\rm(iii)}] $\om \mid_{U_j} =dy_j,\ \ \forall j=1,...,k.$
\end{itemize}

\vglue.05in Since $\gamma(0)=\gamma(1)$, we can assume that:

\begin{itemize}
\item[{\rm(iv)}] $((x_1,y_1),U_1)=((x_k,y_k),U_k)=((x,y),U)$, where $x(p_o)=0\in \mathbb{C}^{n-1}$.
\end{itemize}

\vglue.05in Let us consider the transverse sections
$\Sigma_j=\{(x_j,y_j)\in U_j;\ \
x_j=x_{j}^{o}=x_j(\gamma(t_j))\}\subset U_j,\ j=1,...,k$, where we
shall compute the holonomy in the transverse section
$\Sigma=\Sigma_k$. For sake of  simplicity,  we shall denote the
point $(x_{j}^{o},y_j)\in \Sigma_j$ by $y_j$. Moreover, for the
uniformization of  the notation we shall put $\Sigma_0=\Sigma$ and
$y_0=y=y_k$. \vglue.05in Let us compute the holonomy $f_j \colon
\Sigma_{j-1} \to \Sigma_j,\ j=1,...,k$. This holonomy is of the form
$y_j=f_j(y_{j-1})$. It is enough to prove that
$f_{j}(y_{j-1})=y_{j-1},\ j=1,...,k$. This implies that the holonomy
of $\gamma$, which is the composition  $f_k \circ ...\circ f_1$, is
equal to the  identity of $\Sigma$. \vglue.05in

Indeed, since $\om \mid_{U_{j-1}\cap U_j}=dy_{j-1}=dy_j$, we obtain
that $d(y_j - y_{j-1})=0$ in $U_{j-1}\cap U_j$. This implies that
the difference $y_j - y_{j-1}$ is constant in the  connected
component of $U_{j-1}\cap U_j$ that contains $\gamma(t_{j-1})$, say
$y_j=y_{j-1}+c$. On the other hand, because $U_{j}\cap L=(y_j=0)$
and $U_{j-1}\cap L=(y_{j-1}=0)$, we obtain that $c=0$ completing the
proof.
\end{proof}

\section{Foliations given by closed meromorphic $1$-forms}
Let $M$ be a holomorphic manifold of dimension $n\ge 2$ and $\om$ a
closed meromorphic $1$-form (not holomorphic) on $M$. We shall
denote the polar divisor of $\om$ by $(\om)_{\infty}\subset M$ (see
the definition in \cite{[Ho]}). In the current case we have
$(\om)_{\infty}\ne \emptyset$, since  $\om$  is assumed to be not
holomorphic. Since $\om$ is closed and holomorphic in the open set
$N=M\setminus (\om)_{\infty}$, it  defines a codimension one
foliation of $N$ which we shall denote  $\fa$. Let us see that,
indeed, $\fa$ extends to a foliation of $M$.

\begin{Proposition}
\label{Proposition:15.I}  The foliation $\fa$ extends to $M$ as a
foliation $\fa'$  such that {\it $(\om)_\infty$ is invariant by
$\fa'$}.
\end{Proposition}
\vglue.05in \begin{proof} In the proof of the extension of  $\fa$ we
shall use  the  following fact (see \cite{Gunning 2},
\cite{gunning-rossi}):

\begin{Fact}
The set $L$ of smooth points of $(\om)_{\infty}$ is open and dense
in $(\om)_{\infty}$. Moreover,  the set $S= (\om)_{\infty} \setminus
L$ is an analytic subset of $M$ de codimension $\ge 2$.
\end{Fact}
\noindent The idea is to prove first that  $\fa$ extends to the set $M\setminus S$ and then make use of
Proposition~\ref{Proposition:extensionfoliation}. \vglue.05in With this aim,
let us fix a point  $p\in L$. Since $p$ is a smooth point  of $(\om)_{\infty}$ and
this  set has codimension one, there exists a system of holomorphic
coordinates in  a neighborhood $U$ of $p$, $w=(x,y)\colon U \to
\mathbb{C}^n$, where $w(U)$ is a polydisc,
$x=(x_1,...,x_{n-1})\colon U \to \mathbb{C}^{n-1}$ and such that
$U\cap L=U\cap (\om)_{\infty}=\{(x,y);y=0\}$. \vglue.05in Observe
now  that by  definition of the polar set, there exists  $j>0$ such
that $y^j\,\om$ extends to a holomorphic form on $U$. Let
$$k= \min \{j>0;\ \ y^j\,\om\ \text{extends to a holomorphic form em}\ U \}$$
and let us define $\eta {=} y^k.\om$. We can write $\eta {=} a_n .dy
+ \sum\limits_{j=1}^{n-1}a_j\,dx_j$, where some of the  functions
$a_1,...,a_n$ does not vanish identically  in $U\cap L$. Note that
$\eta$ is integrable in $U\setminus L$ and defines the same
foliation that $\om$ in this set. Henceforth,  $\eta$ is integrable
in $U$, and there it defines a  foliation of $U$ that extends $\fa
\mid_U$. \vglue.05in On the other hand, from $\om=y^{-k}.\eta$ we
obtain
$$d(y^{-k})\wedge \eta + y^{-k}.d\eta=d\om=0\ \Rightarrow\ (*)\ \ dy\wedge \eta =k^{-1}.y.d\eta$$
\vglue.05in Then, from  (*) and from
Proposition~\ref{Proposition:4.I} we get that $(y=0)=L\cap U$ is
invariant by the  foliation defined by $\eta$. This implies then
that $\fa$ extends to $M\setminus S$ in such a way that  $L$ is
invariant by the extension, as desired.
\end{proof}  \vglue.05in

\begin{Remark}\label{Remark:2.I}{\rm Observe that the connected components of  $L$ are leaves
of $\fa'$. Moreover, , given  a leaf $L_0\subset L$ and a coordinate
system  $w=(x,y)\colon U \to \mathbb{C}^{n-1}\times \mathbb{C}$ such
that $U\cap L=U\cap L_0=(y=0)$, we can consider the number
$$k= \min \{j>0;\ \ y^j\,\om\ \text{extends to a holomorphic form on}\ U \}$$
\vglue.05in It is possible to prove that  $k$ only depends on $L_0$,
{\it i.e.}, does not depend on the coordinate system that we
consider. \vglue.05in We shall say then that $\om$ {\it has a pole
of order $k$ in  $L_0$}. \vglue.05in  Next we shall see  a
particular  example of the situation above. }
\end{Remark}

\subsection{Holonomy: meromorphic case}
Let $M$ a complex manifold of dimension $\ge 2$ and $\om$ a closed
meromorphic $1$-form on $M$. According to
Proposition~\ref{Proposition:15.I}, the foliation with singularities
defined by $\omega$ on $M\setminus (\om)_{\infty})$ can be extended
to a foliation of $M$, which we shall  denote by $\fa$, with the
property that $(\om)_{\infty}$ is invariant by $\fa$. In other
words, the smooth part of of the polar set of $\omega$ is a reunion
of  leaves of $\fa$. Let us now see how to compute the holonomy of a
leaf of this foliation. The first remark  is that a leaf which is
not contained in the polar set, has trivial holonomy by
Proposition~\ref{Proposition:trivialholonomy}. Let us then compute
the holonomy of the leaves contained in the polar set
$(\omega)_\infty$.

For this sake we need some notation.  \vglue.05in For $\mathbb N \ni
k\ge 2$ and $a\in \mathbb{C}$, we consider  the following vector
field:
$$ Y^{k,a}= \frac{y^k}{1+a.y^{k-1}} \frac{\partial}{\partial y}$$
defined in the open set  $\{ y\in \mathbb{C};\ 1+a.y^{k-1}\ne 0\}$.
Note that $Y^{k,a}$ generates a  local flow in a neighborhood of
$0\in \mathbb{C}$, which will be denote  by $Y^{k,a}_{z}$. Thus, for
a fixed $z\in \mathbb{C}$, $Y^{k,a}_{z}$ is a biholomorphic map
between  neighborhoods of $0\in \mathbb{C}$, since
$Y^{k,a}_{z}(0)=0$. Let us  denote the germ  of $Y^{k,a}_{z}$ at the
origin by $[Y^{k,a}_{z}]$. \vglue.05in Note that, if $k\ge 3$, then,
$[Y^{k,a}_{z}]$  commutes with the rotation $R_{\lambda}(y)=\lambda
.y$, where $\lambda^{k-1}=1$.  Then, for  every $k\ge 2$ and every
$a\in \mathbb{C}$, the set
$$G_{k,a}= \{ [R_{\lambda}\circ Y^{k,a}_{z}]\ ;\ z\in \mathbb{C}\, \ \lambda^{k-1}=1 \}$$
is an abelian group. \vglue.05in An important particular case is
when  $k=2$ and $a=0$. In this case $G_{2,0}$ is the group of {\it
homographies}\index{homography}  of the form
$$\{y \rightarrow \frac{y}{1+ay}\ ;\ a\in \mathbb{C}\}$$, obtained by straightforward integration of the differential equation
$\frac{dy}{dz}=y^2$.

\vglue.05in Our main result in this section is the following
computation of the holonomy:
\begin{Proposition}
\label{Proposition:18.I} Let $M$ a complex manifold of dimension
$\ge 2$ and $\om$ a closed meromorphic $1$-form on $M$ with polar
set $(\omega)_\infty \subset M$. Denote by $\fa$ the codimension one
foliation of $M$ corresponding to $\omega$.  Let $L$ be a leaf of
$\fa$. Then:

\noindent{\rm(a)} If $L \subset M \setminus (\om)_{\infty}$, then  $\Hol(L)$ is trivial.

\noindent{\rm(b)} If $L \subset (\om)_{\infty}$ and $\om$ has poles of order $1$
in  $L$, then  $\Hol(L)$ is abelian and analytically linearizable,
that is, $\Hol(L)$ is analytically conjugate to a subgroup of linear
maps  $(z\mapsto \lambda z)$ of $\mathbb{C}$,  in some neighborhood
of the origin.

\noindent{\rm(c)} If $L \subset (\om)_{\infty}$ and $\om$ has poles of order  $k\ge 2$ in $L$,
then  $\Hol(L)$ is analytically conjugate to   a subgroup of $G_{k,a}$, for some $a\in \mathbb{C}$.

\end{Proposition}

\begin{proof}  Case (a) has already been addressed.
Suppose now that $L\subset
(\om)_{\infty}$. Fix $p\in L$ and a local foliation chart of $\fa$,
$(x,y)\colon U \to \mathbb D^{n-1}\times \mathbb D \subset
\mathbb{C}^{n-1} \times \mathbb{C}$ such that $(\om)_{\infty}\cap
U=L\cap U=(y=0)$ and the plaques of $\fa$ in $U$ are of the form
$y=c^{te}$. We claim that:
\begin{Claim}
\label{Claim:meromorphicnormalform} We have
$$(*)\ \ \om \mid_U=\frac{g(y)}{y^k}.dy,$$
where $g$ is holomorphic in $\mathbb D$, $g(0)\ne 0$ and $k$ is the
order of the polar set of  $\om$ in $L$.
\end{Claim}
\begin{proof}[Proof of Claim~\ref{Claim:meromorphicnormalform}]
\vglue.05in Indeed, since $\om$ defines $\fa$ on $M\setminus
(\om)_{\infty}$ and the plaques of $\fa$ em $U$ are of the form
$y=c^{te}$, we have $\om \mid_U = h(x,y)dy$, where $h$ is
meromorphic in $U$ with poles in $(y=0)$. On the other hand, because
$\om$ is closed we have $\partial h /\partial x\equiv 0$, {\it
i.e.}, $h=f(y)$, only depends on  $y$. Let $k$ be the order of the
polar set of  $\om$ in $L$. Then $f$ writes as  in (*), as it can be
easily verified.
\end{proof}

Now we need:

\begin{Lemma}
\label{Lemma:7.I}  Let $\alpha$ be a meromorphic $1$-form on a
neighborhood of $0\in \mathbb{C}$. Suppose that $0$ is  a pole of
order $k\ge 1$ for $\alpha$. Then there exists a system of
coordinates $y \colon V \to \mathbb{C}$ with $0\in V$, $y(0)=0$ and
such that $\alpha$ writes in this  system of coordinates as:

\noindent{\rm{\rm(i)}}\,\,\, $\alpha=a.\frac{dy}{y}\ \, a\ne 0$,\,\,\, if\, $k=1$.

\noindent{\rm{\rm(ii)}}\,\,\, $\alpha=\frac{1+a.y^{k-1}}{y^k}.dy$,\,\,\, with\, $a\in \mathbb{C}$,\,\, if $k>1 $.
\end{Lemma}

\vglue.05in \begin{proof} Let us prove the case of simple pole,
$k=1$. The case $k>1$ is similar and is a simple exercise. In the
case of a simple pole  we can write $\alpha= \frac{g(z)}{z}dz$,
where $g$ is holomorphic in neighborhood $W$ of $0$ and $g(0)=a \ne
0$. Note that $a=Res(\alpha,0)$, that is invariant by change of
coordinates (see \cite{[Al]}). We have then $g(z)=a+z.u(z)$, where
$u$ is holomorphic $W$, so that
$$\alpha= a\frac{dz}{z}+u(z)dz.$$
\vglue.05in Let $\varphi(z)$ be  a primitive (integral) of the
holomorphic $1$-form $\frac{u(z)}{a}dz$ in a  neighborhood of $0$.
Let us consider the function $y(z)=z.\exp(\varphi(z))$. Since
$y(0)=0$ and $y^{\prime}(0)\ne 0$, we conclude that  $y$ is a
biholomorphic map between two neighborhoods of $0$. On the other
hand,
$$a\frac{dy}{y}=a\frac{dz}{z}+a d \varphi =a \frac{dz}{z}+u(z)dz=\alpha,$$ completing the proof.
\end{proof}  \vglue.05in

Let us resume the proof of the proposition.

Assume first that  $k=1$. Fix a closed curve $\gamma \colon I \to L$
with $\gamma(0)=\gamma(1)=p_o$.

From Lemma~\ref{Lemma:7.I}  and with similar arguments to the
holomorphic case we have:

\begin{Claim}
There is a  collection  $\mathcal C=\{((x_j,y_j),U_j) \}_{j=1}^{k}$
of foliation charts of  $\fa$,  and a partition $\{0=t_0
<t_1<...<t_k=1 \}$ of $I$, such that:

\noindent{\rm{\rm(i)}\,\, $\cup_{j=1}^{k}U_j =\gamma (I)$.}

\noindent{\rm{\rm(ii)}\,\, $\gamma([t_{j-1},t_{j}])\subset U_j,\ \ \forall j=1,...,k.$}

\noindent{\rm{\rm(iii)}\,\, $\om \mid_{U_j} =a_j\,\frac{dy_j}{y_j},\ \ \forall j=1,...,k.$}

\end{Claim}

Since $\gamma(0)=\gamma(1)$ we can assume that:

\noindent
{\rm{\rm(iv)}\,\, $((x_1,y_1),U_1)=((x_k,y_k),U_k)=((x,y),U)$, where $x(p_o)=0\in \mathbb{C}^{n-1}$.}

Let us consider also sections  $\Sigma_j,j=0,...,k$ as in the
holomorphic case, with  $\Sigma_0=\Sigma_k=\Sigma$. \vglue.05in  Now
we remark that if $A$ is a connected component of $U_{j-1}\cap U_j$
that contains $\gamma(t_{j-1})$, then
$a_{j-1}.\frac{dy_{j-1}}{y_{j-1}}=a_j\,\frac{dy_j}{y_j}$ in $A$.
Comparing the residues of these two forms in $0$, we conclude that
$a_{j-1}=a_j$. We can then conclude that
$\frac{dy_{j-1}}{y_{j-1}}=\frac{dy_j}{y_j}$ in $U_{j-1}\cap U_j$,
for every $j=1,...,k$. This allows us to relate the coordinates
$y_j$ and $y_{j-1}$ in  $A$. Indeed, if  $y_j=f(y_{j-1})$ in $A$,
then we must have $$ \frac{dy_j}{y_j}=
\frac{f^{'}(y_{j-1})}{f(y_{j-1}}dy_{j-1}=\frac{dy_{j-1}}{y_{j-1}}\ \
\Rightarrow\
$$
$$
z.f^{'}(z)=f(z)\ \ \Rightarrow\ \ f(z)=c_j\,z
$$
for some constant $c_j$, as it follows from a straightforward
integration of the differential equation $z.f^{\prime}=f$. This last
implies that the intermediate holonomy maps $f_j\colon
\Sigma_{j-1}\to \Sigma_j$ are linear. Since the composition of
linear maps is also linear, we obtain that the holonomy of $\gamma$
is linear in the coordinate system that we consider. Since this
system only depends on $\om$ (not on the curve $\gamma$), we
conclude  that the holonomy de $L$ is analytically linearizable.
\vglue.05in Now we consider the case where $k\ge 2$. Again we fix a
closed curve $\gamma \colon I \to L$ with $\gamma(0)=\gamma(1)=p_o$.
Similarly to the case $k=1$, from Lemma~\ref{Lemma:7.I} and similar
arguments to the ones in the holomorphic case,  we have:
\begin{Claim}
There is a collection $\mathcal C=\{((x_j,y_j),U_j)
\}_{j=1}^{m}$ of foliation charts
of $\fa$,  and a partition  $\{0=t_0 <t_1<...<t_m=1 \}$ of $I$,
satisfying {\rm(i)},{\rm(ii)},{\rm(iv)} and
\noindent{\rm{\rm(iii)} $\om \mid_{U_j} =\frac{1+a_j\,y_{j}^{k-1}}{y_{j}^{k}}.dy_j,\ \ \forall j=1,...,m.$}
\end{Claim}

\vglue.05in Observe that $a_j=Res(\om,y_j=0)$. Thus, similarly to
the above argumentation,  we conclude  that $a_1=...=a_m=a$.
\vglue.05in As we have seen above, it is enough to relate  $y_j$ and
$y_{j-1}$. Let us do it in the case $k=2$. The case $k>2$ is left to
the reader. For the sake of simplicity of notation let us put
$y_j=w$ and $y_{j-1}=z$, so that $w=f(z)$ is the change of
coordinates. \vglue.05in

\noindent{\bf Case $a=0$.} In this case, in $U_{j-1}\cap U_j$ we
have
$$\frac{dz}{z^2}=\frac{dw}{w^2}=\frac{f'(z)}{(f(z))^2}.dz\ \ \Rightarrow\ \ z^2.f'=f^2\ \ \Rightarrow
\ \ f(z)=\frac{z}{1+c.z}$$ from where we conclude that $f$ is in the
group $G_{2,0}$. Similarly to above, since the composition of
elements in $G_{2,0}$ is also in $G_{2,0}$,  we obtain that
$\Hol(L)$ is analytically conjugate to  a subgroup of $G_{2,0}$.
\vglue.05in

\noindent{\bf Case $a \ne 0$}. The argumentation in the case $k=2$
and $a\ne 0$ is quite similar to the one above: it is enough to
prove that  $y_j$ and $y_{j-1}$ are related  by an  element of
$G_{2,a}$. As it is easily checked, if  $y_j= f(y_{j-1})$, then $f$
satisfies the differential equation
$$(*)\ \ z^2.(1+a.f(z)).f'(z)= (f(z))^2.(1+a.z).$$ Thus it remains to show that  in this case $f$
is in $G_{2,a}$. Let us show how this is done.

\vglue.05in \noindent{\bf Step 1}.  Given  $b\in \mathbb{C}$, there
exists a unique solution $f$ of (*), defined in a  neighborhood of
$0$, such that $f^ \prime(0)=1$ and $f^{\prime}(0)=b$. \vglue.05in
Indeed, if $f$ is a solution of (*) and $g(z)=\frac{f(z)-z}{z^2}$,
then $g$ satisfies the  following differential equation:
$$(**)\ \ g^\prime=g.\frac{a+azg-g}{1+az+azg}=F(z,g).$$
\vglue.05in Since $F$ is holomorphic in a neighborhood of $(0,b/2)$,
(**) has a unique solution $g$, defined in a neighborhood of $0\in
\mathbb{C}$, with $g(0)=b/2$. Setting  $f(z)=z+z^2.g(z)$, we obtain
the required solution.

\noindent{\bf Step 2}. For each $c\in \mathbb{C}$ the function
$f_c=Y^{k,a}_{c}$ is a solution of (*) with initial condition
$f_{c}^{\prime}(0)=1$ and $f_{c}^{\prime \prime}(0)=2c$. \vglue.05in
The proof is a straightforward computation.

Proposition~\ref{Proposition:18.I} then follows from Steps 1 and 2
above.
\end{proof}

\chapter{Reduction of singularities}
\section{Introduction}
\label{section:reduction} The subject of reduction of singularities
is one of the most developed and interesting in the theory of
holomorphic foliations. In the last decades various authors have
given important contributions to the subject. Although the general
problem of the reduction of singularities for a germ of a
holomorphic foliation singularity remains open, several are the
cases where it is already finished. This knowledge has been of great
importance in the development of the theory of holomorphic
foliations with singularities as we shall see in some of the results
to be presented in this text (see
Chapters~\ref{chapter:matteimoussu} and \ref{chapter:limitsets} for
instance). We would like to mention the recent work of Camacho-Lins
Neto-Sad (\cite{C-LN-S1}), Camacho-Sad (\cite{Camacho-Sad}),
Cano(\cite{Cano}) and Cano-Cerveau (\cite{Cn-Ce}) as a  partial
list.

Before going into the reduction of singularities, we recall some
very basic facts from the theory of singularities of holomorphic
vector fields. In short, in the next section we study the analytic
forms of non-degenerate  isolated singularities of  holomorphic
vector fields, according to Poincar\'e and Dulac.
\section{Poincaré and Poincaré-Dulac normal forms}
\label{section:poincaredulac}

In this section we shall consider a holomorphic vector field $X$
defined in a neighborhood $U$  of the origin $0 \in \bc^n$, $n\ge2$,
with a singularity (isolated) at the origin. Since we are interested
in the local analytical description of $X$ or $\fa(X)$ in a
neighborhood of $0$, we shall consider $U$ as small as necessary
without further mention to this.

The {\it eigenvalues\/} of $X$ at $0$ are those of the linear part
$DX(0)$ (here regarded as a  linear map $\bc^n \to \bc^n$). Assume
that $DX(0)$ is non-singular.

\begin{Definition} {\rm We shall say that the singularity $0$ of $X$
{\it is in the Poincaré domain\/}\index{Poincaré domain} if the
convex hull of its eigenvalues in $\mathbb R^2 \simeq \bc$  does not
contain the origin.
Otherwise we say that the singularity $0$ of $X$ {\it is in the Siegel domain}\index{Siegel domain}.
}
\end{Definition}

\begin{Remark}{\rm  If $m=2$ then the singularity is in the Siegel domain if, and only if,
the quotient of its eigenvalues is real negative.
}\end{Remark} A geometrical characterization of singularities in the
Poincaré domain, as well as a rich  dynamical description of them,
is found in the series of works of the  senior Japanese
mathematician, Prof. Toshikazu Ito. An excerpt from those is the
following (\cite{Ito}):

\begin{Theorem}[T. Ito, 1992\index{Theorem! of Ito}]
A singularity of $\bc^m$ of a holomorphic vector field $X$ is in the
Poincaré domain if, and only if,  the foliation $\fa(X)$ is
transverse to every small sphere $S^{2n-1}(r;0)$ of center $0 \in
\bc^n$ and radius $r > 0$. This is the case if $\fa(X)$ is
transverse to {\bf some} sphere contained in its definition domain.
\end{Theorem}

\noindent Let now $m=2$. We shall say that the eigenvalues $\la_1$,
$\la_2$ of $X$ {\it are in resonance\/}  if $\la_1 = k\,\la_2$ or
$\la_2 = k\,\la_1$ for some $k \in \bn$, $k \ge 2$.

For dimension two the complete description of the analytical normal
forms of singularities in  the Poincaré-domain is given below:

\begin{Theorem}[Poincaré linearization theorem, \cite{Arnold,Brjuno,Dulac}\index{Theorem! of Linearization of Poincaré}]

Let $X$ be a holomorphic vector field with a singularity in the
Poincaré domain at the origin $0 \in \bc^2$.  Suppose that the
eigenvalues $\la_1$, $\la_2$ of $X$ are {\tt not} in resonance. Then
there is a unique holomorphic diffeomorphism $\xi$ between
neighborhoods of $0 \in \bc^2$, $\xi(0)=0$, such that:
\begin{itemize}
\item[{\rm(i)}] $\xi'(0) = \Id$
\item[{\rm(ii)}] $\xi_*X = DX(0)$, i.e., $\xi$ takes $X$ into its linear part.
\end{itemize}
\end{Theorem}

\begin{Remark}{\rm (i)\, The above theorem is similar to a result, by Poincaré, for
diffeomorphisms $f\colon (\mathbb C^m,0) \to (\mathbb C^n,0)$. The
proof is based on the  convergence of the formal (power series)
solution to the linearization problem. }
\end{Remark}
For the case of resonances we have:

\begin{Theorem}[Poincaré-Dulac theorem,\cite{Dulac}\index{Theorem! of Poincaré-Dulac}] Let $X$ be a
holomorphic vector field with a singularity in the Poincaré domain
at $0 \in \bc^2$. Suppose that  the eigenvalues are in resonance ,
$\la_1 = k\,\la_2$ for some $k \in \bn$, $k\ge2$. Then these is an
unique holomorphic diffeomorphism $\xi$ between neighborhoods of $0
\in \bc^2$, $\xi(0)=0$, such that:
\begin{itemize}
\item[{\rm(i)}] $\xi'(0) = \Id$
\item[{\rm(ii)}] $\xi_*X = (\la_1\,x + ay^n)\,\dfrac{\po}{\po x} + \left(\la_2\,y\,\dfrac{\po}{\po y}\right)$,
\end{itemize}
for some $a \in \bc$.
\end{Theorem}
If $a=0$ then we are in the {\it analytically linearizable\/} case.
The proof is in the same spirit of the one for the non-resonant
case. We refer to the book of Ilyashenko and Yakovenko
(\cite{Ilyashenko}) for the proofs and a detailed study in these
normal forms.

\section{Blow-up at the origin\,\, (quadratic blow-up)}

The (quadratic) blow-up at the origin is defined as follows: We
consider two copies of $\bc^2$  with coordinates $(a,t)$ and $(u,y)$
and change of coordinates given by $\begin{cases} y=tx&\\
x=uy&\end{cases}$

\noindent In particular $ut = 1$ so that we obtain a complex surface
$\bc_0^2$ that contains a  projective line $\bp \hookrightarrow
\widetilde{\bc}_0^2$ (given by $(x=0)$ and by $(y=0)$). Define the
map $\pi\colon \widetilde{\bc}_0^2 \to \bc^2$ by $\pi(x,t) = (x,tx)$
and $\pi(u,y) = (uy,y)$; then $\pi$ defines a proper holomorphic
projection which is a diffeomorphism between
$\widetilde{\bc}_0^2\setminus \bp$ and $\bc^2\setminus\{0\}$.

The manifold $\widetilde{\bc}_0^2$ is the {\it blow-up of $\bc^2$ at
the origin\/} and it is a  complex (rank 1) vector bundle with basis
$\bc P(1) \approx \bp$, and fiber $\bc$. Roughly speaking, the
``explosion" of $\bc^2$ at the origin ``separates" the lines though
$0$. It is a remarkable fact that the self-intersection $\bp \cdot
\bp$ of the projective line $\bp \hookrightarrow
\widetilde{\bc}_0^2$ (called {\it exceptional divisor}) is negative
equal to $-1$.

By using local coordinates we can introduce the notion of {\it
blow-up of complex surface $M^2$ at  a point $p \in M$\/} is a
natural way and successively by repeating this
process as many times as desired.
\vglue .1in  \begin{figure}[ht]
\begin{center}
\includegraphics[scale=0.65]{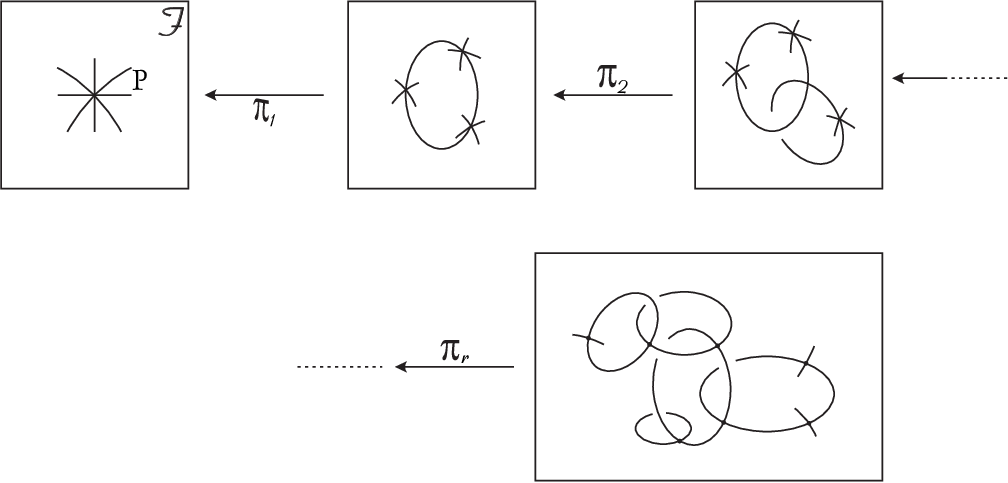}
\caption{ }
\end{center}
\end{figure}

\noindent The final configuration of  the exceptional divisor,
depends on the position of the blow-up centers.

\section{Blow-up on surfaces}
Let us recall the {\it blow-up}\index{blow-up} of $\mathbb{C}^2$ at
$0$. Let us consider two copies of $\mathbb{C}^2$, say $U$ and $V$,
with coordinates  $(t,x)$ and $(s,y)$ respectively. We define a
complex manifold $\tilde{\mathbb{C}^2}$, identifying the point
$(t,x)\in U\setminus (t=0)$ with the point $(s,y)=\alpha
(t,x)=(1/t,tx)\in V\setminus (s=0)$. \vglue.05in The {\it
divisor}\index{blow-up! divisor} of $\tilde {\mathbb {C}^2}$ is
defined as the  submanifold $D$ of $\tilde {\mathbb {C}^2}$ such
that $U \cap D= (x=0)$ and $V\cap D=(y=0)$. Note that, since $y=tx$,
$D$ is well defined and is biholomorphic to the Riemman sphere
$\ov{\mathbb{C}}=\mathbb{C}P(1)$. Moreover,  we can define a
submersion $P \colon \tilde {\mathbb {C}^2}\to D $ by $P
\mid_U(t,x)=t$ and $P\mid_V(s,y)=s$. Then $(\tilde {\mathbb
{C}^2},P,D)$ is a vector fiber space  with basis $D$, projection $P$
and fiber  $\mathbb{C}$, having  $D$ as zero section. \vglue.05in
Let us now consider  the holomorphic map $\pi \colon \tilde {\mathbb
{C}^2}\to \mathbb{C}^2$ defined by $\pi \mid_U(t,x)=(x,tx)$ and $\pi
\mid_V(s,y)=(sy,y)$. Note that $\pi$ is well-defined, since in
$U\cap V$ we have $y=tx$ and $x=sy$. Moreover, , $\pi$ has the
following properties: \noindent{\rm(a)} $\pi^{-1}(0)=D$.
\noindent{\rm(b)} $\pi \mid_{\tilde {\mathbb {C}^2}\setminus D}
\colon \tilde {\mathbb {C}^2}\setminus D \to \mathbb{C}^2\setminus
\{0\}$ is biholomorphic. \noindent{\rm(c)} $\pi$ is a proper map.

\vglue.05in Then, with this structure,
$\tilde {\mathbb {C}^2}$ is
the {\it blow-up}\index{blow-up}  of $\mathbb{C}^2$ at $0$, with
projection\index{blow-up! projection}\index{blow-up! map}
(blow-down) map  $\pi$. \vglue.05in

Let us introduce this same concept in a complex surface. We consider
a complex manifold of dimension two $M$ and a point $q\in M$. The
{\it blow-up}\index{blow-up! of a surface at a point} of $M$ at $q$
is defined as follows: take a
holomorphic local
chart $\varphi \colon A \to B \subset \mathbb{C}^2$ with $q\in A$
and $\varphi(q)=0$. Let $\pi \colon \tilde{\mathbb{C}^2}\to
\mathbb{C}^2$ be the  blow-up map at $0$, with  divisor $D$ and
$\tilde{B}= \pi^{-1}(B)$. In the  disjoint union
$M'=(M\setminus\{q\})\uplus \tilde{B}$ we define an equivalence
relation  $\sim$ by setting $p_o \sim p_1$ if, and only if,
$p_o=p_1$ or, otherwise, $p_o \in A\setminus \{q\}$, $p_1\in
\tilde{B}\setminus D$ and $p_1 = \pi^{-1}(\varphi(p_o))$. The {\it
blow-up}\index{blow-up! of a surface at a point} of $M$ at $q$ is
the quotient $\tilde{M}=M'/\sim$. \vglue.05in Since $\tilde{B}$ is a
manifold
and $\pi^{-1}\circ
\varphi \colon A\setminus \{q\}\to \tilde{B}\setminus D$ is a a
biholomorphic map, it follows that
$\tilde{M}$ is a complex
manifold. Roughly, $\tilde{M}$ is obtained from  $M$ by replacing
the point $q$ by a projective line
$D\simeq
\ov{\mathbb{C}}$. Indeed, the divisor $D$, thanks to the above
procedure, there is a natural embedding of the divisor $D$ into
$\tilde{M}$. \vglue.05in

Given a point $p\in\tilde{M}$, we have three possibilities:
\noindent{\rm(1) The equivalence class of $p$ is in $D$.}
\noindent{\rm(2) The equivalence class of $p$ is in $M\setminus A$.}
\noindent{\rm(3) The equivalence class of $p$ contains two points $p_o\in A\setminus\{q\}$
and $p_1\in \tilde{B}\setminus D$.} \vglue.05in Thus, the points of
$\tilde{M}$ are divided into two classes: those points as in (1),
that will be called   {\it points of the divisor}\index{blow-up!
divisor}, and the  points of $\tilde{M}\setminus D$, that will be
regarded as  points of $M$ (as in  (2) or (3)). \vglue.05in A map of
blow-down $\Pi \colon \tilde{M} \to M$ is defined by $\Pi(p)=q $ in
case (1), $\Pi(p)=p$ in case (2) and $\Pi(p)=p_o$ in case (3). It is
not difficult to see that $\Pi$ has  properties analogous to those
of  $\pi$, {\it i.e.,},
\noindent{\rm(a')$\Pi^{-1}(q)=D$.}
\noindent{\rm(b')$\Pi \mid_{\tilde {M}\setminus D} \colon \tilde {M}\setminus D \to M\setminus
\{q\}$ is a biholomorphic map.
\noindent{\rm(c')$\Pi$ is a proper map.}
\vglue.05in

The above process can be iterated: start with the manifold
$M$ and a point $q_o\in M$. Blowing-up
$M$ at $q_o$, we obtain a manifold  $M_1$ and a blow-up map $\Pi_1
\colon M_1 \to M$ with divisor $D_1=\Pi_{1}^{-1}(p_o)$. Next we
consider any point  $q_1\in M_1$, and perform the blow-up of $M_1$
at $q_1$ obtaining in this way a manifold
$M_2$ and a map
blow-up map  $\Pi_2 \colon M_2 \to M_1$ with divisor $D_2$.
Proceeding in this way, after  $k$ blow-ups, we obtain a manifold
$M_n$ and a blow-up map $\Pi_n \colon M_n \to
M_{n-1}$ with  divisor $D_n$. The composition $\Pi^n=\Pi_n \circ
...\circ \Pi_1 \colon M_n \to M$ is a proper
holomorphic
map that we will call the  {\it blow-up}\index{blow-up! map} map.

\vglue.05in The  divisor $D^n$ of $\Pi^n$ is defined inductively as
follows:
\noindent{\rm(I)$D^1=D_1$.}
\noindent{\rm(II)$D^n=D_n \cup \Pi^{-1}_{n}(D^{n-1})$.}
\vglue.05in Note that $\Pi(D^n)$ is a finite subset
of $M$: indeed, this corresponds to
centers of the explosions.
Moreover, the  map $\Pi^n \mid_{M_n\setminus D^n} \colon
M_n\setminus D^n \to M\setminus \Pi^n(D^n)$ is a biholomorphism.
\vglue.05in The divisor $D^n$ is the  union of $k$ complex curves,
each curve biholomorphic to  $\ov{\mathbb{C}}$. For instance, the
for the second blow-up, if $q_1\in D_1$, then  $D^2=D_2\cup
\Pi_2^{-1}(D_1)$. It follows that  $\Pi_2^{-1}(D_1)\simeq
\ov{\mathbb{C}}$ and that $D_2$ intersects
$\Pi_2^{-1}(D_1)$
transversally at a single point, {\it i.e.},
$D^2$ is the  union of two
embedded projective lines  in $M_2$ with a single common
point. For sake of simplicity we shall
use the same notation for the  $D_i$ and their
successive inverse images
$\Pi_i,...,\Pi_n$. Then, we can say that
$D^n=\cup_{j=1}^n D_j$.
\vglue.05in In the case that for every $j=1,...,n-1$ the j-th
blow-up is centered at  a point of $D^j$, $D^n$ will be a ``graph
free of cycles", that is, for every $i$ the projective
$D_i$ another projective $D_j$ transversally
at a single point, called a {\it corner}\index{corner} of $D^n$,
in such a way that  if
$D_{i_1}\cap D_{i_2}\ne \phi$,...,$D_{i_{m-1}}\cap D_{i_m} \ne
\emptyset$, then $D_{i_1}\ne D_{i_m}$.

Such a process will be called
{\it  blow-up process centered at  $q$}\index{blow-up! process}. \vglue.05in

\subsection{Resolution of curves}

Let us now see in what consists the
``resolution of a  singularity of a curve".
Let us
consider a curve $C=(f(x,y)=0)\subset A \subset \mathbb{C}^2$, where
$f(0,0)=0$, that is, $0\in C$. We assume that the the
Taylor expansion of $f$
is $f=\sum_{j=k}^{\infty}f_j$, where $f_j$ is a homogeneous
polynomial of degree  $j$. Let $\pi\colon \tilde{\mathbb{C}^2}\to
\mathbb{C}^2$ the blow-down of $\mathbb{C}^2$ in $0$. The expression
of  $\pi$ in the chart  $((t,x),U)$ of $\tilde{\mathbb{C}^2}$, we
obtain
\begin{align*}
f\circ \pi (t,x) &=f(x,tx)=\sum_{j=k}^{\infty}f_j(x,tx) \\
&= x^k.\sum_{j=k}^{\infty} x^{j-k}.f_j(1,t)= x^k.f_U(t,x),
\end{align*}
so that $\pi^{-1}(C)\cap U=(x=0)\cup (f_U(t,x)=0) $. Analogously, we
in the chart $((s,y),V)$, we have $\pi^{-1}(C)\cap V= (y=0)\cup
(f_V(s,y)=0)$, where $f_V(s,y)=\sum_{j=k}^{\infty}y^{j-k}.f_j(s,1)$.
Hence, we have $\pi^{-1}(C)=D\cup \tilde{C}$, where
$\tilde{C}=(f_U=0)\cup(f_V=0)$. The curve
$\tilde{C}$ is
called {\it strict transform of}\index{strict transform of a curve}
$C$.

Note that $\tilde{C}\cap D$
is a finite set. Indeed, $\tilde{C}\cap D \cap U = \{(t,0)\ ;\
f_k(1,t)=0\}$, while  $\tilde{C}\cap D \cap V =\{(s,0)\ $;
$f_k(s,1)=0\}$. \vglue.05in

In general, if we consider a blow-up process
$\Pi^n \colon A_n \to A$ with divisor
$D^n=D_1\cup ...\cup D_n$, we obtain $(\Pi^n)^{-1}(C)=D^n\cup C_n$,
where $C_n\cap D^n$ is a finite set. The curve
$C_n$ is called the {\it
strict transform}\index{strict transform of a curve} of $C$ by
$\Pi^n$. \vglue.05in
\begin{Definition}\label{Definition:14capI}{\rm
Let $C$ be a holomorphic curve in a complex surface
$M$.
We say that the  blow-up process, $\Pi^n\colon M_n \to M$, with
divisor $D^n=\cup_{j=1}^n D_j$ is a {\it
resolution}\index{resolution of a curve}
of $C$, if the corresponding strict transform
$C_n$ satisfies the following properties:
\begin{enumerate}[(i)]
\item $C_n$ is regular.

\item $C_n$ meets each $D_j\subset D^n$
transversally.

\item $C_n\cap D^n$ contains no corners.

\end{enumerate}

\vglue.05in } \end{Definition}

Let us illustrate this with an example.

\vglue.05in
\begin{Example}\label{Example:18.I} {\rm  We consider the singular plane curve  $C\subset \mathbb C^2$ given by $f(x,y)=y^2-x^3=0$.
Let $\pi_1\colon M_1=\tilde {\mathbb C^2_0}\to \mathbb{C}^2$ the
blow-up map of $\mathbb{C}^2$ at the origin $0\in \mathbb C^2$. In
the local coordinates $((t,x),U)$ of $M_1$, we obtain
$$ f\circ \pi_1 (t,x)=f(x,tx)=x^2.(t^2-x).$$
Therefore,  $\pi_{1}^{-1}(C)\cap U$ consists of the divisor $(x=0)$
and of the strict transform   $C_1$ of $C$, with equation $x-t^2=0$.
Then clearly  $\pi_{1}^{-1}(C)\subset U$, so that it is not
necessary to consider the other blow-up chart. The strict transform
$C_1$ of $C$  is regular but not transverse to the divisor $D_1$,
since $C_1\cap D=(0,0)\in U$ and $(x-t^2=0)$ is tangent to $(x=0)$
at this very point. In other words, the curve is still not resolved.
\vglue.05in Then a second  blow-up $\pi_2(u,t)=(t,tu)=(t,x)$ is made
at  $(0,0)\in U$. The divisor $D^2$ of this second blow-up  is the
union of two projective lines, $D_1\cup D_2$, and in the chart
$(u,t)$, $D_1$ is represented by $(u=0)$ and $D_2$ by $(t=0)$. We
have then $f\circ \pi_1 \circ \pi_2(u,t)= t^3.u^2.(t-u)$. Thence,
the strict transform  $C_2$ of $C$ is $(t-u=0)$. This last curve
cuts $D^2$ at the corner $(0,0)=D_1\cap D_2$, therefore this curve
is still not resolved. With a final  blow-up $\pi_3$ at the point
$(u,t)=(0,0)$, of the form $t=vu$ (in one of the charts),
we obtain a new  divisor $D_3$,
represented  by $(u=0)$. The strict transform  $C_3$ of $C$ with
equation
$v-1=0$, cuts $D_3$ transversally at the point
point $(v,u)=(1,0)$, which is not a corner point. Henceforth, $C_3$
is a resolution of $C$.



Note that, the original coordinates  $(x,y)$ are related to $(v,u)$
by $$(x,y)=\pi_1\circ \pi_2 \circ \pi_3
(v,u)=(v.u^2,v^2.u^3)=\pi^3(v,u).$$ \vglue.05in Using this and the
parametrization  $u \rightarrow (1,u)$ of $C_3$, we can obtain a
parametrization $u \rightarrow (u^2,u^3)$ of $C$. \vglue.05in }
\end{Example}

\begin{Theorem}[Resolution of singularities for curves
\cite{C-LN-S1}\index{Theorem! of resolution of singularities! for
curves}] \label{Theorem:5.I}
Every holomorphic curve in a complex surface admits
a
resolution.
\end{Theorem}
\vglue.05in
\begin{Corollary}\label{Corollary:5.I}
Let $S$ be a holomorphic curve in a complex surface $M$.
Given a point  $q\in S$, there exist a neighborhood $U$ of $q$ and
holomorphic curves $S_1,...,S_m \subset U$ such that:
\begin{itemize}
\item[{\rm(a)}]  $q\in S_j$ for every $j=1,...,m$.

\item[{\rm(b)}]  $S\cap U \subset S_1\cup ... \cup S_m$. \vglue.05in

\item[{\rm(c)}]  $S_i\cap S_j=\{q\}$, if $i\ne j$.

\item[{\rm(d)}]  Given $j=1,...,m$, there exists
a holomorphic injective map
$\alpha_j \colon \mathbb D_r \to U$,
where $\mathbb D_r=\{z\in \mathbb{C}\ ;\ \mid z \mid <r\}$, such
that $\alpha_j(0)=q$, $\alpha_j(\mathbb D_r)=S_j$ and the
restriction $\alpha_j \mid_{\mathbb D_r \setminus \{0\}}$ is an
embedding.

\vglue.05in In particular, each curve $S_j$ is homeomorphic to the
disc $\mathbb D$.
\end{itemize}
\end{Corollary}
\begin{Definition}\label{Definition:15capI}{\rm
The germs at $q$ of curves $S_1,...,S_m$ are
called the {\it local branches}\index{local branches}\index{branch} of $S$ at $q$.
For each  $j=1,...,m$, the map $\alpha_j$, is called {\it Puiseux
parametrization}\index{Puiseux parametrization} of the branch $S_j$.
}
\end{Definition} \vglue.05in

\begin{proof}[Proof of  Corollary~\ref{Corollary:5.I}]
In the case where  the point $q$ is not a singular point of  $S$ the
result is straightforward. Indeed, in this case, the has only one
branch at  $q$. \vglue.05in Suppose now that $q$ is a singularity
of $S$. Let $\pi \colon \tilde{M} \to M$ a resolution of $S$, with
divisor $D=\cup_{j=1}^{n}D_j$, and $\tilde{S}$ the strict transform
of $S$.  Then $\tilde{S}$ cuts transversally $D$, at non corner
points, forming  a  finite set, say $\{q_1,...,q_m\}$. Because the
curve  $\tilde{S}$ is regular, for each $j=1,...,m$, we can obtain
an embedding  $\beta_j \colon \mathbb D_r \to \tilde{M}$, which is a
parametrization of a   neighborhood of $q_j$ in $\tilde{S}$ such
that $\beta_j(0)=q_j$. Taking the restrictions of the  $\beta_j$ to
a smaller disc if necessary,  we can assume that $\beta_i(\mathbb
D_r)\cap \beta_j(\mathbb D_r) = \phi$ if $i\ne j$. Let us put
$\alpha_j=\pi \circ \beta_j$ and $S_j=\alpha_j(\mathbb D_r)$. It is
not difficult to check that  $S_1,...,S_m$ and
$\alpha_1,...,\alpha_m$ satisfy (a),(c) and (d). We leave the
verification of  (b) as an exercise.
\end{proof}

As a consequence of the above  we obtain:

\begin{Theorem} [cf. \cite{Griffiths-Harris}]
\label{Theorem:6.I}
Let $S$ be a holomorphic  curve in a
complex manifold $M$. Then there exist a
Riemann surface $\tilde{S}$ and a holomorphic map  $\phi \colon
\tilde{S} \to M$ with the following properties:
\begin{itemize}
\item[{\rm(a)}] $\phi(\tilde{S})=S$.
\item[{\rm(b)}] There are discrete subsets
$A\subset \tilde{S}$ and
$B\subset S$ such that $\phi \mid_{\tilde{S}\setminus A} \colon
\tilde{S}\setminus A \to S\setminus B$ is an embedding.

\item[{\rm(c)}]  $\phi^{-1}(B)=A$. Moreover, , $B$ is the
singular set of $S$, and for every $p\in B$, $\phi^{-1}(p)$ is a
subfinite set of $A$.
\end{itemize}

\end{Theorem}

\begin{Definition}\label{Definition:16capI}{\rm
The curve $\tilde{S}$ is called the  {\it normalization
}\index{normalization of a curve}of the curve  $S$.
} \end{Definition}

The theorem above  implies that, given a
singularity $p\in S$, we can define the  {\it
branches}\index{branch} of $S$ at  (through) $p$ in the  following
way: since $\phi^{-1}(p)=\{q_1,...,q_r\}$, is a finite subset of
$\tilde{S}$, we can obtain for each $j=1,...,r$ a disc $\mathbb D_j
\subset \tilde{S}$ such that $p_j\in \mathbb D_j$ and $\mathbb
D_i\cap \mathbb D_j =\phi$ se $i\ne j$. The germs at $p$ of
$\phi(\mathbb D_1),...,\phi(\mathbb D_r)$ are the branches of $S$ by
$p$. The maps $\phi \mid_{\mathbb D_j} \colon \mathbb D_j \to S$,
$j=1,...,r$, are the  {\it Puiseux parametrizations} of these
branches. \vglue.05in

\section{Blow-up of a singular point of a foliation}

Given a foliation $\fa$ of $M$ and  $p \in M$, {\it the blow-up of
$\fa$ at\/} $p$ is the pull-back  foliation $\widetilde{\fa} =
\pi^*(\fa)$ of $\fa$ by the blow-up map $\pi\colon \widetilde{M}_p
\to M$. Both foliations are equivalent of $M\setminus\{p\}$ and
$\widetilde{M}_p\setminus \pi^{-1}(p)$ (recall that $\pi^{-1}(p)$ is
the exceptional divisor $\pi^{-1}(p) \cong \bp$) so that they have
the same leaves. Eventual ``new" singularities are introduced in the
exceptional divisor $\pi^{-1}(p)$. It may occur that $\pi^{-1}(p)$
is invariant or not. If $\pi^{-1}(p)$  is invariant by the foliation
$\fa = \pi^*(\fa)$ we say that the blow-up is {\it
non-dicritical}\index{blow-up! non-dicritical}, otherwise it is
called {\it dicritical}\index{blow-up! dicritical}.
Let us
study this procedure more closely. Let us first see what occurs
with a foliation after  a single blow-up.

Let us consider a holomorphic foliation $\fa$ in a
neighborhood of $0\in \mathbb{C}^2$ with an isolated singularity at
the origin $0$. We assume that  $\fa$ is represented by the vector
field  $X=(P(x,y),Q(x,y))$ or, equivalently, by the dual  $1$-form
$\om=P(x,y)dy-Q(x,y)dx$.
We shall denote by $\fa^*$ the  foliation
with isolated  singularities  $\fa^*=\pi^*(\om)$. Thus $\fa^*$ is
the pull-back of $\fa$ via the blow-up map $\pi\colon \tilde
{\mathbb C^2_0} \to \mathbb C^2$. Let us investigate its expression
in local coordinates. We can write the Taylor expansion
of $\om$ at $0$ as:
$$\om=\sum_{j=k}^{\infty}(P_jdy-Q_jdx),$$
where $P_j$ and $Q_j$ are homogeneous polynomials of degree $j$,
with $P_k \not\equiv 0$ or $Q_k \not\equiv 0$. The $1$-form
$\pi^*(\om)$ writes in the chart  $((t,x),U)$ as:
$$\pi^*(\om)= \sum_{j=k}^{\infty}(P_j(x,tx)d(tx)-Q_j(x,tx)dx)=$$
$$=x^k.\sum_{j=k}^{\infty} x^{j-k}.[(tP_j(1,t)-Q_j(1,t))dx - xP_j(1,t)dt].$$
\vglue.05in Dividing the above $1$-form by
$x^k$ we obtain:
$$(*)\ \ x^{-k}.\pi^*(\om)= (tP_k(1,t)-Q_k(1,t))dx + xP_k(1,t)dt + x.\alpha$$
where
$\alpha=\sum_{j=k+1}^{\infty}x^{j-k-1}.[(tP_j(1,t)-Q_j(1,t))dx+xP_j(1,t)dt]$.
\vglue.05in Set  $R(x,y)=yP_k(x,y)-xQ_k(x,y)$, in such a way that
$x^{-k}.\pi^*(\om)=R(1,t)dx +xP_k(1,t)dt+x.\alpha$.
Analogously, computing the expression
of $\pi^*(\om)$ in the chart $((s,y),V)$, we obtain:
$$(**)\ \ y^{-k}.\pi^*(\om)=R(s,1)dy - yQ_k(s,1)ds +y.\beta.$$
\vglue.05in The  polynomial $R(x,y)$ is the {\it tangent
cone}\index{tangent cone}  of $\om$. We have two cases to consider:

\noindent{\rm(a)} $R\equiv 0$. In this case,
we shall say that the singularity is {\it dicritical}\index{singularity! dicritical}.

\noindent{\rm(b)} $R\not\equiv 0$. In this case,
the singularity is non-dicritical. The tangent cone has then degree
$k+1$.

\vglue.05in

Let us take a closer look at the above cases.

\noindent{Case (a)}. In this case, the forms
in (*) and (**) are still divisible  by $x$ and $y$
respectively. Dividing  (*) by $x$ we obtain
\begin{align*}
\om_1 &=P_k(1,t)dt + \alpha \\
&=P_k(1,t)dt+(tP_{k+1}(1,t)-Q_{k+1}(1,t))dx+x.\alpha_1,
\end{align*}
and this form cannot be divided by
$x$, since  $P_k\not\equiv 0$.

The  foliation $\fa^*$
is then represented in this chart by
$\om_1$ and in the other chart
by the  form $\om_2$, obtained from the division
of (**) by $y$. Note that, at the points of the
divisor $(x=0)$, of the form
$(t_o,0)$ such that $P_k(1,t_o)\ne 0$, the leaves of $\fa^*$
are transversal to the  divisor. The points $(t_o,0)$ such that
$P_k(1,t_o)=0$ will be the singular points
of $\fa^*$, or tangency points of the leaves of $\fa^*$
with the divisor. \vglue.05in

Note also that each  leaf transversal to the
divisor, will originate a  local separatrix  of $\fa$ via
blow-down.
Therefore, {\em a dicritical singularity admits infinitely many
separatrices}. \vglue.05in

\noindent{Case (b)}. In this case the forms in (*)
and (**) cannot be divided anymore. Therefore
they already represent the foliation
$\fa^*$ in their respective charts.
In particular,
{\em  the divisor is invariant by $\fa^*$}. Moreover, the
singularities of $\fa^*$ in the  divisor, are the
points, of the
$(x,t)$ chart, of the form $(0,t_o)$ where $R(1,t_o)=0$, and also
the point $(0,0)$, of the second chart, if $0$ is a zero
of
$R(s,1)=0$. We also have that  $\fa^*$ has
$k+1$ singularities,
counted with  multiplicity, in the  divisor. \vglue.05in

Note that, if some of the singularities of
$\fa^*$ has some separatrix $S$, then $\pi(S)$
is a separatrix of $\fa$ in $0$. \vglue.05in

Let now  us consider the blow-up process at $0\in \mathbb{C}^2$,
consisting of a succession of $ n$ blow-ups, $\Pi\colon M \to
\mathbb{C}^2$, with divisor $D=\cup_{j=1}^{n}D_j$.

The above argumentation proves that
we can obtain  a foliation $\tilde{\fa}$, with
isolated singularities, such  that in
$M\setminus D \simeq \mathbb{C}^2\setminus
\{0\}$ coincides with $\fa$. We shall say that the divisor $D_j$ is
{\it non-dicritical}\index{divisor! non-dicritical},if it is
invariant by $\tilde{\fa}$. Otherwise, we shall say that the divisor
is {\it dicritical}\index{divisor! dicritical}. \vglue.05in

\begin{Example}{\rm  The singularity $xdy-ydx=0$, is not irreducible $(\la = 1 \in \bq_+)$ and the blow-up $y = tx$ gives:
\begin{align*}
xd(tx) &- txdx = 0\\
xtdx + &x^2dt - txdx = 0\\
x^2dt &= 0\\
dt &= 0.
\end{align*}
The blow-up is therefore a non-singular foliation on the surface
$\widetilde{\bc}_0^2$\,, transverse to $\bp$. }
\end{Example}

Later on we shall see what is understood as
``reduction  of a singularity
of a foliation".
Before we shall introduce some notions and notations.

Let $X$ be a holomorphic vector field defined in a
neighborhood of $0\in \mathbb{C}^2$ such that
$0$ is an isolated
singularity of $X$. Let $\lambda_1$ and
$\lambda_2$ be the
eigenvalues of  $DX(0)$.

\begin{Definition}
[simple singularity] \label{Definition:simplesingularity} {\rm We
say that $0$ is a {\it simple singularity}\index{singularity!
simple} of $X$, if:
\begin{itemize}
\item[{\rm(a)}] $\lambda_1\ne 0$ and
$\lambda_2=0$ (or vice-versa). In this case, we shall say that the
singularity is a {\it saddle-node}\index{singularity! saddle-node}.
\item[{\rm(b)}] $\lambda_1,\lambda_2\ne 0$ and
$\lambda_2/\lambda_1\not\in\mathbb Q_+$. The numbers
$\lambda_2/\lambda_1$ and $\lambda_1/\lambda_2$ are then called the
{\it characteristic numbers}\index{characteristic numbers of a
singularity} of the singularity.
\end{itemize}
}
\end{Definition} \vglue.05in

Note that the above  conditions are invariant under holomorphic
changes of coordinates and under multiplication of $X$
by a holomorphic function that does not vanish at
$0$. Thus, the above notion may be extended to
the isolated singularities of
holomorphic foliations on complex surfaces.

\vglue.05in

In few words, the theorem of reduction of singularities
assures that given a foliation $\fa$ with a finite number of
(isolated) singularities
on a complex surface $M$, there exists a (finite) blow-up process,
$\pi
\colon\tilde  M \to M$, such that the pull-back foliation
$\fa^*=\pi^*(\fa)$, which is biholomorphically equivalent
to  $\fa$
outside of the divisor of $\pi$, has only simple singularities.

Indeed, it is possible to say more.

\begin{Definition}\label{Definition:18capI}{\rm
We shall say that the blow-up process  is a {\it reduction of the
singularity}\index{reduction of a singularity} or {\it resolution of
the singularity}\index{resolution of a singularity}, if:

\begin{itemize}
\item[{\rm(i)}]
All singularities of  $\tilde{\fa}$ in $D$ are simple.
\item[{\rm(ii)}] A dicritical divisor $D_j$ contains
no singularities of $\tilde \fa$, and no tangency points of $\tilde
\fa$ with  $D_j$. \end{itemize}
}
\end{Definition}

\begin{Theorem} [Theorem of the reduction of singularities,
Seidenberg \cite{seidenberg}\index{Theorem! of reduction of singularities}
\index{Theorem! of resolution of singularities! for foliations}]
\label{Theorem:7.I} Every isolated singularity of
a  holomorphic foliation of a complex surface
admits a reduction by a blow-up process.
\end{Theorem}

\section{Irreducible singularities}
\label{section:irreduciblesingularities} The reduction of
singularities theorem above mentioned can be made more accurate.

\begin{Definition}{\rm A singularity of a holomorphic vector field $X$ in dimension two is
called {\it irreducible\/} if it belongs to one of the following
categories: (up to a change of coordinates)

\noindent (i)\,\, $X(x,y) = \la x(1+a(x,y))\,\dfrac{\po}{\po x} +
\mu y(1+b(x,y))\,\dfrac{\po}{\po y}$\, $\la/\mu \in
\bc\setminus\bq_+$\,, $a(x,y)$, $b(x,y)$ are holomorphic with
$a(0,0) = b(0,0) = 0$.

\noindent This will be called {\it non-degenerate irreducible case}.

\noindent (ii)\, $X(x,y) = \la\left(x^{k+1}\,\dfrac{\po}{\po x} +
[y(1+\mu x^k)+xb(x,y)]\,\dfrac{\po}{\po y}\right) $ where $\la,\mu
\in \bc$, $\la \ne 0$, $k \in \bn$, $b(x,y)$ is analytic of order
$\ge k+1$ at $0 \in \bc^2$.

\noindent This is the {\it saddle-node case}. A basic model for that
is the {\it formal normal form\/} presented below: }
\end{Definition}

\begin{Theorem} [Martinet-Ramis \cite{martinet-ramisselano}, Hukuara-Kimura-Matuda,
\cite{H-K-M}\index{Theorem! of Hukuara-Kimura-Matuda and
Martinet-Ramis}] A germ of a saddle-node  foliation singularity is
formally equivalent to an unique model $\fa_{\la,k}$ given by
$\begin{cases}\dot x = x^{k+1}&\\ \dot y = y(1+\la x^k)&\end{cases}$

\noindent $\la=$ residue of the saddle-node

\noindent $1+k=$ multiplicity of the saddle-node.
\end{Theorem}

\begin{Remark}{\rm  Not all saddle-node singularities are analytically conjugate to the formal
normal form. Indeed, the formal normal form admits two separatrices
(i.e., two invariant manifolds  of dimension one, through the
origin), this is not the case of the well-know {\it Euler
equation\/} $\begin{cases}\dot x = x^2&\\ \dot y = x+y&\end{cases}$
}\end{Remark}

\begin{Theorem} [Seidenberg 1968, \cite{seidenberg}\index{Theorem! of Seidenberg}
\index{Theorem! of reduction of singularities}]
\label{Theorem:seidenberg} Let $\fa$ be a holomorphic foliation with an isolated
singularity at $0 \in \bc^2$. There is a finite sequence of quadratic blow-ups $\pi(j)\colon M_j \to M_{j-1}$\,, $(j=1,\dots,\ell)$
such that $\pi(1)$ is the blow-up of $\bc^2$ at $0 \in \bc^2$ and $\pi(j)$ is
the blow-up of $M_{j-1}$ at some point $p_{j-1} \in M_{j-1}$\,, with the following properties:

\begin{itemize}
\item[{\rm(a)}]\, The pull-back foliation $\widetilde{\fa} := \pi^*(\fa)$, where
$\pi = \pi_\ell \circ\cdots\circ \pi_1\colon M_\ell \to \bc^2$,
(i.e., is a foliation  with only irreducible singularities, of one
of the two types):
\begin{itemize}
\item[{\rm(i)}] $xdy - \la y\,dx + hot = 0, \,\,\, \la \in \bc\setminus \bq_+
\,\, x^{k+1}dy - (y(1+\la x^k)+ \hot)dx = 0,$
\item[{\rm(ii)}] $ \la \in \bc,\quad k \in \bn$\,\,\, (saddle node), where $\hot$ stands for {\it higher order terms}.
\end{itemize}

\item[{\rm(b)}]\, The exceptional divisor $D = \pi^{-1}(0)$ is a connected union of embedded
projective lines $D = \bigcup\limits_{j=1}^\ell \bp_j$\,, without
triple points, transverse intersections and

\item[{\rm(c)}]\, A component $\bp$ of $D$ is either invariant by the foliation $\widetilde{\fa}$ or
it is transverse to $\widetilde{\fa}$ (without tangent points).
\end{itemize}
\end{Theorem}

The singularity $\fa$ is called {\it non-dicritical\/} if all components of $D$ are invariant by $\widetilde{\fa}$,
and it is called {\it dicritical\/} if $\widetilde{\fa}$ exhibits some {\it dicritical\/} component.

A characterization of dicritical singularities is now possible. This
will be done in the next section.

\section{Separatrices: dicricity and existence}
\label{section:separatrices} We begin with  a definition that comes
from the theory of (Real) ordinary differential equations.

\begin{Definition} [Separatrix\index{separatrix}] {\rm Given a (germ of a) holomorphic singularity $\fa$ at $0 \in \bc^2$,
a {\it separatrix\/} of $\fa$ is (a germ of) an irreducible analytic curve $\Ga \ni 0$ which
is invariant by $\fa$ (i.e., $\Ga\setminus\{0\}$ is contained in a leaf of $\fa$).
}
\end{Definition}

Since a separatrix $\Ga$ is a germ of an  irreducible analytic curve
ate $0 \in \bc^2$, the Newton-Puiseux parametrization theorem gives
a parametrization $\vr\colon (\bc,0) \to (\Ga,0)$ of type
$(t^n,t^m)$, $\langle n,m\rangle=1$ so that $(\Ga,0)$ is
homeomorphic to a disc $(\bd,0)$.

\begin{Example} [Holomorphic first integral] {\rm Given a foliation with an isolated singularity
$\fa$ at $0 \in \bc^2$ we say that a non-constant holomorphic
function $f\colon (\bc^2,0) \to (\bc,0)$ (defined in a neighborhood
of the origin $0 \in \bc^2$) is a {\it holomorphic first integral\/}
\index{holomorphic first integral} for $\fa$ if $f$ is constant
along the leaves of $\fa$. If  $\fa$ is given by the vector field
$X$ (with $\sing(X)=\{0\}$) then this is equivalent to $df(X) \equiv
0$. If we consider the dual $1$-form $\omega$ then this is
equivalent to $\om \wedge df \equiv 0$. In any case, $X$ is parallel
to the Hamiltonian vector field
$$
H_f := -f_y\, \frac{\po}{\po x} + f_x\,\frac{\po}{\po y}
$$
and $\omega$ is of the form $\om=g\cdot df = g(f_x\,dx + f_y\,dy)$
for some meromorphic function  $g$ ($g$ is holomorphic and
non-vanishing if $f$ is chosen to be {\it reduced}). We assume that
$f(0)=0$\,. The separatrices of $\fa$ are the branches of $\{f=0\}$.
It is well-know from the local theory of analytical functions
(\cite{gunning-rossi}) that $f$ can be written (in a small bidisc
centered at $0 \in \bc^2$) as $f = \prod\limits_{j=1}^r f_j^{n_j}$,
with $n_j \in \bn$,\, $f_j$ holomorphic and such that:

(a)\,\, $f_j=0)$ is irreducible

(b)\,\, $f_i$ and $f_j$ are relatively prime

\noindent in the local ring $O_2$ so that $(f_i=0) \cap (f_j=0) =
\{0\}$.

\noindent(*)\,\, If moreover $f$ has irreducible/connected fibers
then $\langle n_1,\dots,n_r\rangle = 1$.

The separatrices are then given by $(f_j=0)$, $j=1,\dots,r$.

The other leaves of $\fa$ are closed in a neighborhood of $0$ and do
not accumulate at $0$.  The foliation is therefore {\it with
finitely many separatrices}.}
\end{Example}

\begin{Example} [Meromorphic first integral] {\rm A natural extension of the above definition
gives as the notion of {\it meromorphic first
integral}\index{meromorphic first integral}, $f\colon M \to
\overline{\bc}$ of a singularity $\fa$ at $0 \in \bc^2$. Writing $f
= \dfrac gh$ for $g,h\colon(\bc^2,0) \to (\bc,0)$ holomorphic
functions which (in the non-trivial case) vanish at $0 \in \bc^2$ we
have that the leaves of $\fa$ are contained in the curves $ag+bh=0$
with $(a,b) \in \bc^2\setminus\{0\}$. In particular {\it an leaves
are contained in separatrices}.}
\end{Example}

\begin{Lemma} Let $\fa$ be a singularity at $0 \in \bc^2$. A leaf $L$ of $\fa$ is contained in
a separatrix if, and only if, $\overline{L}\setminus L=\{0\}$, i.e.,
$L$ accumulates only at the singular point.
\end{Lemma}

\begin{proof} Let $L$ be a leaf of $\fa$ such that $\overline{L}\setminus L=\{0\}$. By Remmert-Stein extension theorem
(Theorem~\ref{Theorem:Remmert-Stein}) the closure $\Ga :=
\overline{L} \subset \bc^2$ is an analytic subset of dimension $1$.
Since $\Ga$ is clearly $\fa$-invariant and irreducible
$(\Ga\setminus\{0\}=L$ is connected).
We conclude that $\Ga = \overline{L} = L \cup \{0\}$ is a separatrix of $\fa$. The converse is clear. \end{proof}

As a corollary of Seidenberg's theorem
(Theorem~\ref{Theorem:seidenberg}) we have.

\begin{Proposition} A foliation singularity $\fa$ at $0 \in \bc^2$ is dicritical if,
and only if, it exhibits infinitely many separatrices.
\end{Proposition}

\begin{proof} First we observe that after the reduction of singularities, a leaf $\widetilde{L}$
of $\widetilde{\fa}$ projects onto a (leaf contained in a)
separatrix of $\fa$ if, and only if,  $\widetilde{L}$ is not
contained in the exceptional divisor $D$.
On the other hand, $\fa$ is dicritical if and only if there is a component $\bp_j$ of $D$ for
which every leaf $\widetilde{L}$ intersecting $\bp_j$ is transverse to $D$ except may be for
those at the {\it corners\/} $\bp_j \cap \bp_i \ne \phi$.

Finally the invariant components of $D$ only originate finitely many
separatrices of $\fa$ (the only possibilities come from separatrices
of singularities which are not at corners, but these singularities
are irreducible and therefore exhibit at most two separatrices).
\end{proof}

As we have already mentioned, a separatrix $\Ga$ admits a
parametrization  $\vr\colon (\bc,0) \to (\Ga,0)$ of type
$(t^n,t^m)$, $\langle n,m\rangle=1$ so that $(\Ga,0)$ is
homeomorphic to a disc $(\bd,0)$.
In particular, the leaf $L = \Ga\setminus\{0\}$ has the topology of a punctured
disc $\bc^* = \bc\setminus\{0\}$. Its fundamental group is cyclic isomorphism to $\bz$, generated by a loop $\ga \simeq S^1$.
The (local) holonomy group of the leaf $L$ is then cyclic generated by a single diffeomorphism $f\colon (\bc,0) \to (\bc,0)$.

\begin{Remark}{\rm  It is here that we notice a drastic difference between the singular
and the non-singular case for foliations: In general it is not so
common to find leaves with fundamental group in the non-singular
case. On the other hand, these leaves are quite common in the
singular foliations framework; thanks to the following theorem.
}\end{Remark}

\begin{Theorem} [Separatrix theorem of Camacho-Sad, \cite{Camacho-Sad}\index{Theorem! existence of separatrix}\index{Theorem! of Camacho-Sad}]
\label{Theorem:separatrix} Every holomorphic foliation singularity
$\fa$ at $0 \in \bc^2$ admits some separatrix.
\end{Theorem}

\begin{Remark}{\rm The above result is typical from the holomorphic case since there are
examples in the real analytic case where the foliation/vector field
admits no separatrix:
take $X$ as the Hamiltonian of $f = x_1^2+x_2^2$ in the real plane $\mathbb R^2 \ni (x_1,x_2)$.
The orbits are concentric circles, no separatrix is allowed.
}\end{Remark}

The Separatrix theorem is a by product of a suitable strategy on the reduction theorem of
Seidenberg (organizing the reduction into {\it linear chains}) and a
residue theorem applied to an {\it index\/} defined
in association with a separatrix of a singularity.
Just to give a few more words about this
important theorem we have:

Let $\fa$ be a germ at $0 \in \bc^2$ of a singularity and $\Ga$ a
smooth separatrix of $\fa$.  In local coordinates we may assume that
$\Ga\colon (y = 0)$. Then we choose $\omega =  A(xy)dx + B(x,y)dy$ a
holomorphic $1$-form with $\sing \omega = \{0\}$, defining $\fa$.
Since $\Ga$ is $\fa$-invariant we can write $\omega =  y\,A_1(x,y)dx
+ B(x,y)dy$ with $A_1(x,y)$ holomorphic. Then we consider the
$1$-form $\eta := \dfrac{-A_1(x,0)}{B(x,0)}\,dx$. It is a
meromorphic $1$-form with no poles off $x=0$ (notice that $y \nmid
B(x,y)$ as a holomorphic function, otherwise $\om$ would have
non-isolated zeros). The {\it Index of $\fa$ relative to $\Ga$ at
$0$\/} is defined as $I(\fa,\Ga,0) := \Res \eta(x=0)$. The Index
admits a geometrical interpretation as follows:

\noindent Given $x$, the ``inclination" of the tangent space
$T_{(x,y)}\fa$  of the leaf $L_{x,y)}$ of $\fa$ through $(x,y)$ is
given by
$$
\te_x(y) = \frac{dy}{dx} = \frac{-y\,A_1(x,y)}{B(x,y)}\,\cdot
$$

The derivative of this function $\te_x$ of $y$ at $y=0$ is then
$\te_x'(0)dx = \dfrac{-A_1(x,0)}{B(x,0)}\,dx$.  The Index
$I(\ga,\Ga,0)$ is then the residue of this $1$-form at $x=0$.

The Camacho-Sad index theorem (\cite{Camacho-Sad}, \cite{LN}) states
that the sum of indexes of a foliation $\fa$  at all the
singularities in a compact analytic smooth invariant curve $\Ga$ on
a complex surface $M^2$ is equal to the self-intersection (first
Chern class) of $\Ga$ in $M$, does not depend therefore on the
foliation $\fa$.
\[
\sum\limits_{p\in\sing(\fa)\cap\Ga}\,\,I(\fa,\Ga_p,p) = \Ga\cdot\Ga
\in \bz.
\]

\begin{Exercise}{\rm  Compute the index of each separatrix in the following cases: linear case,
Poincaré-Dulac normal form and saddle-node case. }
\end{Exercise}

\section{Holonomy and analytic classification}
\label{section:holonomyclassification}

\subsection{Holonomy of irreducible singularities}

As we have already seen there is always a separatrix through a
holomorphic singularity in dimension two.  Such a leaf has a
holonomy map $f\colon \bc,0 \to \bc,0$ and we shall study this map
in some particular cases:

\begin{Example}[linear case] {\rm
$$
\begin{cases}
\dot x = \la x&\\
\dot y = \mu y&
\end{cases}
$$
We fix the separatrix $\Ga\colon (y=0)$.

\noindent Choose a transverse section $\Sigma\colon\{x=0\}$.

\begin{Remark} {\rm $\Sigma$ is a (complex) disc.}
\end{Remark}

We consider the loop $\ga(t) = (x_0e^u,0) \subset 1$; $t \in
[0,2\pi]$. Let $C$ be the product  $\ga \times \Sigma \simeq
S^1\times\bd$, it is a solid torus. In this solid torus $\fa$
induces a real flow given by the ordinary differential equation
\begin{align*}
\frac{dy}{dt} &= \frac{dy/dx}{dx/dt} = \frac{\dfrac{\mu y(t)}{\la x(t)}}{x'(t)} = \frac{\mu}{\la}\, y(t) \, \frac{x'(t)}{x(t)}\\
\frac{dy}{dt} &= \frac{\mu}{\la}\, y(t).i
\end{align*}

\noindent The solutions are $y(t) = y\cdot e^{ t i\mu/\la}$.
Therefore the first return (holonomy)  map is given by $f(y) =
e^{2\pi i\mu/\la}$. This is a linear map.

\noindent In general for a non-degenerate irreducible singularity
$\fa\colon \la x[1+a(x,y)]dy -$
$\mu y[1+b(x,y)]dx=0$ the holonomy map of the separatrix $\Ga\colon (y=0)$
is given by $f(y) = e^{2\pi i\mu/\la}\,y+hot$ in particular its linear part is
$$
f'(0) = e^{e\pi i\mu/\la}.
$$

}
\end{Example}

\begin{Example} [saddle-node normal form case] {\rm We consider a saddle-node in the normal form
$$
y^{k+1}\,dx - x(1+\la y^k)dy = 0.
$$
The strong manifold $\Ga\colon (y=0)$ has holonomy map $h(y)$ given
by a similar procedure.

\begin{align*}
x(t) &= x_0\,e^{it}\\
\frac{dy}{dt} &= \frac{dy/dx}{dx/dt} = \frac{y^{k+1}/x(1+\la
y^k)}{dx/dt} = \frac{i\, y^{k+1}(t)}{1+\la y^k(t)}
\end{align*}
$$
\begin{cases}
y'(t) = i\,\dfrac{y^{k+1}(t)}{1+\la y^k(t)}&\\
y(0) = y&
\end{cases}
$$
$f(y) = y(2\pi)$.

\noindent For instance if $k=1$ and $\la=0$, i.e., for the
saddle-node $\begin{cases} \dot y = y^2&\\ \dot x = x&\end{cases}$
we then have
\begin{align*}
y'(t) &= i\,y^2(t) \Rightarrow \frac{d}{dt}\left(\frac{-1}{y(t)}\right) = i\\
&\Rightarrow \frac{-1}{y(t)} = it + c \Rightarrow c = \frac{-1}{y(0)}\\
&\Rightarrow \frac{-1}{y(t)} = it - \frac{-1}{y(0)} \Rightarrow y(t) = \frac{y(0)}{1-ity(0)}\\
&\Rightarrow f(y) = y(2\pi) = \frac{y}{1-2\pi iy} \quad \text{this
is a homography}
\end{align*}
In general, for a general form saddle-node, the strong manifold
(given by $(y=0)$  in the form $y^{k+1}\,dx - [x(1+\la
y^k)+(\cdots)]dx=0$) the holonomy map of $\Ga$ is given by $f(y) = y
+ a_{k+1}\,y^{k+1} +\cdots$ where $a_{k+1} \ne 0$. It is a map
tangent to the identity.

}
\end{Example}

\begin{Exercise}{\rm  Calculate the holonomy map of the Poincaré-Dulac normal  form
$\fa\colon ydx - (nx+ay^n)dy=0$, $n \ge 2$; for the (only) separatrix $\Ga\colon (y=0)$.
}
\end{Exercise}

\subsection{Holonomy and analytic classification of irreducible singularities}

\vglue .1in

Let us now discuss on of the most important aspects of the concept
of holonomy for singularities of foliations.

It is well-known that for a pair of regular foliations a conjugation
between holonomy groups of diffeomorphic  leaves induces some
conjugation between the foliations in neighborhoods of the given
leaves.

\noindent This is not immediate in case of foliations with
singularities (how to extend the equivalence/conjugation to the
singularities?)

In the local framework we have for irreducible singularities a
precise answer to this question  thanks to the work of
Martinet-Ramis and some other authors.

\begin{Theorem} [Martinet-Ramis, 1983 \cite{martinet-ramisselano}] Let $\fa_1$, $\fa_2$ be two germs  of
saddle-node singularities at $0 \in \bc^2$. We assume that $(y=0)$
is the strong  separatrix of $\fa_1$ and $\fa_2$ and denote by
$f_1,f_2\colon \bc,0 \to \bc,0$ the holonomy map if $\Ga$ with
respect to $\fa_1$\,, $\fa_2$\,. Then there is a germ of a
holomorphic diffeomorphism $\Phi\colon \bc^2,0 \to \bc^2,0$ taking
the leaves of $\fa_1$ onto the leaves of $\fa_2$ (preserving off
course the  strong separatrix $\Ga\colon (y=0)$) if, and only if,
there  is a germ of a holomorphic diffeomorphism $\vr\colon \bc,0
\to \bc,0$ defining a conjugation $\vr \circ f_1 \circ \vr^{-1} =
f_2$ between $f_1$ and $f_2$\,.
\end{Theorem}
The same idea holds for non-degenerate irreducible singularities.
Thanks to Poincaré-Dulac theorem we only need  to consider
singularities in the Siegel domain:

\begin{Theorem} [Mattei-Moussu \cite{mattei-moussu}, Martinet-Ramis \cite{Martinet-Ramisresonant}]
\label{Theorem:analyticconjugation} Let $\fa_1$\,, $\fa_2$ be two
germs of non-degenerate singularities  $\fa_j\colon xdy-\la
y(1+b_j(x,y))dx=0$, with $b_j(x,y)$ holomorphic, $b_j(0,0)=0$, $\la
\in \re_-$\,. Denote by $f_j\colon \bc,0 \to \bc,0$ the holonomy map
of $\Ga\colon (y=0)$ with respect to $\fa_j$\,. Then $\fa_1$ and
$\fa_2$ are analytically conjugate by a holomorphic diffeomorphism
$\Phi\colon \bc^2,0 \to \bc^2,0$ if, and only if, the holonomy maps
$f_1$ and $f_2$ are analytically conjugate in $\Diff(\bc,0)$.
\end{Theorem}

In the Siegel non-degenerate case the idea of the proof is as
follows:

\noindent We choose vector field
$$
X_j(x,y) = x \,\frac{\po}{\po x} + \la y(1+b_j(x,y))\,\frac{\po}{\po
y}
$$
defining $\fa_j$ $(j=1,2)$ in a small bidisc $U_j \subset \bc^2$. We
fix a point  $(x_0,0) \in \Ga\setminus\{0\}$ close enough to origin,
and consider the holonomy maps $h_j\colon (\bc_y,0 \to \bc_y,0)$
defined by $\fa_j$ as the first return map of the (transversely
holomorphic) induced flow $\L_j := \fa_j\big\vert\_C$ where $C$ is
the solid torus $\ga \times \Sigma \cong S^1 \times \bd$ as above.

\noindent Because $\fa_j$ has a saddle like dynamics we know that
for $0 < R < |x_0|$  the leaves of $\fa_j$ intersect transversally
the solid torus $C_R$ given by $|z| = R$, $|y| < \delta$ where
$\delta > 0$ is small enough. \vglue.1in

\begin{Remark}{\rm In the abstract picture of the situation above described,
the leaves of $\fa_j$ are visualized as real curves but actually
they are complex curves, i.e., real surfaces of real dimension two.}
\end{Remark}

The analytic conjugation $\vr\colon (\Sigma,(x_0,0)) \to
(\Sigma(x_0,0))$ between $h_1$ and $h_2$ satisfies
\[
\vr \circ h_1 = h_2 \circ \vr
\]

\noindent First we extend $\vr$ to the solid torus $C = C'_{|x_0|}$
by setting
$$
\vr\big(y_1(t,y)\big) = y_2\big(t,\vr(y)\big)
$$
where $y_j(t,y)$ is the solution of the flow $\L_j =
\fa_j\big\vert\_C$ that starts from $y \in \Sigma$.

\noindent The above definition is consistent/valid because
$y_1(2\pi,y) = h_1(y)$ and
$$
y_2(2\pi,\vr(y)) = h_2(\vr(y))
$$
and by hypothesis we have
$$
\vr \circ h_1(y) = h_2 \circ \vr(y).
$$
Thus we have extended  $\vr$ to the solid torus $C$. Now we show how
to extend $\vr$ ``radially" to the torii $C_R$ with $0 < R < |x_0|$.

We consider the {\it induced radial flow\/} where we consider
$\widetilde{x}(t) = e^{-t}x$ and $\widetilde{y}_j(t,y)$ the solution
of
$$
\begin{cases}
\dfrac{d\widetilde{y}_j}{dt} = X_j(e^{-t}x, \widetilde{y}_j)&\\
\widetilde{y}_j(0) = y&
\end{cases}
$$

\noindent Then we extend $\vr\colon C \to C'$ to the ``interior" of
$C$ by setting
$$
\vr\big(\widetilde{y}_1(t,y)\big) =
\widetilde{y}_2\big(t,\vr(y)\big).
$$
Some estimative shows that
$$
\big\vert \widetilde{y}_1(t,y)\big\vert \le e^{\la.(1+\ve)t}\,|y|
$$
for some constant $0 < \ve \ll |\la|$. (Recall that (important!)
$\la < 0$). Similarly we also have
$$
\big\vert\widetilde{y}_2(t,y)\big\vert \le |y|.e^{|\la|(1+\ve)t}.
$$
This shows that we have
$$
|\widetilde{y}_2(t,y)| \le A.|\widetilde{y}_1(t,y)|
$$
for some constant $A > 0$.

Riemann extension theorem (Theorem~\ref{Theorem:riemannextension})
now shows that $\vr$ extends to the vertical axis $x=0$ for $|y| \le
\delta$ for a certain $0 < \delta$.

\vglue.1in

We shall say that a germ $\fa$ at $0 \in \bc$ is {\it analytically
linearizable\/}  if there is a holomorphic diffeomorphism
$\Phi\colon \bc^2,0 \to \bc^2,0$ taking the leaves of $\fa$ onto
leaves of a linear foliation $\fa_\la\colon xdy-\la y dx=0$, $\la
\in \bc\setminus\{0\}$. In this case $\fa_\la$ is unique and $\fa$
is of the form $\fa\colon xdy-\la y dx + hot = 0$. Similarly a germ
of a holomorphic diffeomorphism $f\colon \bc,0 \to \bc,0$ is {\it
analytically linearizable\/} if there is a germ of a holomorphic
diffeomorphism $\vr\colon \bc,0 \to \bc,0$  such that  that $\vr
\circ f = f_\la \circ \vr$ for some linear map $f_\la(y)=e^{2\pi i
\lambda}.y$.

Then, as a corollary of Theorem~\ref{Theorem:analyticconjugation}
above we have:

\begin{Theorem} {\it A germ of an irreducible non-degenerate singularity $\fa\colon xdy-\la ydx +\cdots = c$
is analytically linearizable if, and only if, its holonomy map of a given separatrix is analytically linearizable.}

\end{Theorem}

\chapter{Holomorphic first integrals}
\label{chapter:matteimoussu}

\section{Mattei-Moussu theorem}
\label{section:matteimoussu}

From the structural viewpoint the singularities with a holomorphic
first  integral are the most simple singularities of a holomorphic
foliation.

\begin{Definition} {\rm Given a singular foliation germ $\fa$ at $0 \in \bc^2$,
a holomorphic function $F\colon \bc^2,0 \to \bc,0$ (i.e., a
holomorphic function $F$  defined in some neighborhood $V$ of $0 \in
\bc^2$ and with $F(0) = 0 \in \bc$), is a {\it holomorphic first
integral\/} of $\fa$ if it is constant along the leaves of $\fa$. }
\end{Definition}

If $\fa$ is given by the vector field $X$ with an isolated
singularity at $0$ then the above condition is equivalent to $df(X)
\equiv 0$. In terms of the dual $1$-form $\omega$, the condition
becomes $\omega \wedge df \equiv 0$. This later, thanks to Saito's
division lemma (\cite{saito}), is equivalent to $\omega =  gdf$.

As we have already seen if $\fa$ admits a (non-constant) holomorphic
first integral $F\colon \bc^2,0 \to \bc,0$ then its leaves satisfy:
\begin{itemize}
\item[{\rm(i)}] the leaves of $\fa$ are closed outside of the origin
\item[{\rm(ii)}] there are only finitely many leaves that accumulate at the origin.
\end{itemize}

In terms of the language of Seidenberg's theorem
(Theorem~\ref{Theorem:seidenberg})  we have: \vglue.1in

\noindent\qquad $\fa$ admits a holomorphic first integral

\noindent\qquad\qquad\qquad\quad{$\Downarrow$}

\begin{itemize}
\item[{\rm(i)}] $\fa$ is non-dicritical.
\item[{\rm(ii)}] The leaves are closed off the singularity.
\end{itemize}

All this for a sufficiently small neighborhood of the singular
points $0 \in \bc^2$.

\begin{Remark}{\rm  Thanks to Remmert-Stein extension theorem (Theorem~\ref{Theorem:Remmert-Stein})
a leaf which is closed off $0 \in \bc^2$ and is not (contained in) a
separatrix, is contained in an analytic curve and is called {\it
analytic leaf}. }\end{Remark}

\noindent The theorem of Mattei-Moussu states a direct converse to
the above and can be stated in dimension two as follows:

\begin{Theorem} [Mattei-Moussu, \cite{mattei-moussu}\index{Theorem! of Mattei-Moussu}]
\label{Theorem:matteimoussu} Let $\fa$ be a holomorphic foliation
singularity at $0 \in \bc^2$.  Assume that for a small neighborhood
$V$ of $0 \in \bc^2$ we have:
\begin{itemize}
\item[{\rm(i)}] The leaves of $\fa$ in $V$ are closed in $V\setminus\{0\}\subset \bc^2$.
\item[{\rm(ii)}] Only a finite number of leaves of $\fa$ in $V$ accumulate at $0 \in \bc^2$.
\end{itemize}

\noindent Then $\fa$ admits a (non-constant) holomorphic first
integral $F\colon W \to \bc$ in  some open subset $0 \in W \subset
V$.
\end{Theorem}

The classical proof relies on the Reduction of Singularities, as
well as on the dynamics  of holomorphic diffeomorphisms $f\colon
\bc,0 \to \bc,0$.

R. Moussu has given an alternative proof based on the classical Reeb
local stability theorem (\cite{Godbillon}) and in the dynamics of
diffeomorphisms $f\colon \bc,0 \to \bc,0$ (\cite{moussu}).

The next example, due to M. Suzuki and Cerveau-Mattei,  shows that
there is no topological criteria for  the existence of a meromorphic
first integral.

\begin{Example}[Suzuki's example, \cite{Suzuki3}]
\label{Example:Suzuki}  {\rm Consider the {\it germ of singular
foliation}\index{germ of a singular foliation} $\underset{=}{\fa}$
at the origin  $0 \in{\mathbb C^2}$  given by: $\Omega = 0$ where
$\Omega = (y^3+y^2-xy)dx-(2xy^2+xy-x^2)dy.$ The germ
$\underset{=}{\fa}$ has the Liouvillian first integral $f(x,y) =
\frac{x}{y} \exp[\frac{y(y+1)}{x}]$ and  the following remarkable
properties:

\begin{itemize}

\item[{\rm (i)}] $\underset{=}{\fa}$ is
$\underline{\mu-\text{simple}}$, that is, it is a dicritical germ
which is desingularized with only one blow-up and the resulting
foliation  has no singularities on the exceptional divisor, it is
transverse to this projective line everywhere except for  (a unique)
tangent point (see \cite{Klughertz}).
\end{itemize}
Therefore it follows that:
\begin{itemize}
\item[{\rm(i)'}] Every leaf of $\underset{=}{\fa}$ is a separatrix and
therefore is given by some equation $(f = 0)$ where $f \in \mathcal
O_2$.

\item[{\rm (ii)}] $\underset{=}{\fa}$ does not admit a meromorphic
first integral in any neighborhood of the origin $0\in{\mathbb C^2}$
(see \cite{Ce-Mt} for a proof, or follow our argumentation).
\end{itemize}

Performing a blow-up  $(y=tx)$ at the origin $0 \in \bc^2$ we obtain
the foliation
\[
\tilde{\fa}\colon t^3dx + (2xt^2+t-1)dt = 0
\]
given by the vector
field
\[
\dot x  = 2xt^2 + t-1, \quad \dot t = t^3.
\]

The initial
foliation has the Liouvillian first integral $f =
\frac{x}{y}\exp(\frac{y(y+1)}{x})$ and therefore the foliation above
has the Liouvillian first integral $f(x,t) =
\frac{1}{t}\,\,e^{t(xt+1)}$. Restricting this function to the
projective line $(x=0)$ we obtain $f(0,t) = \frac1t \,\,e^{\frac
1t}$ which is a Liouvillian function on $\ov{\mathbb C}$. The map
$\sigma\colon(\ov{\mathbb C},1) \to (\ov{\mathbb C},1)$ defined by
mapping the point $p\in(\ov{\mathbb C},1)$ onto the other
intersection point of the leaf $L_p$ of $\tilde{\fa}$ though $p$
with the projective line, is (because of the order-2 tangency) a
germ of involution on $(\ov{\mathbb C},1)$. This germ is given by
the relation $f(0,t)\circ\sigma = f(0,t)$, that is, $\frac1t\, e^t =
\frac{1}{\sigma(t)}\,\,e^{\sigma(t)}$. This defines $\sigma(t)$ as a
nonalgebraic Liouvillian  function on $\ov{\mathbb C}$ and according
to \cite{Klughertz} this is enough to conclude that
$\underset{=}{\fa}$ does not admit a nontrivial meromorphic first
integral.}
\end{Example}

\section{Groups of germs of holomorphic diffeomorphisms}
\label{section:groupsofgerms}

We shall start with some basic facts. We denote by $\Diff(\bc,0)$
the group of germs of  holomorphic diffeomorphisms $f\colon \bc,0
\to \bc,0$. Such a map germ has representatives given by maps
$f_V\colon V \to f_V(V)$, where $0 \in V \subset \bc$ is an open
set, and $f_V\colon V \to f_V(V)$ is a holomorphic diffeomorphism
with $f_V(0)=0$. It can be identified (the germ $f$) with a power
series $f(z) = f'(0)z + \sum\limits_{j\ge2} a_jz^j \in \bc\{z\}$
where $f'(0) \ne 0$.

\begin{Lemma}
\label{Lemma:finitelinearizable}{\it Let $G \subset \Diff(\bc,0)$ be a finite subgroup.
Then $G$ is cyclic and analytically conjugate to a group
$\big\langle z \mapsto e^{\frac{2\pi i}{\nu}}z\big\rangle$ $\nu \in \bn$.}
\end{Lemma}

\begin{proof} We define a map $\Phi\colon \bc,0 \to \bc,0$ by $\Phi(z) = \sum\limits_{g\in G} \dfrac{g(z)}{g'(0)}\,\cdot$
This is a well-defined holomorphic map because $G$ is finite (in
particular all elements of $G$  have a common definition domain
around $0 \in \bc$). Moreover $\Phi'(0)=|G| \ne 0$ so that $\Phi \in
\Diff(\bc,0)$. Given now any element $g_0 \in G$ we have
\begin{align*}
\Phi(g_0(z)) &= \sum_{g\in G} \frac{g(g_0(z))}{g'(0)} = g_0'(0) \sum_{g\in G} \frac{(g\circ g_0)(z)}{g'(0)g_0'(0)}\\
&= g_0'(0) \sum_{g\in G} \frac{(g\circ g_0)(z)}{(g\circ g_0)(0)} =
g_0'(0) \Phi(z).
\end{align*}
Thus $\Phi \in \Diff(\bc,0)$ is an analytic conjugation of $G$ with
a finite subgroup of the  linear group $\GL(1,\bc) = \bc^*$; this
ends the proof.\end{proof}

Given a holomorphic function $\ell(z)$ defined in a neighborhood of
$0 \in \bc$ we consider the {\it invariance group\/} of $\ell(z)$ as
$\Inv(\ell) = \{g \in \Diff(\bc,0); \, \ell \circ g=\ell\}$  in
terms of germs. Assume that $\ell$ is not constant, $\ell(0) = 0$.
Since up to a change of coordinates we have $\ell(z) = z^\nu$ for
some $\nu \in \bn$ we conclude that:

\begin{Lemma}
\label{Lemma:invariancefinite} The invariance group $\Inv(\ell)$ is
a finite cyclic group.
\end{Lemma}
By the two above lemmas we have that the finite subgroups of
$\Diff(\bc,0)$ are the  invariance groups of holomorphic functions
$\ell\colon \bc,0 \to \bc,0$.

A final simple remark concerning the finiteness of subgroups of
$\Diff(\bc,0)$ is:

\begin{Lemma}
\label{Lemma:periodicgroupfinite}{\it Let $G \in \Diff(\bc,0)$ be a {\rm finitely generated}
subgroup such that each element $g \in G$ is periodic $(\exists\, n_g \in \bn$\, {\rm  such that } $g^{n_g} = \Id)$. Then $G$ is finite.}
\end{Lemma}

\begin{proof} First we observe that a non-trivial flat element $g(z) = z+a_{k+1}\,z^{k+1} +\cdots$ \,\, $(a_{k+1}\ne0)$
is not periodic; indeed $g^n(z) = z + na_{k+1}\,z^{k+1}+\cdots$\,.
The commutator of two elements $g_1,g_2 \in G$ is $[g_1,g_2] =
g_1\,g_2\,g_1^{-1}\,g_2^{-1}$  is a map tangent to the identity.
Thus by hypothesis this is the identity and $G$ is abelian. This
implies that $G$ is finite because it is finitely generated.
\end{proof}

Now we get to the main point in our argumentation:

\begin{Proposition} [Finiteness condition]
\label{Proposition:finitenesscondition} {\it A germ of a holomorphic
diffeomorphism $f\colon \bc,0 \to \bc,0$ has finite order $($ i.e.,
$f^n = \Id$  for some $n \in \bn)$ if, and only if, for some
neighborhood $V$ of $0$ all the orbits of $f$ in $V$ are finite.}
\end{Proposition}

The classical proof is as follows:

\begin{Lemma}
\label{Lemma:boundarypoint}{\it Let $K \subset \mathbb R^n$ be a
compact connected neighborhood of $0 \in \mathbb R^n$ and $h\colon K
\to h(K) \in \mathbb R^n$ a diffeomorphism with $h(0)=0$. Then there
exists a boundary point $x \in \po K$ such that the number of
iterates of $x$ by $h$ contained in $K$ is infinite.}
\end{Lemma}
In general, for $x \in K$ we define $\mu_K(x) := \{n; h^n(x) \in
K\}$.

\begin{proof} Assume by contradiction
that $\mu_K$ is bounded in $\po K$,  i.e., $\mu_K\big\vert_{\po K} <
N < \infty$ for some $N \in \bn$. Let us consider
\begin{align*}
A &= \{x \in K; \mu_K(x) < N\} \supset \po K\\
B &= \{x \in \Tni(K); \mu_{\overset{\circ}{K}}(x) \ge N\}, \,\,\,
\overset{\circ}{K} = \Tni(K)
\end{align*}
Then $A$, $B$ are open subset of $K$, $A \supset \po K$, $0 \in B$
and $A \cap B = \phi$.  Therefore, since $K$ is connected,
$\exists\, x_0 \in K$  such that  $x_0 \notin A \cup B$, i.e.,
$\mu_K(x_0) \ge N > \mu_{\overset{\circ}{K}}(x_0)$. Thus there
exists $n_0$ such that $y_0 = h^{n_0}(x_0) \in \po K$. Then we have
$\mu_K(y_0) = \mu_K(x_0) \ge N$; contradiction. This proves
Lemma~\ref{Lemma:boundarypoint}.\end{proof}

Our strategy to prove the Finiteness  condition proposition is to
prove:

\noindent{\it If $h \in \Diff(\bc,0)$ is  not periodic then exists
fundamental system of neighborhoods $\mathcal U$ of $0 \in \bc$ such
that for each neighborhood $U \in \U$ the set of points $x \in U$
where the $U$-orbit is not finite is not countable and contains the
origin $0 \in \bc$ in its closure.}

We start by fixing a compact disc $0 \in D \subset \bc$ such that
$h\colon D \to h(D)$  is a (holomorphic) diffeomorphism. Now we
consider the 3 following sets:
\begin{itemize}
\item $P \subset D$ is the set of periodic points of $h\big\vert_D$\,.
\item $F \subset D$ is the set of non-periodic points of finite $D$-orbit.
\item $I \subset D$ is the set of points with infinite $D$-orbit.
\end{itemize}

Now we define a sequence of compact subsets $A_n \subset D$ as
follows:
\begin{align*}
A_0 &:= D\\
A_1 &:= D \cap h^{-1}(D)\\
\vdots\\
A_n &:= D \cap h^{-1}(D) \cap\dots\cap h^{-n}(D)\\
\vdots
\end{align*}
Also let $C_n$ be the connected component of $A_n$ that contains the
origin $0 \in \bc$ and  $C := \bigcap\limits_{n\ge 0} C_n$\,. Then
by construction a point in $C$ is periodic or has infinite
$D$-orbit, i.e., $C \subset I \cup P$.

\begin{Claim}
\label{Claim:notcountable} If $C$ is countable then $I$ is not
countable.
\end{Claim}

\begin{proof}[Proof of Claim~\ref{Claim:notcountable}] By hypothesis $C$ is countable.
Therefore there is a subdisc $D_r \subset D$ of radius $0 < r <$
radius of $D$;  such that $C \cap \po D_r = \phi$ (otherwise $C$ is
not countable).

\noindent Therefore $\exists\, m \in \bn$ such that $C_m \cap \po
D_r = \phi$.

Let $K$ be a compact, connected neighborhood of $C_m$ that does not
intersect the  other components of $A_m$\,. In particular we have
$\phi = \po K \cap A_m = \po K \cap D \cap h^{-1}(D) \cap\dots\cap
h^{-m}(D)$ so that there must exist for every $x \in \po K$ a $m \in
\bn$  such that  $\bn \ni p < m$ and such that $h^p(x) \notin D$.
Consequently we get $\underline{P \cap \po K = \phi}$.

Now we denote the sets

\noindent $\widetilde{P} :=$ $h$-periodic points in $K$

\noindent $\widetilde{F} :=$ non periodic points but with finite
$K$-orbit

\noindent $\widetilde{I} :=$ points with infinite $K$-orbit.

We stratify $\widetilde{P}$ as $\widetilde{P} =
\bigcup\limits_{n\ge0} \widetilde{P}_n$  where $\widetilde{P}_n
:=\{$ points periodic of period $n \}$.

We observe that the {\it boundary $\po \widetilde{P}_n$ is finite:}
indeed, if $x_0 \in \po \widetilde{P}_n$  is an accumulation point
of points $\po \widetilde{P}_n$ there exists a neighborhood $V$ of
$x_0$ such that $h^n(x) = x$, $\forall\, x \in V$ because $V \cap
\widetilde{P}_n \ne \phi$ (Identity Principle). Since the orbit of
$x_0$ does not intersect $\po K$ we have, for $V$ sufficiently
small, that $h^j(V) \subset \Tni(K)$\, $\forall\,j=0,1,\dots,n-1$.
Thus $V \subset \Tni(\widetilde{P}_n)$, i.e., $x_0 \in
\Tni(\widetilde{P}_n)$. This actually means that points in $\po
\widetilde{P}_n$ are isolated $(\Rightarrow \po \widetilde{P}_n$ is
finite).

\noindent As a consequence {\it the boundary of $\widetilde{P}$ is
countable}.

The compact set $K$ can be decomposed as $K = \Tni(\widetilde{P})
\cup \widetilde{F} \cup (\widetilde{I} \cup \po\widetilde{P})$.
Notice that $\Tni(\widetilde{P})$ is an open subset of $\bc$ since
by hypothesis $\widetilde{P} \cap \po K = \phi$. By its turn,
$\widetilde{F}$ is open in $K$ and may intersect $\po K$. If
$\widetilde{F} = \phi$ then $\po K \subset \widetilde{I} \cup
\po\widetilde{P}$ and since $\po\widetilde{P}$ is countable (as we
have just seen above) and $\po K$ is not countable it follows that
$\widetilde{I}$ {\it is not countable}. If $\widetilde{F} \ne \phi$
and $\Tni(\widetilde{P}) = \phi$ we conclude that for every $r > 0$
small enough, $\widetilde{I} \cup \po \widetilde{P}$ intersects the
boundary $\po D_r$ so that $\widetilde{I}$ is not countable. If now
$\Tni(\widetilde{P}) \ne \phi$ and $\widetilde{F} \ne \phi$ then,
because these are disjoint open subsets, the set
$K-\Tni(\widetilde{P}) \cup \widetilde{F}$ is not countable. Since
$\po\widetilde{P}$ is a countable set this implies that
$\widetilde{I}$ is not countable; proving
Claim~\ref{Claim:notcountable}. \end{proof}

\begin{Lemma}
\label{Lemma:notcountable}{\it Suppose that $C$ is not countable.
Then $I$ is not countable.}
\end{Lemma}

\begin{proof} By Lemma~\ref{Lemma:boundarypoint} for every compact disc $0 \in D_r \subset D$
we have $(P \cup I) \cap \po D_r \ne \phi$. Therefore $P \cup I$ is
not countable.  Suppose by contradict that $I$ is countable; then
$P$ is not countable. Since $C = (C \cap I) \cup (C \cap P)$ (recall
that $C \subset I \cup P$) we have that $C \cap P$ is not countable
(otherwise $C$ would be countable, contradiction). Let us write
$$
C \cap P = \bigcup_{n\ge0} P_n
$$
where $P_n=$ points periodic of period $n$ in $P \cap C$. Then there
exists $n_0$ such that $P_{n_0}$  is not finite with an accumulation
point in $C_{n_0}$\,.  Since $h^{n_0}$ is holomorphic in an open
neighborhood of $C_{n_0}$ we conclude from the Identity Principle
that $h^{n_0} \equiv \Id$, contradicting the non-periodicity of $h$.
\end{proof}

All together, Claim~\ref{Claim:notcountable} and
Lemma~\ref{Lemma:notcountable} prove  that the set $I$ of points $x
\in D$ with infinite $D$-orbit is not countable. This proves
Proposition~\ref{Proposition:finitenesscondition}.

\section[Irreducible singularities]{Irreducible singularities}
\label{section:irreduciblematteimoussu}

We shall now address the irreducible case proving Mattei-Moussu in
this situation:

\begin{Lemma}
\label{Lemma:irreduciblesingintegrable}{\it Let $\fa$ be a germ of an irreducible
singularity at $0 \in \bc^2$. Assume that $\fa$ has closed leaves off the origin
$0 \in \bc^2$. Then $\fa$ admits a holomorphic first integral. Indeed, $\fa$ is
analytically conjugate to $nxdy+mydx=0$, for some $n,m \in \bn$; in a neighborhood of $0 \in \bc^2$.}
\end{Lemma}

\begin{proof} Since by hypothesis $\fa$ is irreducible, we divide the proof in two cases:

\noindent\textbf{Case 1:} $\fa$ is a saddle-node singularity germ.

\noindent In this case we have the strong manifold say $\fa\colon
y^{k+1}\,dx - [x(1+\la y^k)+\cdots]dy=0$ $\Ga\colon (y=0)$ The
holonomy map by of the strong manifold is of the form
$$
h_\ga(y) = y + a_{k+1}\,y^{k+1} +\cdots\quad a_{k+1} \ne 0.
$$
Therefore $h_\ga$ is not periodic. This implies that the orbits of
$h_\ga$ are not all of them closed and thus $\fa$ has some
non-analytic leaves on small neighborhoods of the origin.  Thus,
this case cannot occur.

\begin{Remark}{\rm  Indeed it is possible to say much more about the dynamics of $h_\ga$ is
no orbit is finite except for the fixed point. }\end{Remark}

\noindent\textbf{Case 2:} (non-degenerate case).

\noindent We write $\fa$ as $xdy - \la y dx + hot = 0$\,\, $\la \in
\bn\setminus \bq_+$ with invariant axes.

\noindent Again we consider the holonomy of $\Ga\colon (y=0)$. It is
a map of the form $h_\ga(y) = e^{2\pi i\la}\,y +\cdots$\, If $\la
\notin \re$ then Poincaré Linearization theorem implies that $\fa$
is analytically linearizable and therefore $h_\ga(y)$ is conjugate
to $y \mapsto e^{2\pi i\la}\,y$. This later map only has finite
orbits when $\la \in \bq$ so we would have $\la \in \bq_+$
contradiction. Therefore $\la \in \re_-$ and indeed because $h_\ga$
is periodic we conclude that it is linearizable and the same holds
for $\fa$; moreover (still because $h_\ga$ is periodic) we must have
$\la \in \bq_-$ say, $\la = -n/m$, \, $n,m \in \bn$ and $\fa$ is
analytically linearizable as $nxdy + mydx = 0$. \end{proof}

\section{The case of a single blow-up}
\label{section:oneblowupmatteimoussu}

In order to illustrate the main ideas in the proof we consider the
following situation: $\fa$ can be reduced with a single blow-up
$\pi\colon \widetilde{\mathbb C}_0^2 \to \bc^2$. Since $\fa$ is
non-dicritical the exceptional divisor is
$\widetilde{\fa}$-invariant, $\widetilde{\fa} = \pi^*(\fa)$. Because
the leaves of $\fa$ are closed off the origin, the non-separatrices
are analytic leaves and the same holds for the leaves of $\fa$ which
are contained in $\bp$ or the separatrices of $\widetilde{\fa}$
transverse to $\bp$.

\begin{Remark}{\rm  $\bp\setminus\sing(\fa)$ is a leaf of $\fa$.
}\end{Remark} Write $\Sing(\widetilde{\fa}) =
\{\widetilde{p}_1,\dots,\widetilde{p}_r\} \subset \bp$. From the
irreducible case above we conclude that for each $j \in
\{1,\dots,r\}$ there is a neighborhood $\widetilde{V}_j$ of
$\widetilde{p}_j$ in $\widetilde{\bc}_0^2$ where $\widetilde{\fa}$
admits a (non-constant) holomorphic first integral say $F_j\colon
\widetilde V_j \to \bc$. We may assume that $\widetilde V_j$ is a
product $\widetilde V_j = D_j \times \bd_\ve$ of a disc
$\widetilde{p}_j \in D_j \subset \bp$ and a small disc $0 \in
\mathbb D_\ve$ of radius $\ve > 0$. Assume also that $D_i \cap D_j =
\phi$,\, $\forall\, i \ne j$. \noindent Fix now a point
$\widetilde{p}_0 \in \bp\setminus\bigcup\limits_{j=1}^r D_j$\,.
Since $\bp$ is a 2-sphere we may choose a simply-connected domain
$A_j \subset \bp$ such that
$$
A_j \cap \{\widetilde{p}_0,\widetilde{p}_1,\dots,\widetilde{p}_r\} =
\{\widetilde{p}_0,\widetilde{p}_j\}.
$$
\begin{Claim} {\it We may extend $F_j$ to a holomorphic first
integral $\widetilde{F}_j$ for $\widetilde{\fa}$ in a neighborhood
$U_j$ of $D_j \cup A_j$\,.} \end{Claim}

\begin{proof} Indeed this is the classical  {\it holonomy
extension\/} which is possible because $A_j$ is simply-connected
therefore it has no further holonomy that the one already in
$D_j$\,. Let us be more precise:
\noindent Fix a point $\widetilde{a}_j \in
(D_j\setminus\{\widetilde{p}_j\}) \cap A_j$ and a transverse disc
$\Sigma _{\widetilde{q}_j}$ to $\widetilde{\fa}$ with
$\{\widetilde{q}_j\} = \Sigma_{\widetilde{q}_j} \cap\bp$.

\noindent We consider a simple path $\delta_j\colon [0,1] \to t_j$
with $\delta_j(0) = \widetilde{p}_0$\,. $\delta_j(1) =
\widetilde{q}_j$\,.

We consider the holonomy map $h_{\delta_j}\colon \Sigma \to
\Sigma\widetilde{q}_j$ of the path $\delta_j$ in the leaf\linebreak
$\bp\setminus\{\widetilde{p}_1,\dots,\widetilde{p}_r\}$ of
$\widetilde{\fa}$.

\noindent Then we define $\widetilde{F}_j(y) :=
F_j(h_{\delta_j}(y))$ for $y \in \Sigma$ close enough to
$\widetilde{p}_0$\,.

\noindent This is well-defined because $A_j$ is simply-connected and
$F_j$ is already invariant by the local holonomy map associated to a
small loop $\ga_j \subset D_j\setminus\{\widetilde{p}_j\}$.

\noindent Thus we can construct the (holonomy)  extension
$\widetilde{F}_j$ as above. We may also that $V_j \subset U_j$\,.
\end{proof}

For each extension $\widetilde{F}_j$ we consider the restriction
$\widetilde{F}_j\big\vert_\Sigma$ and the {\it invariance group}
$$
\Inv\big(\widetilde{F}_j\big\vert_\Sigma\big) = \{h \in
\Diff(\Sigma,\widetilde{p}_0)\colon \widetilde{F}_j \circ h =
\widetilde{F}_k\}.
$$
Finally we consider the {\it global invariance group\/}
$\Inv(\widetilde{\fa},\Sigma) :=$ subgroup of
$\Diff(\Sigma,\widetilde{p}_0)$  generated by the invariance groups
$\Inv(\widetilde{F}_j\big\vert_\Sigma)$\, $j=1,\dots,r$.

\begin{Claim} $\Inv(\widetilde{\fa},\Sigma)$ is a finite group.
\end{Claim}

\begin{proof} Indeed, each map
in $h \in\Inv(\widetilde{F}_j\big\vert_\Sigma)$ takes leaves of
$\widetilde{\fa}$ into leaves of $\widetilde{\fa}$ so the same holds
for the maps $h \in \Inv(\widetilde{\fa},\Sigma)$. Since the leaves
of $\widetilde{\fa}$ are analytic curves (except for those contained
in separatrices) we conclude that each map in
$\Inv(\widetilde{\fa},\Sigma)$ has closed orbits. This implies that
each map in $\Inv(\widetilde{\fa},\Sigma)$ is periodic (Finiteness
condition Proposition\ref{Proposition:finitenesscondition}). Because
$\Inv(\widetilde{\fa},\Sigma)$ is finitely generated this implies
that $\Inv(\widetilde{\fa},\Sigma)$ is a finite group, indeed cyclic
generated by a periodic rotation $z \mapsto e^{\frac{e\pi
i}{\nu}}\,z$, $(\nu \in \bn)$.

\end{proof}

Thus there is a holomorphic function $F$ on $\Sigma$ such that every
map in $\Inv(\widetilde{\fa},\Sigma)$ leaves invariant the map $F$;
i.e., $F \circ h = F$, $\forall\, h \in
\Inv(\widetilde{\fa},\Sigma)$. Because $\Inv(\widetilde{\fa},\Sigma)
\supset \Inv(\widetilde{\fa}_j\big\vert_{\Sigma})$ this implies that
$F$ {\it is also constant along the levels of $\widetilde{F}_j$ in
$U_j$}.

Thus we have constructed a holomorphic first integral $F$ for
$\widetilde{\fa}$ in $\bigcup\limits_{j=1}^r U_j$\,.

Now the complementary part of $\bigcup\limits_{j=1}^r A_j$ in $\bp$
is topologically a disc, i.e., it is simply-connected.

This complement then has trivial holonomy and therefore $F$ admits a
holonomy and  therefore $F$ admits a holonomy extension as a
holomorphic first integral of $\widetilde{\fa}$ to a neighborhood
$\widetilde{U}$ of $\bp$ in $\widetilde{\bc}_0^2$\,. This projects
into a holomorphic first integral for $\fa$ in a neighborhood of  $0
\in \bc^2$.

\section{The general case}
\label{section:generalmateimoussu}

As a final word about the proof in the general case, where several
blow-ups may be needed to finish the reduction of singularities, we
have a similar argumentation as above by the Induction on the number
of blow-ups in the reduction of singularities.

\chapter{Dynamics of a local diffeomorphism}
\label{chapter:dynamics}

We consider $f$ a germ of holomorphic diffeomorphism at $0 \in \bc$
fixing $0$; say $f(z) = \la z +a_{k+1} z^{k+1} + \dots$, $k \geq 1,
\, \la \in \bc^*$. We shall describe its dynamics (i.e., the
dynamics of its pseudo-orbits) according to the multiplier
$\lambda=f^\prime(0)$. For a  quick and  effective reference in this
subject we recommend \cite{bracci}.

\section{Hyperbolic case}
\noindent (i)\, \textbf{Hyperbolic Case:}\quad $|\la| \ne 0,1$.

By Poincaré-K\"onigs theorem (1884) there is a unique holomorphic
diffeomorphism $\vr$ with $\vr'(0)=1$, that conjugates $f$ to the
linear map $z \mapsto \la z$.

Then we have if $|\la| < 1$ that $f^n(z) \to 0$, as $n \to \infty$,
for each $z \approx 0$, following a spiralling or linear path,
depending on whether $\la$ is complex on pure real.

\section{Parabolic case}

\noindent (ii) \textbf{Parabolic Case:}\quad $|\la|=1$, $\la^k=1$
for some $k \in \bn$.

\begin{Proposition} Let $h\colon (\bc,0) \to (\bc,0)$
be a germ of a holomorphic diffeomorphism tangent to the identity
$h(z) = z + \sum\limits_{j\ge 2} a_jz^j$, $a_2 \ne 0$. Then there
exist sectors $S^+$ and $S^-$ with vertex at $0 \in \bc$, angles
$\pi-\te_0$ (where $0 < \te_0 < \pi/2)$ and opposite bisectrices in
such a way that:
\begin{itemize}
\item[{\rm(i)}] $h(S^+) \subset S^+$, $\lim\limits_{h\to+\infty} h^n(z)=0$, $\forall\, z \in S^+$
\item[{\rm(ii)}] $h^{-1}(S^{-1}) \subset S^-$, $\lim\limits_{n\to+\infty} h^{-n}(z)=0$, $\forall\, z \in S^-$
\end{itemize}
\end{Proposition}
\begin{proof} A linear change of coordinates allows us to write
$h(z) = z-z^2 + \sum\limits_{j\ge3} a_jz^j$. For $|z| < \delta$
small enough we consider $g(z)$ defined by
$$
\begin{cases}
g(0) = 1&\\
h(z) = \dfrac{z}{1+zg(z)}&
\end{cases}
$$
Then $g(z)$ is analytic (for $|z| < \delta$ small enough) and we
consider the coordinate $\xi = \dfrac 1z$\,, $z \ne 0$. Under these
coordinates we have $h$ as $\widetilde{h}(\xi) = \xi + g(1/\xi)$,
$\widetilde{h}(\infty) = \infty$, a diffeomorphism defined in a
neighborhood of $\xi = \infty$. Because $g(0)=1$ we have $g(1/\xi) =
1 + \dfrac{1}{\xi}\,\widetilde{g}(\xi)$ with $\widetilde{g}$ bounded
as $\xi \to \infty$. Therefore $\widetilde{h}(\xi) = \xi + 1 +
\dfrac{1}{\xi}\,\widetilde{g}(\xi)$ is a map which is close to the
translation $T(\xi) = \xi+1$ for $|\xi|$ big enough.

Hence we can find a sector $S^+$ of horizontal
bisectrix and angle
$\pi-\te_0$ $(0 < \te_0 < \pi/2)$ and a disc $D^+ = (\xi_1-R)^2 +
\xi_2^2 \le (2R)^2$ so that if $\widehat{S}^+ := \widetilde{S}^+
\cap (CD^+)$ then $\widetilde{h}(\widehat{S}^+) \subset
\widehat{S}^+$ and $\lim\limits_{n\to\infty} \widetilde{h}^n(\xi) =
\infty$. Similarly we obtain $\widehat{S}^-$  such that
$\widetilde{h}^{-1}(\widehat{S}^-) \subset \widehat{S}^-$ and
$\lim\limits_{n\to\infty} \widetilde{h}^{-n}(\xi) = \infty$ on
$\widehat{S}^-$.

The sector $S^+$ and $S^-$ are then the images of $\widehat{S}^+$
and $\widehat{S}^-$ respectively, by the change of coordinates $z =
1/\xi$. \end{proof}

Similarly, for the case, $h\colon \bc,0 \to \bc,0$\,\, $h(z) = z +
\sum\limits_{j\le k+1}^\infty a_jz^j$, $a_{k+1} \ne 0$; we have
sectors $S_1^+,\dots,S_k^+$ and sectors $S_1^-,\dots,S_k^-$ so that
$S_j^+$ alternates with $S_j^-$ which alternates with $S^+_{j+1}$\,;
having angles $\ge \pi/k$, in such a way that
\begin{itemize}
\item[{\rm(a)}] $S_j^+ \cap S_j^{-1} \ne \phi$
\item[{\rm(b)}] $h(S_j^+) \subset S_j^+$\,, \, $h_j^{-1}(S_j^-) \subset S_j^-$ \newline
$\lim\limits_{n\to\infty} h^n(z)=0$, \, $\forall\,z \in
\bigcup\limits_{j=1}^k S_j^+$ \newline $\lim\limits_{n\to\infty}
h^{-n}(z) = 0$, \, $\forall\, z \in \bigcup\limits_{j=1}^k S_j^-$
\end{itemize}
\noindent Indeed, more can be said:

\begin{Theorem} [Camacho, \cite{camacho}\index{Theorem! of Camacho}]
Let $f(z) = \la z + O(z^2)$ be a holomorphic germ of a complex
diffeomorphism, $\la^n=1$ for some $n \in \bn$ (with $n$ minimal).
Then:
\begin{itemize}
\item[{\rm(i)}] Either $f^n(z) = z$,
\item[{\rm(ii)}] or there exists $k \in \bn$
such that $f$ is topologically conjugate
to $T_{k,\la,n}\colon z \mapsto \la z(1-z^{nk})$.
\end{itemize}
If moreover $f(z) = z + a_{k+1}\,z^{k+1} + O(z^{k+2})$\,\,
$a_{k+1}\ne0$, then $f$ is topologically conjugate to $T_k\colon z
\mapsto z+z^{k+1}$.
\end{Theorem}

And also:

\begin{Theorem} [Fatou-Leau flower theorem, \cite{abate}\index{Theorem! of Fatou-Leau}]

Let $f(z) = z + z^{k+1} + O(z^{k+2})$, $k \in \bn$. Then there exist
$k$ domains called {\it petals}, $P_j$\,, symmetric with respect to
the $k$ directions $\arg(z) = \dfrac{2\pi q}{k}$\,, $q=0,\dots,k-1$,
such that:
\begin{itemize} \item[{\rm(i)}] $P_j \cap P_k = \phi$ for $j \ne k$;
$0 \in \po P_j$ and each petal $P_j$ is holomorphic to the
right-half plane $\bh \subset \mathbb R^2$ \item[{\rm(ii)}] for each
$z \in P_j$ we have $f^m(z) \to 0$ as $m \to \infty$ moreover.
\item[{\rm(iii)}] For each $j$ the map $f\big\vert_{P_j}$ is
holomorphically conjugate to the parabolic automorphism $z \mapsto
z+i$ on $\bh$. \end{itemize}
If $f(z) = z+z^{k+1} + O(z^{k+2})$ then $f^{-1}(z) = z-z^{k+1} +
O(z^{k+2})$.
\end{Theorem}
Therefore we get, from Leau-Fatou theorem, $k$ attracting $Q_j$ for
$f^{-1}$ symmetric with respect to the $k$ directions $\arg(z) =
\dfrac{(zq+1)...}{k}$, $q=0,\dots,k-1$. These directions are the
bisectrices of the angles between two consecutive attracting
directions for $f$. The $Q_j$'s are repelling petals for $f$,
intersecting the $P_j$'s and $\bigcup\limits_j P_j \cup Q_j \cup
\{0\}$ is an open neighborhood of $0 \in \bc$ in $\bc$.  We have
therefore a pretty clear description of the dynamics of $f$.

\section{Elliptic case}

\noindent\textbf{Elliptic Case:}\quad $|\la|=1$, \,\, $\la = e^{2\pi
i\te}$ for some $\te \in \re\setminus\bq$.

The case is pretty rich and it is  a subject of deep research. The
main question first posed is whether $f = \la z+\dots$ is always
analytically linearizable. It was Cremer who first gave an example
of an elliptic   map which is not analytically linearizable. Indeed
Cremer introduced the following: \vglue.1in

\noindent\textbf{Cremer condition, \cite{cremer}:}\quad for $\te \in
\re\setminus\bq$ {\it if}
$$
\limsup_{n\to\infty} |\{n\te\}|^{-1/n}=\infty
$$
{\it then\/} there exist an elliptic germ $f(z) = e^{2\pi i\te}\,z +
O(z^2)$ which is not linearizable.

In the above statement, for a positive real number $x \in \mathbb
R$, we have
\[
\{x\} := x-[x]\quad [x] = \text{the integral part of\/}\, x.
\]

\noindent A number $\te$ satisfying the Cremer condition above is
called a {\it Cremer number\/}.  Cremer numbers form a dense subset
of $\re$ of zero Lebesgue measure \cite{bracci}.

On the other hand there are  arithmetical conditions by originally
by Siegel and recently J.C. Yoccoz and Bryuno (\cite{Brjuno}) for
assuring that if $\te \in \re\setminus\bq$ satisfies this
arithmetical condition then $f$ is always analytically linearizable.
These $\te$ form a full Lebesgue measure subset of $\re$.

Another remarkable result is:
\begin{Theorem}[Sigel-Bryuno-Yoccoz, \cite{perezmarco}]
Let $\te \in \re\setminus\bq$ and $\lambda = e^{2 \pi i \theta}$. If
the germ $p_\lambda(z)= \lambda z + z^2$ is analytically
linearizable then every germ $f\in \Diff(\mathbb C,0)$ with
$f^\prime(0)=\lambda$, is also analytically linearizable.

\end{Theorem}

Regarding the dynamics we have the  remarkable work of Pérez-Marco
(\cite{perezmarco,P6,P7}). He introduces the  following concept:

\begin{Definition} [Small cycles property]
{\rm A {\it small cycle\/} for $f$ is a finite  orbit of $f$ (a
subset $\{p_1,\dots,p_n\} \subset \bc\setminus\{0\}$ such that $p_i
\ne p_j$ and $f(p_i) = p_{j+1}$ (mod $n$)). We say that $f$ has the
{\it small cycles property\/} if for any open neighborhood $U$ of
$0$ then exists a small cycle for $f$ contained  in $U$. }
\end{Definition}

In this case the small cycles accumulate at $0$. In particular the
germ $f$ is not linearizable.

\begin{Theorem} [Pérez-Marco\index{Theorem! of Pérez-Marco}]

{\it There exist elliptic germs with the small cycles property. Not
all non-linearizable elliptic germs have the small cycles property.}
\end{Theorem}

\noindent Pérez-Marco gives an arithmetic condition on $\te$ in
order to decide whether the non-linearizable germ has the small
cycles property.

\noindent

Actually Pérez-Marco work goes much further with the introduction of
the {\it Hedge-Hogs}. He also concludes

\begin{Theorem} [Pérez-Marco]

{\it A non-linearizable elliptic map always has arbitrarily close to
the origin some orbit which accumulates at the origin.}
\end{Theorem}
Such an orbit cannot be closed.

For an elliptic germ $f\colon \bc,0 \to \bc,0$ we can add: Choose an
open connected subset $0 \in U \subset \bc$ where $f$ is univalued
the {\it stable set\/} of $f$ in $U$ is $K(U,f) =
\bigcap\limits_{j=0}^\infty f^{-j}(U)$.

\begin{Theorem} [Pérez-Marco] Let $f\colon \mathbb C,0 \to \mathbb
C,0$ be an elliptic map germ with stable set $K(f,U)$. Then:
\begin{itemize} \item[{\rm(i)}] $K(f,U)$ is compact, connected, {\it
full\/} (i.e., $U\setminus K(f,U)$ is connected), contains the
origin but is not restricted to the origin, i.e.,  $0 \in K(f,U) \ne
\{0\}$. Moreover, $K(f,U)$ is not locally connected at any point
distinct from the origin. \item[{\rm(ii)}] Any point of
$K(f,U)\setminus\{0\}$ is {\it recurrent\/}\index{recurrent point}
(that is, it is a limit point of its orbit). \item[{\rm(iii)}] There
is an orbit in $K(f,U)$ that accumulates at the origin.
\item[{\rm(iv)}] No non-trivial orbit converges to the origin.
\end{itemize}

\end{Theorem}
The stable set  $K(f,U)$ is called a {\it
hedge-hog}\index{hedge-hog}.

\chapter{Foliations on complex projective spaces}
\label{chapter:foliationsprojective}

\section{The complex projective plane and foliations}
\label{section:complexprojetiveplane}

The {\it complex projective plane\/} $\bc P(2)$ is the quotient
space $\bc^3\setminus 0$ by the equivalence relation $p, q \in
{\bc^3}\setminus 0$,\, $p \sim q \Leftrightarrow p= \la. q$ for some
$\la \in \bc\setminus\{0\}$. Thus $\bc P(2)$ is the space of lines
through the origin of $\bc^3$. By introducing homogeneous
coordinates $[(x_1;x_2;x_3)]$ on $\bc P(2)$ we conclude that $\bc
P(2)$ is equipped with an atlas $u$ consisting of three affine
charts $(x,y)$, $(u,v)$, $(r,s)$ with the following changes of
coordinates:

\centerline{$u=1/x, v=y/x$\,; \quad $r=1/y, s=x/y$\,; \quad $s =1/v$
\quad $r=u/v$}

\noindent A well-known fact  is then:

\begin{Proposition} The complex projective plane $\bc P(2)$ is a compact, connected and simply-connected complex surface.
\end{Proposition}

Let us now investigate the structure of the space of holomorphic
foliations with singularities (foliation of dimension one are the
ones who are interesting) on the complex projective plane. Recall
that a holomorphic foliation with singularities on $\bc P(2)$ is
given by an open cover $\bc P(2) = \bigcup\limits_j U_j$ such that
on each open subset $U_j$ the plaques of $\fa$ are given by a
holomorphic vector field $X_j$ in $U_j$ and if $U_i \cap U_j \ne
\phi$ then $X_i\big\vert_{U_i\cap U_j} = g_{ij}$
$X_j\big\vert_{U_i\cap U_j}$ for some non-vanishing holomorphic
function $g_{ij}\colon U_i \cap U_j \to \bc$. Finally we assume that
$\sing(X_j) = \phi$ or consists of a single point $p_j \in U_j$\,.
In particular we have the following example:

\begin{Example} [Polynomial Vector Fields on $\bc^2$]{\rm \, Let
$X(x,y) = P(x,y)\,\dfrac{\po}{\po x} + Q(x,y)\,\dfrac{\po}{\po y}$
be a polynomial vector field on $\bc^2$.\,\, $(P,Q \in \bc[x,y])$.
For our purposes we may assume that $P$, $Q$ have no common factor
on $\bc[x,y]$ and (equivalently) $(P=0) \cap (Q=0)$ is a finite
subset of $\bc^2$. Thus $X$ {\it has isolated singularities on\/}
$\bc^2$. Let us show that the foliation with singularities induced
by $X$ on $\bc^2$ extends to $\bc P(2)$. For this we consider the
change of coordinates $u=1/x$, $v=y/x$\,. We have:
$$
\begin{cases}
\dot x = P(x,y)&\\
\dot y = Q(x,y)&
\end{cases}
$$
Therefore $\dot u = \dfrac{-\dot x}{x^2} = -u^2.P\left(\dfrac 1u,
\dfrac vu\right)$.

\noindent Similarly\, $v = yu \Rightarrow \dot v = \dot y u + y\dot
u$.

\noindent Thus $\dot v = uQ\left(\dfrac 1u, \dfrac vu\right) +
\dfrac vu \cdot (-u^2)\cdot P\left(\dfrac 1u, \dfrac vu\right)$. Now
we write $P\left(\dfrac 1u, \dfrac vu\right) = \dfrac{1}{u^n} \cdot
\widetilde{P}(u,v)$ for a polynomial $\widetilde{P}(u,v)$ and some
$n \in \bn$, such that $u \nmid \widetilde{P}(u,v)$. Thus
$$
\dot u = \frac{-u^2}{u^n}\cdot \widetilde{P}(u,v) =
\frac{-\widetilde{P}(u,v)}{u^{n-2}}\,\cdot
$$
Similarly we write
$$
Q\left(\frac 1u, \frac vu\right) = \frac{1}{u^m}\,\widetilde{Q}(u,v)
$$
with $m \in \bn$ and $u  \nmid \widetilde{Q}(u,v) \in \bc[u,v]$.
Then
$$
\dot v = \frac{u\cdot\widetilde{Q}(u,v)}{u^m} - \frac{uv}{u^n}\,
\widetilde{P}(u,v).
$$
$$
\begin{cases}
\dot u = \dfrac{-\widetilde{P}(u,v)}{u^{n-2}}&\\
\dot v = \dfrac{\widetilde{Q}(u,v)}{u^{m-1}} -
\dfrac{v\widetilde{P}(u,v)}{u^{n-1}}\qquad
\end{cases}
$$

\noindent Thus we have a polynomial vector field
$\widetilde{X}(u,v)$ with isolated singularities in $\bc^2_{(u,v)}$
such that in the intersection of the spaces $\bc^2_{(u,v)}$ and
$\bc^2_{(u,v)}$ we have $X = \dfrac{1}{u^\ell}\,\widetilde{X}$ for
some $\ell \ge 0$.}
\end{Example}
\begin{Remark}{\rm  In any case we conclude that $\widetilde{X}$
defines a foliation on $\bc^2_{(u,v)}$ and this foliation coincides
with the one induced by $X$ on $\bc^2_{(x,y)}$\,. }\end{Remark}

\noindent Similarly for the coordinates $r = \dfrac 1y$\,, $s =
\dfrac xy$ we conclude that there is a polynomial vector field
$\overline{X}(r,s)$, with isolated singularities on the plane
$\bc^2_{(r,s)}$\,, such that in the intersection $\bc^2_{(x,y)} \cap
\bc^2_{(r,s)}$
$$
X = \frac{1}{r^k}\,\overline{X}
$$
for some $0 \le k \in \bn$.

\noindent Therefore, just applying the definition, we conclude that
$X$ defines a foliation $\fa(X)$ on the complex projective plane in
a natural way.

\noindent Summarizing the above example we have: {\it A polynomial
vector field $X$ on $\bc^2$  defines/induces a foliation $\fa(X)$ on
$\bc P(2)$.} Conversely we have:

\begin{Proposition} {\it Every holomorphic foliation with singularities on $\bc P(2)$ is the one induced
by  a certain polynomial vector field on $\bc^2$.}
\end{Proposition}

\begin{proof} Let $\fa$ be a foliation on $\bc P(2)$ with (finite) singular set $\sing(\fa) \subset \bc P(2)$.
Denote by $\pi\colon \bc^3\setminus\{0\} \to \bc P(2)$ the canonical
projection.

By definition we have $\bc P(2) = \bigcup\limits_j U_j$ a (finite
because $\bc P(2)$ is compact)  finite open cover where $\fa$ has
its plaques given by holomorphic vector fields $X_j$ on $U_j$ and
having isolated singularities. We can assume that the intersections
$U_i \cap U_j$\,\, $(i\ne j)$ contain no singularities of $\fa$ and
that if $U_i \cap U_j \ne \phi$ then $X_i = g_{ij}\,X_j$ with
$g_{ij} \in O^*(U_i \cap U_j)$.

\begin{Remark}{\rm  The idea if to ``lift" $\fa$ to $\bc^3\setminus\{0\}$  by
$\pi\colon \bc^3\setminus 0 \to \bc P(2)$ and then prove that $\fa$
(because it comes  from a foliation on $\bc P(2)$) can be given by a
polynomial system. Nevertheless, the lifted foliation $\pi^*\fa$
will have codimension one i.e., dimension two. So we must pass to
differential forms instead of vector fields. }\end{Remark} Let us
then choose dual $1$-forms $\om_j$ on the $U_j$ such that on each
$U_i \cap U_j \ne \phi$ we have  $\om_i = g_{ij}\,\om_j$ and
$\fa\big\vert_{U_j}$ is given by $\om_j=0$. Then we lift $\{\om_j\}$
and $\{g_{ij}\}$ by $\pi\colon \bc^3\setminus 0 \to \bc P(2)$
obtaining in this way $1$-forms $\widetilde{\om}_j$ in the open sets
$\widetilde{U}_j = \pi^{-1}(U_j) \subset \bc^3\setminus 0$ and such
on each intersection $\widetilde{U}_i \cap \widetilde U_j \ne \phi$
we have $\widetilde{\om}_i =
\widetilde{g_{ij}}\cdot\widetilde{\om_j}$\,. Notice that if $U_i
\cap U_j \cap U_k \ne \phi$ then $g_{ij}\cdot g_{ik} = g_{ik}$ on
$U_i \cap U_j \cap U_k$\,. Thus we have
$\widetilde{g_{ij}}\,\widetilde{g_{ij}}\,\widetilde{g_{jk}} =
\widetilde{g_{ik}}$ on $\widetilde{U_i} \cap \widetilde{U_j} \cap
\widetilde{U_k}$\,. Finally $\bc^3\setminus\{0\} = \bigcup\limits_j
\widetilde{U_j}$ so that the data $\big\{\widetilde{U_j},
\widetilde{g_{ij}}\big\}$ defines a Multiplicative Cocycle on
$\bc^3\setminus\{0\}$. Because $H^1(\bc^3\setminus 0)=0$ and
$\overset{\vee}{H^2}(\bc^3\setminus 0,\bz) = 0$ (Cartan's theorem,
\cite{Cartan}) the second (multiplicative) Cousin Problem has a
solution on $\bc^3\setminus\{0\}$ so that there are holomorphic
functions $\widetilde{g}_j\colon \widetilde{U}_j \to \bc^*$ with the
property that $\widetilde{g_{ij}} =
\dfrac{\widetilde{g_i}}{\widetilde{g_j}}$ on each $\widetilde{U_i}
\cap \widetilde{U_j} \ne \phi$ (\cite{Gunning 1}). Therefore we have
on each
$$
\widetilde{U_i} \cap \widetilde{U_j} \ne \phi\colon
\widetilde{\om_i} = \frac{\widetilde{g_i}}{\widetilde{g_j}}\,
\widetilde{\om_j} \Rightarrow \frac{1}{\widetilde{g_i}}\,
\widetilde{\om_i} = \frac{1}{\widetilde{g_j}}\, \widetilde{\om_j}\,.
$$
In this way we can define a holomorphic $1$-form $\widetilde{\om}$
on $\bc^3\setminus\{0\}$ by setting
$\widetilde{\om}\big\vert_{\widetilde{U_j}} :=
\dfrac{1}{\widetilde{g_j}}\, \widetilde{\om_j}$\,.  Thanks to
Hartogs' extension theorem $\widetilde{\om}$ extends as a
holomorphic $1$-form on $\bc^3$. Now we consider the pull-back
foliation $\widetilde{\fa} = \pi^*(\fa)$ (induced by the pull-back
of $\fa$) on $\bc^3\setminus\{0\}$.
\begin{Claim} The $1$-form $\widetilde\omega$ is integrable (i.e.,
$\widetilde\om \wedge d\widetilde\omega = 0$) and the foliation
induced by $\widetilde\omega$ coincides with $\widetilde{\fa}$.
\end{Claim}

This is quite clear since each $\om_j$ (and therefore each
$\dfrac{1}{\widetilde{g_j}}\,\widetilde{\om_j}$) is integrable (in
dimension two any $1$-form is integrable).

Now we write $\widetilde\om = \widetilde\om_\nu +
\widetilde\om_{\nu+1} +\cdots+ \widetilde\om_j +\cdots$ where
$\widetilde\om_j$ is a polynomial homogeneous $1$-form of degree $j$
(use Taylor/power Series).

Then $d\widetilde\om = d\widetilde\om_\nu + d\widetilde\om_{\nu+1}
+\cdots+ d\widetilde\om_j +\cdots$ and $0 = \widetilde\om \wedge
d\widetilde\omega =  (\widetilde\om_\nu +
\widetilde\om_{\nu+1}+\cdots) \wedge (d\widetilde\om_\nu +
d\widetilde\om_{\nu+1} +\cdots) = \widetilde\om_\nu \wedge
d\widetilde\om_\nu +\cdots$\,. Therefore $\widetilde\om_\nu \wedge
d\widetilde\om_\nu \equiv 0$, i.e., $\widetilde\om_\nu$ is
integrable. The details in the proof of following claim are then
lift to the reader:

\begin{Claim} The $1$-form $\widetilde\om_\nu$ defines the same foliation as $\widetilde\omega$ on $\bc^3$.
\end{Claim}

\begin{proof} Notice that $\widetilde\om_\nu$ is homogeneous so that it is {\it radially saturated}\index{foliation! radially saturated} (if $p \in \bc^3\setminus 0$ then the leaf of $\widetilde\om_\nu=0$ containing $p$, also contains the line $\{\la p; \la \in \bc^*\}$).

By its turn, the $1$-form $\widetilde\omega$ defines
$\widetilde{\fa}$, which is  the pull-back of a foliation on $\bc
P(2)$. Therefore, $\widetilde\omega$ is also radially saturated.
Thus the claim follows: given $p \in \bc^3\setminus\{0\}$ and
$\vec{v} \in \bc^3$ we have $\forall\,t \in \bc$
\begin{align*}
\widetilde\om(tp)\cdot\vec{v} &= \widetilde\om_\nu(tp)\cdot\vec{v} + \widetilde\om_{\nu+1}(tp)\cdot \vec{v} +\cdots\\
&= t^\nu\big[\widetilde\om_\nu(p)\cdot\vec{v} +
t\widetilde\om_{\nu+1}(p)\cdot\vec{v} +\cdots\big].
\end{align*}
Thus for $\vec{v} \in T_p(\widetilde\fa)$ we have
$$
0 = t^\nu\big[\widetilde\om_\nu(p)\cdot\vec{v} +
t\widetilde\om_{\nu+1}(p)\cdot\vec{v} +\cdots\big], \forall\,t
$$
so that $\widetilde\om_\nu(p)\cdot\vec{v}=0$ and hence  $\vec{v} \in
T_p(\{\widetilde\om_\nu=0\})$.  By comparing dimensions we conclude
the proof of the claim. \end{proof}

Since $\widetilde\om_\nu$ is homogeneous polynomial it induces a
polynomial $1$-form $\om(x,y) = P(x,y)dy - Q(x,y)dx$ on
$\bc^2_{(x,y)} \subset \bc P(2)$ so that
$\fa\big\vert_{\bc^2_{(x,y)}}$ is given by $\omega(x,y)=0$. This
ends the proof of the proposition.

\end{proof}

We shall adopt the following convention:

\noindent Given an algebraic (irreducible) (not necessarily  smooth)
curve $C \subset \bc P(2)$ given  on $\bc^2$ by an affine polynomial
equation $f(x,y)=0$ we consider its lift to $\bc^3\setminus\{0\}$
and then to $\bc^3$; denoted by $\widetilde{C} \subset \bc^3$. Then
$\widetilde{C}$ is an algebraic (not necessarily smooth)
hypersurface which has an homogeneous equation
$\widetilde{f}(x_1,x_2,x_3)=0$. We denote by $Z(\widetilde{f})$ the
curve $C$ on $\bc P(2)$ and by $(\widetilde{f}=0)$ the hypersurface
$\widetilde{C}$.

\section{The theorem of Darboux-Jouanolou}
\label{section:darboux}

Given a foliation $\fa$ on $\bc P(2)$ we may ask whether has leaf
which is a closed analytic subset of $\bc P(2)$.  In order to study
this question shall use:

\begin{Proposition} {\it Given a foliation $\fa$ on $\bc P(2)$ and a leaf $L \in \fa$ the following are equivalent:}
\begin{itemize}
\item[{\rm(i)}] {\it $L$ is contained in some algebraic curve $C \subset \bc P(2)$.}
\item[{\rm(ii)}] {\it $\overline{L}$ is an algebraic (invariant) curve $C \subset \bc P(2)$.}
\item[{\rm(iii)}]{\it $\overline{L}\setminus L \subset \sing(\fa)$, i.e., $L$ only accumulate at singular points of $\fa$.}
\item[{\rm(iv)}] {\it $\overline{L} = (\overline{L} \cap \sing(\fa)) \cup L$.}
\end{itemize}
\end{Proposition}

\noindent The above proposition is a consequence of the above
discussion and results. We then define a leaf $L \in \fa$ as an {\it
algebraic leaf\/}\index{leaf! algebraic} if $\overline{L} = C$ is an
algebraic curve (which is necessarily invariant by $\fa$).

\begin{Example}{\rm  Let $R = \dfrac PQ$ be a rational function on $\bc^2$;\, $P(x,y)$, $Q(x,y) \in \bc[x,y]$ have no common factor.

\noindent Then $R$ defines a foliation $\fa$ on $\bc P(2)$ whose
leaves are algebraic contained in the algebraic curves $aP+bQ=0$,
$(a,b) \in \bc^2 \ne 0$. In particular $\fa$ has infinitely many
algebraic leaves. }
\end{Example}
\begin{Theorem}[Theorem of Darboux-Jouanolou,
\cite{Jo}\index{Theorem! of Darboux-Jouanolou}] If a foliation $\fa$
on $\bc P(2)$ admits infinitely many algebraic leaves then $\fa$
admits a rational first integral. In particular, all leaves are
algebraic. \end{Theorem} \begin{proof} Choose a polynomial $1$-form
$\omega(x,y) = P(x,y)dy - Q(x,y)dx$ that defines $\fa$ on $\bc^2$,
with isolated singularities. Given an algebraic curve $C \subset \bc
P(2)$ with $C \cap \bc^2$ having irreducible polynomial equation
$f(x,y)=0$ we put $C^* = C\setminus(C \cap \sing(\fa))$.

\begin{Claim} \label{Claim:darboux}{$C^*$ is an algebraic leaf
of $\fa$ if, and only if, $\dfrac 1f\, \om \wedge df$ is a
polynomial $2$-form on $\bc^2$.}
\end{Claim}

\begin{proof}[Proof of
Claim~\ref{Claim:darboux}] Assume that $C^*$ is invariant by $\fa$.
Choose a point $p \in C^*$ and a local chart $(\widetilde x,
\widetilde y) \in \widetilde{U}$ centered at $p$ such that
$f(\widetilde x,\widetilde y) = \widetilde{y}$ and $C^* \cap
\widetilde{U}\colon (\widetilde y=0)$ and write $\om(\widetilde
x,\widetilde y) = B(\widetilde x,\widetilde y)d\widetilde x -
A(\widetilde x,\widetilde y) d\widetilde y$ where $A$, $B$ are
holomorphic near $(0,0)$. The vector field $\widetilde{X}(\widetilde
x,\widetilde y) := A(\widetilde x,\widetilde y)\,\dfrac{\po}{\po
\widetilde x} + B(\widetilde x,\widetilde y)\,\dfrac{\po}{\po
\widetilde y}$ then defines the foliation of $\widetilde{U}$.

\noindent Since $C^*\cap \widetilde U\colon (\widetilde y=0)$ is
$\widetilde X$-invariant we conclude that $B(\widetilde x,0) \equiv
0$.

\noindent Therefore we may assume that $\widetilde y\big\vert
B(\widetilde x,\widetilde y)$ in local ring of holomorphic
functions. Therefore $\dfrac 1f\,\om \wedge df =
\dfrac{1}{\widetilde y}\, B(\widetilde x,\widetilde y) d\widetilde x
\wedge d\widetilde y$ is holomorphic in $\widetilde{U}$. Thanks  to
Hartogs' extension theorem this shows that $\dfrac 1f\, \om \wedge
df$ is holomorphic in all points of $C \cap \bc^2$ and therefore it
is polynomial (we already have a priori that it is rational).

\noindent Assume now that $\dfrac 1f\,\om \wedge df$ is polynomial.
Then similarly to  above, $\dfrac{1}{\widetilde{y}}\,\om \wedge
d\widetilde y$ is holomorphic in $\widetilde{U}$ and therefore
$\widetilde y\big\vert B(\widetilde x,\widetilde y)$ so that
$B(\widetilde x,\widetilde y) = \widetilde y\cdot B_1(\widetilde
x,\widetilde y)$ for some holomorphic $B_1(\widetilde x,\widetilde
y)$. In particular $B(\widetilde x,0) \equiv 0$ and therefore
$\widetilde{X}$ is tangent to $(\widetilde y = 0)$ in
$\widetilde{U}$. Thus $C^* \cap \widetilde{U}$ is $\fa$-invariant
and by the Identity Principle or by Hartogs $C$ is an algebraic leaf
of $\fa$.\end{proof}

Notice that if we consider homogeneous coordinates $(x_1;x_2;x_3)$
on $\bc P(2)$ then as we have seen before, the pull-back foliation
$\widetilde\fa$ of $\fa$ to $\mathbb C^3$, can be defined by a
homogeneous polynomial $1$-form $\Om(x_1,x_2,x_3)$ of degree $\nu$;
this form satisfying $\Om\cdot\vec{R}\equiv 0$  ($\Om$ is radially
saturated) where $\vec{R} = x_1\,\dfrac{\po}{\po x_1} +
x_2\,\dfrac{\po}{\po x_2} + x_3\,\dfrac{\po}{\po x_3}$ is the radial
vector field. Given an algebraic curve $C \subset \bc P(2)$ we
consider an irreducible homogeneous polynomial $f(x_1,x_2,x_3)$ such
that $\{f=0\}$ is the homogeneous equation of $C = Z(f)$.

Let $k$ be the degree of the coefficients of the homogeneous
$1$-form $\Om$. Then the above claim rewrites as follows:

\begin{Claim} {\it $C$ is an algebraic invariant curve by $\fa$ if, and only if, there exists a
$2$-form $\te(x_1,x_2,x_3)$ such that: {\rm (i)} $df \wedge \Om =
f\te$\quad {\rm (ii)}\,  the coefficients of $\te$ are homogeneous
polynomial of degree $k-1$}.
\end{Claim}
Notice that the conclusion about the degree of (the coefficients of)
$\te$ is immediate since $\Om \wedge df$  is homogeneous of degree
$k+\deg f-1$ while $f\te$ is homogeneous of degree $\deg f + \deg
\te$.

\noindent Let now $E_k = \big\{\te; \te$  is $2$-form with
homogeneous coefficients of degree $k-$ on $\bc^3\big\}$.

\noindent Then $E$ is a finite dimension $\bc$-vector space of
finite dimension say $N(k) = \dim E_k$\,.  Assume that the foliation
$\widetilde{\fa}$ has $N(k)+1$ algebraic solutions given by
$(f_0=0),\dots,(f_{N(k)}=0)$ where $f_j(x_1,x_2,x_3)$ is homogeneous
of degree $k-1$, irreducible and $f_i$, $f_j$ are relatively prime
if $i \ne j$.

We write $\dfrac{df_j}{f_j} \wedge \Om = \te_j$\,, $j=0,\dots,N(k)$
as in the claim above.  Since $\dim E_k = N(k)$ the set
$\{\te_0,\dots,\te_{N(k)}\}$ is linearly dependent. There is
$(a_0,\dots,a_{N(k)}) \in$ $\bc^{N(k)+1}\setminus\}0\}$ such that
$\sum\limits_{j=0}^{N(k)} a_j\,\te_j = 0$ so that
$\left(\sum\limits_{j=0}^{N(k)}a_j\,\dfrac{df_j}{f_j}\right) \wedge
\Om = 0$ and then $\left(\prod\limits_{j=0}^{N(k)} f_j\right)\,\al
\wedge \Om =0$ for $\al := \sum\limits_{j=0}^{N(k)}
a_j\,\dfrac{df_j}{f_j}\,\cdot$ Because cod\,sing$(\Om) \ge 2$ we
must have $f_0\ldots f_N\, \cdot \al = g \Om$ for some homogeneous
polynomial $g(x_1,x_2,x_3)$ and some $N \le N(k)$ (given by the
number of non-zero coefficient $a_j$ in $\al$). Thus since $\al =
\sum\limits_{j=0}^N a_j\,\dfrac{df_j}{f_j}$ is closed, we conclude
that $\widetilde{\fa}$ is given by a closed rational $1$-form. If
$\widetilde{\fa}$ admits another algebraic solution $(f_{N(k)+1} =
0)$ then we write $df_{N(k)+1} \wedge \Om = f_{N(k)+1} \te$ and
because $\{\te,\te_1,\dots,\te_{N(k)}\}$ is linearly independent a
similar argumentation shows that for some
$$
\be := \sum_{j=1}^N b_j \, \frac{df_{i(j)}}{f_{i(j)}}
$$
with $f_{i(1)}\dots f_{i(N')}\be = h\Om$ where $h$ is a homogeneous
polynomial, $b_j \ne 0$, $i(j) \ne 0$, $j=1,\dots,N' \le N(k)$.

Then $\al = F\cdot \be$ where $F = \dfrac{gf_{i(1)}\dots
f_{i(N')}}{hf_0\dots f_N}\,\cdot$ Notice that $F$ is not constant
because $\Res_{(f_0=0)} \al = \al_0 \ne 0$ and $\Res_{(f_0=0)}\be =
0$. Since $\al$ and $\be$ are closed we have $0 = dF \wedge \be$ and
therefore $d F \wedge \Om = 0$. Thus $F$ is a first integral for
$\widetilde{\fa}$. The result follows.

\end{proof}

As a corollary of the above argumentation we have:

\begin{Corollary} For each $k \in \bn$ there exists $N(k) \in \bn$
such that  if a foliation on $\bc P(2)$ has more than $N(k)$
algebraic leaves then it has a rational meromorphic first integral.
\end{Corollary}

\begin{Remark}{\rm  Joaunolou-Darboux theorem is an algebraic parallel to Mattei-Moussu theorem.
}\end{Remark}

\section{Foliations given by closed $1$-forms}
\label{section:closedoneforms}

Since $\bc P(2)$ is compact and simply-connected it does not admit a
non-trivial closed  $1$-form which is holomorphic (indeed, such an
$1$-form $\omega$ would be exact $\om = dF$ for some holomorphic
function $F$ on $\bc P(2)$, but then $F$ must be constant and $\om
\equiv 0$ because $\bc P(2)$ is compact).

Nevertheless there are non-trivial closed meromorphic $1$-forms.
Since any meromorphic function on $\bc P(2)$ is already a rational
function (Liouville's theorem, \cite{gunning-rossi}) we consider the
class of closed rational $1$-forms on $\bc P(2)$ which we will
denote by $\Om(\bc P(2))$.

An interesting example is the Poincaré-Dulac normal form
$(nx+ay^n)dy - ydx=0$ \, $(n\ge2, a \in \bc^*)$  that defines a
foliation $\fa_{a,n}$ on $\bc P(2)$ that is also given by the closed
rational $1$-form $\Om_{a,n} \in \Om(\bc P(2))$ defined by
$\Om_{\,n} = \dfrac{(nx+ay^n)dy-ydx}{y^{n+1}}\,\cdot$ The poles of
$\Om_{a,n}$ in $\bc$ are given by $(y=0)$ which is the polar set of
$\Om_{a,n}$\,. This is a general fact as we shall below. Another
important example is the class of {\it linear logarithmic
foliations\/} (also called {\it Darboux foliations}) given by
$1$-forms as $\om = \big(\prod\limits_{j=1}^r
f_j\big)\cdot\big(\sum\limits_{j=1}^r \la_j
\,\dfrac{df_j}{f_j}\big)$ where $f_j \in \bc[x,y]$, $\la_j \in
\bc^*$. In this case the foliation is given by an element $\Om =
\sum\limits_{j=1}^r \la_j\,\dfrac{df_j}{f_j} \in \Om(\bc P(2))$.
Next we describe the structure of the elements in $\Om(\bc P(2))$.

\begin{Proposition} \label{Proposition:integrationlemma} Let $\om$
be a closed rational $1$-form on $\bc P(2)$ and let $\Om :=
\pi^*(\om)$ be its lift to $\bc^3$ where $\pi\colon \bc^3\setminus 0
\to \bc P(2)$ is the canonical projection. Then we have $$ \Om =
\sum_{j=1}^r \la_j\, \frac{df_j}{f_j} +
d\bigg(\frac{g}{f_1^{n_1-1}\dots f_r^{n_r-1}}\bigg) $$ where
\begin{itemize} \item[{\rm(a)}] $r \ge 2$ and $f_1,\dots,f_r \nmid
g$ are homogeneous polynomials in $\bc^3$; \item[{\rm(b)}]
$f_1,\dots,f_r$ are irreducible and pairwise relatively prime.
\item[{\rm(c)}] If $n_j > 1$ then $f_j \nmid g$ \item[{\rm(d)}]
$\deg(g) = \sum\limits_{j=1}^r \deg(f_j)\,(n_j-1)$, i.e., $\deg(g) =
\deg(f_1^{n_1-1} \dots f_r^{n_r-1})$; \item[{\rm(e)}]
$\la_1,\dots,\la_r \in \bc$ and $\sum\limits_{j=1}^r
\la_j\,\deg(f_j)=0$. \item[{\rm(f)}] If $n_j=1$ then $\la_j \ne 0$.
\newline Moreover: \item[{\rm(g)}] The polar set of $\om$ is given
by $\bigcup\limits_{j=1}^r (f_j=0)$ where $n_j=$ order of $(f_j=0)$
as a polar curve of $\om$ and $\la_j = \Res_{(f_j=0)} \om$.
\end{itemize}
\end{Proposition}

\begin{proof} As we have seen, $\om$ cannot be holomorphic, so that
its polar set $(\om)_\infty$ is not empty. Because $(\om)_\infty$
has codimension one it can be written as $(\om)_\infty
=\bigcup\limits_{j=1}^r (f_j=0)$ where $f_1,\dots,f_r$ are
irreducible polynomials in $\bc^3$, pairwise relatively prime.

\noindent Let $\la_j = \Res_{(f_j=0)} \om \in \bc$ be the residue of
$\omega$ in $(f_j=0) =: Z(f_j) \subset \bc P(2)$ and $n_j=$ order of
$Z(f_j)$ as pole of $\om$.

\begin{Remark}{\rm  $\la_j$ can be calculated/defined as follows.
Choose a point $p_j \in Z(f_j)$ which is not a singular point of the
variety $Z(f_j)$. Take a transverse disc $Z_{p_j}$ centered at
$p_j$\,, transverse to $Z(f_j)$ and such that
$\left(\bigcup\limits_{j=1}^r Z(f_i)\right) \cap \Sigma_{p_j} =
\{p_j\}$. }\end{Remark}

\noindent Choose a small simple loop $\ga_j \subset
\Sigma_{p_j}\setminus\{p_j\}$ positively oriented.

\noindent Then
$$
\la_j := \dfrac{1}{2\pi\sqrt{-1}}  \oint_{\ga_j} \om =
\dfrac{1}{2\pi\sqrt{-1}} \oint_{\ga_j} \om\big\vert_{\Sigma_{p_j}}
$$
By a classical result of Deligne  \cite{Deligne} $\la_j$ is
well-defined (recall that $Z(f_j)$ is irreducible). Now we claim:

\begin{Claim} $\sum\limits_{j=1}^r \la_j \deg(f_j)=0$.
\end{Claim}

\begin{proof} We consider a linear embedding $E\colon \bc P(1) \to
\bc P(2)$ such that the line $E(\bc P(1)) =: \bl$ intersects
$(\om)_\infty$ only at non-singular points of the variety
$(\om)_\infty$ and the intersection is always transverse.

\noindent We consider the restriction, i.e., the induced  $1$-form
$\xi := \om\big\vert_\bl = E^*(\om)$. Then $\xi$ has polar at the
points that correspond to the intersection points $\bl \cap
(\om)_\infty$\,. Moreover, since $\bl$ induces a transverse disc
(like the discs $\Sigma_{p_j}$) at each intersection points $p \in
\bl \cap (\om)_\infty$ we have that if $p \in \bl \cap Z(f_j)$ then
$\Res_p\,\xi = \la_j$\,. By its turn, Bezout's theorem implies that
the number of intersection points $p \in \bl \cap Z(f_j)$ is equal
to $\deg(\bl)\deg(f_j) = \deg(f_j)$. Finally, because $\bl$ is a
Riemann Sphere the theorem of residues applied to $\xi$ says that
$$
0 = \sum_{p\in\bl\cap(\om)_\infty}\, \Res(\xi,p) = \sum_{p=1}^r
\la_j\, \deg(f_j).
$$

\end{proof}

Now we consider the pull-back $\Om = \pi^*(\om)$ which naturally
extends to $\bc^3$.

The polar set of $\Om$ is $(\Om)_\infty = \bigcup\limits_{j=1}^r
(f_j=0) \subset \bc^3$ and $\po\Om=0$ with $\la_j =
\Res_{(f_j=0)}\,\Om$, $n_j =$ order of $(f_j=0)$ in
$(\Om)_\infty$\,. Now we introduce the $1$-form
$$
\al := \sum_{j=1}^r \la_j\, \frac{df_j}{f_j}\,\cdot
$$
This is a closed rational $1$-form on $\bc^3$ such that $\be :=
\Om-\al$ is rational, closed but with all residues equal to zero.
Also $(\be)_\infty \subset (\Om)_\infty$.

\begin{Claim}  $\be$ is exact, i.e., $\be = df$ for some meromorphic function $f$ on $\bc^3$.
\end{Claim}

\begin{proof} Indeed, we start by proving given a closed path
$\ga\colon S^1 \to \bc^3\setminus(\Om)_\infty$ then we have
$\displaystyle\int_\ga \be=0$. Let $\ga\colon S^1 \to
\bc^3\setminus(\Omega)_\infty$ be given and (by approximation
theory) assume that $\ga$ is of class $C^\infty$.
Since $\bc^3$ is simply-connected there is a continuous extension
$F\colon \overline{\bd} \to \bc^3$ ($\overline{\bd}$ is the closed
unit disc $|z| \le 1$ in $\bc$). Such that $F\big\vert_{S^1} = \ga$.
Again we may assume that:
\begin{itemize}
\item[{\rm(i)}] $F$ is of class $C^\infty$
\item[{\rm(ii)}] $F(\overline{\bd})$ avoids the singular set of $(\Om)_\infty$ (which is a finite set).
\item[{\rm(iii)}] $F$ is transverse to (the smooth part of) $(\Om)_\infty$\,.
\end{itemize}
In particular $F(\overline{\bd}) \cap (\Om)_\infty =
\{z_1,\dots,z_m\}$ is a finite set.

\noindent Then by the theorem of residues we have
$$
\int_\ga \be = \int_{S^1} F^*(\be) = \sum_{j=1}^m 2\pi
i\,\Res_{z_j}\,F^*(\be) = 0.
$$

\end{proof}

In order to conclude that $\beta$ is exact on $\bc P(2)$ it is
enough to observe that we have already $\be = df$ for some
meromorphic function $f\colon \bc^3\setminus(\Om)_\infty \to
\overline{\bc}$. Because $(\be)_\infty \subset (\Om)_\infty$\, the
function $f$ is indeed holomorphic in $\bc^3\setminus(\Om)_\infty$.
Once we know that $\be$ is rational we can already conclude from
$df=\be$ that $f$ admits an extension to $\bc^3$ as a rational
function. Thus we have proved that $\beta$ is the differential of a
rational function $f$ on $\bc^3$ with $(f)_\infty \subset
(\Om)_\infty$.

\end{proof}

Thus we can write
$$
\Om = \sum_{j=1}^r \la_j\, \frac{df_j}{f_j} + df \quad\text{on}\quad
\bc^3
$$
or else
$$
\Om = \sum_{j=1}^r \la_j \, \frac{df_j}{f_j} + d\left(\frac
gh\right)\cdot
$$

\section{Riccati foliations}
\label{section:ricattifoliations}

The classical Riccati differential equation, put in terms of complex
variables, is
$$
\begin{cases}
\dot x = p(x)&\\
\dot y = a(x)y^2 + b(x)y + c(x).&
\end{cases}
$$

 We will consider  the case where the coefficients are
polynomials \linebreak $p(x),a(x),b(x), c(x) \in \bc[x]$. The
appropriate space to study the geometry of a {\it Riccati
foliation\/}\index{foliation! Riccati}  is perhaps the surface
$\overline{\bc}\times\overline{\bc} = M$. Let us see why. First we
recall the canonical coordinate changes in $M$.

\noindent On $M$ we have the natural projection $\pi_1\colon M \to
\overline{\bc}$ given by $\pi_1(p_1,p_2) = p_1$\,.  The fibers
$\pi_1^{-1}(p) = \{p\} \times \overline{\bc}$ are Riemann spheres.
Now we consider a Riccati foliation $\fa$ of $M$ obtained by
extending from $\bc^2_{(x,y)}$ to $M$ the foliation induced by the
polynomial vector field $X(x,y) = p(x)\,\dfrac{\po}{\po x} +
(a(x)y^2 + b(x)y + c(x))\,\dfrac{\po}{\po y}\,\cdot$ We claim:

\begin{Claim}  {\it A vertical fiber $\pi_1^{-1}(p_1)$\,\, $p_1\ne \infty$ is invariant by $\fa$ if, and only if, $p_1$ is a zero $p(x)$.}
\end{Claim}
This is quite clear since $\fa$ is given on $\bc_{(x,y)}^2$ by
$\begin{cases}\dot x = p(x)&\\ \dot y =
a(x)y^2+b(x)y+c(x)&\end{cases}$ and $\pi_1(x,y) = x$. Therefore
$\fa$ has a finite number of invariant vertical fibers, that depends
on the degree if $p(x)$. Let us change coordinates on $\fa$:

\begin{align*}
&\quad\dot Y = -\dfrac{\dot y}{y^2} = - Y^2 \left(\dfrac{a(x)}{Y^2} + \dfrac{b(x)}{Y} + c(x)\right)\\
&\begin{cases}
\dot Y = -(c(x)Y^2 + b(x)Y + a(x))\\
\dot x = p(x)
\end{cases}
\end{align*}
Thus once again $\fa$ is given by a Riccati ordinary differential
equation in the $\bc_{(x,Y)}^2$ coordinate system. Therefore, we can
claim:

\begin{Claim}  {\it Given a non-invariant vertical fiber $F = \pi_1^{-1}(p_1)$ the leaves of $\fa$ are transverse to $F$.}
\end{Claim}

\begin{proof} Indeed, we may assume that $p_1 \ne \infty$ and
therefore the non-invariance of $F$ is equivalent to $p(x)$ does not
vanish at $x = p_1$\,.

Since the expressions of $\fa$ in the coordinate systems $(x,y)$ and
$(x,Y)$ are
$$
\begin{cases}
\dot x = p(x)&\\
\dot y = \dots&
\end{cases} \qquad\qquad
\begin{cases}
\dot x = p(x)&\\
\dot Y = \dots&
\end{cases}
$$
we conclude that $F$ is always transverse to $\fa$.\end{proof}

Actually the above transversality is a characterization of Riccati
foliations.
\begin{Claim}  {\it A foliation $\fa$ on $\ov\bc \times \ov\bc = M$
which is transverse to some vertical fiber $\pi_1^{-1}(p_1) = F$ is
a Riccati foliation.} \end{Claim} \begin{proof} We may suppose that
$p_1=x_1 \ne \infty$ and choose a polynomial vector field
$\begin{cases} \dot x = A(x,y)&\\ \dot y = B(x,y)&\end{cases}$ that
defines $\fa$ on $\bc^2_{(x,y)}$\,.

\noindent By the transversality on $\bc^2_{(x,y)}$ we conclude that
$A(x,y) \ne 0$, $\forall\, y \in \bc$. Since $A$ is a polynomial
this implies that $A(x,y)$ depends only on $x$\,, not on $y$. Now,
since the fiber $F \cong \overline{\bc}$ is compact, the foliation
is also transverse to the nearby fibers
$\pi_1^{-1}(\widetilde{x}_1)$ for $\widetilde{x}_1 \approx x_1$\,.
By the same reasoning above we conclude that $A(\widetilde{x}_1,y) =
A(\widetilde{x}_1)$ does not depend on $y$, for every
$\widetilde{x}_1 \approx x_1$\,. Because $A(x,y)$ is a polynomial
this implies that $A(x,y) = A(x)$ depends only on $x$. Thus
$\fa\big\vert_{\bc^2}$ is given by $\begin{cases} \dot x = A(x)&\\
\dot y = B(x,y)&\end{cases}$.

\noindent Now we change coordinates to $(x,Y) = (x,\frac 1y)$. Then
$\dot Y = \frac{-\dot y}{y^2} = -Y^2(B(x,\frac 1Y))$ and then $\dot
Y = \dfrac{-Y\,\widetilde{B}(x,Y)}{Y^n}$ where $\widetilde{B}(x,Y)$
is a polynomial, $Y \nmid \widetilde{B}(x,Y)$ and $n = \deg_yB$.
Since $\fa$ is given by $\begin{cases} \dot x = A(x)&\\ \dot Y =
-\dfrac{\widetilde{B}(x,Y)}{Y^{n-2}} \end{cases}$ and by the same
arguments above, the transversality of $F$ with $\fa$ at the point
implies that $n-2 \le 0$, i.e., $n \le 2$.

\noindent Therefore $B(x,y) = b_0(x)+b_1(x)y+b_2(x)y^2$ and $\fa$ is
a Riccati foliation. \end{proof}

From the structural point of view, Riccati foliations are related to
suspensions of groups of automorphisms of $\bc P(1)$. For this let
us recall the classical concept of Ehresmann.

\begin{Definition} {\rm Given a fiber bundle space $\xi(\pi\colon
\overset{F}{\longrightarrow} B)$ with basis $B$, fiber $F$, total
space $E$ and projection $\pi$; a foliation $\fa$ on $E$ is said to
be {\it transverse to the fibers of} $\xi$ if: \begin{itemize}
\item[{\rm(i)}] $\dim \fa + \dim F = \dim E$; \item[{\rm(ii)}] {\it
$\fa$ is transverse to the fibers $\pi^{-1}(b) \subset E$;}
\item[{\rm(iii)}] {\it Given any leaf $L$ of $\fa$ the restriction
$\pi\big\vert_L\colon L \to B$ is a (surjective) covering map.}
\end{itemize} } \end{Definition}

Recall that in this case we have a natural action of $\pi_1(B)$ on
$\Diff(F)$ given a base point $b_0 \in B$ and a path $\ga \in
\pi_1(B,b_0)$. We define a map $h_\ga\colon F_{b_0} \to F_{b_0}$ as
follows: given $y \in F_{b_0} = \pi^{-1}(b)$ we consider the lifted
path $\widetilde{\ga}_y(t) \subset L_y$ obtained from the covering
map $\pi\big\vert_{L_y}\colon L_y \to B$ (where $y \in L_y$ is a
leaf of $\fa$).

\noindent Then we put $h_\ga(y) := \ga_y(1)$; the final point of the
lifting.

\noindent The image of this group homomorphism $\vr\colon
\pi_1(B,b_0) \to \Diff(F)$ is called the {\it global holonomy\/} of
$\fa$ in $\xi$. Is is well-known that $\fa$ is conjugate to the
suspension of its global holonomy. One important fact is the
following remark by Ehresmann:

\begin{Proposition}[\cite{camacho-linsneto,Godbillon}] {\it Let $\fa$ be a
foliation on the fiber space $\xi(\pi\colon E
\overset{F}{\longrightarrow} B)$. Assume that {\rm(i)}\, $\dim \fa +
\dim F = \dim E$ and {\rm(ii)}\, $\fa$ is transverse to the fibers
$\pi^{-1}(b) \subset E$. Then $\fa$ is transverse to the fibers of
$\xi$ if the fiber $F$ is compact.} \end{Proposition} As a corollary
we obtain:

\begin{Proposition} {\it Let $\fa$ be a Riccati foliation on
$\overline{\bc} \times \overline{\bc}$ and let $F_1,\dots,F_r
\subset \overline{\bc}\times\overline{\bc}$ be the invariant
{\rm(}vertical{\rm)} fibers of $\fa$. Then $\fa_0 =
\fa\big\vert_{(\overline{\bc}\times
\overline{\bc})\setminus\bigcup\limits_{j=1}^r F_j}$ is a foliation
transverse to the fibers of the fiber space $\xi\colon
(\overline{\bc}\times\overline{\bc})\setminus \bigcup\limits_{j=1}^r
F_j \to \overline{\bc}\setminus\bigcup\limits_{j=1}^r F_j \to
\overline{\bc}$\, $(F_j = p_j \times \overline{\bc})$ and in
particular $\fa_0$ is holomorphically conjugate to the suspension of
a finitely generated group of Möebius maps.} \end{Proposition}

\chapter{Foliations with algebraic limit sets}
\label{chapter:limitsets}

\section{Limit sets of foliations}
\label{section:limitsets} Let $\fa$ be a holomorphic foliation with
singularities on a compact complex manifold $M$. Given a leat $L \in
\mathcal F$ we consider an {\it exhaustion}\index{exhaustion} by
compact subsets of $L$, i.e., $L=\bigcup\limits_{j \in \mathbb N}
K_j$, where each $K_j\subset L$ is a compact subset and $K_j \subset
Int(K_{j+1})$ for all $j \in \mathbb N$.

\begin{Definition} {\rm The {\it limit set of the
leaf\/}\index{leaf! limit set}\index{limit set! of a leaf} $L$ is
defined as $\lim(L) = \bigcap\limits_{j\in\bn}(\overline{L\setminus
K_j})$. The {\it limit set of foliation\/}\index{foliation! limit
set}\index{limit set! of a foliation} $\fa$ is $\lim(\fa) =
\bigcup\limits_{L\in\fa} \lim L$. } \end{Definition}

\begin{Remark}{\rm  This notion is clearly motivated by the theory
of real Dynamical Systems and also by the dynamics of groups of
rational maps on the Riemann sphere. }\end{Remark}

The very basic properties of the limit set of a foliation are listed
below:

\begin{Proposition} {\it Let $\fa$ and $M$ {\rm(}compact{\rm)} be as
above, then:} \begin{itemize} \item[{\rm(i)}] $\lim(\fa) \subset M$
is $\fa$-{\it invariant} \item[{\rm(ii)}] $\sing(\fa) \subset
\lim(\fa)$ \item[{\rm(iii)}] {\it If $\dim M=2$ then for a leaf $L
\in \fa$ we have $\lim(L) \subset \sing(\fa) \Leftrightarrow
\overline{L} \subset M$ is an analytic curve.} \item[{\rm(iv)}] {\it
If $M = \bc P(2)$ then $\lim(\fa) \subset \sing(\fa)$ iff $\fa$ has
a rational first integral.} \end{itemize} \end{Proposition}

\begin{proof} Let us prove (iii) since (i) and (ii) are more immediate. Assume that $\dim M=2$ and
that $L \in \fa$ satisfies $\lim(L) \subset \sing(\fa)$. We claim
that $\overline{L}$ is an analytic curve in $M$.  Indeed given a
point $p \in \overline{L}\setminus L$ then necessarily we have $p
\in \sing(\fa)$ (because $\lim(L) \subset \sing(\fa)$). Therefore by
Remmert-Stein extension theorem $\overline{L}\subset M$ is analytic
of dimension one. The converse of (iii) is clear. Let us now prove
(iv). Assume that $M = \bc P(2)$ and $\lim(\fa) \subset \sing(\fa)$.
Then from (iii) every leaf $L$ of $\fa$ is contained in an analytic
curve which is $\fa$-invariant. By Chow's theorem every leaf of
$\fa$ is algebraic. By Darboux-Joaunolou theorem $\fa$ admits a
rational first integral.

\end{proof}

Next we give some examples of limit sets of foliations.

\begin{Example} [Linear foliations] {\rm We consider a {\it linear
foliation}\index{foliation! linear} $\fa$ on $\bc P(2)$, given in an
affine chart by $\begin{cases} \dot x = \la x&\\ \dot y = \mu
y&\end{cases}$, \, $\la,\mu \in \bc\setminus\{0\}$. Then:
\begin{itemize} \item[{\rm(i)}] $\la/\mu \in \bq \Rightarrow$ the
leaves are all algebraic and we have $\lim(\fa) = \sing(\fa)$.
\item[{\rm(ii)}] $\la/\mu \in \re\setminus\bq \Rightarrow \lim(\fa)$
is not algebraic, indeed for each leaf $L \in \fa$ we have $\lim(L)
= M_L^3$ where $M_L^3 \subset \bc P(2)$ is the singular real variety
of dimension $3$ given by $|x|^\mu|y|^{-\la} = c \in \re$ (we assume
that $\begin{cases} \la,\mu \in \re&\\ \la/\mu \notin
\bq&\end{cases}$) for a certain constant $c > 0$. \item[{\rm(iii)}]
$\la/\mu \in \bc\setminus\re \Rightarrow \lim(\fa)$ is algebraic, it
is the union of three projective lines: the compactification of the
axes $(x=0)$ and $(y=0)$, and the line at infinity $\bc P(1)_\infty
= \bc P(2)\setminus\bc^2$. \end{itemize} } \end{Example}

Another important property of the limit set is:

\begin{Lemma} {\it Let $\pi\colon \widetilde{M} \to M$ be a proper holomorphic map,
$\fa$ a foliation of $M$ which is generically transverse to $\pi$
(meaning that the  set of tangent points of $\fa$ with $\pi$ has
codimension $\ge2$ in $M$). Then for a leaf $L \in \fa$ if we denote
by $\widetilde{L}$ the inverse image $\widetilde{L} = \pi^{-1}(L)$
we have that $\widetilde{L}$ is a {\it finite\/} union of leaves of
the pull-back foliation $\widetilde{\fa} := \pi^*(\fa)$; say
$\widetilde{\fa} = \widetilde{L}_1 \cup\dots\cup \widetilde{L}_r$\,.
Moreover we have (since $\pi$ is proper)}
$$
\pi^{-1}(\lim(l)) = \bigcup_{j=1}^r \lim (\widetilde{L}_j).
$$
{\it In particular, we have}
$$
\lim \widetilde{\fa} \subset \pi^{-1}(\lim \fa).
$$
\end{Lemma}

\begin{Remark}{\rm  For the proof, the essential point is that given an
exhaustion $\{K_j\}_{j\in\bn}$ of $L$ by compact subsets, since
$\pi$ is proper, the collection $\{\widetilde{K}_j =
\pi^{-1}(K_j)\}_{j\in\bn}$ defines  a compact exhaustion of $L :=
\pi^{-1}(L)$. }\end{Remark}


As a consequence of this and of the linear case in the above example
we obtain:

\begin{Example}{\rm  A rational pull-back $\pi^*(\fa) = \widetilde{\fa}$ to $\bc P(2)$ of a
linear hyperbolic foliation $\fa\colon xdy-\la ydx = 0$ \, $(\la \in
\bc\setminus\re)$  (by a rational map $\pi\colon \bc P(2) \to \bc
P(2)$) is a foliation with an algebraic limit set of dimension one.
}
\end{Example}

Another example is given by

\begin{Example} [Riccati foliations\index{foliation! Riccati}]{\rm Given a Riccati
foliation $\eR$ on $\overline{\bc}\times\overline{\bc}$ we know that
its dynamics  is strongly related to the dynamics of a finitely
generated group of Möebius transformations. In particular, by
choosing suitable subgroups of $\SL(2,\bc)$ we can obtain Riccati
foliations on $\overline{\bc}\times\overline{\bc}$ with global
holonomy a group $G \subset \SL(2,\bc)$ with one or two fixed points
in the Riemann sphere. This foliation will have an algebraic limit
set, consisting of one or two curves. Any rational map $\pi\colon
\bc P(2) \to \overline{\bc}\times\overline{\bc}$ then induces a
following $\pi^*(\fa)$ on $\bc P(2)$ with algebraic limit set. For
instance, let us take any finitely generated group  of Möebius
transformations $G\subset \SL (2,\co)$. Assume that the limit set of
$G$ on $\mathbb CP(1)$ is a single point, which can be assumed to be
the origin $0\in \overline \co$. The limit point $0$ is a fixed
point of $G$. According to \cite{LN3} we can find a Riccati
foliation $\fa : p(x)dy - (a(x) y^2 + b(x) y + c(x))dx=0$ on
$\overline \co \times \overline \co$, whose holonomy group of the
line $\overline{(y=0)}$ is conjugated to the group $G$. Moreover we
can assume that the singularities of $\fa$ over this horizontal line
are reduced and non degenerate. The line $\overline{(y=0)}$ is
invariant by $\fa$ so that $c(x)=0$, and also it is contained in the
limit set of $\fa$. This example can also be seen in $\co P(2)$
using a birational transformation. This will create a dicritical
singularity for the foliation.}
\end{Example}

\section{Groups of germs of diffeomorphisms with finite limit set}
\label{section:germswithlimitset}

We shall study the dynamics of subgroups of $\Diff(\bc,0)$ which may
be associate to the holonomy groups of foliations with
analytic/algebraic limit sets. We shall start with some very basic
facts:

\begin{Lemma}  [Poincaré linearization lemma\index{Lemma! of Linearization of Poincaré}]
Let $f \in \Diff(\bc,0)$ be such that $|f^\prime(0)| \ne 1$. Then:
\begin{itemize}
\item[{\rm(i)}] $f$  is analytically linearizable: $\exists\,\phi \in \Diff(\bc,0)$ such
that $\phi \circ f(z) = f^\prime(0).\phi(z)$.
\item[{\rm(ii)}]  If $\psi \in \Diff(\bc,0)$ is any map that
linearizes $f$ (i.e., $\psi \circ f(x) = f^\prime(0).\psi(z)$) then
$\phi \circ \psi^{-1}$ is  linear (i.e., $\phi = \mu\psi$ for some
$\mu \in \bc^*$).
\end{itemize}
\end{Lemma}

\begin{proof} (i) is the well-known linearization theorem of Poincaré. Now we claim:

\begin{Claim} {\it Let $g \in \Diff(\bc,0)$ commuting with $f$, i.e., $f \circ g = g \circ f$.
Then $g$ is linear in any coordinate that linearizes $f$.}
\end{Claim}

\begin{proof} Indeed, write $f(z) = \la z$ with $|\la| \ne 1$.
 In particular $\la^n \ne 1$, $\forall\,n\ne0$. Write $g(z) = \sum\limits_{n=1}^\infty g_n\,z^n$.
 Since $g \circ f = f \circ g$ we conclude that
$$
\la.g_n = \la^n.g_n\,, \,\,\, \forall\, n \in \bn.
$$
Since $\la^n \ne \la$, $\forall\,n\ne1$ we get $g_n=0$,
$\forall\,n\ne1$ and therefore $g(z) = g,z$. \end{proof}

If $f_1,f_2 \in \Diff(\bc,0)$ are such that $f_j^{-1}\,f\,f_j = \la
z$\, $j=1,2$ then putting $g = f_1f_2^{-1}$ we conclude that $g
\circ f = f \circ g$ so that $g(z) = \mu z$ and therefore $f_1 = \mu
f_2$ proving (ii).

\end{proof}

The above proof actually gives:

\begin{Lemma}  {\it If $f(z) = \la z$ is linear and $g \in \Diff(\bc,0)$ is such that $f \circ g = g \circ f$, then we have:}
\begin{itemize}
\item[{\rm(i)}] $\la^n \ne 1$\quad $\forall\,n\not\le0 \Rightarrow g(z) = \mu(z)$ {\it (is also linear)}
\item[{\rm(ii)}] {\it $\la^k = 1$ for some $k \in \bn \Rightarrow g(z) = \mu z(1+u(z^k))$ for some holomorphic function $u(z)$ with $u(0)=1$.}
\end{itemize}
\end{Lemma}
We also recall the following (already discussed) result (see
Chapter~\ref{chapter:dynamics}):

\begin{Theorem} Let $f \in \Diff(\bc,0)$ be of the form $f(z) = z + a_{k+1}\,z^{k+1} +\dots,a_{k+1} \ne 0$.
Then $f$ is topologically conjugate to the diffeomorphism $\widehat{f}(z) = \dfrac{z}{(1+a_{k+1}z^k)^{1/k}}$ is a neighborhood of the origin.
\end{Theorem}

In particular we have:
\begin{itemize}
\item[{\rm(1)}] For every point, close enough to the origin, its orbit
is contained in an invariant by $f$ continuous curve that passes through the origin.
\item[{\rm(2)}] For every point $z$, close enough to the origin, $f^n(z)$ or $f^{-n}(z)$ converges to the origin as $n \to +\infty$.
\end{itemize}

\begin{Definition} {\rm Let $G \subset \Diff(\bc,0)$ be a subgroup. Given a
connected neighborhood $V$ of $0$ in $\bc$ and $z \in V$, the ({\it
pseudo-orbit of $z$ by\/}\index{pseudo-orbit} $G$ is defined as
$O(z) = \{f(z); f$  is the representative of some element of $G$ and
$z \in \Dom(f)\}$. Given $z \in V\setminus\{0\}$ we say that {\it
the pseudo-orbit of $z$ is discrete off the origin\/}{\it
pseudo-orbit! discrete off the origin} if $\overline{O(z)}\setminus
O(z) \subset \{0\}$. If this is true for all $z \in V\setminus\{0\}$
then we say that $G$ has {\it discrete pseudo-orbits off the origin
in $V$.}

}
\end{Definition}

A remarkable theorem of Isao Nakai implies the following:

\begin{Theorem} [I. Nakai 1994, \cite{nakai}\index{Theorem! of Nakai}] {\it Let
$G \subset \Diff(\bc,0)$ be a non-solvable finitely generated
subgroup.  Then there is a fundamental system of neighborhoods $0
\in V \subset \bc$ such that on each $V$ the group has no
non-trivial orbit closed or discrete off the origin.}
\end{Theorem}
A particular case of this is proved below:

\begin{Proposition} [key proposition, \cite{C-LN-S2}] {\it Let $G \subset \Diff(\bc,0)$ be a subgroup such that:}
\begin{itemize}
\item[{\rm(1)}] $\exists\, f \in G$ {\it with\/} $|f'(0)| < 1$.
\item[{\rm(2)}] {\it $\exists$\, neighborhood $0 \in V \subset \bc$ such that
the orbits of $G$ in $V$ are discrete off the origin. Then $G$ is abelian.}
\end{itemize}
\end{Proposition}

\begin{proof} We may $f(z) = \la z$ ($|\la|<1$) in some local coordinate; in
some sub-neighborhood of $0$ in $V$. Suppose that $G$ is not
abelian. Then some map  $g \in G$ does not commute with $f$:
otherwise $G$ is linearizable and therefore abelian.

Thus $\exists\, h \in G$, $h = [f,g] \ne \Id$. The dynamics of $h$
is such that $\forall\,z \approx 0$ we have $h^n(z) \to 0$ or
$h^{-n}(z) \to 0$ as $n \to +\infty$. Let now $A \subset \bc$ a {\it
fundamental domain\/} for the attractor $f$, i.e., $A =
\overline{D\setminus f(D)}$ where $0 \in D$ is a small disc centered
at the origin. Notice that for any $z \ne 0$, $z \approx 0$\,
$\exists\, n \in \bz$  such that  $f^n(z) \in A$.

\begin{Claim} $\exists$\, a non-discrete orbit in $A$.
\end{Claim}

\begin{proof} Choose a compact disc $0 \in K \subset D$ so that $K \cap A = \phi$.
For each $z \in A$ there is a minimal $m_1(z) \in \bz$ such that
$h^{m_1}(z) \in K$.  There is also a minimal positive number $n_1(z)
\in \bn$ such that $f^{-n_1}\circ h^{m_1}(z) \in A$. Proceeding in
this way we get a sequence of points $\{z_r\} \subset
\bc\setminus\{0\}$ of the form $z_r = f^{-n_r} \circ h^{m_r}
\circ\cdots\circ f^{-n_1} \circ h^{-m_1}(z) \in A$ such that
$h^{m_r} \circ\cdots\circ f^{-n_1} \circ h^{-m_1}(z) \in K$,
$\forall\, r \in \bn$. Given two sequences of numbers $m =
\left\{m_j\right\}_{j=1}^r$\quad $n = \left\{n_j\right\}_{j=1}^r$ as
above we define the set
$$
V_{m,n} := \{z \in A; f^{-n_r} \circ h^{m_r} \circ\cdots\circ
f^{-n_1} \circ h^{-m_1}(z) = z\}.
$$
{\it Then $V_{m,n}$ is a finite set:} otherwise (since $A$ is
compact) we should have an accumulation point in $A$ and then
$f^{-n_r} \circ h^{m_r} \circ\cdots\circ f^{-n_1} \circ h^{-m_1}(z)
= z$, $\forall\,z$. On the other hand the derivative $\left(f^{-n_r}
\circ h^{m_m} \circ\cdots\circ f^{-n_1} \circ h^{-m_1}\right)'(0)
\ne 1$, contradiction. Thus $\bigcup\limits_{m,n}V_{m,n}$ is
countable so that $A\setminus \bigcup\limits_{m,n} V_{m,n} \ne \phi$
and we have some non-discrete orbit in $A$. \end{proof} The claim
completes the proof of the proposition.
\end{proof}

As a consequence of the above results:

\begin{Proposition} {\it Let $G \subset \Diff(\bc,0)$ be a subgroup such that:}
\begin{itemize}
\item[{\rm(i)}] {\it $G$ contains an attractor.}
\item[{\rm(ii)}] {\it The pseudo-orbits of $G$ are discrete off the origin.}
\end{itemize}

\noindent{\it Then $G$ is abelian and linearizable.}
\end{Proposition}

\section{Virtual holonomy groups}
\label{section:virtualholonomy}

The virtual holonomy group to be introduced below is the geometric
object that measures the accumulations of the leaves around a given
leaf. Let us be more precise:

\noindent Let $\fa$ be a foliation of a complex surface $M$, $L \in
\fa$ a leaf; $q \in L$ a base point $(q \notin \sing(\fa))$ and
$\Sigma$ a transverse disc through $q \in \Sigma \cap L$.

\noindent We consider the holonomy group $\Hol(\fa,L,\Sigma,q)
\hookrightarrow \Diff(\Sigma,q)$ and introduce the {\it virtual
holonomy group}\index{virtual holonomy group}\index{holonomy!
virtual holonomy} as follows (\cite{C-LN-S2}):
$$
\Hol^{\virt}(\fa,L,\Sigma,q) := \{f \in \Diff(\Sigma,q); L_z =
L_{f(z)}\,, \forall\, z \in \Sigma\}
$$

In the above notation $L_y$ is the leaf of $\fa$ (in $M\setminus\sing(\fa))$ that contains $y\in M \setminus \sing(\fa)$.

\noindent Then clearly $\Hol(\fa,L,\Sigma,q) \subset
\Hol^{\virt}(\fa,L,\Sigma,q)$. Then main result of this chapter is
the following:

\begin{Theorem} [Linearization theorem, Camacho-Lins Neto-Sad \cite{C-LN-S2}\index{Theorem! Linearization theorem}]
\label{Theorem:linearization} {\it Let $\fa$ be a holomorphic
foliation on $\bc P(2)$. Assume that the limit set of $\lim(\fa)$ is
algebraic of codimension one and contains an irreducible component
$\La \subset \bc P(2)$ of dimension one such that:}
\begin{itemize}
\item[{\rm(1)}] {\it $\sing(\fa) \cap \La$ is non-dicritical and contains no saddle-node in its reduction of singularities.}
\item[{\rm(2)}] {\it Some component in the reduction of singularities of $\La$ contains an attractor in its virtual holonomy group.}
\end{itemize}
{\it Then there is a rational map $\phi\colon \bc P(2) \to \bc P(2)$
and there is a linear foliation  $L_\la\colon xdy-\la ydx = 0$, $\la
\in (\bc\setminus\re)\cup(\re_-\setminus\bq)$ such that $\fa$ is the
pull-back $\fa = \phi^*(\L_\la)$ of $\L_\la$ by $\phi$.}
\end{Theorem}

\section{Construction of closed meromorphic forms}
\label{section:constructionclosedforms}

We shall now see how to construct closed meromorphic $1$-forms
defining a foliation, based on information an the virtual holonomy
and on the singularities. This is done in a neighborhood of a
compact analytic invariant divisor.

\begin{Proposition}
\label{Proposition:constructionclosedform} \noindent Let  $\fa$ be a
foliation on $M^2$ with $\sing(\fa) \subset \La \subset M$, where
$\La$ is an analytic (compact) invariant curve. Denote by $\pi\colon
(\widetilde{M},D) \to (M,\La)$ the reduction of singularities of
$\fa$ in $\La$ and let $\pi^*\fa = \widetilde{\fa}$. Assume that:
\begin{itemize}
\item[{\rm(1)}] $D$ is invariant and $\sing(\widetilde{\fa}) \subset D$ contains no saddle-node singularity.
\item[{\rm(2)}] Each component $D_j \subset D$ has abelian virtual holonomy group and contains an attractor.
\item[{\rm(3)}] $D$ has no cycles.
\end{itemize}
Then there exists a neighborhood $\widetilde{V}$ of $D$ in
$\widetilde{M}$ where $\widetilde{\fa}$ is defined by a closed
meromorphic $1$-form $\widetilde{\om}$, with simple poles and
$(\widetilde{\om})_\infty \supset D$.
\end{Proposition}

\begin{proof}[Proof of Proposition~\ref{Proposition:constructionclosedform}]
We observe that $D = D_0 \cup D_1 \cup\cdots\cup D_r$ where $D_0$ is
{\it the strict transform of\/} $\La$, i.e., $d_0 =
\overline{\pi^{-1}(\La\setminus\sing(\fa))} \subset \widetilde{M}$.

\noindent Let us first consider the case of a component $D_j \subset
D$.

\begin{Lemma}
\label{Lemma:constructionclosedform} Given a component $D_j \subset
D$  there exists a closed meromorphic $1$-form $\om_j$\,, with
simple poles, defined in a neighborhood $U_j$ of $D_j$\,, such that
$\widetilde{\fa}\big\vert_{U_j}$ is given by $\om_j=0$ off
$(\om_j)_\infty$. The $1$-form $\om_j$ is uniquely determined by the
condition: given $q \in D_j\setminus\sing \widetilde{\fa}$,\,
$\Sigma \ni q$ transverse disc, $\Sigma \cap D_j = q_j$ and a
holomorphic coordinate $z$ in $\Sigma$\,\, $(z(q)=0)$ that
linearizes the virtual holonomy, then $\om_j\big\vert_\Sigma =
\dfrac{dz}{z}\,\cdot$
\end{Lemma}

\begin{proof} [Proof of Lemma~\ref{Lemma:constructionclosedform}] Given a
 point $p \in D_j\setminus\sing\widetilde{\fa}$, choose a holomorphic
 chart $\phi = (x,u)\colon U \to \phi(U) \subset \bc^2$ with $p \in U$, $\phi(p)=(0,0)$, $\phi(U) = \{(x,y)\colon |x|<2, |y|<2\}$ and:
\begin{itemize}
\item[{\rm(1)}] $\widetilde{\fa}\big\vert_U$ is given by $dy=0$;
\item[{\rm(2)}] $D_j \cap U \subset (y=0)$
\item[{\rm(3)}] $\Sigma\colon (x=0)$ is transverse disc to $\widetilde{\fa}$
and $g\big\vert_\Sigma$ a local chart that linearizes $\Hol^{\virt}(\widetilde{\fa},D_j,\Sigma)$.
\end{itemize}

\begin{Remark}{\rm  The existence of $\phi\colon U \to \phi(U)$ is obtained
by extending a local transverse coordinate that linearizes the
virtual holonomy  from the transverse section to a neighborhood (of
bidisc type) as constant along local plaques of $\widetilde{\fa}$.
}\end{Remark}

\noindent We obtain then an open cover
$$
\U = \{(U_\al), (x_\al,y_\al)\colon U_\al \to \bc^2\}_{\al \in
\sqrt\al}
$$
of $D_j\setminus\sing(\widetilde{\fa})$ satisfying (1), (2) and (3)
above. We may also assume that
\begin{itemize}
\item[{\rm(4)}] If $U_\al \cap U_\be \ne \phi$ then $U_\al \cap U_\be$ is connected.
\end{itemize}
Let us the study this situation $U_\al \cap U_\be \ne \phi$.

\begin{Claim} We have $y_\al = c_{\al\be}\cdot y_\be$ for some constant $c_{\al,\be} \in \bc^*$.
\end{Claim}
The claim is an easy consequence of (3) and the fact that
$\Hol^{\virt}(\fa,L_j,\Sigma,q)$ contains an attractor, where $L_j =
D_j\setminus\sing(\widetilde{\fa})$.

\noindent We then conclude that the closed meromorphic $1$-forms
$\dfrac{dy_\al}{y_\al}$ and $\dfrac{y_\be}{y_\be}$ coincide on
$U_\al \cap U_\be$\,.

\noindent This gives a closed meromorphic $1$-form $\om_j$ in $V_j =
\bigcup\limits_{\al\in\be} U_\al$\,, defined by\linebreak
$\omega_j\big\vert_{U_\al} = \dfrac{dy_\al}{y_\al}\,\cdot$

\begin{Claim}
\label{Claim:extensionsingularity} Given any singular point  $p_j
\in \sing \widetilde{\fa} \cap D_j$ the $1$-form $\om_j$ extends
meromorphic to a neighborhood of $p_j$ in $\widetilde{M}$.
\end{Claim}

\begin{proof}[Proof of Claim~\ref{Claim:extensionsingularity}]
The point $p_j$ is not a saddle-node by hypothesis. The local
holonomy of the separatrix $p_j \in \Ga_j$ of $\widetilde{\fa}$
contained in $D_j$ is analytically linearizable because it is in the
virtual holonomy. Therefore, by  Mattei-Moussu
Theorem~\ref{Theorem:matteimoussu}  we conclude that
$\widetilde{\fa}$ is analytically linearizable in a neighborhood of
$p_j$\,. Let then $(x,y)\colon U \to \bc^2$ be a local chart such
that $p_j \in U_j$\,, $x(p_j) = y(p_j) = 0$, $D_j \cap U \subset
(y=0)$ and $\widetilde{\fa}\big\vert_U$ is given by $xdy-\la ydx =
0$, $\la \in \bc^*$ (indeed $\la \in \bc\setminus\bq_+$ because
$p_j$ is irreducible).

\noindent The local holonomy of $\Ga_j$ at the transverse disc
$\Sigma_j\colon (x=1)$ is given by $h(y) = e^{2\pi i\la}\,y$. We set
$\om_{p_j} := \dfrac{dy}{y} - \dfrac{\la dx}{x}$ in $U$.

\noindent Since $\om_{p_j}$ and $\om_j$ define $\fa$ in $V_j \cap U$
($V_j \cap U$ contains a neighborhood of the loop $\ga \subset
\Ga_j$ given by $\ga = \{(x,0); |x|=1\}$) we have $\om_j =
f\cdot\om_{p_j}$ for some meromorphic function $f$ in $V_j \cap U$.
Since

\noindent Given $\ve > 0$ let $V_\ve := \{(x,y)\colon 1-\ve < |x| <
1+\ve; |y| < \ve\}$ and $V_\ve \subset V_j \cap U$ for $\ve > 0$
small enough. Since $\om_j$ and $\om_{p_j}$ have simple poles in $V$
we have that $f$ is holomorphic in $V$ and therefore it can be
represented (in $V_\ve$ and therefore) in $V$ by a Laurent Series
like $f(x,y) = \sum\limits_{i\in\bz,j\ge0} h_{ij}\,x^iy^j$.

\noindent Since $\om_j = f\om_{p_j}$ and $d\om_j=0=d\om_{p_j}$ we
obtain $df \wedge \om_{p_j}=0$. This last relation can be rewritten
as
$$
(*)\qquad x\,f_x + \la y\,f_y = 0
$$
and in terms of the Laurent Series of $f$ as
$$
(**)\qquad (i+\la_j)f_{ij}=0, \,\, \forall\, i \in \bz, \, j \ge 0.
$$

\noindent\textbf{Case 1:} If $\underline{\la\notin\bq}$ then
$f_{ij}=0$\,\, $\forall(i,j) \ne (0,0)$ and therefore $f$ is
constant. Thus $\om_j = c.\om_{p_j}$ in the common domain, for some
$c \in \bc^*$. Comparing residues along $(y=0)$ we obtain $c=1$ and
therefore $\om_j = \om_{p_j}$ in $V_j \cap U$. In particular,
$\om_j$ extends as $\om_{p_j}$ to $I \ni p_j$\,.

\noindent\textbf{Case 2:} If $\la \in \bq_-$ say, $\la = -\dfrac mn$
with $n,m \in \bn$\,\, $\lg m,n\rg=1$. Then (**) implies that
$n_i-m_k=0$ if $f_{ij}\ne0$. That is, $f_{ij}\ne0 \Rightarrow (i,j)
= (km,kn)$ for some $k \ge 0$ (notice that $j\ge 0$). Therefore
$f(x,y) = \ell(x''''\,y''')$ for some holomorphic function $\ell(z)$
and then $f$ admits a holomorphic extension to a neighborhood of
$p_j$\,. The $1$-form $\om_j$ then extends to a neighborhood of
$p_j$ as $f.\om_{p_j}$\,. This ends the proof of the  claim.
\end{proof}

\noindent This and the fact that two coordinates linearizing an
attractor differ up to a multiplication imply
Lemma~\ref{Lemma:constructionclosedform}. \end{proof}

Now we pass from (each) $D_j$ to a neighborhood of $D = D_0 \cup
\dots \cup D_r$\,. Let then $\om_j$ a closed meromorphic $1$-form
with simple poles in a neighborhood $U_j$ of $D_j$\,\,
$(j=0,1,\dots,r)$ in $\widetilde{M}$, such that $\om_j=0$ defines
$\widetilde{\fa}$ in $U_j\setminus(\om_j)_\infty$\,.

\noindent Assume that $D_i \cap D_j \ne \phi$ say $q = D_i \cap
D_j$\,. Then

\noindent in a neighborhood $U_{ij} \subset U_i \cap U_j$ we have
$\om_i = f \om_j$ for some meromorphic\linebreak function $f$.

\begin{Claim}
\label{Claim:constant}\,\,\, $f$ is constant.
\end{Claim}

\begin{proof}[Proof of Claim~\ref{Claim:constant}] Since $D_i$ and $D_j$ are
$\widetilde{\fa}$-invariant the corner $q$ is a singularity of
$\widetilde{\fa}$.  We have seen that is analytically linearizable
say, $\widetilde{\fa}\big\vert_U\colon xdy-\la ydx=0$, in some
neighborhood $U \ni (x,y)$ of $q$, with $D_j \subset (y=0)$, $D_i
\subset (x=0)$.

We have $df \wedge (xdy-\la y dx) \equiv 0$.

\noindent If $\la \notin \bq$ then we have seen that $f$ is
constant. Assume now that $\la = -m/n \in \bq_-$, $m,n \in \bn$,
$\lg m,n\rg = 1$. The virtual holonomy of $D_j$ is linearized by
$y\big\vert_{\Sigma_j}$. {\it We may choose $(x,y) \in U$ such that
also $\Hol^{\virt}(\widetilde{\fa},D_i,\Sigma_i)$ is linear in the
coordinate\/} $x \mapsto (x,1) \in \Sigma_i\colon (y=1)$, let us see
why: Indeed, by hypothesis the virtual holonomy of $D_i$ contains an
attractor say $g$, with $g'(0)=\mu$, $|\mu| < 1$. The local holonomy
of the separatrix $q \in \Ga_i \subset D_i$ is given by $h(x) =
e^{-1\pi i\,\frac nm}\,x$. Since $g$ and $h$ commute we have $g(x) =
\mu x\,\widetilde{g}(x^m)$ for some $\widetilde{g} \in O_1$ with
$\widetilde{g}(0)=1$. Let now $\widetilde{x} = \phi(x)$ be a change
of coordinates valid at $0 \in \Sigma_i$\,, such that $\phi \circ g
\circ \phi^{-1}$ is linear (recall that
$\Hol^{\virt}(\widetilde{\fa},D_i,\Sigma_i,q)$ is abelian
linearizable). Then $\phi(g(x)) = \mu\phi(x)$. Since $\phi$
linearizes $g$ and it also linearizes the local holonomy $h$ of
$\Ga_i$ we conclude that $\phi(x)=x.\phi_1(x^m)$ where
$\widetilde{\phi}_1 \in O_1$, $\widetilde{\phi}_1(0) \ne 0$. We then
consider the change of coordinates $(\widetilde{x},\widetilde{y}) =
\psi(x,y) = (x\widetilde{\phi}(x^my^n),y)$. Then
$\psi\big\vert_{\Sigma_i} \equiv \phi$ so that $\psi$ preserves the
linear foliation $mydx+nxdy=0$. Indeed, $\psi^*(mydx+nxdy) =
u.(m\widetilde{y}d\widetilde{x} + n\widetilde{x}d\widetilde{y})$ for
some $u$ holomorphic with $u(0) \ne 0$.

\noindent Therefore we may assume that $g$ is linear in the
coordinate $(x,y)$.

Now given the linear singularity $q\colon mydx + nxdy = 0$

\noindent We claim:

\begin{Claim}
\label{Claim:virtualattractor} Define $k\colon \Sigma_j \to
\Sigma_j$ by $k(y) := \mu_1\,y$  where $\mu_1^m = \mu^n$. Then $k
\in \Hol^{\virt}(\widetilde{\fa},D_j,\Sigma_j)$ and it is an
attractor.
\end{Claim}

\begin{proof}[Proof of Claim~\ref{Claim:virtualattractor}] We consider the first integral
$\xi = x^my^n$ for $\widetilde{\fa}$ in a neighborhood of $q$. Then
$k(y)$ preserves the  leaves of
$\widetilde{\fa}\big\vert_{\Sigma_j}$ which are given by $g^{n/m} =$
constant.

\noindent It is enough to observe that $\left(k(y)\right)^{m/n} =
y^{m/n}$. \end{proof}

The above implies that $\Hol^{\virt}(\widetilde{\fa},D_j,\Sigma_j)$
is linear in the coordinate $y \mapsto (1,y)$ of $\Sigma_j$ (because
this coordinate linearizes the attractor $k(y)$). But from the
characterization/definition/construction of $\om_i$ we have
$$
\om_i\big\vert_{\Sigma_i} = \frac{dx}{x} \quad\text{and}\quad \om_i
= \frac{dx}{x} + \frac nm\, \frac{dy}{y}
$$
and analogously
$$
\om_j = \frac{dy}{y} + \frac mn\, \frac{dx}{x}
$$
so that $\om_i = \dfrac nm\, \om_j$\,. This proves
Claim~\ref{Claim:constant} in the  resonant case.\end{proof}

Let us proceed with the proof of
Proposition~\ref{Proposition:constructionclosedform}.

\noindent{\it Since $D$ has no cycles\/} we can construct a closed
$1$-form $\omega$ in a neighborhood $\bigcup\limits_{j=0}^r U_j$ of
$D$ by choosing $\om\big\vert_{D_j} = c_j\cdot\om_j$ for some
suitable choice of the constants $c_j \cdot\om_j$ for some suitable
choice of the constants $c_j \in \bc^*$, $j=0,\dots,r$. \end{proof}

\noindent For the case of foliations on $\bc P(2)$ we obtain:

\begin{Proposition}
\label{Proposition:constructionformprojective}{\it Let $\fa$ be a
foliation on $\bc P(2)$ with  an algebraic invariant curve $\La
\subset \bc P(2)$, having a reduction of the singularities for
$\sing(\fa) \cap \La$, \, $\pi\colon (M,D) \to (\bc P(2),\La)$ as
follows:}
\begin{itemize}
\item[{\rm(1)}] {\it $\Sing \widetilde{\fa} \cap D$ contains no saddle-node and is non-dicritical.}
\item[{\rm(2)}] {\it Each irreducible component $D_j$ of $D$ has abelian virtual holonomy containing
an attractor. Then $\fa$ is given by a logarithmic (rational) one form on $\bc P(2)$.}
\end{itemize}
\end{Proposition}

\begin{proof} The point here is that $D = D_0 \cup\dots\cup D_r$ may have cycles.
Nevertheless we have another argument as follows: Since $\fa$ is
defined on $\bc P(2)$ there is a rational $1$-form $\Om$ that
defines $\fa$ and we consider $\widetilde{\Om} = \pi^*(\Om)$ which
is a rational $1$-form defining $\widetilde{\fa}$ on $M$. By the
preceding proposition for each $j=0,\dots,r$ there is a closed
meromorphic $1$-form $\om_j$ in a neighborhood $U_j$ of $D_j$ in $M$
such that $\widetilde{\fa}\big\vert_{U_j}$ is defined by $\om_j$ and
if $U_{ij} = U_i \cap U_j \ne \phi$ then in $U_{ij}$ we have $\om_i
= c_{ij}\,\om_j$ for some $c_{ij} \in \bc^*$. Since
$\widetilde{\Om}$ defines $\widetilde{\fa}$, on each $U_j$ we have
$\widetilde{\Om} = h_i\om_j$ for some meromorphic function $h_j$ in
$U_j$\,. Thus in $U_{ij} \ne \phi$ we have $h_i\om_i = h_j\om_j
\Rightarrow \om_i = \dfrac{h_j}{h_i}\,\om_j$ so that
$\dfrac{h_j}{h_i} = c_{ij}$ and then $\dfrac{dh_i}{h_i} =
\dfrac{dh_j}{h_j}\,\cdot$ Define then $\widetilde{\eta}$ in
$\widetilde{U} \equiv \bigcup\limits_{j=0}^r U_j$ as
$\widetilde{\eta}\big\vert_{U_j} := \dfrac{dh_j}{h_j}\,\cdot$ There
is a closed meromorphic $1$-form $\eta$ in $U = \pi(\widetilde{U})$
such that $\widetilde{\eta} = \pi^*(\eta)$.

\noindent By Levi's extension theorem (Theorem~\ref{Theorem:Levi})
the $1$-form $\eta$ extends as a closed rational one form on $\bc
P(2)$.

\begin{Claim}
\label{Claim:integration}\,\, $\eta = \dfrac{dh}{h}$ for some
rational function $h$ on $\bc P(2)$.
\end{Claim}

\begin{proof}[Proof of Claim~\ref{Claim:integration}] By the description of the closed
rational $1$-forms on $\bc P(2)$ (Proposition~\ref{Proposition:integrationlemma}) it is enough to observe that:
\begin{itemize}
\item[{\rm(1)}] $\eta$ has simple poles on $\bc P(2)$
\item[{\rm(2)}] for each component $H$ of the polar set of $\eta$ we have $\Res_H\,\eta \in \bz$.
\end{itemize}
The proof of (1) and (2) is a consequence of the local description
as $\dfrac{dh_j}{h_j}$ of $\widetilde{\eta}$ in a neighborhood of $D
= \pi^{-1}(\La)$ and of Bézout's theorem (every component of
$(\eta)_\infty$ must intersect $\La$). This proves
Claim~\ref{Claim:integration}.\end{proof}

Finally, let $\widetilde{h} = \pi^*(h) = h \circ \pi$. Then,
$\widetilde{\eta} = \dfrac{dh}{\widetilde h}$ and $\widetilde{h} =$
const. $h_j$ so that $d\left(\dfrac{\widetilde\Om}{\widetilde
h}\right) = 0$. Therefore $d\left(\dfrac{\Om}{h}\right) = 0$.

\noindent Thus $\om := \dfrac{\Om}{h}$ is a closed rational $1$-form
which defines $\fa$ on $\bc P(2)$. We claim (as it is easy to see)
that $\om$ has simple poles so that $\om$ is logarithmic on $\bc
P(2)$. This ends the proof of
Proposition~\ref{Proposition:constructionformprojective}.
\end{proof}

\section{The Linearization theorem}
\label{section:linearizationtheorem}

Now we proceed to prove Theorem~\ref{Theorem:linearization}, the
main result of this chapter.

\begin{proof}[Proof of Theorem~\ref{Theorem:linearization}\index{Theorem! Linearization theorem}]
For the first part it is enough to prove that every component $D_j$
of the exceptional divisor  $D$ in the reduction of singularities
$\pi\colon (M,D) \to (\bc P(2),\La)$ of $\sing(\fa) \cup \La$,
exhibits a hyperbolic attractor or a non periodic linearizable map
in its virtual holonomy group. We consider a component $D_{j_0}$
which has a hyperbolic virtual holonomy map say $h_{j_0}$\,. Then if
$D_i$ is an adjacent component, $D_i \cap D_{j_0} \ne \phi$ we have
two possibilities for the corner singularity $q = D_i \cap D_{j_0}$
\begin{itemize}
\item[{\rm(1)}] $q\colon xdy-\la ydx=0$,\, $\la \notin \bq$, \, $D_i\colon (y=0)$, \, $D_{j_0}\colon (x=0)$. \newline
In this case the virtual holonomy (actually the holonomy) of the
component $D_i$ contains a linearizable non-resonant map $h_i\colon
y \mapsto e^{2\pi i\la}\,y$.
\item[{\rm(2)}] $q\colon xdy-\la ydx=0$, \, $\la = -m/n \in \bq_-$\,, \, $\lg m,n\rg=1$. \newline
In this case we have seen in the above proofs that the first
integral $x^my^n$ permits the ``passage" of the attractor $h_{j_0}$
to an attractor $k$ in the virtual holonomy of $D_i$\,.
\end{itemize}

Thus $\fa$ is given by a logarithmic $1$-form
$\omega\big\vert_{\bc^2} = \sum\limits_{j=1}^\ell \la_j\,
\dfrac{df_j}{f_j}$ where $\la_j \in \bc$, \, $f_j$ is an irreducible
polynomial, $(\om)_\infty = \bigcup\limits_{j=1}^\ell \Ga_j$ where
$\Ga_j = \overline{(f_j=0)}$, for a suitable choice of the
coordinate system $\bc^2 \subset \bc P(2)$ with the line $\bc
P(2)\setminus\bc^2$ not invariant by $\fa$.

We then know that $\sum\limits_{j=1}^\ell d_j\la_j=0$. where $d_j =
\deg(f_j)$. We sketch the main steps. Let $\Ga_1$ be the component
of $(\om)_\infty$ with a virtual holonomy attractor.

\begin{Claim} Fix a point $p_1 \in \Ga_1\setminus\sing(\fa)$, \, $\Sigma_1 \ni p_1$ a
transverse disc and $z \in (\Sigma_1,p_1)$, $z(p_1)=0$ a local
coordinate such that $\Hol^{\virt}(\fa,\Ga_1,\Sigma_1)$  is linear
in the coordinate $z$. Then for each $j\ge2$ the map $h_j(z) =
\left(e^{2\pi\sqrt{-1}}\,\la_j/\la_1\right).z$ belongs to
$\Hol^{\virt}(\fa,\Ga_1,\Sigma_1)$.
\end{Claim}

We also observe that $G = \Hol^{\virt}(\fa,\Ga_1,\Sigma_1)$ is
abelian, contains an attractor, linearizable and has discrete orbits
off the origin. Thus $G$ is generated by an attractor $z \mapsto
e^{2\pi i\la}\,z$ and a rational rotation $z \mapsto e^{\frac{2\pi
i}{m}}\,z$.

From the above claim we then conclude that for each $j\ge2$\,
$\exists\,k_j,\ell_j \in \bz$ such that
$$
\frac{\la_j}{\la_1} = \frac{k_j}{m} + \ell_j\,\la.
$$
Therefore $m\,\dfrac{\la_j}{\la_1} = v_j-u_j\,\la$ for some $u_j,v_j
\in \bz$. We define
$$
\begin{cases}
F_1 = f_1^m\,f_2^{v_2} \dots f_\ell^{v_\ell}\\
F_2 = f_2^{v_2} \dots f_\ell^{u_\ell}
\end{cases}
$$
then we obtain
\begin{align*}
\frac{dF_1}{F_1} - \la\,\frac{dF_2}{F_2} &= m\,\frac{df_1}{f_1} = \sum_{j=2}^\ell (v_j-\la u_j)\, \frac{df_j}{f_j}\\
&= m\left(\frac{df_1}{f_1} + \sum\,
\frac{\la_j}{\la_1}\,\frac{df_j}{f_j}\right) = \frac{m}{\la_1} \cdot
\om.
\end{align*}

Thus the rational map $\phi = (F_1,F_2)\colon \bc P(2) \to \bc P(2)$
is such that $\fa = \phi^*(\L)$ for $\L\colon \dfrac{dx}{x} -
\la\,\dfrac{dy}{y} = 0$. This ends the proof of the Linearization
Theorem. \end{proof}

Theorem~\ref{Theorem:linearization} is proved in a more general
setting in \cite{C-LN-S2}. Indeed it can be stated without the
assumption of absence of saddle-nodes. For this we assume that each
irreducible component of $\Lambda$ contains some hyperbolic
attractor in its virtual holonomy. The precise statement is:

\begin{Theorem} [\cite{C-LN-S2}, page 431]

\noindent{\it Let $\fa$ be a holomorphic foliation on $\bc P(2)$.
Assume that the limit set $\lim(\fa)$ is algebraic. Denote by
$\lim_1(\fa)$ the pure codimension one component of $\lim(\fa)$.
Assume that $\lim_1(\fa)\ne \emptyset$ and:}
\begin{itemize}
\item[{\rm(1)}] {\it $\sing(\fa) \cap \lim_1(\fa)$ and $\lim_1(\fa)$ contains all separatrices of its singularities. }
\item[{\rm(2)}] {\it Any irreducible component of $\lim_1(\fa)$ contains an attractor in its  virtual holonomy group.}
\end{itemize}

\noindent{\it Then there are a rational map $\phi\colon \bc P(2) \to
\bc P(2)$ and  a linear foliation  $L_\la\colon xdy-\la ydx = 0$,
$\la \in (\bc\setminus\re)\cup(\re_-\setminus\bq)$ such that $\fa$
is the pull-back $\fa = \phi^*(\L_\la)$ of $\L_\la$ by $\phi$.}
\end{Theorem}

The proof is a little more elaborate as we have to prove that there
are no saddle-nodes in the reduction of singularities. Nevertheless,
following the same line of reasoning presented in this chapter one
may be able to give an alternative geometrical proof of this fact.

\chapter{Some modern questions}
\label{chapter:modernquestions} We introduce and comment now some
modern questions in the theory of holomorphic foliations with
singularities.

\section{Holomorphic flows on Stein spaces}
\label{section:suzuki}

We shall now discuss several aspects of the dynamics, topology and
analytic classification of the foliation with singularities defined
by an action of the groups $\bc$ or $\bc^*$ on a Stein variety
usually under the presence of a (singular) fixed point of the
action. We consider an action $\vr\colon \bc \times N \to N$, i.e.,
a holomorphic map such that:
\begin{itemize}
\item[{\rm(i)}] $\vr(t,\vr(t_2,x)) = \vr(t_1+t_2,x)$,\, $\forall\, t_1,t_2 \in \bc$, \, $\forall\,x \in N$
\item[{\rm(ii)}] $\vr(0,x) = x$, \, $\forall\, x \in \bc$.
\end{itemize}
The action is {\it periodic\/} of period $\tau \in \bc^*$ if
$\vr(\tau,x)=x$, \, $\forall\, x \in \bn$.

The action induces a group of homomorphism
\begin{align*}
&\bc \to \Aut(N) = \{\text{group of holomorphic diffeomorphisms of\/}\, N\}\\
&t \mapsto \vr_t := \vr\big\vert\_{\{t\}\times N}\colon N \to N
\end{align*}
If the action is periodic of period $\tau$ than $\vr_\tau = \Id$ and
we may induce an action $\psi\colon \bc^*\times N \to N$ by setting
$\psi(a,x) := \vr_{\frac{\tau}{2\pi i}\ell n\,u}(x)$.

Conversely, any action $\psi\colon \bc^*\times N \to N$ defines a
periodic action $\vr\colon \bc\times N \to N$ by $\vr(t,x) =
\psi(e^t,x)$, $\forall\, t \in \bc$, \, $\forall\,x \in N$; of
period $2\pi i$.

Thanks to this, we will focus on actions $\vr\colon \bc \times N \to
N$. Given a point $x \in N$ the {\it orbit\/}\index{orbit} of $x$ is
by definition the subset $O(x) = \{\vr(t,x); t \in \bc\} \subset N$.
Since we are dealing with an action we know that the orbits can
contain important information provided some regularity is required:
Let $Z$ be the holomorphic vector field on $N$ defined by $Z(x) :=
\dfrac{\po\vr}{\po t} (t,x)\bigg\vert_{t=0}$. Then the integral
curves of $Z$ are the orbits of $\vr$ and the fixed points of $\vr$
are the singularities of $Z$, if we assume (and we shall) that $Z$
has isolated singularities on $N$. The vector field $Z$ is {\it
complete\/}\index{vector field! complete} and $\vr$ is its (globally
defined) flow.

\begin{Lemma} The orbits of the action are biholomorphic to
$\bc$, \, $\bc^* \cong \dfrac{\bc}{\bz}$ or a torus $\dfrac{\bc}{\bz\oplus\bz}\,\cdot$
\end{Lemma}
The above is a straightforward consequence of the fact that the
orbit $O(p)$ is through $p$ biholomorphic to the quotient
$\dfrac{\bc}{G_p}$ where $G_p \subset \bc$ is the {\it isotropy\/}
subgroup $G_p = \{t \in \bc; \vr(t,p)=p\}$.

\noindent Since $G_p \subset (\bc,+)$ is a discrete subgroup we
conclude that $G_p \simeq \{0\}$, $\bz$ or $\bz \oplus \bz$.

Assume now that $N$ is a Stein space. Then $N$ contains no positive
dimension compact analytic subset. Therefore we get:

\begin{Corollary} The orbits of $\vr$ on a Stein space are biholomorphic to the
plane $\bc$ or to the cylinder $\bc^* \cong \dfrac{\bc}{\bz}\,\cdot$
\end{Corollary}


\subsection{Suzuki's theory}
\label{section:suzukitheory}

A fundamental contribution to the study of holomorphic flows and
foliations on  Stein surfaces, was made by M. Suzuki who introduced
on this subject the use of techniques from potential theory and the
theory of analytic spaces (cf. \cite{Suzuki1} and \cite{Suzuki2}).
M. Suzuki's work is from  the middle 70's. Some of his main results
are collected below:

\begin{Theorem} [M. Suzuki, \cite{Suzuki1}\index{Theorem! of Suzuki}]
\label{Theorem:Suzuki1}
Given a $\mathbb C$-action $\vr$ on
a normal Stein analytic space $V$ of dimension $n \geq 2$:
\begin{itemize}

\item[{\rm(i)}] There is a subset $\mathcal E\subset V$ of logarithmic capacity
zero such that $\vr_{t}(\mathcal E)=\mathcal E$, for any $t\in \mathbb C$, and all orbits
of $\vr$ in $V\setminus \mathcal E$ are of the same topological type.

\item[{\rm(ii)}] Any leaf of $\fa_{\vr}$ containing an orbit of $\vr$
isomorphic to $\mathbb C^*$ is closed in $V\backslash sing {\fa_\vr}$.

\item[{\rm(iii)}] If $n=2$ and the leaves of $\fa_\vr$ are
properly embedded in $V\setminus \sing(\fa_\vr)$ then there is a
meromorphic first integral of $\fa_{\vr}$ on $V$, not constant, and one
can find a Riemann surface $S$ and a surjective holomorphic map
$p:V\backslash sing { \fa_{\vr}}\to S$, such that:
\begin{itemize}
\item[{\rm(1)}] The irreducible components of the fibers $\{p^{-1}(w); w\in S\}$ of $p$
are the leaves of $\fa_{\vr}$.
\item[{\rm(2)}] The subset $E\subset V$ union of all the reducible levels $p^{-1}(w)$, \,  $w\in S$,
has zero logarithmic capacity.
\end{itemize}

\item[{\rm(iv)}] If $n=2$ and the generic leaf is isomorphic to $\mathbb C^*$,
then any leaf of $\fa_{\vr}$ is closed in $V\backslash sing {
\fa_{\vr}}$ and {\rm(}therefore{\rm)} there is a meromorphic first
integral as in {\rm(iii)}.

\end{itemize}
\end{Theorem}

For the case of actions with isolated singular points we promptly
obtain:

\begin{Theorem} [Suzuki, \cite{Suzuki1,Suzuki2,Suzuki3}\index{Theorem! of Suzuki}]
\label{Theorem:Suzuki2} Given a $\bc$-action $\vr$ with isolated
singularities on a normal Stein analytic space $N$ of dimension
$n\ge2$ we have:
\begin{itemize}
\item[{\rm(i)}] There is a subset $\E \subset N$ of zero logarithmic capacity such that:
\begin{itemize}
\item[{\rm(a)}] $\E$ is invariant by $\vr$\,: $\vr_t(\ve) \subset \ve$, $\forall\, t \in \bc$
\item[{\rm(b)}] all the orbits in $N\setminus\ve$ are diffeomorphic.
\end{itemize}
\item[{\rm(ii)}] A periodic orbit of $\vr$ (i.e., diffeomorphic to $\bc^*$) is closed in $N\setminus\Fix(\vr)$.
where $\Fix(\vr) = \{p \in N; \vr*t,p) = p, \forall\,t \in \bc\}$ is
the set of (singular fixed points of $\vr$. In particular, because
of Remmert-Stein Theorem, such a periodic orbit is contained in an
analytic curve in $N$.
\item[{\rm(iii)}] Assume that $n = \dim N = 2$. Then:
\begin{itemize}
\item[{\rm(a)}] If all the orbits are property embedded in $N\setminus\Fix(\vr)$
then $\tau = \dfrac{\po\vr}{\po t}\bigg\vert_{t=0}$ admits a
meromorphic first  integral $f\colon M \to \overline{\bc}$.
\item[{\rm(b)}] If the {\tt generic} orbit of $\vr$ (i.e., the one in $N \setminus\ve$)
is diffeomorphic to $\bc^*$ then $Z$ admits a meromorphic first integral.
\end{itemize}
\end{itemize}
\end{Theorem}

For the case of analytic actions of $\bc$ on $\bc^2$, in another
remarkable work, M. Suzuki proves:

\begin{Theorem}[Suzuki, \cite{Suzuki2}]
Any $\bc^*$-action on $\bc^2$ is analytically linearizable, {\it
i.e.}, analytically equivalent to an operation of the form $s\circ
(x,y)=(s^n x, s^m y), \, s \in \bc^*, (x,y)\in \bc^2$, for some $n,m
\in \mathbb N$.
\end{Theorem}

The classification of holomorphic $\bc$-actions with proper orbits on $\bc^2$ is
the following:

\begin{Theorem}[\cite{Suzuki2}, Theorem 4]
Every holomorphic action $\vr$ of $\bc$ with proper orbits on
$\bc^2$ is analytically equivalent to one of the following
operations:

\begin{itemize}

\item[{\rm(i)}]Degenerate operations of the form $(\alpha): \, t\circ
(x,y)=(x, y + a(x)t)$ or  \, $(\beta): \, t \circ (x,y)=
(x,e^{\lambda(x)t}(y - b(x)) + b(x))$, where $a(x), \, \lambda(x)$
are entire functions of one variable $x$ and $b(x)$ is a meromorphic
function of $x$ such that $\lambda (x). b(x)$ is holomorphic on
$x\in \bc$.

\item[{\rm(ii)}] Exponential type operations of the form $(\gamma): \, t
\circ (x,y)=(xe^{n\lambda t}, y e^{m\lambda t})$ with $\lambda\in
\bc^*, \, n, m \in \mathbb N$.

\item[{\rm(iii)}] Exponential type operations of the form $(\gamma^\prime):
\,t \circ (x,y)=(xe^{n\lambda(u) t}, y e^{-m\lambda(u) t})$ with
$\lambda$ is an entire function of one variable, $ n, m \in \mathbb
N$ and $u=x^m y^n$.

\item[{\rm(iv)}] Operations of the form $\alpha^{-1}\circ \rho _t \circ
\alpha$ where $\alpha(x,y)=(x, x^\ell y + P_\ell(x)), \, \ell \in
\mathbb N, \, P_\ell$ is a polynomial of degree $\leq \ell -1$ such
that $P_\ell(0) \ne 0$ and $\rho$ is an operation of the form
$(\gamma^\prime)$ above, where $\lambda(z)$ has a zero of order
$\geq \ell/m$ at $z=0$.

\end{itemize}

\end{Theorem}

Using the above M. Suzuki was already able to prove that:
\begin{itemize}
\item[{\rm(iv)}] Any $\bc^*$-action on $\bc^2$ is analytically linearizable, i.e., analytically equivalent/conjugate to
$s \circ (x,y) = (s^nx, s^my), \,\, s \in \bc^*, \,\, (x,y) \in
\bc^2$ for some $n,m \in \bn$.
\end{itemize}
Already present in some of Suzuki's works is the viewpoint of {\it
foliation with singularities.} Let us therefore denote by $\fa_\vr$
the one-dimensional holomorphic foliation with singularities defined
by $\vr$ on $N$:\, $\fa_\vr$ is the foliation $\fa(Z)$ where $\tau =
\dfrac{\po\vr}{\po t}\bigg\vert_{t=0}$\,; the leaves of $\fa_\vr$
are the orbits of $\vr$ and $\sing(\fa_\vr) = \sing(\tau) =
\Fix(\vr)$.

\begin{Example}{\rm  Let $(x,y) \in \bc^2$ be affine coordinates and define
$Z(x,y) = x\,\dfrac{\po}{\po x} + \la y\, \dfrac{\po}{\po y}$ where
$\la \in \bq_+$\,;  say $\la = \dfrac nm$\,, $\lg n,m\rg=1$,\, $n,m
\in \bn$. Then $Z$ is complete and the flow of $Z$ defines a
holomorphic action $\vr\colon \bc\times\bc^2 \to \bc^2$ given by
$\vr_t(x,y) = \big(xe^t, ye^{\frac nm\,t}\big)$. The orbits are all
periodic (diffeomorphic to $\bc^*$) and the origin is a singularity
which is non-degenerate and dicritical for the foliation $\fa(Z)$. }
\end{Example}
\noindent A converse of this example is as follows:

\begin{Theorem} [Global linearization theorem, \cite{scarduaactions}]
\label{Theorem:flows} Let $N$ be a connected Stein surface with
$\overset{\vee}{H^2}(N,\bz) = 0$ equipped with a holomorphic action
$\vr\colon \bc\times N \to N$, with isolated singularities, having a
non-degenerate dicritical singularity $p_0 \in N$. Then $N$ is
biholomorphic to $\bc^2$, indeed there is a biholomorphic map
$\Phi\colon N \to \bc^2$ that conjugates $\vr$ to an action
$\bc\times\bc^2\to\bc^2$,\, $(t,(x,y)) \mapsto \big(xe^{\la_1t},
ye^{\la_2T}\big)$ for some $\la_1,\la_2 \in \bc^*$ with $\la_1/\la_2
\in \bq_+$.
\end{Theorem}

\subsection{Proof of the Global linearization theorem}
Let us give an idea of the proof of Theorem~\ref{Theorem:flows}.
First we  make a simple remark:

\begin{Claim}  The generic leaf is diffeomorphic to $\bc^*$.
\end{Claim}

\begin{proof} Take a leaf $L$ that accumulates at $p_0$ and is contained in
a separatrix $\Ga_{p_0} \ni p_0$\,. We have two possibilities for
$L$:\,\, $L \simeq \bc$ or $L \simeq \bc^*$. If $L \simeq \bc$ then
in a suitable neighborhood of $p_0$\,\, $L \supset \Ga_{p_0} \cup
\{p_0\}$ so we obtain a holomorphic map $p\colon \overline{\bc} \to
N$ holomorphic and non-constant such that $p\big\vert_\bc(\bc)
\subset L$. This is not possible because $N$ cannot contain a
compact holomorphic curve

\noindent Thus every leaf containing a separatrix is a periodic
orbit. Since $p_0 \in \sing(\fa_\vr)$ is dicritical its reduction of
singularities exhibits some non-invariant projective line and
therefore we have an invariant sector $p_0 \in \eS_{p_0} \subset N$
with vertex at $p_0$\,, such that every orbit in $\eS_{p_)}$ is
dicritical and therefore periodic.

\noindent Because $\eS_{p_0}$ has non-empty interior we conclude
that the orbits of $\fa_\vr$ in $\eS_{p_0}$ are not contained in the
exceptional set $\E \subset N$, which has zero logarithmic capacity.
Therefore the generic orbit of $\vr$ is periodic. \quad \end{proof}

Thanks to Suzuki's theorems~\ref{Theorem:Suzuki1} and
~\ref{Theorem:Suzuki2} we conclude that $\fa_\vr$ has a meromorphic
first integral say $f\colon M \to \bc$. If $f$ is holomorphic at
$p_0$ then $\fa_\vr$ has only finitely many separatrices at $p_0$\,,
contradiction. Thus $p_0$ belongs to the set of points of
indetermination of $f$:\, the germ of $f$ at $p_0$ is the quotient
$f = h/g$, \, $h,g \in O_2(p_0)$ with $\lg h,g\rg=1$ and $h(p_0) =
g(p_0) = 0$. Therefore all the orbits of $\vr$ that approach $p_0$
are of type $\bc^*$ and contain separatrices of $\fa_\vr$ at
$p_0$\,.

\noindent Indeed, since $Z$ has a non-degenerate singularity with a
meromorphic first integral at $p_0$ we conclude that:

\noindent $\bullet$\, $DZ(p_0)$ is in the Poincaré-domain,
analytically linearizable and

\noindent $\bullet$\, $Z = \la_1\,\widetilde{x}\,
\dfrac{\po}{\po\widetilde{x}} + \la_2\,\widetilde{y}\,
\dfrac{\po}{\po \widetilde{y}}$ in some local coordinates
$(\widetilde{x},\widetilde{y}) \in V \ni p_0$\,; for some
$\la_1,\la_2 \in \bc^*$ with $\la_1/\la_2 \in \bq_+$\,.

\noindent We than consider the attraction basin of $p_0$ denoted by
$B_{p_0}(Z) = B_{p_0}$ as in the usual Dynamical Systems context.

\noindent Because of the above local linearization we have:

\noindent $\bullet$\, $B_{p_0}$ is an open subset that contains a
neighborhood of $p_0$ in $N$.

\noindent $\bullet$\, The flow of $Z$ is analytically conjugate to
the linear flow (on $B_{p_0})$ given by
$$
Z_{\la_1,\la_2} = \la_1\,x\,\frac{\po}{\po x} +
\la_2\,y\,\frac{\po}{\po y}
$$

\noindent The theorem will be proved if we prove that $B_{p_0} = N$:
indeed, the attraction basin $B_0(Z_{\la_1,\la_2})$ of the origin $0
\in \bc^2$ for the linear vector field $Z_{\la_1,\la_2}$ is $\bc^2$.

\noindent For this we shall prove:

\begin{Lemma}
\label{Lemma:boundarydescription} The boundary $\po B_{p_0}$ is a (possibly empty)
union of isolated singular points and if invariant  analytic curves,
each curve accumulating at a unique non-dicritical singularity of
$Z$.
\end{Lemma}

\noindent In order to prove this lemma we go step by step:

\begin{Claim}  A leaf $L_0 \subset \po B_{p_0}$ diffeomorphic to $\bc^*$ cannot be closed.
\end{Claim}

\begin{proof} Let $L_0 \subset \po B_{p_0}$ be a closed leaf in $N$ which is periodic
(diffeomorphic to $\bc^*$). Then $L_0$ is an analytic smooth curve
in $N$ (recall that, since  $\fa_\vr$ has a meromorphic first
integral, all leaves have analytic closure of dimension one in $N$).
Since $N$ is Stein and $\overset{\vee}{H^2}(N,\bz)=0$ we can take a
holomorphic reduced equation $L_0 := \{h=0\}$, where $h\colon N \to
\bc$ is holomorphic. Because $L_0 \simeq \bc^*$ which is
homeomorphic to the cylinder $\eS^1\times\re$ we can take $\ga\colon
S^1 \to L_0$ generator of the homology of $L_0$\,; and a holomorphic
$1$-form $\al$ on $L_0$ such that $\displaystyle\int_\ga \al=1$

\noindent Since $N$ is Stein by a theorem of Cartan there is an
extension $\widetilde{\al}$ of $\al$ to $N$. Then $\widetilde{\al}$
is a holomorphic $1$-form on $N$ such that
$\widetilde{\al}\big\vert_{L_0} = \al$ and therefore
$\displaystyle\int_\ga \widetilde{\al}=1$. Choose now a transverse
disc $\Sigma$ centered at some point $q_0 \in \ga(S^1) \subset L_0$
and consider the holonomy map $f_\ga\colon (\Sigma,q_0) \to
(\Sigma,q_0)$ induced by $\ga$. Since $\fa_\vr$ admits a meromorphic
first integral this holonomy map $f_\ga$ is periodic. This means
that there $n \in \bn$ such that the lift $\widetilde{\ga}_z^n$ of
$\ga^n$ to the leaf $L_z \ni z$ for $z \in \Sigma$, $z \approx
q_0$\,, is such that $\widetilde{\ga}_z^n$ is closed. Let us assume
for simplicity that $n=1$.

Then the lifts $\widetilde{\ga}_z$ of $\ga$ one closed paths. Since
$L_0 \subset \po B_{p_0}$ we have $\Sigma \cap B_{p_0} \ni q_0$\,,
i.e., we may choose $z \in \Sigma$ arbitrarily close to $q_0$ such
that $L_z \subset B_{p_0}$\,.

\noindent Since $L_z \simeq \bc^*$ and $\overline{L}_z = L_z \cup
\{p_0\} \cong \bc$ is simply-connected we must have
$\displaystyle\int_\delta \om = 0$ for closed $1$-form $\omega$ in
$L_z$ and all closed path $\delta$ in $L_z$\,. Since
$\widetilde{\al}$ is holomorphic its restriction $\om :=
\al\big\vert_{L_z}$ is closed. Thus we conclude (because
$\widetilde{\ga}_z$ is closed) that
$\displaystyle\int_{\widetilde{\ga}_z} \widetilde{\al} = 0$,
$\forall\, z \in \Sigma \cap B_{p_0}$\,. On the other hand for $z
\in \Sigma$ close enough to $q_0$ we have
$\left\vert\displaystyle\int_{\widetilde{\ga}_z} \widetilde{\al} -
\displaystyle\int_\ga \al\right\vert < \dfrac 12$ so that
$\displaystyle\int_{\widetilde{\ga}_z} \ne 0$, contradiction. \quad
\end{proof}

\begin{Claim}  All leaves of $\fa_\vr$ are biholomorphic to $\bc^*$.
\end{Claim}

\begin{proof} Indeed, in $B_{p_0}$ the flow is conjugate to a periodic flow so
that in $B_{p_0}$ it has a certain period say
$\vr_\tau\big\vert_{B_{p_0}} \equiv \Id$.  The Identity Principle
implies that $\vr_\tau \equiv \Id$ in $N$. \quad\end{proof}

\begin{Claim}  Each leaf $L_0 \subset \po B_{p_0}$ accumulates a unique singularity of $Z_\vr$ and this singularity is non-dicritical.
\end{Claim}

\begin{proof} We already know that $L_0$ is not closed so that it accumulates some singularity
say $\widetilde{p}_0 \in \overline{L}_0\setminus L_0$\,,
$\widetilde{p}_0 \in \sing(\fa_\vr)$.  If $\widetilde{p}_0$ is
dicritical then the meromorphic first integral $f$ has an
indeterminacy point at $\widetilde{p}_0$ and by writing $f = g/h$\,
$g,h \in O_2(\widetilde{p}_0)$ with $g(\widetilde{p}_0) =
h(\widetilde{p}_0) = 0$ we conclude that there is an open subset
$V_{\widetilde{p}_0}$ (like $B_{\widetilde{p}_0}$) such that every
leaf accumulating at $\widetilde{p}_0$ is contained in a separatrix
of $\fa_\vr$ through $\widetilde{p}_0$\,.

\noindent In particular, because $L_0 \subset \po B_{p_0}$\,, there
are leaves $L \subset B_{p_0}$ that contain separatrices through
$p_0$ and through $\widetilde{p}_0$\,, such a leaf cannot exist
because $\overline{L}$ would be compact (a Riemann Sphere)

\noindent Thus $\widetilde{p}_0$ is non-dicritical. The same type of
argument shows that $\overline{L}_0\setminus L_0$ is a single point.

\end{proof}
\begin{Claim} \,\, $\po B_{p_0}$ contains no isolated points.
\end{Claim}

\begin{proof} This is for topological reasons: If $q_0 \in \po B_{p_0}$ is an isolated point,
then since $B_{p_0}$ is diffeomorphic to $\bc^2$ (the basin of
$Z_{\la_1,\la_2}$ in $\bc^2$ is $\bc^2$)  we conclude that $\po
B_{p_0} = \{q_0\}$ and therefore $N$ is compact homeomorphic to $S^4
= \mathbb R^4 \cup \{\infty\}$, contradiction because $N$ is
compact.
\end{proof}

Thus we have proved that $\po B_{p_0}$ {\it is a discrete union of
analytic curves} and proved Lemma~\ref{Lemma:boundarydescription}.

\noindent Let us now proceed:

\begin{Claim}
\label{Claim:union} \,\, $N = B_{p_0} \cup \po B_{p_0}$\,.
\end{Claim}

\begin{proof} Put $A = N\setminus \po B_{p_0}$ and $B = B_{p_0}$\,.

\noindent Then $A$ and $B$ are open subsets of $N$; \, $B$ is
connected (diffeomorphic to $\bc^2$) and so is $A$ ($\po B_{p_0}$ is
analytic of dimension $\le 1$). Moreover $A \supset B$ and $\po A =
\po B$. Then since $N$ is connected we have $A=B$, i.e., $N =
B_{p_0} \cup \po B_{p_0}$\,. \quad \end{proof}

\begin{Claim}
\label{Claim:boundaryempty} \,\, $\po B_{p_0} = \emptyset$.
\end{Claim}

\begin{proof} Suppose that $\po B_{p_0} \supset \overline{L}_0$ for some leaf $L_0$\,.
 Let $L_0$ be given by the reduced equation $\{h=0\}$, $h\colon N \to \bc$ holomorphic.
 Put $\al := \dfrac{dh}{h}$ then this a closed meromorphic $1$-form on $N$ with polar set $(\al)_\infty = \overline{L}_0$\,.

\noindent Choose a transverse disc $\Sigma$ to $\fa_\vr$ centered at
a point $q_0 \in L_0 \cap \Sigma$ and a loop $\ga\colon S^1 \to
\Sigma\setminus\{q_0\}$

\noindent that avoids the real codimension $\ge2$ subset $\po
B_{p_0}$\,; $\ga(S^1) \subset B_{p_0}$\,. Then because $B_{p_0}
\simeq \bc^2$ we have $\displaystyle\int_\ga \al=0$. On the other
hand a straightforward computation shows that $\displaystyle\int_\ga
\al = \displaystyle\int_\ga \dfrac{dh}{h} = 2\pi i$, contradiction.
\quad\end{proof}

From Claims~\ref{Claim:union} and ~\ref{Claim:boundaryempty} we
conclude that $N=B_{P_0}$. The proof of  Theorem~\ref{Theorem:flows}
is now complete.

\section{Real transverse sections of holomorphic foliations}
\label{section:ito}

We discuss some problems related to the interplay between geometric
theory of foliations and Holomorphic foliations with singularities.
We start by recalling that a foliation $\fa$ on a manifold $M$ is
{\it transverse\/}\index{manifold transverse to a holomorphic
foliation} to a submanifold $N \subset M$ if for every point $p \in
N$ we have $T_p(\fa) + T_p(N) = T_p(M)$ where $T_p(\fa) \subset
T_p(M)$ is the {\it tangent space\/} of $\fa$ at $p$, defined by
$T_p(\fa) = T_p(L_p)$ where $L_p \ni p$ is the leaf of $\fa$ that
contains $p$. In particular, if $\fa$ is singular, then $\sing(\fa)
\cap N = \phi$.

\noindent In the real codimension one case the existence of a
compact transverse submanifold to the foliation is an important
object in the study of the foliation dynamics. As example of this is
given by

\begin{Theorem} [Haefliger's Theorem, \cite{camacho-linsneto,Godbillon}\index{Theorem! of Haefliger}]
Let $\fa$ be a {\tt real codimension one} foliation of class $C^2$
on a manifold $M$. Suppose that  there exists a an immersed curve
$\ga\colon S^1 \hookrightarrow M$ such that:
\begin{itemize}
\item[{\rm(i)}] $\ga$ is homotopic to zero in $M$.
\item[{\rm(ii)}] $\ga$ is transverse to $\fa$.
\end{itemize}
Then there exists a leaf $L_0 \in \fa$ and a closed path
$\delta\colon S^1 \to L_0$ such that the holonomy map
$h_\delta\colon (\re,0) \to (\re,0)$ corresponding to $\delta$ is a
{\it one-sided\/} map: $\exists\, \ve > 0$ such that, up to change
of orientation on $\ga$ transverse section $\Sigma$ to $L_0$ at
$\delta(0) = \delta(1)$, we have
$h_\delta\big\vert_{\Sigma^+\setminus\{0\}} \equiv \Id$ and
$h_\ga\big\vert_{\Sigma^-\setminus\{0\}}$ is a contraction.
\end{Theorem}

\noindent $\Sigma = \Sigma^+ \cup \{\delta\{0\}\} \cup \Sigma^-$
\begin{figure}[ht]
\begin{center}
\includegraphics[scale=0.6]{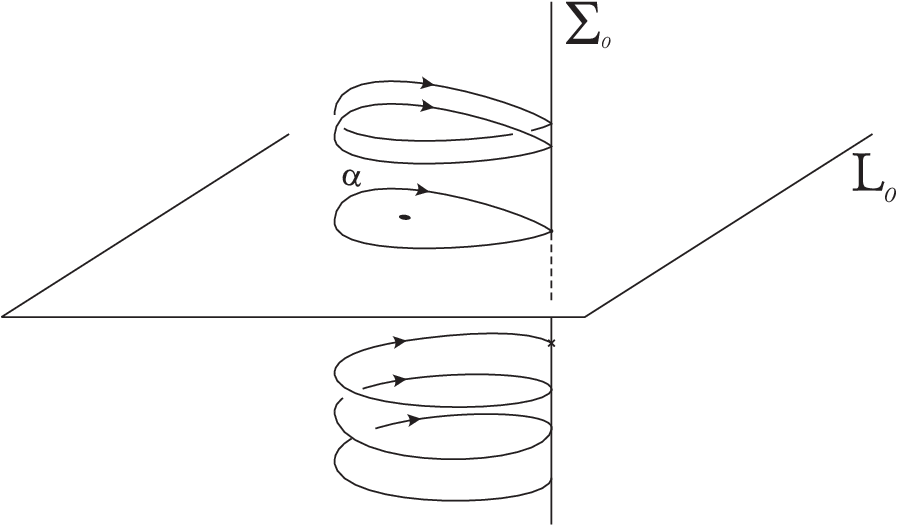}
\caption{}
\end{center}
\end{figure}

\noindent As an important corollary we get:

\begin{Corollary} If $M$ is a compact (real) manifold of dimension $m \ge 2$ with finite
fundamental group, then $M$ does not admit a {\rm(}{\tt non-singular}{\rm)} {\tt real analytic} foliation of {\tt real codimension one}.
\end{Corollary}

\noindent Indeed the above results rely on the following remark:

\begin{Lemma} Let $\fa$ be a real codimension one foliation of a manifold $M$
and $L$ a leaf of $\fa$ such that $\overline{L}\setminus L \ni p$. Then there is a
transverse closed curve $\ga\colon S^1 \to M$, to $\fa$ that intersects the leaf $L$.
\end{Lemma}

\noindent In particular in a compact manifold $M$ any real
codimension one has some closed transverse curve.

Another important feature in the above results of Haefliger is the
following:

Let $X$ be a $C^1$ (real) vector field defined in a neighborhood $U$
of the closed disc $\overline{D}^2 \subset \mathbb R^2$,\, $D^2 =
\{(x,y) \in \mathbb R^2\colon x^2+y^2 \le 1\}$. Assume that: (i)\,
the singularities of $X$ inside the disc are of {\it Morse type\/}
(i.e., centers $x_1^2+x_2^2=c$ or saddles $x_1^2-x_2^2=c$) and
(ii)\, $X$ is transverse (pointing inwards) to the boundary $\eS^1 =
\po \overline{D}^2$. \noindent Then $X$ has some limit cycle or some
graph $\Ga \subset D^2$ with one sided (holonomy) Poincaré map.

\begin{figure}[ht]
\begin{center}
\includegraphics[scale=0.4]{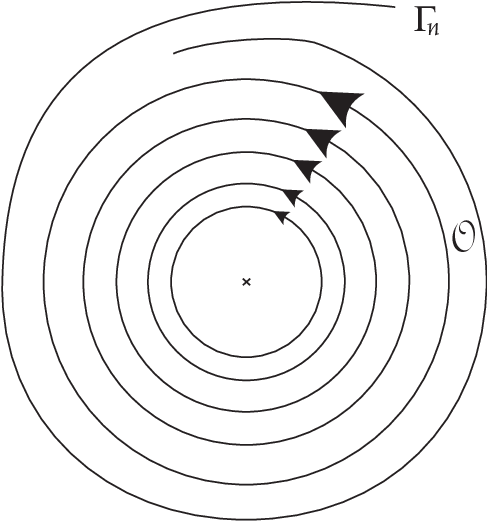}
\caption{}
\end{center}
\end{figure}

\noindent These are fundamental steps in the way to the proof of the
celebrated

\begin{Theorem}[Novikov compact leaf theorem, \cite{camacho-linsneto,Godbillon}\index{Theorem! of Novikov}]
{\it Let $M^3$ be a compact real $3$-manifold with finite
fundamental group. Any foliation of real codimension one on $M^3$
has a compact leaf diffeomorphic to the $2$-torus $S^1\times S^1$ or
has all its leaves compact.}
\end{Theorem}

\noindent Indeed, if not all leaves are compact then $\fa$ has some
{\it Reeb component\/} which is a region on $M$, invariant by $\fa$,
diffeomorphic to the solid torus $S^1 \times \overline{D}^2$ where
the boundary $\po(S^1\times\overline{D}^2) = S^1 \times S^1$ is a
leaf and the interior leaves are diffeomorphic to $\mathbb R^2$.
\noindent Resuming the study of holomorphic foliations we are
therefore interested in the consequences of the existence of compact
transverse sections on the dynamics of the foliation. For some
reasons we allow the transverse manifold to have codimension
different from the dimension of the foliation. We start with the
simplest case: Let $Z$ be a holomorphic vector field on a
neighborhood $0 \in U \subset \bc^n$, $n \ge 2$. Let $0 < r$ be such
that $\overline{B(0,r)} = \{z \in \bc^n; |z| \le r\} \subset U$.

\begin{Definition} {\rm We say that $Z$ is {\it transverse\/}\index{vector field! holomorphic ! transverse to a sphere}
to $S^{2n-1}(0;r) = \po B(0;r)$ at a point $z \in S^{2m-1}(0;r)$ if
the $2$dimensional real foliation induced by $Z$ on
$U\setminus\sing(Z)$ is transverse to the real submanifold
$S^{n-1}(0;r) \subset U \subset \mathbb R^{2n}$ at the point $z$. It
this happens for all $z \in S^{2n-1}(0;r)$ then we shall say that
$Z$ is transverse to $S^{2n-1}(0;r)$  and write $Z
\overline{\pitchfork} S^{2n-1}(0;r)$. }
\end{Definition}
\noindent A characterization of this situation in terms of
holomorphic coordinates is:

\begin{Lemma} Write $Z = \sum\limits_{j=1}^n A_j \, \dfrac{\po}{\po z_j}\,\cdot$
Then $Z \overline{\pitchfork} S^{2n-1}(0;r) \Leftrightarrow \sum\limits_{j=1}^n \bar{z}_j\,A_j \ne 0$,\linebreak $\forall\, z \in S^{2n-1}(0;r)$.
\end{Lemma}

\noindent The geometric interpretation is as follows: Let
$\overset{\rightarrow}{R} = \sum\limits_{j=1}^n z_j\,
\dfrac{\po}{\po z_j}$ be the complex radial vector field on $\bc^n$.
Then $Z \overline{\pitchfork} S^{2n-1}(0;r) \Leftrightarrow \lg
Z,\overset{\rightarrow}{R}\rg \ne 0\,\, \forall\, z \in
S^{2n-1}(0;r)$, where $\lg\,,\,\rg$ denotes the usual hermitian
product on $\bc^n$. In particular a vector field $Z$ with a
non-degenerate singularity at origin say $\Spec(D\,Z(0))=
\{\la_1,\dots,\la_n\} \subset \mathbb R^2$, is in the
Poincaré-domain if and only if $Z \overline{\pitchfork}
S^{2n-1}(0;\ve)$ for every $0 < \ve$ small enough.

\noindent\textbf{Question:} What happens if $Z$ is transverse to a
``big sphere"\,?

\noindent The answer to this question was given by Toshikazu Ito in
1992:

\begin{Theorem} [Ito,\cite{Ito}\index{Theorem! of Ito}] {\it Let $Z$ be a holomorphic vector field on
a neighborhood $U$ of $\overline{B(0,R)}$ in $\bc^n$, $n \ge 2$. Suppose that
$Z \overline{\pitchfork} \po\overline{B(0,R)} = S^{2n-1}(0,R)$. Then:}
\begin{itemize}
\item[{\rm(i)}] {\it $Z$ has exactly one singular point $o \in \sing(Z)$ in the ball $\overline{B(0.R)}$.}
\item[{\rm(ii)}] {\it The singularity $o \in \sing(Z)$ is in the Poincaré domain.}
\end{itemize}
{\it By a Möbius transformation we may assume that $o$ is the origin
$o = 0 \in \bc^n$. In this case:}
\begin{itemize}
\item[{\rm(iii)}] {\it $Z$ is transverse to {\tt all} spheres $S^{2n-1}(0;r)$, $0 < r \le R$.}
\item[{\rm(iv)}] {\it There is a real analytic conjugation between the flow of $Z$ in $B(0;r)\setminus\{0\}$
and the product $[0,\infty)\times\L$ where $\L$ is the real (transversely holomorphic) flow induced by $Z$ in $S^{2n-1}(0;r)$.} \newline
{\it In particular}
\item[{\rm(v)}] {\it Each orbit of $Z$ in $B(0;R)$ accumulates at the singular point $o$.}
\end{itemize}
\end{Theorem}
\noindent Thus the one-dimensional case is somehow well-understood.
On the other hand  the codimension one case still remains open:

\begin{Question} {\it Let $n\ge3$. Is there a codimension one holomorphic foliation $\fa$ of
a neighborhood $B(0;r) \subset 0 \in U \subset \bc^n$, with
$S^{2n-1}(0,R) \overline{\pitchfork} \fa$?}
\end{Question}

\noindent There is a number of partial results, all of them
suggesting that the answer to the above question is NO. One of the
most important is:

\begin{Theorem} [\cite{Ito-Scardua}]
Let $\Omega$ be a holomorphic $1$-form on a neighborhood $0 \in
\overline{B(0,R)} \subset U \subset \bn^n$, $n\ge3$. Suppose that
the distribution $\Ker(\Omega)$ is transverse to $S^{2n-1}(0;R)$.
Then:
\begin{itemize}
\item[{\rm(i)}] $n$ is even (therefore $n\ge4$);
\item[{\rm(ii)}] there exists exactly one singular point $q \in \sing(\Om) \cap B(0;R)$ and this point is {\tt simple}.
\end{itemize}
\end{Theorem}

\noindent If we write $\Om = \sum\limits_{j=1}^n f_j(z)dz_j$ with
$f_j\colon U \to \bc$ holomorphic, then $\sing(\Om) =$\linebreak
$\{z_0 \in U; f_j(z_0)=0, \forall\,j\}$. A singular point $z_0 \in
U$ is {\it simple\/}\index{singularity! simple} if
$$
\Det\left(\frac{\po f_i(z_0)}{\po z_j}\right)_{i,j=1,\dots,n} \ne 0.
$$

\begin{Remark}{\rm

\noindent (1)\, We do not assume that $\Om$ is integrable.

\noindent (2)\, The distribution $\Ker(\Om)$ is defined by: given $p
\in U$ then $\Ker(\Om)(p) := \left\{v \in T_p(\bc^n); \Om(p)\cdot
v=0\right\}$.

\noindent (3)\, By a classical result of Malgrange
(\cite{Malgrange}) a holomorphic foliation of codimension-one admits
a holomorphic first integral in a small neighborhood of a point
where the singular set of the foliation has codimension $\ge3$.
Therefore, in the above statement if $\Om$ is integrable and $q =
\sing(\Om) \cap B(0;R)$ then the foliation $\fa\colon (\Ker(\Om)):
\Om = 0$ admits a holomorphic $f$ first integral: $f\colon V_p \to
\bc$,\, $p \in V_p \subset U$; of {\it Morse type\/} so that $f =
\sum\limits_{j=1}^n \widetilde{z}_j$ in some local coordinates
$(\widetilde{z}_1,\dots,\widetilde{z}_n)$. }\end{Remark}

\noindent Nevertheless, because of the Maximum Principle, $\fa$
cannot be transverse to any sphere in $V_p$\,. This already shows
that the codimension one situation is (in view of Ito's theorem
above) pretty different from the one-dimensional case. Just to
mention an (non-integrable) example we take $\Om =
\sum\limits_{j=1}^n \left(z_{2j}\,dz_{2j-1} - z_{2j} +
dz_{z_j}\right)$ then $\sing(\Om) = \{0\}$ and
$$
\Ker(\Om) \cap S^{2n-1}(0;R), \,\,\, \forall\, R > 0.
$$
Indeed $\Om$ is a symplectic $1$-form associated to a contact
structure on the sphere $S^{2n-1}(0;R)$, $\forall\, R > 0$.

Despite some advance in the subject, we still do not know the answer
to the following question: \vglue.1in \noindent{\bf Question}: {\sl
Is there a codimension one holomorphic foliation $\fa$ defined in a
neighborhood of a closed ball $B\subset \mathbb C^n, n \geq 3$, such
that the foliation is transverse to the boundary sphere $\partial
B$?}

\section{Non-trivial minimal sets of holomorphic foliations}
\label{section:minimalsets} Given a foliation $\fa$ regular, of any
class of differentiability, of $M$, we call a {\it minimal
set}\index{minimal set! of a foliation} of $\fa$ a closed invariant
subset ${\mathcal M} \subset M$ such that ${\mathcal M}$ is minimal
with this property.

\begin{Remark}
{\rm If $M$ is compact, then it is fairly well-known that any
(regular)
foliation exhibits minimal sets.}
\end{Remark}

Let now $\fa$ be a foliation with singular  set $\sing(\fa)$ on a
compact manifold $M$.

\begin{Definition}
A non-trivial minimal set of $\fa$ is a minimal set of
$\fa^1=\fa|_{M\setminus\sing(\fa)}$ on $M\setminus\sing(\fa)$.
\end{Definition}

If $\fa$ has some leaf $L_0$ such that
$\bar{L}_{0}\cap\sing(\fa)\neq\emptyset$ then $m=\bar{L}_{0}$ is a
non-trivial minimal set of $\fa$.  Thus, the problem of existence of
non-trivial minimal sets in the case of singular foliations is
equivalent to know whether all leaves must accumulate singularities
or not.

\begin{Problem}[Camacho-Lins Neto-Sad, \cite{C-LN-S3}]
Is there any non-trivial minimal set for a foliation $\fa$ on ${\co}
P(2)$?
\end{Problem}

\begin{Theorem}[Camacho-Lins Neto-Sad, \cite{C-LN-S3}]
Let $\fa$ be a foliation with a non-trivial minimal set $\mathcal M
$ on ${\co} P(2)$. Then
\begin{enumerate}
\item $\mathcal M$ is unique.
\item If $C\subset {\co} P(2)$ is an algebraic curve then
${\mathcal M}
\cap C\ne\emptyset$.
\item $\fa$ has no algebraic invariant curve.
\item There exists an hermitian metric on ${\co} P(2)\setminus\sing(\fa)$
such
that it is complete and induces negative curvature
$K\leq-\epsilon^2<0$ on the leaves of $\fa|_{\mathcal M}$.  In particular
each leaf $L\subset {\mathcal M}$ is covered by the disc ${\mathbb D}$ and the
family of uniformizations is normal.
\item All leaves $L \subset {\mathcal M}$ have exponential growth.
\end{enumerate}
\end{Theorem}

Given a $3$-dimensional real analytic hypersurface $N^3\subset{\co}
P(2)$, we assume that $N^3 \subset{\co} P(2)$ is Levi-flat.  This
vanishing of its Levi-form is equivalent to the following: for any
$p\in N^3$, there exists a local chart at $p$, $(x,y)\in U$ for
${\co} P(2)$ such that $N \cap U = \{{\rm Im.}y=0\}$. This implies
the following:

\begin{Theorem}[Lins Neto, \cite{LN2}]
There exists a holomorphic foliation $\fa$ on ${\co} P(2)$ such that
$N^3$ is $\fa$-invariant and $\sing(\fa)\cap N^3=\emptyset$. In
particular, $N^3$ gives a non-trivial minimal set of $\fa$.
\end{Theorem}

Thus, the problem of the existence of non-trivial minimal sets for
foliations on ${\co} P(2)$ is also related to the existence of
Levi-flat real hypersurfaces (or submanifolds) $M^3 \subset {\co}
P(2)$. It is still not known whether a codimension one holomorphic
foliation on $\mathbb CP(n)$ can have a non-trivial  minimal  set.

\section{Transversely homogeneous holomorphic foliations}
\label{section:homogeneous} From the structural point of view the
simplest foliations are those with an homogeneous transverse
structure (cf. \cite{Godbillon}). A  holomorphic foliation $\fa$ on
a smooth manifold $M$ has a {\it holomorphic homogeneous transverse
strucutre\/}\index{transverse structure! homogeneous} if there are a
complex Lie group $G$, a connected closed subgroup $H < G$ such that
$\fa$ admits an atlas of submersions $y_j\colon U_j \subset M \to
G/H$ satisfying $y_i = g_{ij}\circ y_j$ for some locally constant
map $g_{ij}\colon U_i \cap U_j \to G$ for each $U_i \cap U_j \ne
\emptyset$. In other words, the transversely holomorphic atlas of
submersions for $\fa$ has transition maps given by left translations
on $G$ and submersions taking values on the homogeneous space $G/H$.
We shall say that $\fa$ is transversely homogeneous {\it of model\/}
$G/H$. Some important properties of transversely homogeneous
holomorphic foliations are listed below:
\begin{enumerate}

\item Any transversely homogeneous holomorphic foliation is
a transversely holomorphic foliation with a holomorphic homogeneous
transverse structure.

\item  Given a foliation $\fa$ of $M$ as in (1) with model
$G/H$ then any real submanifold $M \subset M$ transverse to $\fa$ is
equipped with a transversely holomorphic foliation $\fa_1 = \fa|_M$
with holomorphic homogeneous transverse structure of model $G/H$.

\item  Let $F = G/H$ be an homogeneous space of a
complex Lie group $G$ ($H\triangleleft G$ is a closed Lie subgroup).
Any homomorphism representation $\vr\colon \pi_1(N) \to \Aut(F)$
gives rise to a transversely holomorphic foliation $\fa_\vr$ on
$({\widetilde N\times F})/{\Phi} = M_\vr$ which is holomorphically
transversely homogeneous of model $G/H$.

\item For the case $G=P\SL(2,\mathbb C)$ and $H\subset G$ is the affine group $H=\Aff(\mathbb C)$
(isotropy group of the point at infinity $\infty \in \mathbb C
P(1)$), we have that the quotient $G/H\simeq \mathbb CP(1)$  is the
Riemann sphere and the foliations with this transverse model are
called {\it transversely projective}\index{foliation! transversely
projective}.
\end{enumerate}

Adapting the above notion for the case of holomorphic foliation with
singularities we have:
\begin{Definition}[transversely homogeneous holomorphic foliation with singularities]
\label{Definition:transprojfolnonsing} {\rm A holomorphic foliation
$\fa$ with singularities on a manifold $M$. We shall say  that $\fa$
is {\it transversely homogeneous}\index{foliation! transversely
homogeneous} of model $G/H$ if the underlying non-singular foliation
is transversely homogeneous of model $G/H$ on $M\setminus
\sing(\fa)$. In particular, $\fa$  is called {\it transversely
projective} if there is an open cover $\bigcup\limits_{j\in J} U_j =
M\setminus \sing(\fa)$ such that in each $U_j$ the foliation is
given by a submersion $f_j\colon U_j \to \ov{\mathbb C}$ and if $U_i
\cap U_j \ne \emptyset$ then we have $f_i = f_{ij}\circ f_j$ in $U_i
\cap U_j$ where $f_{ij}\colon U_i \cap U_j \to P\SL(2,{\mathbb C})$
is locally constant. Thus,  on each intersection $U_i \cap U_j \ne
\emptyset$, we have $f_i =
\frac{a_{ij}f_j+b_{ij}}{c_{ij}f_j+d_{ij}}$ for some locally constant
functions $a_{ij}, b_{ij}, c_{ij}, d_{ij}$ with $a_{ij}d_{ij} -
b_{ij}c_{ij} = 1$. The data $\mathcal P = \{U_j, f_j, f_{ij}, j \in
J\}$ is called a {\it projective transverse structure} for $\fa$. }
\end{Definition}

Basic references for transversely affine and
transversely projective foliations (in the nonsingular case) are
found in \cite{Godbillon}.

\begin{itemize}

\item[{\rm(5)}]  Based on the Rieman-Koebe uniformization theorem we have:
\begin{Proposition}[\cite{Scardua1}, Theorem 6.1 page 203).] Let $\fa$ be a holomorphic
singular transversely homogeneous foliation of codimension one on
$M^n$. Then $\fa$ is a transversely projective foliation on $M^n$.
\end{Proposition}
\end{itemize}
\begin{proof} We know that $G/ H$ is a simply-connected complex manifold of dimension one.  By the Riemann-Koebe
uniformization theorem we have a conformal equivalence  $G/H \equiv
\ov{\mathbb C}, \mathbb C$ or $\mathbb D$ the unitary disc. This
implies that either $G\subset \Aut(\ov{\mathbb C})=\mathbb P
SL(2,\mathbb C), G\subset \Aut(\mathbb C)=\Aff(\mathbb C)$ or
$G\subset \Aut(\mathbb D)\cong \mathbb P SL(2,\mathbb R)$.  The
proposition follows.
\end{proof}

\subsection{Transversely Lie foliations}

Let $\fa$ be a codimension $\ell$ foliation of a manifold $M$. If
$\fa$  admits a Lie group transverse structure of model $G$, or a
{\it $G$-transverse structure} for short, then  we shall  call $\fa$
a {\it $G$-foliation} or, simply, {\it Lie foliation}. The
characterization of $G$-foliations in terms of differential forms is
given below. Let $\{\omega_1,...,\omega_\ell\}$ be a basis of the
Lie algebra of  $G$. Then we have $d\omega_k = \sum\limits_{i<j}
c_{ij} ^k \omega_i \wedge \omega_j$ for a family constants
$\{c_{ij}^k\}$ called the {\it structure constants} of the Lie
algebra in the given basis.

\begin{Theorem}[Darboux-Lie, \cite{Godbillon}\index{Theorem! of Darboux-Lie}]
\label{Theorem:Darboux-Lie} Let  $G$ be a complex Lie group of
dimension $\ell$. Let  $\{\omega_1,...,\omega_\ell\}$ be a basis of
the Lie algebra of  $G$ with  structure constants $\{c_{ij}^k\}$.
Suppose that a  complex manifold $V^m$ of dimension $m \geq \ell$
admits a system of  $1$-forms $\Omega_1,...,\Omega_\ell$ on $M$ such
that:
\begin{itemize}
\item[{\rm(i)}]  $\{\Omega_1,...,\Omega_\ell\}$ is a rank $\ell$
integrable system which defines $\fa$.

\item[{\rm (ii)}] $d\Omega_k =\sum_{i<j} c_{ij}^k \Omega_i
\wedge \Omega_j$.
\end{itemize}
Then:
\begin{itemize}

\item[{\rm (iii)}] For each point $p\in M$ there is a neighborhood
$p\in U_p \subseteq M$ equipped with a submersion $f_p\colon U_p \to
G$ which defines $\fa$ in $U_p$ such that $f_p^*
(\omega_j)=\Omega_j$ in $U_p$, for all $j\in \{1,...,q\}$.

\item[{\rm(iv)}]  If $U_p \cap U_q \ne \emptyset$ then in the
intersection we have $f_q = L_{g_{pq}}(f_p)$ for some locally
constant left translation $L_{g_{pq}}$ in $G$.
\item[{\rm(v)}] If $M$ is simply-connected we can take $U_p = M$.

\end{itemize}

\end{Theorem}

\section{Transverselly affine foliations}
Let $\fa$ be a codimension one  holomorphic foliation with singularities of $M$. We say that $\fa$ is {\it transversely
additive}\index{foliation! transversely additive} when the maps
$g_{ij}$ in the definition of holomorphic homogeneous transverse
structure are of the form $g_{ij}(z) = z+b_{ij}$\,, \, $b_{ij} \in
\bc$ locally constant in $U_i \cap U_j$\,. If $g_{ij}(z) = a_{ij}\,z
+ b_{ij}$\,, for locally constant $a_{ij} \in \bc-\{0\}$ and $b_{ij}
\in \bc$ we say that $\fa$ is {\it  transversely
affine\/}\index{foliation! transversely affine} and it is {\it
transversely projective\/}\index{foliation! transversely projective}
if $g_{ij}(z) = \dfrac{a_{ij}\,z + b_{ij}}{c_{ij}\,z
+ d_{ij}}$ with locally constant $\begin{pmatrix} a_{ij} &b_{ij}\\
c_{ij} &d_{ij}\end{pmatrix} \in \SL(2,\bc)$.

The problem of deciding wether there exist affine transverse
structures for a given foliation is equivalent to a problem on
differential forms as stated below:

\begin{Proposition}[\cite{Scardua1}]
\label{Propositio:transvaffineforms} The possible holomorphic affine
transverse structures for $\fa$ in $M$ are classified by the
collections $(\Om_j, \eta_j)$ of differential $1$-forms defined in
the open sets $U_j \subset M$ such that:

{\rm(i)}\,\,$\Om_j$ and $\eta_j$ are transversely holomorphic,
$\Om_j$ is integrable and defines $\fa$ in $U_j$\,, \,\, $d\Om_j =
\eta_j \wedge \Om_j$ and $d\eta_j=0$ in $U_j$\,, if $U_i \cap U_j
\ne \emptyset$ then $\Om_i = g_{ij}\,\Om_j$ and $\eta_i = \eta_j +
\dfrac{dg_{ij}}{g_{ij}}$ for non-vanishing transversely holomorphic
function $g_{ij}\colon U_i \cap U_j \to \bc-\{0\}$.

{\rm(ii)} Two such collections $(\Om_j, \eta_j)$ and
$(\Om_j',\eta_j')$ define the same affine transverse structure for
$\fa$ in $M$ if and only if \, $\Om_j' = g_j\,\Om_j$ and $\eta_j' =
\eta_j + \dfrac{dg_j}{g_j}$ for some transversely holomorphic
non-vanishing functions $g_j\colon U_j \to \bc-\{0\}$.
\end{Proposition}

\begin{proof}
First we prove (i). Assume that $\fa$ is transversely affine with
transversely holomorphic atlas of submersions $y_j\colon U_j \to
\bc$. Given any transversely holomorphic non-singular $1$-form
$\Om_j$ defining $\fa$ in $U_j$ we have $\Om_j = g_j\,dy_j$ for some
transversely holomorphic function $g_j\colon U_j \to \bc-\{0\}$ and
we define $\eta_j = \dfrac{dg_j}{g_j}\,\cdot$ If $U_i \cap U_j \ne
\emptyset$ then $\Om_i = g_{ij}\,\Om_j$ and $y_i = a_{ij}\,y_j +
b_{ij}$ imply $dy_i = a_{ij}\,dy_j$ and therefore $a_{ij}\,g_i =
g_j\,g_{ij}$\,. Thus $\dfrac{dg_i}{g_i} = \dfrac{dg_j}{g_j} +
\dfrac{dg_{ij}}{g_{ij}}$ in $U_i \cap U_j$\,. Clearly $d\eta_j=0$,
\, $d\Om_j = \eta_j \wedge \Om_j$ and $\eta_i = \eta_j +
\dfrac{dg_{ij}}{g_{ij}}\,\cdot$ This proves the first part of (i).
Let us prove the second part of (i), i.e., the converse part. For
this we assume that we have a unique $1$-form  $\Om$ which is
integrable meromorphic on $M$ and defines $\fa$ outside the polar
divisor $(\Om)_\infty$. Assume then that $\Omega$ and $\eta$ are as
in the statement. Since $\eta$ is holomorphic and closed in
$M\backslash(\Om)_\infty$, there exists an open cover $\{U_i\}_{i
\in I}$ of $M\backslash(\Om)_\infty$ and there are holomorphic
functions $h_i \in \Hol(U_i)$ such that $\eta\big|_{U_i} = dh_i$. We
define $g_i = \exp(h_i)$, $g_i \in \O(U_i)^*$ to obtain
$\eta\big|_{U_i} = dg_i/g_i$. From $d\Om = \eta\wedge\Om$ we obtain
$d\left(\frac{\Om}{g_i}\right) = 0$, and therefore $\Om = g_i\,dy_i$
for some holomorphic function $y_i \in \O(U_i)$. This can be done in
$M\backslash(\Om)_\infty$. Given a point $p_i \in (\Om)_\infty$ we
can choose a local chart $(x,y) \in U_i$ such that $p_i =(0,0)$,
$(\Om)_\infty \cap U = \{y=0\}$ and $\eta(x,y) = -n \frac{dy}{y}
+\frac{df}{f}$ where $n =$ order of $(\Om)_\infty$ and $f \in
\O(U_i)^*$. Therefore we have $\eta = \frac{d(f.y^{-n})}{f.y^{-n}} =
\frac{dg_i}{g_i}$\,, $g_i = f.y^{-n}$. The $1$-form
$\frac{\Om}{g_i}$ is closed and holomorphic so that it can be
written as $\frac{\Om}{g_i} = dy_i$ for some holomorphic $y_i$. We
have covered $M\backslash s(\fa)$ with open sets $U_i$ where we have
the relations $\Om =g_i\,dy_i$, $\eta = \frac{dg_i}{g_i}$\,. In each
$U_i \cap U_j \ne \phi$ we have $\frac{dg_i}{g_i} = \eta =
\frac{dg_j}{g_j}$ and $g_i\,dy_i =\Om =g_j\,dy_j$. The first
equality implies $g_j = a_{ij}.g_i$ for some locally constant
$a_{ij}$ and it follows from the second equality that $dy_i
=a_{ij}\,dy_j$ and then $y_i = a_{ij}\,y_j + b_{ij}$ with $b_{ij}$
locally constant in $U_i \cap U_j$. This shows that $\fa$ is
transversely affine in $M$.

Now we prove (ii). For this sake it is enough to prove the following:
\begin{Claim}
\label{Claim:sameaffine} Two pairs $(\Om,\eta)$ and
$(\Om^\prime,\eta^\prime)$ define the same affine structure for
$\fa$ in $M$ if and only if there exists a meromorphic map $g$ on
$M$  satisfying $\Om^\prime = g\Om$ and $\eta^\prime = \eta +
\frac{dg}{g}$\,.

\end{Claim}
\begin{proof}
Let $(\Om,\eta)$ be
given and let $g\colon M \to \ov\co$ be a meromorphic function. We
define $\Om^\prime = g\Om$ and $\eta^\prime = \eta +
\frac{dg}{g}$\,. Using the same notation above we have
$\eta^\prime\big|_{U_i} = \eta\big|_{U_i} + \frac{dg}{g} =
\frac{dg_i}{g_i} + \frac{dg}{g} = \frac{d(g_ig)}{(g_i.g)}$ and
$\Om^\prime\big|_{U_i} = g.\Om\big|_{U_i} = (gg_i)dy_i$, and this
shows that: $g_i^\prime = a_{ij}\,g_j^\prime \quad\text{and}\quad
y_i^\prime = y_i \quad \text{so that}\quad a_{ij}^\prime = a_{ij}$
and $b_{ij}^\prime = b_{ij}$. Hence, the pairs $(\Om,\eta)$ and
$(\Om^\prime,\eta^\prime)$ define the same transversal structure for
$\fa$ in $M$. Finally, suppose that $(\Om,\eta)$ and
$(\Om^\prime,\eta^\prime)$ define the same transversal structure for
$\fa$ in $M$. Since $\Om$ and $\Om^\prime$ define $\fa$, we have
$\Om^\prime = g\Om$ for some $g\colon M \to \ov\co$ meromorphic.
Using the same notation above we write (locally) $\Om = g_i\,dy_i$,
$\Om^\prime = g_i^\prime\,dy_i$, $\eta = dg_i/g_i$ and $\eta^\prime
= dg_i^\prime/g_i^\prime$; but $g_i^\prime = gg_i$ so $\eta^\prime =
\eta + dg/g$ completing the proof of the claim.
\end{proof}
This ends the proof of
Proposition~\ref{Propositio:transvaffineforms}.\end{proof}

Using the above and some other techniques like the index theorem of
Camacho-Sad (cf.\cite{Camacho-Sad})  and some linearization results
like those in Chapter~\ref{chapter:limitsets}, it is possible to
prove that:

\begin{Theorem}[\cite{Scardua1}] Let $\fa$ be a codimension one
foliation on $\co P(2)$ which is transversely affine outside an
algebraic codimension one invariant subset $S \subset \co P(2)$.
Suppose that $\fa$ has reduced non-degenerate singularities in $S$.
Then $\fa$ is a logarithmic foliation.
\end{Theorem}

\section{Transversely projective foliations}
Let $M$ be  a complex manifold and $\fa$ be a codimension one
holomorphic foliation with singularities of $M$. Recall that $\fa$
is called {\it transversely projective} if the underlying
``non-singular" foliation $\fa_0=:\fa\big|_{M\setminus \sing(\fa)}$
is transversely projective. This means that there is an open cover
$\bigcup\limits_{j\in J} U_j = M\setminus \sing(\fa)$ such that in
each $U_j$ the foliation is given by a submersion $f_j\colon U_j \to
\ov{\mathbb C}$ and if $U_i \cap U_j \ne \emptyset$ then we have
$f_i = f_{ij}\circ f_j$ in $U_i \cap U_j$ where $f_{ij}\colon U_i
\cap U_j \to P\SL(2,{\mathbb C})$ is locally constant. Thus,  on
each intersection $U_i \cap U_j \ne \emptyset$, we have $f_i =
\frac{a_{ij}f_j+b_{ij}}{c_{ij}f_j+d_{ij}}$ for some locally constant
functions $a_{ij}, b_{ij}, c_{ij}, d_{ij}$ with $a_{ij}d_{ij} -
b_{ij}c_{ij} = 1$.

\subsection{Development of a transversely projective foliation - Touzet's work}
\label{subsection:development}

We recall the notion of development of a transversely projective
foliation, first mentioned in the beginning of this section, already
adapting it to  our current framework. Let $\G$ be a (non-singular)
holomorphic foliation on a complex surface $N$. Suppose  that $\G$
is transversely projective in $N$. There is a Galoisian (i.e., a
transitive) covering $\pi\colon P \to N$ where $\pi$ is holomorphic,
a homomorphism $h\colon \pi_1(N) \to P\SL(2,{\mathbb C})$ and a
holomorphic submersion $\Phi\colon P \to {{\mathbb C}  P}^1$ such
that:
\begin{itemize}
\item[(i)] $\Phi$ is $h$-{\it equivariant}. This means that
for any homotopy class $[\gamma]\in \pi_1(N)$, we have
\[
h([\gamma]) (\Phi(x)) = \Phi(\widetilde {[\gamma]}(x)), \, \forall x
\in M\setminus S\, ,
\]
where by $\widetilde {[\gamma]}\colon P \to P$ we denote the
covering map induced by $[\gamma]$ in the Galoisian covering
$p\colon P \to N$.

\item[(ii)] $\pi^*\big(\G\big|_{N}\big)$ is the
foliation defined by the submersion $\Phi$.
\end{itemize}

In the above construction of the development, we may take $P$ as the
universal covering $\pi\colon \widetilde{N}\to N$ of $N$. We shall
refer to the submersion  $\Theta\colon \widetilde{N} \to {{\mathbb
C}P(1)}$ as a {\it multiform first integral} \index{first integral!
multiform} of $\G$ given by the projective structure in $N$.   Given
a homotopy class $[\gamma] \in \pi_1(M\setminus S)$, the
corresponding {\it monodromy map} is the image $h([\gamma])\subset
P\SL(2,{\mathbb C})$.

\begin{Definition}
\label{definition:monodromy} {\rm The {\it global monodromy} of the
foliation, with respect to this development,  is the image
$\Mon(\G)=h(\pi_1(N))\subset P\SL(2,{\mathbb C})$. }
\end{Definition}

\medskip

\begin{Remark}{\rm Some remarks about the above construction are as
follows. The construction of the development in \cite{Godbillon}
requires the foliation to be nonsingular. Assume now that $\fa$ is a
foliation with singular set of codimension $\geq 2$ on a complex
manifold $M$. Then $N=M\setminus \sing(\fa)$ is a complex manifold
and $\mathcal G:= \fa\big|_{N}$ is non-singular. By definition $\fa$
is transversely projective if and only if $\mathcal G$ is
transversely projective. Moreover, since $\sing(\fa)\subset M$ has
real codimension $\geq 4$, we conclude that there is a natural
isomorphism $\pi_1(N)\cong \pi_1(M)$. In particular, we can assume
in the above construction that $M=N$, i.e., the notion of
development above introduced can be introduced for foliations with
singularities.   Finally, thanks to Hartogs' extension theorem
(\cite{Gunning 1}), any holomorphic map  from $M\backslash
\sing(\fa) $ to ${\mathbb C}P(1)$  extends uniquely to a holomorphic
map from  $M$ to ${\mathbb C}P(1)$. }
\end{Remark}

Using the notion of monodromy and its properties, F. Touzet has been
able to study the analytic  classification of irreducible
singularities which have a suitable projective transverse structure
off its set of separatrices. He calls such a projective structure
{\it of moderate growth}, meaning that the foliation admits a
meromorphic projective triple defined in a neighborhood of the
singularity. He proves the following:
\begin{Theorem}[cf. \cite{Touzet}, Theorem II.4.2]
Let  $\fa$ a germ of irreducible singularity  of resonant type or of
saddle-node type.  Then  the foliation admits a meromorphic
projective triple near  the singularity if and only if on a
neighborhood of $0\in \mathbb C^2$, $\fa$ is the pull-back of a
Riccati foliation on $\ov \bc \times \ov \bc$ by a meromorphic map.
\end{Theorem}

The proof of this theorem is based on the study and classification
of the Martinet-Ramis cocycles of the singularity expressed in terms
of some classifying holonomy map of a separatrix of the singularity

For the non-resonant case, without the need of the moderate growth
hypothesis he proves:
\begin{Theorem}[\cite{Touzet},
Theorem II.3.1\index{Theorem! of Touzet}] A nondegenerate
nonresonant singularity $xdy - \lambda y dx + \omega_2(x,y)=0, \,
\lambda \in \mathbb C \setminus \mathbb Q_+$,  is analytically
linearizable if and only if the corresponding foliation $\fa$ is
transversely projective in $U\setminus \sep(\fa,U)$ for some
neighborhood $U$ of the singularity.
\end{Theorem}

Another interesting work in this direction is
\cite{Berthier-Touzet}, where the authors study the case of
irreducible singularities with Liouvillian first integral, in the
sense of M. Singer \cite{singer}.

\subsection{Projective structures and differential forms}

Let $\fa$ be  a codimension one holomorphic  foliation with singular
set $\sing(\fa)$ of codimension $\ge 2$ on a complex manifold $M$.
The existence of a projective transverse structure for $\fa$ is
equivalent to the existence of suitable triples of differential
forms as follows:

\begin{Proposition} [\cite{Scardua1}, Proposition 1.1 page 190]
\label{Proposition:forms}  Assume that $\fa$ is given by an
integrable holomorphic $1$-form $\Om$ on $M$ and suppose that there
exists a holomorphic $1$-form $\eta$ on $M$ such that
$\text{\rm{\it(Proj.1)} } d\Om = \eta \wedge \Om$. Then $\fa$ is
transversely projective of $M$ if and only if there exists a
holomorphic $1$-form $\xi$ on $M$ such that $\text{\rm{\it(Proj.2)}
} d\eta = \Om \wedge \xi$ and $\text{\rm{\it(Proj.3)} }d\xi = \xi
\wedge \eta$.
\end{Proposition}

The above proposition helps in the description of some examples of
transversely projective foliations:

\begin{Example}
{\rm
Let $\alpha$ be a closed meromorphic $1$-form on $M$ and let
$f\colon M \to \ov{\mathbb C}$ be a meromorphic function. Define
$(\Omega,\eta,\xi)$ by: $\Omega = df - f^2\alpha,\quad \eta =
2f\alpha \quad\text{and}\quad \xi = 2\alpha.$ Then
$(\Omega,\eta,\xi)$ is a projective triple and therefore $\Omega$
defines a holomorphic foliation of $M$, transversely projective in
the complement of the analytic invariant codimension one set $S
\subset M$, $S = (\alpha)_\infty \cup (f)_\infty$. The same
conclusion holds for $\Omega_\lambda = \Omega + \lambda\alpha$,
where $\lambda\in\mathbb C$. The foliation $\fa(\Omega_\lambda)$ is
also transversely affine in some smaller open set of the form
$M\backslash S^\prime$ where $S^\prime \supset S$, $S^\prime = S
\cup (f^2-\lambda=0)$. (In fact $\frac{\Om_\lambda}{f^2-\lambda} =
\frac{df}{f^2-\lambda} - \alpha$ is closed and holomorphic in
$M\backslash S^\prime$). }
\end{Example}
\begin{Example} {\rm Let $h\colon M\to{\mathbb C^*}$ be holomorphic such that $d\xi = -
\frac{dh}{2h} \wedge \xi$ where $\xi$ is holomorphic. (We can write
this condition as $d(\sqrt h.\xi)=0)$. Let $F$ be any holomorphic
function and write (for $\lambda\in\mathbb C$) $\Omega =
F\cdot\left(\frac{dF}{F} - \frac 12 \frac{dh}{h}\right) -
\left(\frac{F^2}{2} - \frac{\lambda}{2}h\right).\xi,\, \eta = \frac
12 \frac{dh}{h} + F\cdot\xi.$ The triple $(\Omega,\eta,\xi)$
satisfies the conditions of Proposition~\ref{Proposition:forms} and
then $\fa = \fa(\Omega)$ is a transversely projective foliation of
$M$.}
\end{Example}

\subsubsection{Proof of Proposition~\ref{Proposition:forms}}
Let us now give a proof for Proposition~\ref{Proposition:forms}. We
start with a remark about its need.
\begin{Remark}
{\rm Proposition~\ref{Proposition:forms} is stated (for the real
non-singular case)
with an idea of its proof, in \cite{Godbillon} (see Prop. 3.20, pp. 262).
However, it seems that the suggested proof uses some triviality
hypothesis on principal fiber-bundles of structural group
$\Aff(\mathbb C)$, over the manifold $M$ (see \cite{Godbillon} Prop. 3.6 pp. 249-250).
In our case  this is replaced by the existence of the form $\eta$ in
the statement. On the other hand, since some  of its elements
will be useful later, we supply a proof for Proposition~\ref{Proposition:forms}.

}
\end{Remark}

We will use the two following lemmas whose proofs are
straightforward computations or consequence of Darboux-Lie theorem,
Theorem~\ref{Theorem:Darboux-Lie}, therefore left to the reader:
\begin{Lemma}
\label{Lemma:1.1} Let $x,y,\widetilde x,\widetilde y\colon
U\subset\mathbb C^n\to\ov{\mathbb C}$ be meromorphic functions
satisfying:
(i) $ydx-xdy = \widetilde y d\widetilde x - \widetilde x d\widetilde y$;
(ii)$\frac{\widetilde x}{\widetilde y} =
\frac{ax+by}{cx+dy}$\,,\quad $\begin{pmatrix} a &b\\ c
&d\end{pmatrix} \in P\SL(2,\mathbb C)$.

\noindent Then $\widetilde x = \ve.(ax+by)$ and $\widetilde y =
\ve.(cx+dy)$ for some $\ve \in \mathbb C$, $\ve^2 = 1$.
\end{Lemma}

\begin{Lemma}\label{Lemma:1.2} Let $x,y,\widetilde x,
\widetilde y\colon U \subset \mathbb C^n \to \ov{\mathbb C}$ be
meromorphic functions satisfying $\widetilde x = ax+by$, $\widetilde
y = cx+dy$ for some $\begin{pmatrix} a &b\\ c &d\end{pmatrix} \in
P\SL(2,\mathbb C)$. Then $xdy-ydx = \widetilde x d\widetilde y -
\widetilde y d\widetilde x$.
\end{Lemma}
\begin{proof}[Proof of Proposition~\ref{Proposition:forms}] Suppose $\fa$ is transversely projective in
$M^n$, say, $\{f_i\colon U_i \to \mathbb C\}$ is a projective
transverse structure for $\fa$ in $M\backslash s(\fa)$. In each
$U_i$  we have $\Om = -g_i\,df_i$ for some holomorphic $g_i \in
\O(U_i)^*$. In each $U_i \cap U_j \ne \phi$ we have: $g_i\,df_i =
g_j\,df_j$ and $ (1) \, f_i = \frac{a_{ij}f_j +
b_{ij}}{c_{ij}f_j+d_{ij}}$ as in
Definition~\ref{Definition:transprojfolnonsing}.  Since $d\Om =
d(-g_i\,df_i) = \frac{dg_i}{g_i} \wedge \Om$ we have $\eta =
\frac{dg_i}{g_i} - h_i\Om$ for some holomorphic $h_i$ in $U_i$. We
define $x_i,y_i,u_i,v_i\colon U_i \to \mathbb C$ in the following
way: $ (2) \, y_i^2 = g_i, \quad \frac{x_i}{y_i} = f_i, \quad h_i =
\frac{2v_i}{y_i}\quad\text{and}\quad x_iv_i-y_iu_i = 1. $ Thus we
have: $\Om = x_i\,dy_i - y_i\,dx_i$ and $(3)\, \eta = 2(v_i\,dx_i -
u_i\,dy_i)$. This motivates us to define local models (see
\cite{Godbillon} Section 3.18 pp. 261): $\xi_i = 2(v_i\,du_i -
u_i\,dv_i) \quad\text{in}\quad U_i.$ It is easy to check that we
have $d\xi_i = \xi_i \wedge \eta, \quad d\eta = \Om\wedge\xi_i
\quad\text{in}\quad U_i.$ We can assume that $dx_i$ and $dy_i$ are
independent for all $i \in I$. In fact $dx_i \wedge dy_i = 0
\Rightarrow d\Om\big|_{U_i} = 2 \,dx_i\wedge dy_i = 0 \Rightarrow
d\Om = 0$ in $M$ (we can assume $M$ to be connected) $\Rightarrow$
we have $0 = d\Om = \eta\wedge\Om$ so that $\eta = h\Om$ for some
holomorphic function $h\colon M \to \mathbb C \Rightarrow$ we can
choose $\xi = \frac{h^2\Om}{2} + h\eta + dh$ which satisfies the
relations $d\eta = \Om\wedge\xi$ and $d\xi = \xi\wedge\eta$.
\smallskip
\begin{Claim}\label{Claim:1} $\xi_i=\xi_j$ in each $U_i \cap U_j \ne \phi$ and therefore
the $\xi_i$'s can be glued into a holomorphic $1$-form $\xi$ in
$M\backslash s(\fa)$ satisfying the conditions of the statement.
\end{Claim}
\begin{proof} From (1) and (2) we obtain
$\frac{x_i}{y_i} = \frac{a_{ij}x_j+b_{ij}y_j}{c_{ij}x_j +
d_{ij}y_j}$\,. Therefore according to Lemma~\ref{Lemma:1.1} we have
$(4)\, x_i = \ve.(a_{ij}x_j + b_{ij}x_j),\, y_i = \ve.(c_{ij}x_j +
d_{ij}y_j)\, \ve^2 = 1.$ Using (3) and (4) we obtain: $(a_{ij}v_i -
c_{ij}u_i)dx_j + (b_{ij}v_i - d_{ij}u_i)dy_j = \ve.(v_j\,dx_j -
u_j\,dy_j)$ and therefore: $(5) \, v_j =\epsilon( a_{ij}\,v_i -
c_{ij}\,u_i), \, \,u_j = \epsilon (-b_{ij}\,v_i + d_{ij}\,u_j)$. It
follows form (5) and Lemma~\ref{Lemma:1.2} that $v_i\,du_i -
u_i\,dv_i = v_j\,du_j - u_j\,dv_j$ which proves the claim.
\end{proof}

\begin{Claim}\label{Claim:2} We have $\xi = \xi_i = h_i^2 \frac{\Om}{2} + h_i\eta + dh_i$
in each $U_i$.
\end{Claim}
\begin{proof} We have
$h_i^2\Om = \frac{4v_i^2}{y_i^2}\,(x_i\,dy_i - y_i\,dx_i)$, $h_i\eta
= \frac{4v_i}{y_i}\,(v_i\,dx_i - u_i\,dy_i)$, $dh_i =
2d\left(\frac{v_i}{y_i}\right)$. \noindent  Hence
$\frac{h_i^2\Om}{4} + \frac{h_i\eta}{2} + \frac{dh_i}{2} =
\frac{v_i^2}{y_i}\,dx_i - \frac{v_i}{y_i^2}\,(x_iv_i-1)dy_i +
\frac{dv_i}{y_i}$\,.

\noindent  On the other hand a straightforward calculation shows
that $\frac{\xi_i}{2} = v_i\,du_i - u_i\,dv_i =
\frac{v_i^2}{y_i}\,dx_i - \frac{v_i}{y_i}(x_iv_i-1)dy_i +
\frac{dv_i}{y_i}\,.$ And thus Claim~\ref{Claim:2} is proved.
\end{proof}

Since $\codim\,s(\fa) \ge 2$ it follows that $\xi$ extends
holomorphically to $M$. This proves the first part. Now we assume
that $(\Om,\eta,\xi)$ is {\it holomorphic} as in the statement of
the proposition:
\begin{Claim}\label{Claim:3} Given any $p \in M\backslash s(\fa)$ there exist
holomorphic functions $x,y,u,v\colon U \to \mathbb C$ defined in an
open neighborhood $U\ni p$ such that: $\Om = xdy-ydx$, $\eta =
2(vdx-udy)$ and $\xi = 2(vdu-udv)$.
\end{Claim}
\begin{proof} This claim is a consequence of Darboux's Theorem~\ref{Theorem:Darboux-Lie} (see also
\cite{Godbillon} pp. 230), but we can give an alternative proof as
follows: We write locally $\Om = -gdf = xdy-ydx$ and $\eta =
\frac{dg}{g} - h\Om = 2(vdx-udy)$ as in the proof of the first part.
Using Claim~\ref{Claim:2} and the last part of
Proposition~\ref{Proposition: 2.I} below we obtain locally $\xi =
\frac{h^2\Om}{2} + h\eta + dh + \ell.\Om$; for some holomorphic
function $\ell$ satisfying $\frac{d\ell}{-2\ell} \wedge \Om = d\Om$.
This last equality implies that $d(\sqrt\ell.\Om) = 0$ and then
$\ell = \frac{r(f)}{g^2}$ for some holomorphic function $r(z)$. Now
we look for holomorphic functions $\widetilde f$, $\widetilde g$ and
$\widetilde h$ satisfying: $\Om = -\widetilde g d\widetilde f, \quad
\eta = \frac{d\widetilde g}{\widetilde g} - \widetilde h\Om$ and
$\xi = \frac{\widetilde h^2\Om}{2} + \widetilde h\eta + d\widetilde
h$. We try $\widetilde f = U(f)$ for some holomorphic non-vanishing
$U(z)$. Using $\Om = gdf = -\widetilde g d\widetilde f$ we get
$\widetilde g = \frac{g}{U^\prime(f)}$\,. Using $\eta = \frac{dg}{g}
- d\Om = \frac{d\widetilde g}{\widetilde g} - \widetilde h\Om$ we
get $\widetilde h = h- \frac{U^{\prime\prime}}{gU^\prime}$\,. Using
$\xi = \frac{h^2\Om}{2} + h\eta + dh + \ell\Om = \frac{\widetilde
h^2\Om}{2} + \widetilde h\eta + d\widetilde h$ we get
$d\left(\frac{U^{\prime \prime}(f)}{U^\prime(f)}\right) = r(f)df$.

\noindent  Therefore it is possible to write $\Om$, $\eta$ and $\xi$
as in the statement of the claim: define $x = \widetilde f y$, $y =
\sqrt{\widetilde g}$, $v = \frac{\widetilde hy}{2}$ and $u =
\frac{xv-1}{y}$ as in the first part of the proof. This proves Claim
3.
\end{proof}

\noindent  Using Claim~\ref{Claim:3}  we prove that $\fa$ is
transversely projective in $M\backslash s(\fa)$, that is in $M$. The
last part of Proposition~\ref{Proposition:forms} can be proved using
the relation stated above between the projective structure and the
local trivializations for $\Om$, $\eta$ and $\xi$. For instance we
prove the following.

\begin{Claim}\label{Claim:4} $(\Om,\eta,\xi)$ and $(f\Om, \eta + \frac{df}{f}, \frac
1f\,\xi)$ define the same projective structure for $\fa$, for any
holomorphic $f\colon M \to \mathbb C^*$.
\end{Claim}

\begin{proof} Using the notation of the first part we define $\hat x_i = \sqrt
f.\,x_i$, $\hat y_i = \sqrt f.\,y_i$, $\hat u_i = \frac{1}{\sqrt
f}\,.\,u_i$ and $\hat v_i = \frac{1}{\sqrt f}\,.\,v_i$. Then: $f\Om
= \hat x_i\,d\hat y_i - \hat y_i\,d\hat x_i$, $\eta+\frac{df}{f} =
2(\hat v_i\,d\hat x_i - \hat u_i\,d\hat y_i)$ and $\frac 1f\,\xi =
2(\hat v_i\,d\hat u_i - \hat u_i\,d\hat v_i)$. Furthermore we have
$\frac{\hat x_i}{\hat y_i} = \frac{x_i}{y_i} = \frac{a_{ij}x_j +
b_{ij}y_j}{c_{ij}x_j + d_{ij}y_j} = \frac{a_{ij}\hat x_j +
b_{ij}\hat y_j}{c_{ij}\hat x_j + d_{ij}\hat y_j},$ and this proves
the claim and finishes the holomorphic part of the proof.
\end{proof}
Now we only have to observe that if $(\Om, \eta$) is a pair of
meromorphic $1$-forms and if $\fa$ is transversely projective in
$M$, then the same steps of the first part of the proof apply to
construct a meromorphic $1$-form $\xi$ satisfying the relations of
the statement.
\end{proof}

Let $\fa$ be  a codimension one holomorphic  foliation with singular
set $\sing(\fa)$ of codimension $\ge 2$ on a complex manifold $M$.
As mentioned in the Introduction, the existence of a projective
transverse structure for $\fa$ is equivalent to the existence of
suitable triples of differential forms (cf.
Proposition~\ref{Proposition:forms}, see \cite{Scardua1} Section 3,
page 193):

\noindent This motivates the following definition:

\begin{Definition}[projective triple]
\label{Definition:projectivetriple} {\rm Given holomorphic $1$-forms
(respectively, meromorphic $1$-forms) $\Om$, $\eta$ and $\xi$ on $M$
we shall say that $(\Om,\eta, \xi)$ is a {\it holomorphic projective
triple\/} (respectively, a {\it meromorphic projective triple\/}) if
they satisfy relations {\it(Proj.1)}, {\it(Proj.2)} and
{\it(Proj.3)} above. The foliation $\fa^\perp$ defined by the
$1$-form $\xi$ is called {\it transverse foliation}\index{foliation!
transverse} corresponding to the projective triple. If $\eta$ is not
identically zero then $\fa^\perp$ is really a foliation of $M$ which
is transverse to $\fa$ outside of a proper analytic subset. }

\end{Definition}

The following definition  plays a fundamental role in the theory of
transversely projective foliations.

\begin{Definition}[moderate growth (transversely projective foliations)]
\label{Definition:moderategrowth} {\rm A foliation $\fa$ of $M$ will
be called {\it transversely projective of moderate growth} if it
admits a meromorphic projective triple defined in $M$. This means
that $\fa$ is transversely projective in some the complementar of
some analytic subset $M\subset M$ of codimension one.}
\end{Definition}
The termonilogy {\it foliation with moderate growth} has already
been introduced in \cite{Touzet}. With the above definitions,
Proposition~\ref{Proposition:forms} says that $\fa$ is transversely
projective on $M$ if and only if the holomorphic pair ($\Om$,
$\eta$) may be completed to a holomorphic projective triple.
Moreover, a foliation $\fa$ which is transversely projective of
moderate growth exhibits a projective transverse structure $\mathcal
P$ in the complement of some codimension divisor $D\subset M$ ($D$
contained in the polar set of the projective triple). One question
then is whether the projective transverse structure $\mathcal P$
extends to the divisor $D$. The other question, apparently simpler,
is whether the foliation $\fa$ is actually projective of moderate
growth. According to \cite{Scardua1} we may  perform modifications
in a  projective triple as follows:

\begin{Proposition}[\cite{Scardua1}]
\label{Proposition:modificationforms} Let $M$ be a connected complex
manifold.
\begin{itemize}
\item[\rm(i)] Given a meromorphic projective triple $(\Om, \eta,
\xi)$ and meromorphic functions $g$, $h$ on $M$ we can define a new
meromorphic projective triple as follows:

{\rm {\it(Mod.1)}}\,\, $\Om' = g\,\Om$

{\rm {\it(Mod.2)}}\,\, $\eta' = \eta + \frac{dg}{g} + h\,\Om$

{\rm {\it(Mod.3)}}\,\, $\xi' = \frac 1g\,\big(\xi - dh - h\eta -
\frac{h^2}{2}\,\Om\big)$

\item[\rm(ii)] Two holomorphic projective triples $(\Om,\eta,\xi)$
and $(\Om', \eta', \xi')$ define the same projective transverse
structure for a given foliation $\fa$ if and only if  we have {\rm
{\it(Mod.1)}, {\it(Mod.2)}} and {\rm {\it(Mod.3)}} for some
holomorphic functions $g$, $h$ with $g$ non-vanishing.

\item[{\rm(iii)}]
Let $(\Om,\eta,\xi)$ and $(\Om, \eta, \xi')$ be meromorphic
projective triples. Then $\xi' = \xi +F\,\Om$ for some meromorphic
function $F$ in $M$ with $d\,\Om = -\frac 12\, \frac{dF}{F} \wedge
\Om$.

\end{itemize}

\end{Proposition}

\smallskip

\noindent This last proposition   implies that  suitable meromorphic
projective triples also define projective transverse structures.
\noindent We can rewrite condition (iii) on $F$ as $d(\sqrt {F}
\,\Om) = 0$. This implies that if the projective triples $(\Omega,
\eta, \xi)$ and $(\Omega, \eta, \xi ^\prime)$ are not identical then
the foliation defined by $\Omega$ is transversely affine outside the
codimension one analytical invariant subset $S=\{F=0\}\cup
\{F=\infty\}$. (\cite{Scardua1}).

This approach is useful because of the following result:
\begin{Theorem} [\cite{Scardua1} Theorem 4.1 page 197]
\label{Theorem:xifirstintegral} Let $\fa$ be a foliation of $M$
where $M$ is a polydisc $M\subset {\mathbb C}^m$ or a projective
manifold over $\mathbb C$ of dimension $m \geq 2$. Assume that $\fa$
admits a meromorphic projective triple $(\Omega, \eta, \xi)$ defined
in $M$. If $\xi$ admits a meromorphic first integral in $U$ then
$\fa$ is a meromorphic pull-back of a Riccati foliation.
\end{Theorem}

\begin{proof}

By hypothesis, $\xi$ defines a foliation which admits a meromorphic
first integral. Since we are either on a projective manifold or in a
polydisc centered at the origin, we  can write $\xi = g\,dR$ for
some meromorphic  functions $g$  and $R$ (these functions are
rational in the case of a projective surface). Then we may replace
the meromorphic  triple $(\Om, \eta,\xi)$ by $(\Om',\eta',\xi')$
where $\Om' = g\Om$, \, $\eta' = \eta + \frac{dg}{g}$ \,and\, $\xi'
= \frac 1g\,\xi = dR$.  The relations $d\Om' = \eta' \wedge \xi'$,
\,\, $d\eta' = \Om' \wedge \xi'$, \,\, $d\xi' = \xi \wedge \eta'$
imply that $\eta' = HdR$ for some meromorphic function $H$.  Now we
define $\om := \frac{H^2}{2}\,\xi' - H\eta' + dH = \frac 12\, H^2dR
+ dH$ $1$-form such that $d\om = -HdH \wedge dR$. On the other hand
$\eta' \wedge \om = HdR \wedge dH = -HdH \wedge dR$. Thus $d\om =
\eta' \wedge \om$. We also have $d\eta' = dH \wedge dR = (-\frac
12\, H^2dR + dH) \wedge dR = \om \wedge \xi'$. The meromorphic
triple $(\om,\eta',\xi')$ satisfies the projective relations $d\om =
\eta' \wedge \om$, \, $d\eta' = \om \wedge \xi'$, \, $d\xi' = \xi'
\wedge \eta'$ and therefore by
Proposition~\ref{Proposition:modificationforms} (iii)  we conclude
that $\Om' = \om + F.\xi'$ for some meromorphic function $F$ such
that $d\xi' = \xi' \wedge \frac{1}{2}\frac{dF}{F}\,\cdot$ This
implies $dF \wedge dR \equiv 0$. By the classical Stein
Factorization theorem(\cite{G-R})  we may assume from the beginning
that $R$ has connected fibers and therefore $dF \wedge dR \equiv 0$
implies $F = \vr(R)$ for some one-variable meromorphic function
$\vr(z) \in {\mathbb C}(z)$. {\em In the case where  $M$ is a
projective manifold all the meromorphic objects are rational and
therefore $\vr(z)$ is also a rational function.} We obtain therefore
$\Om' = -\frac 12\, H^2dR + dH + \vr(R)dR = = dH - (\frac 12\, H^2 -
\vr(R))dR.$ If we define a meromorphic  map $\sigma\colon
M\dashrightarrow \ov{\mathbb C}\times\ov{\mathbb C}$ \, by \,
$\sigma(p) = \big(R(p), H(p)\big)$ then clearly $\Om' =
\sigma^*(dy-(\frac 12\, y^2 - \vr(x))dx)$ and therefore $\fa$ is the
pull-back $\fa = \sigma^*(\eR)$ of the Riccati foliation $\eR$ given
on $\ov{\mathbb C}\times\ov{\mathbb C}$ by the meromorphic (rational
if $M$ is a projective manifold) $1$-form $\Omega_\vr:=dy -(\frac
12\, y^2 - \vr(x))dx$.
\end{proof}

\begin{Definition} {\rm A meromorphic projective triple
$(\Omega ', \eta ', \xi ')$ is {\it geometric}  if it can be written
locally as in  {\it(Mod.1)}, {\it(Mod.2)} and {\it(Mod.3)} for some
(locally defined) holomorphic projective triple $(\Om, \eta, \xi)$
and some (locally defined) meromorphic functions.}
\end{Definition}

\smallskip

\noindent As an immediate consequence we obtain:

\begin{Proposition}
\label{Proposition:true} A geometric projective triple $(\Om',
\eta', \xi')$ defines a transversely projective foliation $\fa$
given by $\Om'$ on $M$.
\end{Proposition}

\subsubsection{Classification of projective foliations: moderate growth on
projective manifolds} In \cite{Loray-Touzet-Vitorio} we find the
following definition of transversely projective foliation of a
smooth projective manifold. {\sl Let M be a smooth projective
manifold over $\mathbb C$. A (holomorphic singular) codimension one
foliation $\fa$  of $M$. The foliation is said to be {\it
transversely projective} if given a  non zero rational $1$-form
$\omega$  defining $\fa$ (and therefore satisfying the Frobenius
integrability condition $\omega \wedge d \omega  = 0$) we have that
there are {\em rational} $1$-forms $\alpha$ and $\beta$ on $M$ such
that the $sl_2$-connection on the rank $2$ trivial vector bundle
defined by $\Delta = d + \begin{pmatrix} \alpha & \beta \\ \omega &
- \alpha \end{pmatrix}$ is flat. }

Let us compare the above definition with the one we have been using
so far in this survey. Indeed, compared to
Definition~\ref{Definition:transprojfolnonsing}  there is a
difference, quite easy to explain. In the above definition, we
already assume that the foliation admits a rational  projective
triple, i.e., a projective triple meromorphic defined everywhere in
the manifold $M$. This is not necessarily the case if we just start
with a foliation which is (according to our definition
Definition~\ref{Definition:transprojfolnonsing}) transversely
projective in $M\setminus S$ for some algebraic curve $S\subset M$.
Nevertheless, often we cannot extend the projective transverse
structure to the curve $S$ (for instance,  in the case of Riccati
foliations or logarithmic foliations). Thus what is considered in
\cite{Loray-Touzet-Vitorio} are what we have called {\em
transversely projective foliations with moderate growth} (cf.
Definition~\ref{Definition:projectivetriple}). projective structure
in $M\setminus S$.

The authors also introduce the following notion:
\begin{Definition}[\cite{Loray-Touzet-Vitorio}]
{\rm A {\it Riccati foliation}\index{foliation! Riccati} over a
projective manifold $M$ consists of a pair $(\pi\colon  P \to M;H) =
(P;H)$ where $\pi \colon  P \to M$ is a locally trivial $\mathbb
P(1)$ fiber bundle in the Zariski topology, this means that  $P$ is
the projectivization of the total space of a rank two vector bundle
$E$,  and $H$ is a codimension one foliation on $P$ which is
transverse to a general fiber of $\pi$. In the case of a clear
context,  the $\mathbb P(1)$-bundle $P$ is omitted from the
notation. Then  $H$ is called a {\it Riccati foliation}. The
foliation $H$ is defined by the projectivization of horizontal
sections of a (non unique) at meromorphic connection $r$ on $E$. The
connection $r$ is uniquely determined by $H$ and its trace on
$\det(E)$. We say that the Riccati foliation $H$ is {\it regular} if
it lifts to a meromorphic connection $r$ with at worst regular
singularities (see \cite{Deligne}), and irregular if not. It is said
that a Riccati foliation $(P;H)$ over $M$ {\it factors} through a
projective manifold $M^\prime$  if there exists a Riccati foliation
$(\pi^\prime \colon P^\prime \to M^\prime, H ^\prime)$ over
$M^\prime$, and rational maps $\phi \colon M \dashrightarrow M
^\prime$ and $\Phi \colon P \dashrightarrow P ^\prime$,   such that
$\pi^\prime \circ\Phi = \phi \circ \pi$, and  $\Phi$ has degree one
when restricted to a general fiber of $P$, and $H = \Phi^*H
^\prime$.

}
\end{Definition}

Using the notion above, alternatively, in
\cite{Loray-Touzet-Vitorio} the authors state that a foliation
$\mathcal F$ of $M$ is transversely projective if there exists a
triple $\mathcal P = (P;H; \sigma)$ satisfying \begin{enumerate}
\item  $(P;H)$ is a Riccati
foliation over $M$; and

\item $\sigma : M \dashrightarrow P$ is a rational section
generically transverse to $H$ such that $\fa= \sigma^*H$.
\end{enumerate}

After making the conversion between the notions of transversely
projective foliation in \cite{Loray-Touzet-Vitorio} and the one we
consider in this text, we can state the main classification result
of \cite{Loray-Touzet-Vitorio} as follows:

\begin{Theorem}[cf. \cite{Loray-Touzet-Vitorio}, Theorem D]
Let $\mathcal F$ be a codimension one transversely projective
foliation \underline{of moderate growth} on a projective manifold
$M$. Then at least one of the following assertions holds true.

\begin{enumerate}

\item  There exists a generically finite Galois morphism
$f \colon Y \to M$ such that $f^*\fa$ is defined by a closed
rational $1$-form.

\item  There exists a rational map $f \colon
M \dashrightarrow S$ to a ruled surface $S$, and a Riccati foliation
$\mathcal R$ on $S$ such that $\fa = f^*R$.

\item  The transverse projective structure for $\fa$
has at worst regular singularities, and the monodromy representation
of $\fa$ factors through one of the tautological representations of
a polydisc Shimura modular orbifold $\mathcal H$.

\end{enumerate}
\end{Theorem}

There is still a number of interesting questions, on the local and
on the global framework, about the classification and the
description of foliations with projective transverse structure.

\chapter*{Exercises}

\begin{itemize}
\item[{\bf 1.}] Given the linear vector
field $X(x,y) = x\,\dfrac{\po}{\po x} - y\,\dfrac{\po}{\po y}$ on
$\bc^2$  describe the global picture of the corresponding foliation
on $\bc P(2)$ (the foliation exhibits three singularities on $\bc
P(2)$, two of which are dicritical and require one blow-up).
\item[{\bf 2.}] Let $\fa_1$ and $\fa_2$ be two holomorphic foliations with singularities
on $\bc P(2)$. Show that if $\fa_1 \ne \fa_2$ and they have a common
leaf the this leaf  is algebraic (hint.:\, Let $\fa_j$ be given by
the polynomial vector field $X_j$ on $\bc^2$. Then a local
parametrization $\big(x(z),y(z)\big)$, $z \in \bd$ of the common
leaf must satisfy $\big(x'(z),y'(z)\big) = \la_j(z)\cdot
X_j\big(x(z),y(z)\big)$, $\forall\, z \in \bd$ for some holomorphic
$\la_j(z)$. Then write $X_j = (P_j,Q_j)$ where $P_j$\,, $Q_j$ are
polynomials to conclude that $\dfrac{P_1}{Q_1} = \dfrac{P_2}{Q_2}$
on the leaf).
\item[{\bf 3.}] State and prove the following local form of closed meromorphic $1$-forms: \newline {\it Let $\omega$
be a closed meromorphic $1$-form on a neighborhood of the
origin\linebreak $0 \in \bc^n$, $n\ge2$.  Then there is a
neighborhood $0 \in U \subset \bc^n$ where $\omega$ is defined and
writes
$$
\omega\big\vert_U = \sum_{j=1}^r \la_j\, \frac{df_j}{f_j} +
d(g/\prod\limits_{j=1}^r f_j^{n_j-1})
$$
for some holomorphic $f_j,g\colon U \to \bc$,\, $g \ne 0$,\, $\la_j
\in \bc$,\, $n_j \in \bn$.}
\item[{\bf 4.}] Prove that an isolated singularity of a holomorphic vector field $X$ on $\bc^n$,
say $0 \in \bc^n$, which is in the Poincaré domain; it is
necessarily {\it transverse\/} to  the small spheres
$S^{2n-1}(0;\ve)$ (of radius $\ve > 0$) centered at the singular
point. {\rm Where, by {\it transverse\/} we mean transversality
between the leaves of the foliation and the sphere, as real
submanifolds of $\mathbb R^{2n}$. }

\item[{\bf 5.}] Let $\fa$ be a germ of foliation singularity at $0 \in \bc^2$. Assume
that $\fa$ is given by a closed meromorphic one form $\omega$ with
simple poles in a  neighborhood of $0 \in \bc^2$. Prove that $\fa$
is analytically linearizable. \newline Sug.:\, First consider the
case where $\fa$ is not a saddle-node. Show that the holonomy of a
separatrix of $\fa$ is analytically linearizable. In order to do
this, show that in a neighborhood of a point $p \ne 0$ belonging to
a separatrix $\Ga$ we can choose local coordinates such that
$\Ga\colon (y=0)$, $p\colon (x=y=0)$ and $\om(x,y) =
a\,\dfrac{dy}{y}$ where $(a=\Res_\Ga\,\om)$. Then conclude that if
$(\widetilde{x},\widetilde{y})$ are similar coordinates then we have
$\widetilde{y}=$ const. $y$. Using then the Martinet-Ramis formal
normal form (\cite{martinet-ramisselano}), get rid of the
saddle-node case.
\item[{\bf 6.}] Let $G \subset \Diff(\bc,0)$ be an abelian analytically linearizable subgroup
containing an attractor say $f \in G$ with $|f'(0)| < 1$. Suppose
that for every point $z \in (\bc,0)$  we have
$\overline{O(z)}\setminus O(z) \subset \{0\}$. Prove that $G$ is
generated by $f$ and some rational rotation $g(z) = e^{\frac{2\pi
i}{\nu}}\,z$,\, $\nu \in \bz$.
\item[{\bf 7.}] Prove that a non-dicritical germ of a holomorphic foliation admitting a meromorphic
first integral, necessarily admits a holomorphic first integral.
\item[{\bf 8.}] Let $\fa$ be a foliation on $\bc P(2)$ (holomorphic with singularities). Assume that the
limit set $\lim(\fa)$ is algebraic, consisting of points and a
finite number of (invariant) algebraic  curves $\La_j \subset \bc
P(2)$, $j=1,\dots,r$. Given a point $p \in \La_j\setminus\sing(\fa)$
and a transverse disc $p \in \Sigma \cap \La_j$ show that the
virtual holonomy group $G = \Hol^{\virt}(\fa, \La_j, \Sigma, p)$
satisfies the following property: $\forall\, z \in (\Sigma,p)$, \,
$\overline{O(z)}\setminus O(z) \subset \{p\}$.
\item[{\bf 9.}] In the situation of Exercise 8 above assume that the
reduction of singularities of a component $\La_{j_0} \subset
\lim(\fa)$ exhibits only invariant components (i.e., $\sing(\fa)
\cap \La_{j_0}$ is non-dicritical) and each such component has an
attractor on its virtual holonomy. Prove that there are no
saddle-nodes in the reduction of singularities of $\La_{j_0}$\,.
\item[{\bf 10.}] Let $\fa$ be a holomorphic foliation on $\bc P(2)$,
given on $\bc^2$ by a closed meromorphic $1$-form $\omega$ (not
necessarily rational $1$-form). Show that: \begin{itemize}
\item[{\rm(i)}] If the line at infinity $\ell_\infty := \bc
P(2)\setminus \bc^2$ is not $\fa$-invariant then $\om$ is rational,
i.e., $\om$ admits an extension to $\bc P(2)$. \item[{\rm(ii)}] If
$\ell_\infty$ is $\fa$-invariant and contains some irreducible
singularity of type $xdy - \la ydx +\cdots = 0$\,\, $\la \notin
\bq_+$\,, then $\om$ also admits an extension to $\bc P(2)$.
\end{itemize} \item[{\bf 11.}] Complete the details in
Example~\ref{Example:implicit}. \item[{\bf 12.}] Show that there is
no holomorphic foliation of dimension $k$ without singularities on
the complex projective space $\mathbb CP(n)$ for $1\leq k \leq n-1$.
Is it also true for smooth foliations of even dimension? \item[{\bf
13.}] Is it true that a one-dimensional holomorphic foliation with
isolated singularities $\fa$ on the complex projective space
$\mathbb CP(n)$ is given by a rational vector field? \end{itemize}

\chapter*{Some open questions}

Here are some open questions that the author thinks are relevant in
the framework of holomorphic foliations with singularities.

\begin{Question}

Let $\fa$ be a polynomial vector field on the complex affine space
$\mathbb C^3$. Assume that for infinitely many of its orbits they
are complete intersection of two algebraic surfaces on $\mathbb
C^3$. Is is true that the vector field admits a strong rational
first integral, i.e., a rational map $R\colon \mathbb C^3
\dashrightarrow \mathbb C^2$ such that $R$ is constant on each orbit
of $X$.
\end{Question}

\begin{Question}
Let $X$ be a germ of a holomorphic vector field at the origin $0 \in
\mathbb C^3$. Assume that $X$ admits a formal strong  first
integral, i.e., a pair of formal  functions $\hat f, \hat g \in \hat
{\mathcal O}_3$ such that $d\hat f (X)=0$ and $d\hat g(X)=0$. Is is
true that there is a convergent strong first integral $F=(f,g)$ with
$f, g \in \mathcal O_3$?

\end{Question}

\begin{Question}
Let be given a polynomial vector field $X$ on $\mathbb C^n$. Assume
that the set of algebraic orbits of $X$ has positive measure on
$\mathbb C^n$. Is it true that $X$ admits some type of algebraic
first integral?
\end{Question}

\begin{Question}
Let $G\subset \Diff(\mathbb C^2,0)$ be a subgroup of germs of
complex diffeomorphisms at $0 \in \mathbb C^2$. Assume that $G$ has
the origin as an stable fixed point, in the sense of Lyapunov. Is it
true that $G$ is analytically linearizable?

\end{Question}

\begin{Question}
Let $X$ be a polynomial vector field on $\mathbb C^n$ and denote by
$\fa$ the corresponding one-dimensional foliation of the complex
projective space $\mathbb CP(n)$. Assume that there is an algebraic
curve $\Lambda\subset \mathbb CP(n)$ which is irreducible, invariant
by $\fa$  and stable in the sense of Lyapunov. What is the normal
form of $\fa$? Is it true that if the singularities in $\Lambda$ are
hyperbolic then $\fa$ is given by a linear vector field in some
coordinate chart?

\end{Question}

\begin{Question}
Let $\fa$ be a germ of one-dimensional holomorphic foliation at
$0\in \mathbb C^n$. Assume that $\fa$ is induced by a vector field
with a non-resonant singularity at the origin. Suppose that $\fa$ is
transversely homogeneous in the complement of some invariant
analytic hypersurface germ at the origin. Is it true that the vector
field germ of $X$ at $0$ is analytically conjugate to its formal
normal form?
\end{Question}

\begin{Question}
Let $\fa$ be a one-dimensional holomorphic foliation of the complex
projective space $\mathbb CP(n)$. Assume that: (i) there is an
analytic codimension one subset $\Lambda \subset \mathbb CP(n)$ such
that $\fa$ is transversely homogeneous in the complement $\mathbb
CP(n)\setminus \Lambda$; (ii) the singularities of $\fa$ in
$\Lambda$ are {\it generic}. What is the classification of $\fa$?
Does $\fa$ admit some sort of Liouvillian first integral?

\end{Question}

\bibliographystyle{amsalpha}

\printindex

\end{document}